\documentclass[11pt]{amsart}
\usepackage{amssymb}
\usepackage{lscape}
\usepackage[all]{xy}
\usepackage{color}
\usepackage{epic,eepic}

\newcommand{\red}{\color{red}}
\newcommand{\blue}{\color{blue}}

\newtheorem{Thm}{Theorem}[section]
\newtheorem{Lem}[Thm]{Lemma}
\newtheorem{Cor}[Thm]{Corollary}
\newtheorem{Prop}[Thm]{Proposition}
\newtheorem{Conj}[Thm]{Conjecture}
\newtheorem{Rem}[Thm]{Remark}

\newcommand{\A}{\mathbb{A}}
\newcommand{\D}{\mathbb{D}}
\newcommand{\E}{\mathbb{E}}

\newcommand{\Z}{\mathbb{Z}}

\newcommand{\N}{\mathbb{N}}
\newcommand{\C}{\mathbb{C}}

\newcommand{\df}{\colon}

\newcommand{\id}{{\rm id}}

\newcommand{\F}{{\mathcal F}}
\newcommand{\U}{{\mathcal U}}
\newcommand{\G}{{\mathcal G}}
\newcommand{\T}{{\mathcal T}}
\newcommand{\X}{{\mathcal X}}
\newcommand{\DD}{{\mathcal D}}
\newcommand{\RR}{{\mathcal R}}
\newcommand{\II}{{\mathcal I}}
\newcommand{\MM}{{\mathcal M}}
\newcommand{\cA}{{\mathcal A}}
\newcommand{\cS}{{\mathcal S}}
\newcommand{\cP}{{\mathcal P}}

\newcommand{\CC}{{\mathcal C}}
\newcommand{\stCC}{{\underline{\mathcal C}}}

\newcommand{\g}{\mathfrak{g}}
\newcommand{\n}{\mathfrak{n}}
\newcommand{\h}{\mathfrak{h}}

\newcommand{\mm}{\mathbf{m}}
\newcommand{\nn}{\mathbf{n}}
\newcommand{\ee}{\mathbf{e}}
\newcommand{\aaa}{\mathbf{a}}

\newcommand{\ii}{{\mathbf i}}
\newcommand{\jj}{{\mathbf j}}
\newcommand{\kk}{{\mathbf k}}
\newcommand{\yy}{{\mathbf y}}

\newcommand{\hI}{\widehat{I}}

\newcommand{\LL}{\Lambda}
\newcommand{\GG}{\Gamma}
\newcommand{\si}{\sigma}
\newcommand{\vpi}{\varpi}
\newcommand{\vph}{\varphi}
\newcommand{\la}{\lambda}

\newcommand{\ad}{\operatorname{ad}}

\newcommand{\rk}{\operatorname{rk}}
\newcommand{\orb}{\mathcal O}

\newcommand{\md}{\operatorname{mod}}
\newcommand{\nil}{\operatorname{nil}}
\newcommand{\rep}{\operatorname{rep}}

\newcommand{\add}{\operatorname{add}}
\newcommand{\Gen}{\operatorname{Fac}}
\newcommand{\Cogen}{\operatorname{Sub}}

\newcommand{\gldim}{\operatorname{gl.dim}}
\newcommand{\pdim}{\operatorname{proj.dim}}
\newcommand{\soc}{\operatorname{soc}}

\newcommand{\Hom}{\operatorname{Hom}}
\newcommand{\Ext}{\operatorname{Ext}}
\newcommand{\End}{\operatorname{End}}

\newcommand{\Ima}{\operatorname{Im}}
\newcommand{\Ker}{\operatorname{Ker}}
\newcommand{\Coker}{\operatorname{Coker}}

\newcommand{\dm}{{\rm dim}\,}
\newcommand{\dimv}{\underline{\dim}}

\newcommand{\Span}{\operatorname{Span}}

\newcommand{\bil}[1]{\langle #1\rangle}
\newcommand{\abs}[1]{\left| #1\right|}
\newcommand{\ebrace}[1]{\langle #1 \rangle}

\newcommand{\op}{{\rm op}}
\newcommand{\Res}{\operatorname{Res}}

\newcommand{\irr}{\operatorname{Irr}}

\newcommand{\bsm}{\begin{smallmatrix}}
\newcommand{\esm}{\end{smallmatrix}}
\newcommand{\bbm}{\begin{matrix}}
\newcommand{\ebm}{\end{matrix}}

\newcommand{\GL}{\operatorname{GL}}

\newcommand{\Gmin}{G^{\rm min}}
\newcommand{\Nmin}{N^{\rm min}}
\newcommand{\Norm}{\operatorname{Norm}}
\newcommand{\Lie}{\operatorname{Lie}}

\def\shuff#1#2{\mathbin{
\hbox{\vbox{ \hbox{\vrule \hskip#2 \vrule height#1 width 0pt}%
\hrule}%
\vbox{ \hbox{\vrule \hskip#2 \vrule height#1 width 0pt
\vrule}%
\hrule}%
}}}

\def\SHUF{{\mathchoice{\shuff{7pt}{3.5pt}}%
{\shuff{6pt}{3pt}}%
{\shuff{4pt}{2pt}}%
{\shuff{3pt}{1.5pt}}}}%

\def\shuffle{\,\SHUF\,\,}

\makeindex

\begin{document}

\date{April 5, 2008}

\bigskip
\title[Cluster algebra structures and semicanonical bases]
{Cluster algebra structures and semicanonical bases for
unipotent groups}

\author{Christof Gei{\ss}}
\address{Christof Gei{ss}\newline
Instituto de Matem\'aticas\newline
Universidad Nacional Aut{\'o}noma de M{\'e}xico\newline
Ciudad Universitaria\newline
04510 M{\'e}xico D.F.\newline
M{\'e}xico}
\email{christof@math.unam.mx}

\author{Bernard Leclerc}
\address{Bernard Leclerc\newline
Laboratoire LMNO\newline
Universit\'e de Caen\newline
F-14032 Caen Cedex\newline
France}
\email{leclerc@math.unicaen.fr}

\author{Jan Schr\"oer}
\address{Jan Schr\"oer\newline
Mathematisches Institut\newline
Universit\"at Bonn\newline
Beringstr. 1\newline
D-53115 Bonn\newline
Germany}
\email{schroer@math.uni-bonn.de}



\begin{abstract}
Let $Q$ be a finite quiver without oriented cycles, and let
$\LL$ be the associated preprojective algebra.
To each terminal $\C Q$-module $M$ (these are certain
preinjective $\C Q$-modules), we attach
a natural subcategory $\CC_M$ of $\md(\LL)$.
We show that $\CC_M$ is a Frobenius category,
and that
its stable category $\underline{\CC}_M$ is a 
Calabi-Yau category of dimension two.

Then we develop a theory of mutations of maximal rigid
objects of $\CC_M$, analogous to the mutations of clusters
in Fomin and Zelevinsky's theory of cluster algebras.
We also provide an explicit quasi-hereditary structure
on the endomorphism algebra of a distinguished maximal 
rigid object of $\CC_M$, and we use it to describe the 
combinatorics of mutations. 

Next, we show that $\CC_M$ yields a categorification 
of a cluster algebra $\cA(\CC_M)$, which is
not acyclic in general.
We give a realization of $\cA(\CC_M)$ as a
subalgebra of the graded dual of the enveloping algebra
$U(\n)$, where $\n$ is a maximal nilpotent subalgebra of the 
symmetric Kac-Moody Lie algebra $\g$ associated to the quiver $Q$.

Let $\cS^*$ be the dual of Lusztig's semicanonical
basis $\cS$ of $U(\n)$.
We show that all cluster monomials of $\cA(\CC_M)$
belong to $\cS^*$, and that
$\cS^* \cap \cA(\CC_M)$ is a $\C$-basis 
of $\cA(\CC_M)$.

Next, we prove that $\cA(\CC_M)$ is naturally isomorphic
to the coordinate ring 
$\C[N(w)]$ of the finite-dimensional unipotent
subgroup $N(w)$ of the Kac-Moody group $G$ attached to $\g$.
Here $w = w(M)$ is the adaptable element of the Weyl group of $\g$
which we associate to each terminal $\C Q$-module $M$.

Moreover, we show that the cluster algebra obtained from 
$\cA(\CC_M)$ by formally inverting the generators
of the coefficient ring is isomorphic to the algebra 
$\C[N^w]$ of regular
functions on the unipotent cell $N^w :=  N \cap (B_-wB_-)$ 
of $G$.
We obtain a corresponding dual semicanonical 
basis of $\C[N^w]$.

Finally, by ``specializing coefficients''
we obtain a dual semicanonical basis for
a coefficient free cluster algebra $\cA_w$ associated to $w$.
As a special case,
we obtain a dual semicanonical basis
of the
(coefficient free) acyclic cluster algebras $\cA_Q$
associated to $Q$, which naturally extends the set of cluster monomials
in $\cA_Q$.
\end{abstract}

\maketitle

\setcounter{tocdepth}{1}

\tableofcontents

\parskip2mm



{\Large\part{Introduction and main results}}



{\Large\section{Introduction}}


\subsection{Introduction}
This is the continuation of an extensive project to obtain a 
better understanding
of the relations between the following topics:
\begin{itemize}

\item[(i)]
Representation theory of quivers,

\item[(ii)]
Representation theory of preprojective algebras,

\item[(iii)]
Lusztig's (semi)canonical basis of universal enveloping algebras,

\item[(iv)]
Fomin and Zelevinsky's theory of cluster algebras,

\item[(v)]
Frobenius categories and 2-Calabi-Yau categories, 

\item[(vi)]
Cluster algebra structures on coordinate algebras of 
unipotent groups, Bruhat cells and flag varieties.

\end{itemize}
The topics (i) and (iii) are closely related.
The numerous connections have been studied by many authors.
Let us just mention Lusztig's work on canonical bases of quantum
groups, and Ringel's Hall algebra approach to quantum groups.
An important link between (ii) and (iii), due to Lusztig
\cite{Lu1,Lu2} and Kashiwara and Saito \cite{KSa}
is that the elements of the (semi)canonical basis are naturally
parametrized by the irreducible components of the varieties of
nilpotent representations of a preprojective algebra.

Cluster algebras were invented by Fomin and Zelevinsky
\cite{BFZ,FZ1,FZ2},
with the aim of providing a new algebraic and combinatorial setting
for canonical bases and total positivity.
One important breakthrough was the insight that 
the class of acyclic cluster algebras with a skew-symmetric exchange matrix
can be categorified using the so-called cluster categories.
Cluster categories were introduced by Buan, Marsh, Reineke, Reiten and
Todorov \cite{BMRRT}.
In a series of papers by some of these authors and also by Caldero and
Keller \cite{CK,CK2}, it was established that cluster categories have
all necessary properties to provide the mentioned categorification.
We refer to the nice overview article \cite{BM} for more details on
the development of this beautiful theory which established a
strong connection between the topics (i), (iv) and (v).

In \cite{GLSRigid} we showed that the representation theory of
preprojective algebras $\LL$ of Dynkin type 
(i.e. type $\A$, $\D$ or $\E$)
is also closely related to cluster algebras.
We proved that $\md(\LL)$ can be regarded as a categorification
of a natural cluster structure on the polynomial algebra $\C[N]$.
Here $N$ is a maximal unipotent subgroup of a complex Lie group 
of the same type as $\LL$. Let $\n$ be its Lie algebra,
and $U(\n)$ be the universal enveloping algebra of $\n$.
The graded dual $U(\n)^*_{\rm gr}$ can be identified with the coordinate
algebra $\C[N]$.
By means of our categorification, we were able to prove that 
all the cluster monomials of $\C[N]$ belong to the dual
of Lusztig's semicanonical basis of $U(\n)$.
Note that the cluster algebra $\C[N]$ is in general not acyclic.

The aim of this article is to extend these results to the more
general setting of Kac-Moody groups and their unipotent cells.
We also provide additional tools for studying these
categories and cluster structures.
For example we show that the endomorphism algebras of 
certain maximal rigid modules are quasi-hereditary and deduce from
this a new combinatorial algorithm for mutations.

More precisely, we consider preprojective algebras 
$\LL=\LL_Q$ attached to quivers $Q$ which are not necessarily of Dynkin type. 
These algebras are therefore infinite-dimensional in general.
The category $\nil(\LL)$ of all finite-dimensional nilpotent
representations of $\LL$ is then too large to be related to a 
cluster algebra of finite rank.
Moreover, it does not have projective or injective objects,
and it lacks an Auslander-Reiten translation. 
However, we give a general procedure to attach to certain 
preinjective representations $M$ of $Q$ a natural 
subcategory $\CC_M$ of $\nil(\LL)$.
We show that these subcategories $\CC_M$ are Frobenius categories
and that the corresponding stable categories ${\stCC}_M$ are 
Calabi-Yau categories of dimension two.
Each subcategory $\CC_M$ comes with two distinguished maximal rigid modules
$T_M$ and $T_M^\vee$ described combinatorially.
In the special case where $Q$ is of Dynkin type and $M$ is the
sum of all indecomposable representations of $Q$ (up to isomorphism)
we have $\CC_M = \md(\LL)$, the modules $T_M$ and $T_M^\vee$
are those constructed in \cite{GLSAus}, and we recover the
setting of \cite{GLSRigid}.
In another direction, if $Q$ is an arbitrary (acyclic) quiver
and $M = I \oplus \tau(I)$, where $I$ is the sum of the indecomposable
injective representations of $Q$ and $\tau$ is the Auslander-Reiten
translation, it follows from a result of Keller and Reiten \cite{KR}
that the stable category $\underline{\CC}_M$ is triangle equivalent
to the cluster category $\CC_Q$ of \cite{BMRRT}. 
We provide in this case a natural functor $G\df \CC_M \to \CC_Q$
inducing an equivalence  $\underline{G}\df \underline{\CC}_M \to \CC_Q$
of additive categories.

We then develop, as in \cite{GLSRigid}, a theory 
of mutations for maximal rigid objects $T$ in $\CC_M$,
and we study their endomorphism algebras $\End_\LL(T)$.
We show that these algebras have global dimension 3 and
that their quiver has neither loops nor 2-cycles.
Special attention is given to the algebra $B:=\End_\LL(T_M)$
for which we provide an explicit quasi-hereditary structure.
We prove that $\CC_M$ is anti-equivalent to the category
of $\Delta$-filtered $B$-modules.
This allows us to describe the mutations of maximal
rigid $\LL$-modules in terms of the $\Delta$-dimension
vectors of the corresponding $\End_\LL(T_M)$-modules.
We also exhibit a simple sequence of mutations between $T_M$
and $T_M^\vee$ and describe all the maximal rigid modules 
arising from this sequence.

In the last part we associate to the subcategory $\CC_M$ 
a cluster algebra $\cA(\CC_M)$ which in general is
not acyclic,
and we show that 
$\CC_M$ can be seen as a categorification of
$\cA(\CC_M)$. 
(As a very special case, we also obtain in this way a new 
categorification of every acyclic cluster algebra with a
skew-symmetric exchange matrix and a certain choice of coefficients.)
The proof relies on the fact that
the algebra $\cA(\CC_M)$ has a natural realization
as a certain subalgebra of the graded dual $U(\n)^*_{\rm gr}$, where 
$\n$ is now the positive part of the Kac-Moody Lie algebra 
$\g = \n_- \oplus \h \oplus \n$ 
of the same type as $\LL$. We show that again all the
cluster monomials belong to the dual
of Lusztig's semicanonical basis of $U(\n)$. 
Next, we prove that $\cA(\CC_M)$ has a 
simple monomial basis
coming from the objects of the additive
closure $\add(M)$ of $M$. 
We call it the {\it dual PBW-basis} of $\cA(\CC_M)$,
and regard it as a generalization (in the dual setting) of the bases 
of $U(\n)$ constructed by Ringel in terms of quiver representations,
when $\g$ is finite-dimensional \cite {Ri6}. 
We use this to prove that $\cA(\CC_M)$ is spanned
by a subset of the dual semicanonical basis of $U(\n)^*_{\rm gr}$.
Thus, we obtain another natural basis of $\cA(\CC_M)$
containing all the cluster monomials.
We call it the {\it dual  semicanonical basis} of $\cA(\CC_M)$.
Finally, we prove that $\cA(\CC_M)$ is isomorphic
to the coordinate ring of a finite-dimensional unipotent
subgroup of the Kac-Moody group $G$ attached to $\g$.
Moreover, we show that the cluster algebra obtained from 
$\cA(\CC_M)$ by formally inverting the generators
of the coefficient ring is isomorphic to the algebra of regular
functions on a certain unipotent cell of $G$.    
This solves Conjecture III.3.1 of \cite{BIRS} for all unipotent
cells attached to an {\it adaptable} element $w$ of the Weyl
group of $G$ (the definition of adaptable is given below, see \S~\ref{wadapt}).
Note also that in the Dynkin case, we recover a result of
\cite{BFZ} for the double Bruhat cells of type $(e,w)$
with $w$ adaptable, but our proof is different and shows
that the coordinate ring of the cell is not only an
upper cluster algebra but a genuine cluster algebra. 
In the last section, we explain how the results of this
paper are related to those of \cite{GLSFlag}, in which
a cluster algebra structure on the coordinate ring of the unipotent radical
$N_K$ of a parabolic subgroup of a complex simple algebraic group of
type $\A, \D, \E$ was introduced. We give a proof of 
Conjecture~9.6 of \cite{GLSFlag} in the case where 
the Weyl group element $w_0w_0^K$ is adaptable.

Our results have some overlap with the recent work of 
Buan, Iyama, Reiten and Scott \cite{BIRS}.
Up to a simple duality, our categories $\CC_M$ coincide
with the categories $\CC_w$ introduced in \cite{BIRS}, but
only for adaptable Weyl group elements $w$, and our maximal
rigid modules $T_M^\vee$ are some of the cluster-tilting objects
of the categories $\CC_w$ described in \cite{BIRS}. 
However our methods are very different, and for our smaller
class of categories we can prove stronger results,
like the existence of quasi-hereditary endomorphism algebras or
the existence of semicanonical bases for the corresponding cluster algebras.

\subsection{Plan of the paper}
This article is organized as follows.

In Sections~\ref{section2}, \ref{section3}, we give a more detailed 
presentation of our results.

Part~\ref{part2} is devoted to the study of the subcategories
$\CC_M$.
Some known results on quiver representations and preprojective algebras
are collected in Section~\ref{quivers}.
In Section~\ref{section4} we introduce the important concept of a 
{\it selfinjective torsion class} of $\md(\LL)$. 
Some technical but crucial results on the lifting of certain $KQ$-module
homomorphisms to $\LL$-module homomorphisms are proved in 
Section~\ref{liftsection}.
These results are used in Section~\ref{completerigid} to construct a
$\CC_M$-complete rigid module $T_M$ and to compute the quiver of its
endomorphism algebra.
Then we show in Section~\ref{section7} that $\CC_M$ is a Frobenius 
category whose stable category $\stCC_M$ is a 2-Calabi-Yau category. 
In particular, it turns out that $\CC_M$ is a selfinjective torsion
class of $\md(\LL)$.
In Sections~\ref{section8} and \ref{section9} 
we prove some basic properties of $\CC$-maximal rigid modules, where
$\CC$ is now an arbitrary selfinjective torsion
class of $\md(\LL)$.
In Section~\ref{clustercat}
we show that for every quiver $Q$ without oriented cycles 
there exists a terminal $KQ$-module $M$ such that the
stable category $\underline{\CC}_M$ is triangle equivalent
to the cluster category $\CC_Q$ as defined in \cite{BMRRT}.
This uses a recent result by Keller and Reiten \cite{KR}.
We also construct an explicit functor
$\CC_M \to \CC_Q$ which then yields an equivalence 
of additive categories $\stCC_M \to \CC_Q$.

Part~\ref{part3} develops the theory of mutations of rigid objects
of $\CC_M$.
Sections~\ref{section11}, \ref{endosection} and \ref{graphmutation} 
contain the adaptation of the results in \cite[Sections 5,\,6,\,7]{GLSRigid} 
to our more general situation of selfinjective torsion classes.
In Section~\ref{section15} we prove that cluster variables
are determined by their dimension vector, in the appropriate sense,
and that one can describe the mutation of clusters in terms of
these dimension vectors.
We also obtain a characterization
of all short exact sequences of $\LL$-modules which become split
exact sequences of $KQ$-modules after applying the restriction
functor $\pi_Q$.
We show in Section~\ref{quasihereditary} that 
$\End_\LL(T_M)$ is a quasi-hereditary algebra.
This can be used to reformulate the 
results in Section~\ref{section15} in
terms of $\Delta$-dimension vectors, see Section~\ref{section16}.
In Section~\ref{section17} we construct a sequence of mutations
starting with our module $T_M$ which yields generalizations of
classical determinantal identities.

Part~\ref{part4} contains the applications of the previous 
constructions to cluster algebras.
In Section~\ref{semican} we repeat several known results
on Kac-Moody Lie algebras, and also recall our results 
about the 
multiplicative behaviour of the functions $\delta_X$.
One of the central parts of our theory is the construction
of dual PBW- and dual semicanonical bases for the
cluster algebras $\RR(\CC_M)$.
This is done in Section~\ref{PBWsection}.
The special case of acyclic cluster algebras is discussed in 
Section~\ref{acycliccase}.
In Section~\ref{coordinaterings} we 
prove all our results on cluster algebra structures of coordinate
rings.
Finally, we present some open problems in
Section~\ref{openproblems}.

\subsection{Notation}
Throughout let $K$ be an algebraically closed field.
For a $K$-algebra $A$ let $\md(A)$ be the category of
finite-dimensional left $A$-modules.
By a {\it module} we always mean a finite-dimensional left module.
Often we do not distinguish between a module and its isomorphism 
class.
Let ${\rm D} = \Hom_K(-,K)\df \md(A) \to \md(A^\op)$ be the usual
duality. 

For a quiver $Q$ let $\rep(Q)$ be the category of finite-dimensional
representations of $Q$ over $K$.
It is well known that
we can identify $\rep(Q)$ and $\md(KQ)$.

By a {\it subcategory} we always mean a full subcategory.
For an $A$-module $M$ let $\add(M)$ be the subcategory of all
$A$-modules
which are isomorphic to finite direct sums of direct summands of $M$.
A subcategory $\U$ of $\md(A)$ 
is an {\it additive subcategory} if any finite direct
sum of modules in $\U$ is again in $\U$.
By $\Gen(M)$ 
\index{$\Gen(M)$}
(resp. $\Cogen(M)$) 
\index{$\Cogen(M)$}
we denote the subcategory of
all $A$-modules $X$ such that there exists some $t \ge 1$ and some
epimorphism $M^t \to X$ (resp. monomorphism $X \to M^t$).

For an $A$-module $M$ let $\Sigma(M)$ be the number of isomorphism classes
of indecomposable direct summands of $M$.
An $A$-module is called 
{\it basic} 
\index{basic module}
if it can be written as a direct sum
of pairwise non-isomorphic indecomposable modules.

For a vector space $V$ we sometimes write
$|V|$ for the dimension of $V$.
If $f\df X \to Y$ and $g\df Y \to Z$ are maps, then the composition
is denoted by $gf = g \circ f\df X \to Z$.

If $U$ is a subset of a $K$-vector space $V$, then let
$\Span_K\ebrace{U}$ be the subspace of $V$ generated by
$U$.

By $K(X_1,\ldots,X_r)$ (resp. $K[X_1,\ldots,X_r]$)
we denote the field of rational functions (resp. the polynomial ring)
in the variables $X_1,\ldots,X_r$ with coefficients in $K$.

Let $\C$ be the field of complex numbers, and let $\N = \{0,1,2,\ldots\}$
be the natural numbers, including $0$.
Set $\N_1 := \N \setminus \{ 0 \}$.
For natural numbers $a \le b$ let $[a,b] = \{ i \in \N \mid a \le i \le b\}$.

Recommended introductions to representation theory of finite-dimensional
algebras and Auslander-Reiten theory are the books \cite{ARS,ASS,GR,Ri1}.


{\Large\section{Main results: Rigid modules over preprojective algebras}
\label{section2}}


\subsection{Preprojective algebras}
Throughout, let $Q$ be a finite quiver without oriented cycles, and let
$$
\LL = \LL_Q = K\overline{Q}/(c)
$$ 
be the associated 
{\it preprojective algebra}.
\index{preprojective algebra}
We assume that $Q$ is connected and has vertices
$\{1,\ldots,n\}$ with $n$ at least two.
Here $K$ is an algebraically closed field,  
$K\overline{Q}$ is the path algebra of the 
{\it double quiver } $\overline{Q}$ 
\index{double quiver}
\index{$\overline{Q}$}
of $Q$ which is obtained
from $Q$ by adding to each arrow $a\df i \to j$ in $Q$ an arrow
$a^*\df j \to i$ pointing in the opposite direction, and $(c)$ is
the ideal generated by the element
$$
c = \sum_{a \in Q_1} (a^*a - aa^*)
$$
where $Q_1$ is the set of arrows of $Q$.
Clearly, the path algebra $KQ$ is a subalgebra of $\LL$. 
We denote by 
$$
\pi_Q\df \md(\LL) \to \md(KQ)
$$
the corresponding restriction functor.

\subsection{Terminal $KQ$-modules}\label{terminal}
Let $\tau = \tau_Q$ be the Auslander-Reiten translation of $KQ$, and let
$I_1,\ldots,I_n$ be the indecomposable injective $KQ$-modules.
A $KQ$-module $M$ is called 
{\it preinjective} 
\index{preinjective $KQ$-module}
if $M$ is isomorphic
to a direct sum of modules of the form
$
\tau^j(I_i)
$
where $j \ge 0$ and $1 \le i \le n$.
There is the dual notion of a preprojective module.

A $KQ$-module  $M = M_1 \oplus \cdots \oplus M_r$ with 
$M_i$ indecomposable and $M_i \not\cong M_j$ for all $i \not= j$
is called a 
{\it terminal $KQ$-module} 
\index{terminal $KQ$-module}
if the following hold:
\begin{itemize}

\item[(i)]
$M$ is preinjective;

\item[(ii)]
If $X$ is an indecomposable $KQ$-module with $\Hom_{KQ}(M,X) \not= 0$, then
$X \in \add(M)$;

\item[(iii)]
$I_i \in \add(M)$ for all indecomposable injective $KQ$-modules $I_i$.

\end{itemize}
In other words, the indecomposable direct summands of $M$ 
are the vertices of a subgraph of the preinjective component of the
Auslander-Reiten quiver of $KQ$ which is closed under successor.
We define
$$
t_i := t_i(M) := \max\left\{ j \ge 0 \mid \tau^j(I_i) \in \add(M) \setminus 
\{0\}\right\}.
$$
\index{$t_i(M)$}\noindent

\subsection{The subcategory $\CC_M$}
Let $M$ be a terminal $KQ$-module, and let
$$
\CC_M := \pi_Q^{-1}(\add(M))
$$
\index{$\CC_M$}\noindent
be the subcategory of all $\LL$-modules $X$ with
$\pi_Q(X) \in \add(M)$.
Notice that if $Q$ is a Dynkin quiver and 
$M$ is the sum of all indecomposable representations of $Q$
then $\CC_M=\md(\LL)$.

\begin{Thm}\label{main1}
Let $M$ be a terminal $KQ$-module.
Then the following hold:
\begin{itemize}

\item[(i)]
$\CC_M$ is a Frobenius category
with $n$ indecomposable $\CC_M$-projective-injectives;

\item[(ii)]
The stable category ${\stCC}_M$ is a 2-Calabi-Yau category;

\item[(iii)]
If $t_i(M) = 1$ for all $i$, then ${\stCC}_M$
is triangle equivalent to the cluster category 
$\CC_Q$ associated to $Q$.

\end{itemize}
\end{Thm}

Part (i) and (ii) of Theorem~\ref{main1} are proved in Section~\ref{section7}.
Based on results in Section~\ref{completerigid}, Part (iii) is shown
in Section~\ref{section9}.

\subsection{An example of type $\A_3$}\label{subsectexampleA3}
Let $Q$ be the quiver
{\small
$$
\xymatrix@-0.8pc{
1 \ar[dr]&&3\ar[dl]\\
&2
}
$$
}

\begin{figure}
$$
\xymatrix@-1.0pc{
&\blue{\fbox{$\bsm1&&0\\&1&\esm$}} \ar[dr] && 
\blue{\fbox{$\bsm0&&1\\&0&\esm$}}\\
{\bsm0&&0\\&1&\esm} \ar[ur]\ar[dr] && 
\blue{\fbox{$\bsm1&&1\\&1&\esm$}} \ar[ur]\ar[dr]\\
&\blue{\fbox{$\bsm0&&1\\&1&\esm$}} \ar[ur] && 
\blue{\fbox{$\bsm1&&0\\&0&\esm$}}
}
$$
\caption{A terminal module in type $\A_3$}
\label{fig0}
\end{figure}

Figure~\ref{fig0} shows the Auslander-Reiten quiver of $KQ$.
The indecomposable direct summands of a terminal $KQ$-module $M$
are marked in blue colour.
In Figure~\ref{fig1} we show the Auslander-Reiten quiver of
the preprojective algebra $\LL$ of type $\A_3$.
We display the graded dimension vectors of
the indecomposable $\LL$-modules.
(There is a Galois covering of $\LL$, and the associated push-down
functor is dense, a property which only holds for Dynkin
types $\A_n$, $n \le 4$.
In this case, all $\LL$-modules are uniquely 
determined (up to $\Z$-shift) by their
dimension vector in the covering.
See \cite{GLSSemi1} for more details.)
Note that one has to identify the objects on the two dotted vertical
lines.
The indecomposable $\CC_M$-projective-injective modules  
are marked in red colour, all other indecomposable modules in $\CC_M$
are marked in blue.
Observe that $\CC_M$ contains 7 indecomposable modules,
and three of these are $\CC_M$-projective-injective.
The stable category ${\stCC}_M$ is triangle equivalent to the
product $\CC_{\A_1} \times \CC_{\A_1}$ of two
cluster categories of type $\A_1$.

\begin{figure}
\[
\xymatrix@-1.0pc{
&&&\\
\red{\fbox{$\bsm0&&0\\&0\\1&&0\\&1\\0&&1\esm$}} \ar[dr]\ar@{.}[dd] &&&&&&
\red{\fbox{$\bsm0&&1\\&1\\1&&0\\&0\\0&&0\esm$}}\ar@{.}[dd]\\
&\blue{\fbox{$\bsm0&&0\\&0\\1&&0\\&1\\0&&0\esm$}} \ar[dr] &&
\blue{\fbox{$\bsm0&&0\\&0\\0&&1\\&0\\0&&0\esm$}} \ar[dr] &&
{\bsm0&&0\\&1\\1&&0\\&0\\0&&0\esm} \ar[dr]\ar[ur]\\
{\bsm0&&0\\&0\\0&&0\\&1\\0&&0\esm} \ar[dr]\ar[ur]\ar@{.}[dd] &&
\red{\fbox{$\bsm0&&0\\&0\\1&&1\\&1\\0&&0\esm$}} \ar[dr]\ar[r]\ar[ur] &
{\bsm0&&0\\&1\\1&&1\\&1\\0&&0\esm} \ar[r] &
{\bsm0&&0\\&1\\1&&1\\&0\\0&&0\esm} \ar[dr]\ar[ur] &&
{\bsm0&&0\\&1\\0&&0\\&0\\0&&0\esm}\ar@{.}[dd]\\
&\blue{\fbox{$\bsm0&&0\\&0\\0&&1\\&1\\0&&0\esm$}} \ar[ur] &&
\blue{\fbox{$\bsm0&&0\\&0\\1&&0\\&0\\0&&0\esm$}} \ar[ur] &&
{\bsm0&&0\\&1\\0&&1\\&0\\0&&0\esm} \ar[dr]\ar[ur]\\
\red{\fbox{$\bsm0&&0\\&0\\0&&1\\&1\\1&&0\esm$}} \ar[ur] &&&&&&
\red{\fbox{$\bsm1&&0\\&1\\0&&1\\&0\\0&&0\esm$}}
}
\]
\caption{
A category $\CC_M \subset \md(\LL)$ 
with $\stCC_M$ triangle equivalent to
$\CC_{\A_1 \times \A_1}$
}\label{fig1}
\end{figure}

\subsection{Maximal rigid modules and their endomorphism algebras}
A $\LL$-module $T$  is 
{\it rigid} 
\index{rigid module}
if $\Ext_\LL^1(T,T) = 0$.
For a module $X$ let 
$\Sigma(X)$ 
\index{$\Sigma(X)$}
be the number of isomorphism classes
of indecomposable direct summands of $X$.

Let $\CC$ be a subcategory of $\md(\LL)$.
Define the 
{\it rank} 
\index{rank of a subcategory}
of $\CC$ as
$$
\rk(\CC) = \max\{ \Sigma(T) \mid T \text{ rigid in } \CC \}
$$
\index{$\rk(\CC)$}\noindent
if such a maximum exists, and set $\rk(\CC) = \infty$, otherwise.

The category $\CC$ is called $Q$-{\it finite} if there exists some
$M \in \md(KQ)$ such that 
$$
\pi_Q(\CC) \subseteq \add(M).
$$
In this case, if $\CC$ is additive, 
one can imitate the proof of \cite[Theorem 1.1]{GSExt} 
to show that $\rk(\CC) \le \Sigma(M)$.
If $M$ is a terminal $KQ$-module, then we prove that $\rk(\CC_M) =
\Sigma(M)$,
see Corollary~\ref{completerigidT_M}.

Recall that for all $X, Y \in \md(\LL)$ we have
$\dm \Ext_\LL^1(X,Y) = \dm \Ext_\LL^1(Y,X)$, see \cite{CB} and also 
\cite{GLSSemi2}.
Assume that $T$ is a rigid $\LL$-module in an
additive subcategory $\CC$ of $\md(\LL)$ with $\rk(\CC) < \infty$.
We need the following definitions:
\begin{itemize}

\item
$T$ is $\CC$-{\it complete rigid} if $\Sigma(T) = \rk(\CC)$;
\index{complete rigid $\LL$-module}\noindent

\item
$T$ is $\CC$-{\it maximal rigid} if 
\index{maximal rigid $\LL$-module}\noindent
$
\Ext_\LL^1(T \oplus X,X) = 0
$
with $X \in \CC$ implies $X \in \add(T)$;

\item
$T$ is $\CC$-{\it maximal 1-orthogonal} if 
\index{maximal 1-orthogonal $\LL$-module}\noindent
$
\Ext_\LL^1(T,X) = 0
$
with $X \in \CC$
implies $X \in \add(T)$.

\end{itemize}
The notion of a maximal 1-orthogonal module is due to Iyama \cite{I1}.
These modules are also called {\it cluster tilting objects}.
\index{cluster tilting object}\noindent

\begin{Thm}\label{main9}
Let $M$ be a terminal $KQ$-module.
For a $\LL$-module $T$ in $\CC_M$ the following are equivalent:
\begin{itemize}

\item[(i)]
$T$ is $\CC_M$-complete rigid;

\item[(ii)]
$T$ is $\CC_M$-maximal rigid;

\item[(iii)]
$T$ is $\CC_M$-maximal 1-orthogonal.

\end{itemize}
If $T$ satisfies one of the above equivalent conditions, then
the following hold:
\begin{itemize}

\item
$\gldim(\End_\LL(T)) = 3$;

\item
The quiver $\GG_T$ of $\End_\LL(T)$ has no loops 
and no 2-cycles.

\end{itemize}
\end{Thm}

The proof of
Theorem~\ref{main9} can be found in Section~\ref{endosection}.

\subsection{The complete rigid modules $T_M$ and $T_M^\vee$}\label{existrigid}
Let $M = M_1 \oplus \cdots \oplus M_r$ be a terminal $KQ$-module.
Without loss of generality assume that $M_{r-n+1},\ldots,M_r$ are injective.
Let $\GG_M$ 
\index{$\GG_M$ (where $M$ is a terminal $KQ$-module)}\noindent
be the quiver of $\End_{KQ}(M)$.
Its vertices $1,\ldots,r$ correspond to $M_1,\ldots,M_r$, and the
number of arrows $i \to j$ equals
the number of arrows in the Auslander-Reiten quiver
of $\md(KQ)$ which start in $M_i$ and end in $M_j$.

Let $\GG_M^*$ 
\index{$\GG_M^*$ (where $M$ is a terminal $KQ$-module)}
be the quiver which is obtained from
$\GG_M$ by adding an arrow $i \to j$ whenever $M_j = \tau(M_i)$.
Our results in Sections~\ref{liftsection} and \ref{completerigid}
yield the following theorem:

\begin{Thm}\label{main6}
There exist two $\CC_M$-complete rigid $\LL$-modules $T_M$ and $T_M^\vee$
such that
$$\GG_{T_M} = \GG_{T_M^\vee} = \GG_M^*.$$
\end{Thm}

For the explicit description of $T_M$ and $T_M^\vee$ see
Section~\ref{completerigid}.
Here we just note that the $\CC_M$-projective direct summands of $T_M$
correspond to the rightmost vertices of $ \GG_M^*$, 
whereas the  $\CC_M$-projective direct summands of $T_M^\vee$
correspond to the leftmost vertices of $ \GG_M^*$.

\subsection{A quasi-hereditary algebra}

Now consider $B := \End_\LL(T_M)$.
We prove in Section~\ref{quasihereditary} the following theorem:

\begin{Thm}\label{main6'}
\begin{itemize}

\item[(i)]
$B$ is a quasi-hereditary algebra;
 
\item[(ii)]
The restriction of the contravariant functor 
$\Hom_\LL(-,T_M)\df \md(\LL) \to \md(B)$ induces an anti-equivalence
$
F\df \CC_M \to \F(\Delta)
$
where $\F(\Delta)$ is the category of $\Delta$-filtered $B$-modules
and 
$$
\Delta := \{ F(M_i) \mid 1 \le i \le r \}
$$
is the set of standard modules.
(We interpret $M_i$ as a $\LL$-module using the obvious embedding
functor.);

\item[(iii)] 
For a short exact sequence $0 \to X \to Y \to Z \to 0$ in $\CC_M$
the following are equivalent:
\begin{itemize}

\item[(a)]
The short exact sequence
$0 \to \pi_Q(X) \to \pi_Q(Y) \to \pi_Q(Z) \to 0$
splits;

\item[(b)]
The sequence 
$0 \to F(Z) \to F(Y) \to F(X) \to 0$
is exact.

\end{itemize}

\end{itemize}
\end{Thm}

It turns out that $\Hom_\LL(T_M^\vee,T_M)$ is the characteristic
tilting module over $B$.In particular, $\End_\LL(T_M^\vee)$ is also quasi-hereditary.


{\Large\section{Main results: Cluster algebras and semicanonical bases}
\label{section3}}


\subsection{The cluster algebra $\cA(\CC_M)$}\label{clustintro}
We refer to \cite{FZSurv} for an excellent survey on cluster algebras.  
Here we only recall the main definitions and introduce 
a cluster algebra $\cA(\CC_M,T)$ associated to
a terminal $KQ$-module $M$ and any $\CC_M$-maximal rigid module $T$.

If $\widetilde{B} = (b_{ij})$ is any $r \times (r-n)$-matrix
with integer entries, then
the {\it principal part} $B$ of 
$\widetilde{B}$ is obtained from 
$\widetilde{B}$ by deleting the last $n$ rows.
Given some $k \in [1,r-n]$ define a new $r \times (r-n)$-matrix 
$\mu_k(\widetilde{B}) = (b_{ij}')$ by
$$
b_{ij}' =
\begin{cases}
-b_{ij} & \text{if $i=k$ or $j=k$},\\
b_{ij} + \dfrac{|b_{ik}|b_{kj} + b_{ik}|b_{kj}|}{2} & \text{otherwise},
\end{cases}
$$
where $i \in [1,r]$ and $j \in [1,r-n]$.
One calls $\mu_k(\widetilde{B})$ a {\it mutation} 
\index{mutation of a matrix}
of $\widetilde{B}$.
If $\widetilde{B}$ is an integer matrix whose principal part is
skew-symmetric, then it is 
easy to check that $\mu_k(\widetilde{B})$ is also an integer matrix 
with skew-symmetric principal part.
In this case, Fomin and Zelevinsky define a cluster algebra
$\cA(\widetilde{B})$ as follows.
Let $\F = \C(y_1,\ldots,y_r)$ be the field of rational
functions in $r$ commuting variables 
$\yy = (y_1,\ldots,y_r)$. 
One calls $(\yy,\widetilde{B})$ the {\it initial seed} 
\index{initial seed}\noindent
of
$\cA(\widetilde{B})$.
For $1 \le k \le r-n$ define 
\begin{equation}\label{mutationformula}
y_k^* = 
\frac{\prod_{b_{ik}> 0} y_i^{b_{ik}} + \prod_{b_{ik}< 0} y_i^{-b_{ik}}}{y_k}.
\end{equation}
The pair 
$
(\mu_k(\yy),\mu_k(\widetilde{B}))
$, 
where
$\mu_k(\yy)$ is obtained from $\yy$ by replacing $y_k$ by 
$y_k^*$,
is the 
{\it mutation in direction} 
\index{mutation of a seed}
$k$ of the seed
$(\yy,\widetilde{B})$. 

Now one can iterate this process of mutation and obtain inductively
a set of seeds.
Thus
each seed consists of an $r$-tuple of algebraically independent 
elements of $\F$
called a 
{\it cluster} 
\index{cluster}
and of a matrix called the 
{\it exchange matrix}.
\index{exchange matrix} 
The elements of a cluster are its 
{\it  cluster variables}.
\index{cluster variable}
A seed has $r-n$ neighbours obtained by mutation in direction
$1 \le k \le r-n$.
One does not mutate the last $n$ elements of a cluster, they serve
as ''coefficients'' and belong to every cluster.
The {\it cluster algebra}
\index{cluster algebra} 
$
\cA(\widetilde{B})
$
\index{$\cA(\widetilde{B})$}
is by definition the subalgebra of $\F$ generated by the
set of all cluster variables appearing in all seeds obtained 
by iterated mutation starting with the initial seed. 

It is often convenient to define a cluster algebra using an
oriented graph, as follows.
Let $\GG$ be a quiver without loops or 2-cycles
with vertices $\{ 1,\ldots, r\}$.
We can define an $r \times r$-matrix $B(\GG) = (b_{ij})$ 
\index{$B(\GG)$}
by setting
$$
b_{ij} = 
(\text{number of arrows $j \to i$ in $\GG$}) 
-(\text{number of arrows $i \to j$ in $\GG$}).
$$
Let $B(\GG)^\circ$ 
\index{$B(\GG^\circ)$}
be the $r \times (r-n)$-matrix obtained by deleting
the last $n$ columns of $B(\GG)$.
The principal part of $B(\GG)^\circ$ is skew-symmetric, hence
this yields a cluster algebra 
$\cA(B(\GG)^\circ)$.
\index{$\cA(B(\GG)^\circ)$}

We apply this procedure to our subcategory $\CC_M$.
Let $T = T_1 \oplus \cdots \oplus T_r$ be a basic 
$\CC_M$-maximal rigid $\LL$-module
with $T_i$ indecomposable for all $i$.
Without loss of generality assume that $T_{r-n+1}, \ldots, T_r$ are
$\CC_M$-projective. 
By $\GG_T$ 
\index{$\GG_T$ (where $T$ is a maximal rigid $\LL$-module)}
we denote the quiver of the endomorphism algebra $\End_\LL(T)$.
We then define the cluster algebra
$$
\cA(\CC_M,T) := \cA(B(\GG_T)^\circ).
$$
\index{$\cA(\CC_M,T)$}\noindent
In particular, we denote by
$\cA(\CC_M)$ the cluster algebra 
$\cA(\CC_M,T_M)$ attached to the complete rigid
module $T_M$ of Section~\ref{existrigid}. 
Thus $\cA(\CC_M) := \cA(B(\GG_M^*)^\circ)$.

\subsection{Mutation of  rigid modules}\label{intro1.6}
Let $T$ be a basic 
$\CC_M$-maximal rigid $\LL$-module.
Write $B(T) = B(\GG_T) = (t_{ij})_{1 \le i,j \le r}$.
\index{$B(T)$ (where $T$ is a maximal rigid $\LL$-module)}
For $k \in [1,r-n]$ there is a short exact sequence
$$
0 \to T_k \xrightarrow{f} \bigoplus_{t_{ik} > 0} T_i^{t_{ik}}
\to T_k^* \to 0
$$
where $f$ is a minimal left 
$\add(T/T_k)$-approximation of $T_k$, 
i.e. the map $\Hom_\LL(f,T)$
is surjective, and every morphism $h$ with $hf = f$ is an 
isomorphism.
Set 
$$
\mu_{T_k}(T) = T_k^* \oplus T/T_k.
$$
We show that $\mu_{T_k}(T)$ is again a basic $\CC_M$-maximal rigid module.
In particular, $T_k^*$ is indecomposable.
We call $\mu_{T_k}(T)$ the 
{\it mutation of $T$ in direction $T_k$}.
\index{mutation of a maximal rigid $\LL$-module}

There is also a short exact sequence 
$$
0 \to T_k^* \to \bigoplus_{t_{ik}<0} T_i^{-t_{ik}} \xrightarrow{g} T_k \to 0
$$
where $g$ is now a minimal right $\add(T/T_k)$-approximation of $T_k$.

It turns out that the quivers of the endomorphism
algebras $\End_\LL(T)$ and $\End_\LL(\mu_{T_k}(T))$ are related
via Fomin and Zelevinsky's mutation rule:

\begin{Thm}\label{main4}
Let $M$ be a terminal $KQ$-module.
For a basic $\CC_M$-maximal rigid $\LL$-module $T$ as above and $k \in [1,r-n]$
we have
$$
B(\mu_{T_k}(T))^\circ = \mu_k(B(T)^\circ). 
$$
\end{Thm}

The proof can be found in Section~\ref{graphmutation}.

We conjecture that, given two  basic $\CC_M$-maximal rigid
$\LL$-modules $T$ and $T'$, there always exists a sequence of 
mutations changing $T$ into $T'$.
Using Theorem~\ref{main4}, this would imply that the cluster 
algebras $\cA(\CC_M,T)$ do not depend on $T$.
The following weaker result is given and illustrated with
examples in Section~\ref{section17}:
\begin{Thm}\label{main4'}
There is a sequence of mutations changing $T_M$ into $T_M^\vee$.
Therefore 
\[
\cA(\CC_M,T_M^\vee)=\cA(\CC_M,T_M)= \cA(\CC_M).
\]
\end{Thm}

\subsection{The dual semicanonical basis}\label{sect_dualsemi}
We recall the definition of the dual semicanonical basis and its
multiplicative properties, following \cite{Lu1, Lu2, GLSSemi1,GLSSemi2}

From now on, assume that $K = \C$.
Let $\LL_d$ be the affine
variety of nilpotent $\LL$-modules with dimension vector 
$d\in \N^n$.
For a module $X \in \LL_d$ and an $m$-tuple 
$\ii = (i_1,\ldots,i_m)$ with $1 \le i_j \le n$ for
all $j$, let
$\F_{\ii,X}$ denote the projective variety of composition series of $X$
of type $\ii$.
Thus an element in 
$\F_{\ii,X}$ 
\index{$\F_{\ii,X}$}
is a flag
$$
(0 = X_0 \subset X_1 \subset \cdots \subset X_m = X)
$$
of submodules $X_j$ of $X$
such that for all $1 \le j \le m$ the subfactor 
$X_j/X_{j-1}$ is isomorphic to the simple $\LL$-module $S_{i_j}$
associated to the vertex $i_j$ of $Q$.
Let 
$$
d_\ii\df \LL_d \to \C
$$
be the map which
sends $X \in \LL_d$ to $\chi_c(\F_{\ii,X})$, the topological
Euler characteristic of $\F_{\ii,X}$ with respect to cohomology with 
compact support.
Let $\MM_d$ be the $\C$-vector space spanned by the maps
$d_\ii$ and set
$$
\MM := \bigoplus_{d \in \N^n} \MM_d.
$$
\index{$\MM$}\noindent
Lusztig \cite{Lu1,Lu2} has introduced a ``convolution product''
$\star \df \MM \times \MM \to \MM$ such that $d_\ii\star d_\jj = d_\kk$,
where $\kk := (i_1,\ldots ,i_m,j_1,\ldots ,j_l)$ is the concatenation
of $\ii$ and $\jj$.
(The definition of $\star$ is recalled in Section~\ref{section_semican}.) 
He proved that $\MM$ equipped with this product is isomorphic to the 
enveloping algebra $U(\n)$ of the maximal nilpotent subalgebra
$\n$ of the Kac-Moody Lie algebra $\g$ associated to $Q$.

Since $U(\n)$ is a cocommutative Hopf algebra, the graded dual 
$$
U(\n)^*_{\rm gr} \equiv \MM^* := \bigoplus_{d \in \N^n} \MM_d^*
$$
is a commutative  $\C$-algebra.
For $X \in \LL_d$ we have an evaluation map $\delta_X \in \MM_d^*
\subset \MM^*$, given by 
$$
\delta_X(f) := f(X)
$$ 
\index{$\delta_X\df \MM_d \to \C$}\noindent
for $f\in\MM_d$.
In particular, $\delta_X(d_\ii) := \chi_c(\F_{\ii,X})$.
It is shown in \cite{GLSSemi1} that 
$$
\delta_X\delta_Y = \delta_{X\oplus Y}.
$$ 
In \cite{GLSSemi2} a more complicated formula 
is given, expressing
$\delta_X\delta_Y$ as a linear combination of $\delta_Z$ where $Z$ runs over 
all possible middle terms of non-split short exact sequences with end terms 
$X$ and $Y$.
The formula is especially useful when 
$\dm\Ext_\LL^1(X,Y)=1$ (see Section~\ref{subsectclusterchar}).

Let $\irr(\LL_d)$ be the set of irreducible components of $\LL_d$.
For each $Z \in \irr(\LL_d)$
there exists a dense open subset $U$ in $Z$ such that
$\delta_X = \delta_Y$ for all $X,Y \in U$.
If $X\in U$ we say that $X$ is a {\it generic} point of $Z$.
Define 
$\rho_Z := \delta_X$
\index{$\rho_Z$}
for some $X \in U$.
By \cite{Lu1,Lu2},
$$
\cS^* := \left\{ \rho_Z \mid Z \in \irr(\LL_d), d \in \N^n \right\}
$$
\index{$\cS^*$}\noindent
is a basis of $\MM^*$, the dual of Lusztig's semicanonical
basis $\cS$ of $\MM$.
In case $X$ is a rigid $\LL$-module, the orbit of $X$ in $\LL_d$ is open,
its closure is an irreducible component $Z$, and $\delta_X=\rho_Z$ belongs to
$\cS^*$.

\subsection{The cluster algebra $\cA(\CC_M)$          
as a subalgebra of $\MM^*\equiv U(\n)^*_{\rm gr}$}

For a terminal $\C Q$-module 
$M = M_1 \oplus \cdots \oplus M_r$
let $\T(\CC_M)$ 
be the graph with vertices the isomorphism
classes of basic $\CC_M$-maximal rigid $\LL$-modules and with edges given 
by mutations. 
Let $T=T_1\oplus \cdots \oplus T_r$ be a vertex of $\T(\CC_M)$, and
let $\T(\CC_M,T)$ 
denote the connected 
component of $\T(\CC_M)$ containing $T$.
Denote by 
$\RR(\CC_M,T)$ 
\index{$\RR(\CC_M,T)$}
the subalgebra of $\MM^*$
generated by the $\delta_{R_i}\ (1\le i\le r)$ where 
$R = R_1 \oplus \cdots \oplus R_r$ runs
over all vertices of $\T(\CC_M,T)$.

\begin{Thm}\label{main7}
Let $M = M_1 \oplus \cdots \oplus M_r$ be a terminal
$\C Q$-module.
Then the following hold:
\begin{itemize}

\item[(i)]
There is a unique isomorphism $\iota \df \cA(\CC_M,T)\to
\RR(\CC_M,T)$ such that 
\[
\iota(y_i)=\delta_{T_i} \qquad (1\le i\le r);
\]

\item[(ii)]
If we identify the two algebras $\cA(\CC_M,T)$ and
$\RR(\CC_M,T)$ via $\iota$, then the clusters of
$\cA(\CC_M,T)$ are identified with the $r$-tuples
$\delta(R)=(\delta_{R_1},\ldots,\delta_{R_r})$, where $R$
runs over the vertices of the graph $\T(\CC_M,T)$.
Moreover, all cluster monomials belong to 
the dual semicanonical basis 
$\cS^*$ of $\MM^*\equiv U(\n)^*_{\rm gr}$.

\end{itemize}
\end{Thm}

The proof relying on Theorem~\ref{main4} and the multiplication formula
of \cite{GLSSemi2} is given in Section~\ref{clustercattoclusteralg}.

We call $(\CC_M,T)$ a {\it categorification} of the
cluster algebra $\cA(\CC_M,T) = \cA(B(\GG_T)^\circ)$.

Write $\RR(\CC_M) := \RR(\CC_M,T_M)$.
\index{$\RR(\CC_M)$}
By Theorem~\ref{main4'}, we also have 
$\RR(\CC_M)=\RR(\CC_M,T_M^\vee)$.
Theorem~\ref{main7} shows that the cluster algebra 
$\cA(\CC_M)$ is canonically isomorphic to the subalgebra
$\RR(\CC_M)$ of $U(\n)^*_{\rm gr}$.

\subsection{Which cluster algebras did we categorify?}\label{preinjective}

The reader who is not familiar with representation theory of quivers will
ask which cluster algebras are now categorified by our approach.
We explain this in purely combinatorial terms.

As before,
let $Q = (Q_0,Q_1,s,t)$ 
\index{quiver}
be a
finite quiver without oriented cycles.
Here $Q_0 = \{1,\ldots,n\}$ 
denotes the set of vertices and $Q_1$ the set of arrows
of $Q$.
An arrow $a \in Q_1$ starts in a vertex $s(a)$ and terminates in a
vertex $t(a)$.
Let $Q^\op$ be the 
{\it opposite quiver}
of $Q$.
This is obtained from $Q$ by just reversing the direction of all
arrows.

Assume that $Q$ is not a Dynkin quiver.
Then the preinjective component $\II_Q$ 
of the Auslander-Reiten quiver of $KQ$ can be identified
with the translation quiver $\N Q^\op$ which is defined as follows.
The vertices of $\N Q^\op$ are $(i,z)$ with $1 \le i \le n$ and
$z \in \N = \{0,1,2,3,\ldots\}$.
For each arrow $a^*\df j \to i$ in $Q^\op$ there are arrows
$(a^*,z)\df (j,z) \to (i,z)$ and $(a,z)\df (i,z+1) \to (j,z)$
for all $z \in \N$.
Let $\tau(i,z) := (i,z+1)$ be the {\it translation} in $\N Q^\op$.
The vertices $(1,0),\ldots,(n,0)$ are the {\it injective vertices}
of $\N Q^\op$.

Now take any finite successor closed full subquiver $\GG$ of
$\N Q^\op$ such that $\GG$ contains all $n$ injective vertices.
Define a new quiver $\GG^*$ which is obtained from $\GG$ by adding an arrow 
from $(i,z)$ to $(i,z+1)$ whenever these vertices are both in $\GG$.
Then Theorem~\ref{main7} provides a categorification 
of the cluster algebra $\cA(B(\GG^*)^\circ)$.

For example,
if $\GG$ is the full subquiver of $\N Q^\op$ with vertices
$$
\{(1,1),(1,0),(2,1),(2,0),\ldots,(n,1),(n,0)\}, 
$$
then
$\cA(B(\GG^*)^\circ)$ is the acyclic cluster algebra associated
to the quiver $Q$.
But note that this cluster algebra comes along with $n$ coefficients
labelled by the vertices $(1,0),\ldots,(n,0)$.

If $Q$ is a Dynkin quiver, then one obtains categorifications of
cluster algebras
$\cA(B(\GG^*)^\circ)$
in a similar way. 
The only difference is that now we have to work with successor closed
full subquivers $\GG$ of the \emph{finite} Auslander-Reiten quiver of $KQ$.
We will not repeat here how to construct this quiver in this case,
but see e.g. \cite{GLSAus}.

Let us discuss another example.
If $Q$ is the quiver
{\small
$$
\xymatrix{
1 \ar@<0.5ex>[r] \ar@<-0.5ex>[r] & 2 \ar[r] & 3
}
$$
}\noindent 
then the quiver
$\N Q^\op$ looks as indicated in the following picture:
{\small
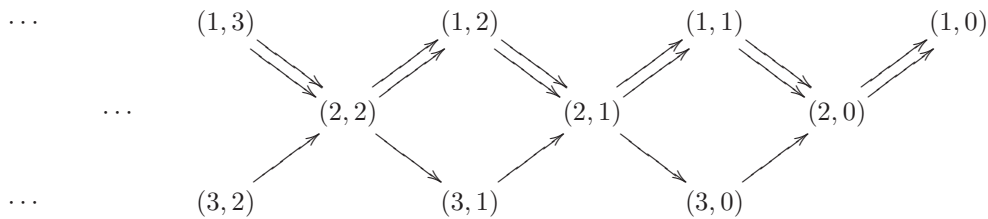
$$
\xymatrix@-0.5pc{
\cdots && (1,3) \ar@<0.5ex>[dr]\ar@<-0.5ex>[dr] && 
(1,2) \ar@<0.5ex>[dr]\ar@<-0.5ex>[dr] &&
(1,1) \ar@<0.5ex>[dr]\ar@<-0.5ex>[dr] && (1,0)
\\
&\cdots && (2,2) \ar@<0.5ex>[ur]\ar@<-0.5ex>[ur]\ar[dr] && 
(2,1) \ar@<0.5ex>[ur]\ar@<-0.5ex>[ur]\ar[dr] && 
(2,0) \ar@<0.5ex>[ur]\ar@<-0.5ex>[ur]
\\
\cdots && (3,2) \ar[ur] && 
(3,1) \ar[ur] && 
(3,0) \ar[ur]
}
$$
}\noindent 
Now let $\GG$ be the full subquiver with vertices 
$$
\{(1,2),(1,1),(1,0),(2,1),(2,0),(3,1),(3,0)\}.
$$
Clearly, $\GG$ is successor closed.
Then $\GG^*$ looks as follows:
{\small
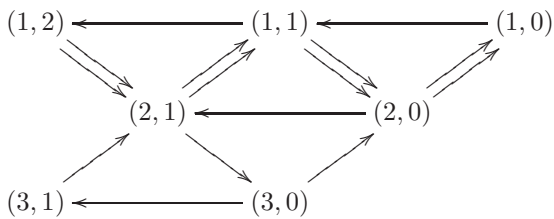
$$
\xymatrix@-0.5pc{ 
(1,2) \ar@<0.5ex>[dr]\ar@<-0.5ex>[dr] &&
(1,1) \ar[ll]\ar@<0.5ex>[dr]\ar@<-0.5ex>[dr] && (1,0) \ar[ll]
\\
&(2,1) \ar@<0.5ex>[ur]\ar@<-0.5ex>[ur]\ar[dr] && 
(2,0) \ar@<0.5ex>[ur]\ar@<-0.5ex>[ur]\ar[ll]
\\
(3,1) \ar[ur] && 
(3,0) \ar[ur]\ar[ll]
}
$$
}

\subsection{Dual PBW-bases and dual semicanonical bases}

In the spirit of Ringel's 
construction of PBW-bases for quantum groups \cite{Ri6},
we construct dual PBW-bases for our cluster algebras 
$\cA(\CC_M)$.
This yields the following result, which we prove in Section~\ref{PBWsection}.

\begin{Thm}\label{basesthm}
Let $M = M_1 \oplus \cdots \oplus M_r$ be a terminal
$\C Q$-module.
\begin{itemize}

\item[(i)]
The cluster algebra $\RR(\CC_M)$ is a polynomial ring in
$r$ variables.
More precisely, we have
\[
\RR(\CC_M) = \C[\delta_{M_1},\ldots,\delta_{M_r}] 
= \Span_\C\ebrace{\delta_X \mid X \in \CC_M},
\]
where we interpret $M_i$ as a $\LL$-module using the obvious embedding
functor;

\item[(ii)]
The set 
$
\cP_M^* := \left\{ \delta_{M'} \mid M' \in \add(M) \right\}
$ 
\index{$\cP_M^*$}
is a $\C$-basis of $\RR(\CC_M)$;

\item[(iii)]
The subset of the dual semicanonical basis
$\cS_M^* := \cS^*\cap\RR(\CC_M)$
\index{$\cS_M^*$}
is a $\C$-basis of $\RR(\CC_M)$ containing all cluster
monomials.
\end{itemize}
\end{Thm}

Let 
$\widetilde{\RR}(\CC_M)$
be the algebra
obtained from $\RR(\CC_M)$ by formally inverting
the elements $\delta_P$ for all $\CC_M$-projectives $P$.
In other words, $\widetilde{\RR}(\CC_M)$ is the cluster algebra
obtained from $\RR(\CC_M)$ by inverting the generators of
its coefficient ring. 
Similarly, let $\underline{\RR}(\CC_M)$ be
the cluster algebra obtained from $\RR(\CC_M)$ by specializing
the elements $\delta_P$ to 1.
For both cluster algebras  $\widetilde{\RR}(\CC_M)$
and $\underline{\RR}(\CC_M)$ we get a $\C$-basis which is
easily obtained from the dual semicanonical basis $\cS^*_M$
and again contains all cluster monomials, see Sections~\ref{generalbase} and
\ref{sectinvertspecial}.

\subsection{Adaptable elements of $W$}\label{wadapt}
Let $W$ be the Weyl group of $\g$, with Coxeter generators
$s_1,\ldots,s_n$.
We say that $w\in W$ is {\it $Q$-adaptable} if there exists a
reduced decomposition $w = s_{i_t} \cdots s_{i_2}s_{i_1}$ such that
$i_1$ is a sink of $Q$, and $i_{k+1}$ is a sink of
$\si_{i_k}\cdots\si_{i_2}\si_{i_1}(Q)$ for every $1 \le k \le t-1$.
Here $\si_i$ is the operation on quivers which changes the
orientation of all the arrows incident to the vertex $i$.
In this case $(i_t,\ldots ,i_1)$ is called a 
{\it $Q$-adapted reduced expression} of $w$.
We say that $w$ is {\it adaptable} if it is $Q$-adaptable for
some (acyclic) orientation $Q$ of the Dynkin diagram of $\g$.

For example if $W$ is finite, the longest element $w_0$ of $W$
is always adaptable. If $Q$ has only two vertices, for
instance if $Q$ is a (generalized) Kronecker quiver,
then every $w$ in $W$ is adaptable.
On the other hand, if $Q$ is a Dynkin quiver of type $\D_4$
with central node labelled $3$, then $w = s_3s_1s_2s_3$
is not adaptable.

It is easy to associate to a terminal $\C Q$-module $M$ a
$Q^\op$-adaptable element $w$.
Indeed, let $\Delta^+_M := \{\dimv(M_1), \ldots, \dimv(M_r)\}$
be the set of dimension vectors of the indecomposable direct
summands of $M$.
It is well known that $\Delta^+_M$ is a subset of the set
$\Delta^+$ of positive real roots of $\g$.  
In fact $\Delta^+_M= \{ \alpha \in \Delta^+ \mid
w(\alpha) < 0 \}$ for a unique $w=w(M)\in W$, and $w(M)$ is 
$Q^\op$-adaptable, see Lemma~\ref{bracketclosed2}.
Conversely any $Q^\op$-adaptable
$w$ (not contained in a proper parabolic subgroup of $W$)
comes from a unique terminal $\C Q$-module $M$. 
Moreover, if $Q'$ is a quiver obtained from $Q$ by changing
the orientation and if $w$ is also ${Q'}^\op$-adaptable,
then $\CC_M=\CC_{M'}$ where $M'$ is the terminal $\C Q'$-module attached 
to $w$ (see Section~\ref{subcatCw}). 
This implies that $\CC_M$ depends only on the adaptable
element $w$ of $W$, and we sometimes write $\CC_M=\CC_w$.

\subsection{Unipotent subgroups and cells}    
Let $M$ be a terminal $\C Q$-module, and let $w=w(M)$ be the associated
Weyl group element.
Let
\[
\n_M=\n(w) = \bigoplus_{\alpha\in\Delta^+_M} \n_\alpha
\]
be the corresponding sum of root subspaces of $\n$.
This is a finite-dimensional nilpotent Lie algebra.
Let $N_M=N(w)$ be the corresponding finite-dimensional unipotent group.

If $G$ is the Kac-Moody group attached to $\g$ as in \cite[Chapter 6]{Ku},
which comes with a pair of subgroups $N$ and $N_-$ (denoted by
$\U$ and $\U_-$ in \cite{Ku}), then
\[
N(w) = N \cap (w^{-1}N_-w).
\]
We also define the unipotent cell
\[
N^w =  N\cap(B_-wB_-)
\]
where $B_-$ is the standard negative Borel subgroup of $G$. 

The following theorem, proved in Section~\ref{coordinaterings}, 
is one of our main
results.
It shows that $\CC_M$ can be regarded as a categorification
of both $N(w)$ and $N^w$.

\begin{Thm}\label{wident}
The algebras of regular functions on $N(w)$ and $N^w$
have a cluster algebra structure.
More precisely, we have natural isomorphisms
\begin{align*}
&\C[N(w)] \cong \RR(\CC_M),\\
&\C[N^w] \cong \widetilde{\RR}(\CC_M).
\end{align*}
\end{Thm}
Note that in the Dynkin case the cluster algebra structure 
on $\C[N^w]$ was already known by work of Berenstein, Fomin and
Zelevinsky \cite{BFZ}, but our proof is different and yields
the additional result that the upper cluster algebra
is in fact a cluster algebra.

\subsection{An example}
We are going to illustrate some of the previous results on an example.
Let $Q$ be the quiver
{\small
$$
\xymatrix@-0.8pc{
1 \ar[dr]\\&2\ar[dr]&&4\ar[dl]\\&&3
}
$$
}
\begin{figure}
$$
\xymatrix@-1.0pc{
&&\blue{\fbox{$\bsm1&&&\\&1&&0\\&&1&\esm$}} \ar[dr] && 
\blue{\fbox{$\bsm0&&&\\&0&&1\\&&0&\esm$}}\\
&\blue{\fbox{$\bsm0&&&\\&1&&0\\&&1&\esm$}} \ar[ur]\ar[dr] && 
\blue{\fbox{$\bsm1&&&\\&1&&1\\&&1&\esm$}} \ar[ur]\ar[dr]\\
{\bsm0&&&\\&0&&0\\&&1&\esm} \ar[ur]\ar[dr] && 
\blue{\fbox{$\bsm0&&&\\&1&&1\\&&1&\esm$}} \ar[ur]\ar[dr] &&
\blue{\fbox{$\bsm1&&&\\&1&&0\\&&0&\esm$}} \ar[dr]\\
&{\bsm0&&&\\&0&&1\\&&1&\esm} \ar[ur] &&
\blue{\fbox{$\bsm0&&&\\&1&&0\\&&0&\esm$}} \ar[ur] && 
\blue{\fbox{$\bsm1&&&\\&0&&0\\&&0&\esm$}}
}
$$
\caption{The Auslander-Reiten quiver of a path algebra of type $\A_4$}
\label{fig2}
\end{figure}
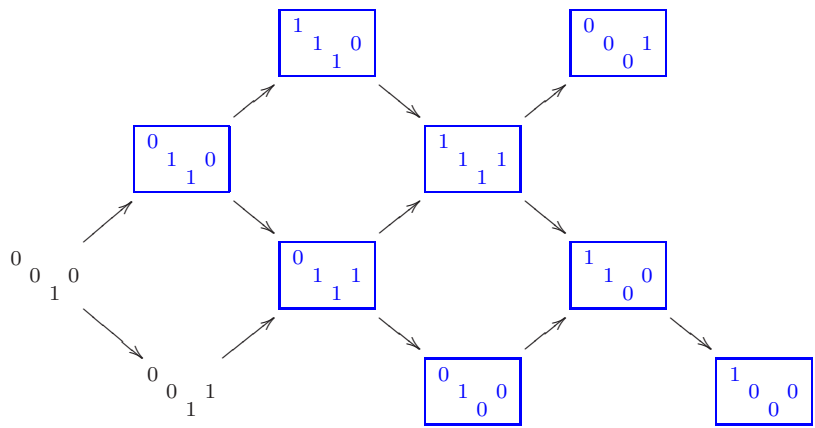

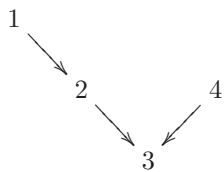
Figure~\ref{fig2} shows the Auslander-Reiten quiver of $\C Q$.
The indecomposable direct summands of a terminal $\C Q$-module $M$
are marked in blue colour.
In Figure~\ref{fig3} we show the Auslander-Reiten quiver of
the preprojective algebra $\LL$ of type $\A_4$.
As in Section~\ref{subsectexampleA3}, we display the graded dimension vectors of
the indecomposable $\LL$-modules.
One has to identify the objects on the two dotted vertical lines.
The indecomposable $\CC_M$-projective-injective modules  
are marked in red colour, all other indecomposable modules in $\CC_M$
are marked in blue.
Observe that $\CC_M$ contains 18 indecomposable modules,
and 4 of these are $\CC_M$-projective-injective.
The stable category ${\stCC}_M$ is triangle equivalent to the
cluster category~$\CC_Q$.

The maximal rigid module $T_M$ has 8 indecomposable direct summands,
namely, the 4 indecomposable $\CC_M$-projective-injective modules
$$
I_4 = \bsm1&&0\\&1&&0\\0&&1\\&0&&1\esm\qquad
I_3 = \bsm&1&&0\\1&&1\\&1&&1\\0&&1\esm\qquad
I_2 = \bsm&1&&1\\1&&1\\&1&&0\esm\qquad
I_1 = \bsm0&&0\\&1&&0\\1&&0\esm\qquad
$$
and the modules
$$
T_1=\bsm1&&0\\&1&&0\\0&&1\esm\qquad
T_2=\bsm0&&0\\&1&&0\\0&&1\esm\qquad
T_3=\bsm0&&0\\&1&&1\\0&&1\esm\qquad
T_4=\bsm0&&0\\&1&&0\\0&&0\esm.
$$
Similarly, $T_M^\vee$ has 4 non-projective indecomposable direct summands,
namely, 
$$
T_1^\vee = \bsm0&&0\\&0&&1\\0&&0\esm\qquad
T_2^\vee = \bsm1&&0\\&1&&1\\0&&1\esm\qquad
T_3^\vee = \bsm1&&0\\&1&&0\\0&&0\esm\qquad
T_4^\vee = \bsm1&&0\\&0&&0\\0&&0\esm.
$$
Here, the group $N$ can be taken to be the group of upper
unitriangular $5\times 5$ matrices with complex coefficients.
Given two subsets $I$ and $J$ of $\{1,2,\ldots,5\}$
with $|I|=|J|$, we denote 
by $D_{IJ}\in\C[N]$ the regular function mapping an element
$x\in N$ to its minor $D_{IJ}(x)$ with row subset $I$ and 
column subset $J$.  
Every $\LL$-module $X$ in $\CC_M$ gives rise to a linear form
$\delta_X \in U(\n)_{\rm gr}^*$ and by means of the isomorphism 
$U(\n)_{\rm gr}^* \cong \C[N]$ to a regular function $\vph_X$.
For example,
$$
\vph_{I_4} = D_{1234,2345},\quad
\vph_{I_3} = D_{123,345},\quad
\vph_{I_2} = D_{12,35},\quad
\vph_{I_1} = D_{1,3},\quad
$$
$$
\vph_{T_1} = D_{123,234},\quad
\vph_{T_2} = D_{123,134},\quad
\vph_{T_3} = D_{123,135},\quad
\vph_{T_4} = D_{12,13},\quad
$$
$$
\vph_{T^\vee_1} = D_{1234,1235},\quad
\vph_{T^\vee_2} = D_{123,235},\quad
\vph_{T^\vee_3} = D_{12,23},\quad
\vph_{T^\vee_4} = D_{1,2}.\quad
$$
The Weyl group element attached to $\CC_M$ is 
$w=s_3s_4s_2s_1s_3s_4s_2s_1$.
The corresponding unipotent subgroup $N(w)$ consists of all
$5\times 5$ matrices of the form
$$
\left[
\begin{matrix}
1&u_1&u_2&u_7&u_4\\
0&1&u_5&u_8&u_6\\
0&0&1&0&0\\
0&0&0&1&u_3\\
0&0&0&0&1
\end{matrix}
\right],
\qquad
(u_1,\ldots,u_8\in\C).
$$
The unipotent cell $N^w$ is a locally closed subset of $N$
defined by the following equations and inequalities:
\begin{multline*}
N^w = \{ x \in N \mid
D_{1,4}(x) = D_{1,5}(x) = D_{12,45}(x) = 0,\\
D_{1234,2345}(x) \not= 0,\
D_{123,345}(x) \not= 0,\
D_{12,35}(x) \not= 0,\
D_{1,3}(x) \not= 0 \}
\end{multline*}
Note that the 4 inequalities are given by the non-vanishing
of the 4 regular functions 
$\vph_{I_k}\ (k = 1,\ldots,4)$ 
attached to the indecomposable $\CC_M$-projective-injective modules.

Our results show that the polynomial algebra $\C[N(w)]$
has a cluster algebra structure, of which
$(\vph_{T_1},\vph_{T_2},\vph_{T_3}, \vph_{T_4},
\vph_{I_1},\vph_{I_2},\vph_{I_3},\vph_{I_4})$
and
$(\vph_{T^\vee_1},\vph_{T^\vee_2},\vph_{T^\vee_3},\vph_{T^\vee_4},
\vph_{I_1},\vph_{I_2},\vph_{I_3},\vph_{I_4})$
are two distinguished clusters.
Its coefficient ring is the polynomial ring in the four variables
$(\vph_{I_1},\vph_{I_2},\vph_{I_3},\vph_{I_4})$.
The cluster mutations of this algebra come from mutations of
maximal rigid modules in $\CC_M$.
Moreover, if we replace the coefficient ring by the ring of
Laurent polynomials in the four variables
$(\vph_{I_1},\vph_{I_2},\vph_{I_3},\vph_{I_4})$,
we obtain the coordinate ring $\C[N^w]$.

{
\begin{landscape}
\begin{figure}
\[
\xymatrix@-1.3pc{
\ar@{.}[d]&&&{\bsm0&&0\\&0&&1\\0&&1\\&1&&0\\1&&0\esm}\ar[dr] &&&&&& 
\red{\fbox{$\bsm1&&0\\&1&&0\\0&&1\\&0&&1\\0&&0\esm$}}\ar[dr] &&&\ar@{.}[d]\\
\blue{\fbox{$\bsm0&&0\\&0&&0\\0&&0\\&0&&1\\0&&0\esm$}}\ar[dr]\ar@{.}[dd] && 
{\bsm0&&0\\&0&&0\\0&&1\\&1&&0\\1&&0\esm}\ar[ur]\ar[dr] && 
{\bsm0&&0\\&0&&1\\0&&1\\&1&&0\\0&&0\esm}\ar[dr] && 
\blue{\fbox{$\bsm0&&0\\&0&&0\\1&&0\\&0&&0\\0&&0\esm$}}\ar[dr] &&
\blue{\fbox{$\bsm0&&0\\&1&&0\\0&&1\\&0&&1\\0&&0\esm$}}\ar[dr]\ar[ur] && 
\blue{\fbox{$\bsm1&&0\\&1&&0\\0&&1\\&0&&0\\0&&0\esm$}}\ar[dr] &&
\blue{\fbox{$\bsm0&&0\\&0&&1\\0&&0\\&0&&0\\0&&0\esm$}} \ar@{.}[dd]\\
&{\bsm0&&0\\&0&&0\\0&&1\\&1&&1\\1&&0\esm}\ar[ur]\ar[dr]&& 
{\bsm0&&0\\&0&&0\\0&&1\\&1&&0\\0&&0\esm}\ar[ur]\ar[dr] && 
{\bsm0&&0\\&0&&1\\1&&1\\&1&&0\\0&&0\esm}\ar[ur]\ar[dr] && 
\blue{\fbox{$\bsm0&&0\\&1&&0\\1&&1\\&0&&1\\0&&0\esm$}}\ar[ur]\ar[dr] &&
\blue{\fbox{$\bsm0&&0\\&1&&0\\0&&1\\&0&&0\\0&&0\esm$}}\ar[ur]\ar[dr]&& 
\blue{\fbox{$\bsm1&&0\\&1&&1\\0&&1\\&0&&0\\0&&0\esm$}}\ar[ur]\ar[dr]&\\
{\bsm0&&0\\&0&&0\\1&&1\\&2&&1\\1&&1\esm} \ar[dr] \ar[ur]\ar@{.}[d] && 
{\bsm0&&0\\&0&&0\\0&&1\\&1&&1\\0&&0\esm} \ar[dr] \ar[ur] &&
{\bsm0&&0\\&0&&0\\1&&1\\&1&&0\\0&&0\esm} \ar[dr] \ar[ur] && 
{\bsm0&&0\\&1&&1\\1&&2\\&1&&1\\0&&0\esm} \ar[dr] \ar[ur] &&
\blue{\fbox{$\bsm0&&0\\&1&&0\\1&&1\\&0&&0\\0&&0\esm$}} \ar[dr] \ar[ur] && 
\blue{\fbox{$\bsm0&&0\\&1&&1\\0&&1\\&0&&0\\0&&0\esm$}} \ar[dr] \ar[ur] &&
{\bsm1&&1\\&2&&1\\1&&1\\&0&&0\\0&&0\esm}\ar@{.}[d]\\
\blue{\fbox{$\bsm0&&0\\&0&&0\\0&&0\\&1&&0\\0&&0\esm$}} \ar[r]\ar@{.}[d] & 
{\bsm0&&0\\&0&&0\\1&&1\\&2&&1\\0&&1\esm}\ar[r]\ar[dr]\ar[ur]&
{\bsm0&&0\\&0&&0\\1&&1\\&1&&1\\0&&1\esm} \ar[r]\ar[ddr] & 
{\bsm0&&0\\&0&&0\\1&&1\\&1&&1\\0&&0\esm}\ar[r]\ar[dr]\ar[ur]&
\blue{\fbox{$\bsm0&&0\\&1&&0\\1&&1\\&1&&1\\0&&0\esm$}}\ar[r] & 
{\bsm0&&0\\&1&&0\\1&&2\\&1&&1\\0&&0\esm}\ar[r]\ar[dr]\ar[ur] &
{\bsm0&&0\\&0&&0\\0&&1\\&0&&0\\0&&0\esm} \ar[r] & 
{\bsm0&&0\\&1&&1\\1&&2\\&1&&0\\0&&0\esm}\ar[r]\ar[dr]\ar[ur]&
\red{\fbox{$\bsm0&&0\\&1&&1\\1&&1\\&1&&0\\0&&0\esm$}}\ar[r]\ar[ddr] & 
\blue{\fbox{$\bsm0&&0\\&1&&1\\1&&1\\&0&&0\\0&&0\esm$}}\ar[r]\ar[dr]\ar[ur] &
{\bsm0&&1\\&1&&1\\1&&1\\&0&&0\\0&&0\esm} \ar[r] & 
{\bsm0&&1\\&2&&1\\1&&1\\&0&&0\\0&&0\esm}\ar[r]\ar[dr]\ar[ur]&
\blue{\fbox{$\bsm0&&0\\&1&&0\\0&&0\\&0&&0\\0&&0\esm$}}\ar@{.}[d]\\
\ar@{.}[d]{\bsm0&&0\\&0&&0\\0&&1\\&1&&1\\0&&1\esm} \ar[ur] && 
\blue{\fbox{$\bsm0&&0\\&0&&0\\1&&0\\&1&&0\\0&&0\esm$}} \ar[ur] &&
{\bsm0&&0\\&0&&0\\0&&1\\&0&&1\\0&&0\esm} \ar[ur] && 
\blue{\fbox{$\bsm0&&0\\&1&&0\\1&&1\\&1&&0\\0&&0\esm$}} \ar[ur]&&
{\bsm0&&0\\&0&&1\\0&&1\\&0&&0\\0&&0\esm}\ar[ur] && 
\red{\fbox{$\bsm0&&0\\&1&&0\\1&&0\\&0&&0\\0&&0\esm$}}\ar[ur] &&
{\bsm0&&1\\&1&&1\\0&&1\\&0&&0\\0&&0\esm}\ar@{.}[d]\\
&&&\red{\fbox{$\bsm0&&0\\&1&&0\\1&&1\\&1&&1\\0&&1\esm$}} \ar[uur] &&&&&&
{\bsm0&&1\\&1&&1\\1&&1\\&1&&0\\0&&0\esm} \ar[uur] &&&
}
\]
\caption{
A category $\CC_M \subset \md(\LL)$ with $\stCC_M$ triangle
equivalent to a cluster category of type $\A_4$
}\label{fig3}
\end{figure}
\end{landscape}
}

\newpage


{\Large\part{The category $\CC_M$}
\label{part2}}



{\Large\section{Representations of quivers and preprojective algebras}\label{quivers}}


\subsection{Nilpotent varieties}
A $\LL$-module $M$ is called {\it nilpotent} 
\index{nilpotent $\LL$-module}
if a composition
series of $M$ contains only the simple modules $S_1,\ldots,S_n$ 
associated to the vertices of $Q$.
Let $\nil(\LL)$ 
\index{$\nil(\LL)$}
be the abelian category of finite-dimensional
nilpotent $\LL$-modules.

Let $d = (d_1, \ldots, d_n) \in \N^n$.
By 
$$
\rep(Q,d) = \prod_{a \in Q_1} 
\Hom_{K}(K^{d_{s(a)}},
K^{d_{t(a)}}) 
$$
\index{$\rep(Q,d)$}\noindent
we denote the affine space of representations of $Q$ with dimension
vector $d$.
Furthermore, let $\md(\LL,d)$ 
be the affine variety of elements
$$
(f_a,f_{a^*})_{a \in Q_1} \in
\prod_{a \in Q_1} \left( \Hom_K(K^{d_{s(a)}},
K^{d_{t(a)}}) \times \Hom_K(K^{d_{t(a)}},
K^{d_{s(a)}}) \right) 
$$
such that the following holds:
\begin{itemize}

\item[(i)]
For all $i \in Q_0$ we have
$$ 
\sum_{a \in Q_1: s(a) = i} f_{a^*}f_a =
\sum_{a \in Q_1: t(a) = i} f_a f_{a^*}.
$$

\end{itemize}
By $\LL_d = \nil(\LL,d)$ 
\index{$\LL_d$}
\index{$\nil(\LL,d)$}
we denote the variety
of all 
$(f_a,f_{a^*})_{a \in Q_1} \in \md(\LL,d)$ such that
the following condition holds:
\begin{itemize}

\item[(ii)]
There exists some $N$ such that for each path $a_1a_2 \cdots a_N$ of
length $N$ in the double quiver $\overline{Q}$ of $Q$
we have $f_{a_1}f_{a_2} \cdots f_{a_N} = 0$.

\end{itemize}
If $Q$ is a Dynkin quiver, then (ii) follows already from condition
(i).
One can regard (ii) as a nilpotency condition, which explains why
the varieties $\LL_d$ are often called {\it nilpotent varieties}.
\index{nilpotent variety}
Note that $\rep(Q,d)$ can be considered as a subvariety of $\LL_d$.
In fact $\rep(Q,d)$ forms an irreducible component of $\LL_d$.
Lusztig \cite[Section 12]{Lu1} proved that 
all irreducible components of $\LL_d$ have the same dimension, namely
$$
\dm \rep(Q,d) = \sum_{a \in Q_1} d_{s(a)}d_{t(a)}.
$$
One can interpret $\LL_d$ as the variety of nilpotent
$\LL$-modules with dimension vector $d$.
The group 
$$
\GL_d = \prod_{i=1}^n \GL(d_i,K)
$$ 
acts on $\md(\LL,d)$, 
$\LL_d$ and $\rep(Q,d)$ by conjugation.
Namely, for 
$g = (g_1,\ldots,g_n) \in \GL_d$ and
$x = (f_a,f_{a^*})_{a \in Q_1} \in \md(\LL,d)$ define
$$
g \cdot x := (g_{t(a)}f_a g_{s(a)}^{-1},
g_{s(a)}f_{a^*} g_{t(a)}^{-1})_{a \in Q_1}.
$$
The action on $\LL_d$ and $\rep(Q,d)$ is obtained via restriction.
The isomorphism classes of $\LL$-modules
in $\md(\LL,d)$ and $\LL_d$, and $KQ$-modules in $\rep(Q,d)$,
respectively, correspond to the orbits of these actions.
For a module $M$ with dimension vector $d$ over $\LL$ or over $KQ$ 
let $\orb_M$ 
be its $\GL_d$-orbit
in $\md(\LL,d)$, $\LL_d$ or $\rep(Q,d)$, respectively.

\subsection{Dimension formulas for nilpotent varieties}
The restriction functor 
$$
\pi_Q\df \md(\LL) \to \md(KQ)
$$ 
\index{$\pi_Q$}\noindent
induces a surjective morphism of varieties
$
\pi_{Q,d}\df \md(\LL,d) \to \rep(Q,d).
$
\index{$\pi_{Q,d}$}\noindent
There is a bilinear form $\bil{-,-}\df \Z^n \times \Z^n \to \Z$
associated to $Q$ defined by
$$
\bil{d,e} = \sum_{i \in Q_0} d_ie_i -
\sum_{a \in Q_1} d_{s(a)}e_{t(a)}.
$$
The dimension vector of a $KQ$-module $M$ is denoted by $\dimv(M)$.
For $KQ$-modules $M$ and $N$ set
$$
\bil{M,N} := 
\dm \Hom_{KQ}(M,N) - \dm \Ext_{KQ}^1(M,N).
$$
It is known that $\bil{M,N} = \bil{\dimv(M),\dimv(N)}$.
Let $(-,-)\df \Z^n \times \Z^n \to \Z$
be the symmetrization of the bilinear form $\bil{-,-}$,
i.e. $(d,e) := \bil{d,e} + \bil{e,d}$.
Thus for $\LL$-modules $X$ and $Y$ we have
$$
(\dimv(X),\dimv(Y)) = \bil{\pi_Q(X),\pi_Q(Y)} + \bil{\pi_Q(Y),\pi_Q(X)}.
$$

\begin{Lem}[{\cite[Lemma 1]{CB}}]\label{extandhom}
For any $\LL$-modules $X$ and $Y$ we have
$$
\dm \Ext_\LL^1(X,Y) = \dm \Hom_\LL(X,Y) + \dm \Hom_\LL(Y,X) - 
(\dimv(X),\dimv(Y)).
$$
\end{Lem}

\begin{Cor}\label{extandhomcor}
$\dm \Ext_\LL^1(X,X)$ is even, and 
$\dm \Ext_\LL^1(X,Y) = \dm \Ext_\LL^1(Y,X)$.
\end{Cor}

For a $\GL_d$-orbit $\orb$ in $\LL_d$ let 
${\rm codim}\, \orb = \dm \LL_d - \dm \orb$ be its codimension.

\begin{Lem}\label{corlu1}
For any nilpotent $\LL$-module $M$ we have
$\dm \Ext_\LL^1(M,M) = 2\, {\rm codim}\, \orb_M$.
\end{Lem}

\begin{proof}
Set $d = \dimv(M)$.
By Lemma~\ref{extandhom} we have
$$
\dm \Ext_\LL^1(M,M) = 2\, \dm \End_\LL(M) - (d,d).
$$
Furthermore, 
$
\dm \orb_M = \dm \GL_d - \dm \End_\LL(M)$.
Thus
\[
{\rm codim}\, \orb_M = \dm \LL_d - \dm \orb_M =
\sum_{\alpha \in Q_1} d_{s(\alpha)}d_{t(\alpha)} - 
\sum_{i=1}^n d_i^2
+ \dm \End_\LL(M).
\]
Combining these equations yields the result.
\end{proof}

\begin{Cor}\label{corlu2}
For a nilpotent $\LL$-module $M$ with dimension vector
$d$ the following are equivalent:
\begin{itemize}

\item
The closure $\overline{\orb_M}$ of $\orb_M$ is an irreducible
component of $\LL_d$;

\item
The orbit $\orb_M$ is open in $\LL_d$;

\item
$\Ext_\LL^1(M,M) = 0$.

\end{itemize}
\end{Cor}

\begin{Lem}[{\cite[Theorem 3.3]{CB2}}]\label{CBfibre}
For each $N \in \rep(Q,d)$ the fibre $\pi_{Q,d}^{-1}(N)$ is isomorphic
to ${\rm D}\Ext_{KQ}^1(N,N)$.
\end{Lem}

\begin{Cor}[{\cite[Lemma 3.4]{CB2}}]\label{CBfibrecor1}
For $N \in \rep(Q,d)$ we have
\begin{align*}
\dm \pi_{Q,d}^{-1}(\orb_N) &= \dm \orb_N + \dm \Ext_{KQ}^1(N,N)\\
&= \sum_{i \in Q_0} d_i^2 - \dm \End_{KQ}(N) + \dm \Ext_{KQ}^1(N,N)\\
&= \sum_{i \in Q_0} d_i^2 - \bil{d,d}
= \sum_{a \in Q_1} d_{s(a)}d_{t(a)} = \dm \rep(Q,d).
\end{align*}
Furthermore,
$
\pi_{Q,d}^{-1}(\orb_N)
$
is locally closed and irreducible in $\md(\LL,d)$.
\end{Cor}

\begin{Cor}\label{CBfibrecor2}
Let $N \in \rep(Q,d)$ such that
$\pi_{Q,d}^{-1}(N) \subseteq \LL_d$.
Then 
$$
\overline{\pi_{Q,d}^{-1}(\orb_N)}
$$ 
is an irreducible component of $\LL_d$.
\end{Cor}

\subsection{Terminal $KQ$-modules and irreducible components}
Following Ringel \cite{Ri2} 
we define a $K$-category $\CC(1,\tau)$ as follows:
The objects are of the form $(X,f)$ where $X$ 
is in $\md(KQ)$ and $f\df X \to \tau(X)$ is a $KQ$-module homomorphism.
Here $\tau = \tau_Q$ denotes the Auslander-Reiten translation in $\md(KQ)$.
The morphisms from $(X,f)$ to $(Y,g)$ are just the $KQ$-module 
homomorphisms $h\df X \to Y$
such that the diagram
$$
\xymatrix{
X \ar[d]^f \ar[r]^h & Y \ar[d]^g\\
\tau(X) \ar[r]^{\tau(h)} & \tau(Y) 
}
$$
commutes.
Then the categories $\md(\LL)$ and $\CC(1,\tau)$ are
isomorphic \cite[Theorem B]{Ri2}.
More precisely, there exists an isomorphism of categories
$$
\Psi\df \md(\LL) \to \CC(1,\tau)
$$ 
such that $\Psi(X) = (Y,f)$ implies $\pi_Q(X) = Y$ for all $X \in \md(\LL)$.

\begin{Lem}
Let $M$ be a terminal $KQ$-module.
Assume $N \in \add(M)$.
Then 
$$
\overline{\pi_{Q,d}^{-1}(\orb_N)} 
$$
is an irreducible component of $\LL_d$.
In particular, $\CC_M \subseteq \nil(\LL)$.
\end{Lem}

\begin{proof}
Since $M$ is a terminal $KQ$-module and $N \in \add(M)$, 
we know that $\pi_{Q,d}^{-1}(N)$ is contained in $\LL_d$.
Indeed, by \cite{Ri2} for every $KQ$-module $X$ the intersection 
$$
\pi_{Q,d}^{-1}(X) \cap \LL_d
$$ 
can be identified with the space of 
$KQ$-module homomorphisms 
$f\df X \to \tau(X)$ such that the composition
$$
X \xrightarrow{f} 
\tau(X) \xrightarrow{\tau(f)} 
\tau^2(X) \xrightarrow{\tau^2(f)} \cdots
\xrightarrow{\tau^{m-1}(f)}
\tau^m(X)
$$ 
is zero for some $m \ge 1$.
Since $N$ is preinjective, such an $m$ always exists, namely we
have $\Hom_{KQ}(N,\tau^m(N)) = 0$ for $m$ large enough.
Then use Corollary~\ref{CBfibrecor2}.
\end{proof}


{\Large\section{Selfinjective torsion classes in $\nil(\LL)$}\label{section4}}


\subsection{Tilting modules and torsion classes}\label{perp}
We need to recall some facts on torsion theories and
tilting modules.
Let $A$ be a $K$-algebra, and let $\U$ be a subcategory of $\md(A)$.

A module $C$ in $\U$ is a 
{\it generator} 
\index{generator of a subcategory}
(resp. {\it cogenerator}) 
\index{cogenerator of a subcategory}
of
$\U$ if for each $X \in \U$ there exists some $t \ge 1$
and an epimorphism $C^t \to X$ (resp. a monomorphism $X \to C^t$).

Let
\begin{align*}
\U^\perp &= \{ X \in \md(A) \mid \Hom_A(U,X) = 0 \text{ for all } U \in \U\},\\
{^\perp\U} &= \{ X \in \md(A) \mid \Hom_A(X,U) = 0 \text{ for all } U \in \U\}.
\end{align*}
The following lemma is well known:

\begin{Lem}[{\cite[Section 1.1]{Bo2}}]\label{torsionlemma}
For a subcategory $\T$ of $\md(A)$
the following are equivalent:
\begin{itemize}

\item[(i)]
$\T = {^\perp(\T^\perp)}$;

\item[(ii)]
$\T$ is closed under extensions and factor modules.

\end{itemize}
\end{Lem}

A pair $(\F,\T)$ of subcategories of $\md(A)$ is called a {\it torsion theory}
\index{torsion theory}
in $\md(A)$ if $\T^\perp = \F$ and $\T = {^\perp\F}$.
The modules in $\F$ are called 
{\it torsion-free modules} 
and the ones in
$\T$ are called {\it torsion modules}.
A subcategory $\T$ of $\md(A)$ is a {\it torsion class}
\index{torsion class}
if it satisfies one of the equivalent conditions in Lemma~\ref{torsionlemma}.

An $A$-module $T$ 
is a {\it tilting module} 
\index{tilting module}
if there exists some $d \ge 1$ such that
the following three conditions hold:
\begin{itemize}

\item[(1)] $\pdim(T) \le d$;

\item[(2)]
$\Ext_A^i(T,T) = 0$ for all $i \ge 1$;

\item[(3)] There exists a short exact sequence
$
0 \to {_A}A \to T_0 \to T_1 \to \cdots \to T_d \to 0
$
with $T_i \in \add(T)$ for all $i \ge 0$.

\end{itemize}
Such a module $T$ is a 
{\it classical tilting module}
\index{classical tilting module}
if we can take $d=1$.
Note that over path algebras $KQ$ every tilting module is a classical
tilting module.

Any classical tilting module $T$ over an algebra
$A$ yields a torsion theory
$(\F,\T)$ where
\begin{align*}
\F &= \{ N \in \md(A) \mid \Hom_A(T,N) = 0\},\\
\T &= \{ N \in \md(A) \mid \Ext_A^1(T,N) = 0\}.
\end{align*}
It is a well known result from tilting theory that
$
\T = \Gen(T).
$
As a reference for tilting theory we recommend
\cite{ASS,Bo1,HR,Ri1}.

Let $A$ and $B$ be finite-dimensional $K$-algebras.
The algebras $A$ and $B$ are {\it derived equivalent} if their
derived categories
$D^b(\md(A))$ and $D^b(\md(B))$ 
are equivalent as triangulated categories, see for example 
\cite[Section 0]{H1}.
We need the following results:

\begin{Thm}[{\cite[Section 1.7]{H1}}]\label{Happel}
If $T$ is a tilting module over $A$, then
$A$ and $\End_A(T)^\op$ are derived equivalent.
\end{Thm}

\begin{Thm}[{\cite[Section 1.4]{H1}}]
If $A$ and $B$ are derived equivalent,
then $\gldim(A) < \infty$ if and only if $\gldim(B) < \infty$.
\end{Thm}

\subsection{Approximations of modules}
Let $A$ be a $K$-algebra, and let $M$ be an $A$-module.
A homomorphism $f\df X \to M'$ in $\md(A)$ is 
a 
{\it left $\add(M)$-approximation} 
\index{left approximation}
of $X$ if $M' \in \add(M)$ 
and the induced map
$$
\Hom_A(f,M)\df \Hom_A(M',M) \to \Hom_A(X,M) 
$$
is surjective.
A morphism $f\df V \to W$ is called
{\it left minimal} 
\index{left minimal morphism}
if every morphism $g\df W \to W$ with $gf = f$ is an 
isomorphism.
Dually, one defines right $\add(M)$-approximations and
right minimal morphisms.
\index{right approximation}
\index{right minimal morphism}
Some well known basic properties of approximations can be found in 
\cite[Section 3.1]{GLSRigid}.

\subsection{Frobenius subcategories}
Let $\CC$ be a subcategory of $\md(\LL)$ 
which is closed under extensions.
A $\LL$-module $C$ in $\CC$ is called $\CC$-{\it projective}  
\index{$\CC$-projective}
(resp. $\CC$-{\it injective})
\index{$\CC$-injective}
if $\Ext_\LL^1(C,X) = 0$ (resp. $\Ext_\LL^1(X,C) = 0$) for all $X \in \CC$.
If $C$ is $\CC$-projective and $\CC$-injective, then $C$ is also called
$\CC$-{\it projective-injective}.
\index{$\CC$-projective-injective}

We say that $\CC$ has {\it enough projectives} 
(resp. {\it enough injectives})
if for each $X \in \CC$ there exists a short exact sequence
$
0 \to Y \to C \to X \to 0
$
(resp. $0 \to X \to C \to Y \to 0$)
where $C$ is $\CC$-projective (resp. $\CC$-injective)
and $Y \in \CC$.

\begin{Lem}
For a $\LL$-module $C$ in $\CC$ the following are equivalent:
\begin{itemize}

\item
$C$ is $\CC$-projective;

\item
$C$ is $\CC$-injective.

\end{itemize}
\end{Lem}

\begin{proof}
This follows immediately from Corollary~\ref{extandhomcor}.
\end{proof}

If $\CC$ has enough projectives and enough injectives, then
$\CC$ is called a {\it Frobenius subcategory} 
\index{Frobenius subcategory of $\md(\LL)$}
of $\md(\LL)$.
In particular, $\CC$ is a Frobenius category in the sense of
Happel \cite{H2}.

\subsection{Cluster torsion classes}\label{defclustertorsionclass}
Let $\CC$ be a subcategory of $\nil(\LL)$.
We call $\CC$ a {\it selfinjective torsion class} 
\index{selfinjective torsion class}
if the following hold:
\begin{itemize}

\item[(i)]
$\CC$ is closed under extensions;

\item[(ii)]
$\CC$ is closed under factor modules;

\item[(iii)]
There exists a generator-cogenerator $I_\CC$ of $\CC$ which is
$\CC$-projective-injective;

\item[(iv)]
$\rk(\CC) < \infty$.

\end{itemize}
It follows from the definitions and Lemma~\ref{torsionlemma} that
for a selfinjective torsion class $\CC$ of $\nil(\LL)$ we have
$
{^\perp(\CC^\perp)} = \CC.
$
In particular, $(\CC^\perp,\CC)$ is a torsion theory in $\md(\LL)$.

We will show that each selfinjective torsion class $\CC$ of $\nil(\LL)$ can
be interpreted as a categorification of a certain cluster algebra,
provided the following holds:
\begin{itemize}

\item[($\star$)]
There exists a $\CC$-complete rigid module $T_\CC$ such that
the quiver $\GG_{T_\CC}$ of $\End_\LL(T_\CC)$ has no loops.

\end{itemize}
A selfinjective torsion class satisfying $(\star)$ will be called
a {\it cluster torsion class}.
\index{cluster torsion class}

We will prove in Proposition~\ref{CMtorsion} that for every terminal 
$KQ$-module
$M$, the subcategory $\CC_M$ is a cluster torsion class.
For simplicity, in the introduction,  
Theorems~\ref{main4} and \ref{main7} were only stated for
subcategories of the form $\CC_M$.
But the proofs (Sections~\ref{graphmutation} 
and \ref{clustercattoclusteralg}) 
are carried out more generally for 
cluster torsion classes.

\subsection{$\CC_M$ is a torsion class}
Let 
$M = M_1 \oplus \cdots \oplus M_r$ 
be a terminal $KQ$-module as defined in Section~\ref{terminal}.
Set
$$
T := \bigoplus_{i=1}^n \tau^{t_i(M)}(I_i).
$$
Note that the $KQ$-module $T$ is a basic tilting module
with
$\Gen(T) = \add(M)$.
We can identify $\CC_M$ with the category of pairs $(X,f)$
with $X \in \add(M)$ and
$f\df X \to \tau(X)$ a $KQ$-module homomorphism.
Clearly, $\CC_M$ is an additive subcategory.

\begin{Lem}\label{extensionclosed}
$\CC_M$ is closed under extensions.
\end{Lem}

\begin{proof}
Let
$
0 \to (X,f) \to (Y,g) \to (Z,h) \to 0
$ 
be a short exact
sequence of $\LL$-modules with $(X,f),(Z,h) \in \CC_M$.
Applying the functor $\pi_Q$ we get
a short exact sequence
$$
0 \to X \to Y \to Z \to 0
$$ 
in $\md(KQ)$ with $X,Z \in \add(M)$.
For each indecomposable direct summand $Y_i$ of $Y$ there exists
a non-zero map $X \to Y_i$, which implies that $Y_i \in \add(M)$, 
or $Y_i$ is a non-zero direct summand of $Z$, which also
implies $Y_i \in \add(M)$.
Thus $(Y,g) \in \CC_M$.
\end{proof}

\begin{Lem}\label{factorclosed}
$\CC_M$ is closed under factor modules.
\end{Lem}

\begin{proof}
Let $(Y,g)$ be a factor module of some $(X,f) \in \CC_M$.
Then for every indecomposable direct summand $Y_i$ of $Y$ there
exists a non-zero map $X \to Y_i$, which implies
$Y_i \in \add(M)$.
It follows that $(Y,g) \in \CC_M$.
\end{proof}

\begin{Cor}
For a terminal $KQ$-module $M$ the following hold:
\begin{itemize}

\item[(i)]
$\CC_M = {{^\perp}(\CC_M^\perp)}$;

\item[(ii)]
For each $(X,f) \in \CC_M$ there exists a short exact sequence
$$
0 \to (X_1,f_1) \to (X,f) \to (X_2,f_2) \to 0
$$
with $(X_1,f_1) \in \CC_M$ and $(X_2,f_2) \in \CC_M^\perp$;

\item[(iii)]
$\CC_M^\perp = \{ (Y,g) \in \md(\LL) \mid Y \cap \add(M) = 0\}$.

\end{itemize}
\end{Cor}

\begin{proof}
Part (i) follows from Lemma~\ref{torsionlemma}, Lemma
\ref{extensionclosed} and Lemma~\ref{factorclosed}.
To prove (ii),
let $(X,f)$ be a $\LL$-module.
We can write 
$$
X = X_1 \oplus X_2
$$ 
where
$X_1$ is a maximal direct summand of $X$ such that
$X_1 \in \add(M)$, and $X_2$ is some complement.
Note that $X_1$ and $X_2$ are uniquely determined up to
isomorphism.
By $f_1$ we denote the restriction of $f$ to $X_1$.
Since $M$ is a terminal $KQ$-module, we get $\Hom_{KQ}(X_1,\tau(X_2)) = 0$.
Thus, the image of $f_1$ is contained in $\tau(X_1)$.
So we can regard $(X_1,f_1)$ as a $\LL$-module.
In particular,
$(X_1,f_1)$ is a submodule of $(X,f)$.
We get a short exact sequence of the form
$$
0 \to (X_1,f_1) \to (X,f) \to (X_2,f_2) \to 0
$$
with $(X_1,f_1) \in \CC_M$ and $(X_2,f_2) \in \CC_M^\perp$.
Also (iii) follows easily from these considerations.
\end{proof}


{\Large\section{Lifting homomorphisms from $\md(KQ)$ to $\md(\LL)$}
\label{liftsection}}


As before, let $I_1,\ldots,I_n$ be the indecomposable
injective $KQ$-modules.
For natural numbers $a \le b$ define
$$
I_{i,[a,b]} = \bigoplus_{j=a}^b \tau^j(I_i),
$$
\index{$(I_{i,[a,b]},e_{i,[a,b]})$}\noindent
and let
$$
e_{i,[a,b]} = \left(\bsm 
0&1&&\\
 & \ddots &\ddots  &\\
 &   &       0&1\\
 &   &        &0
\esm\right)
\df I_{i,[a,b]} \to \tau(I_{i,[a,b]})
$$
be the $KQ$-module homomorphism with $e_{i,[a,b]}(\tau^a(I_i)) = 0$ and 
whose restriction to $\tau^j(I_i)$ is
the identity for $a+1 \le j \le b$.
The $\LL$-modules of the form
$
(I_{i,[a,b]},e_{i,[a,b]})
$
are crucial for our theory.

Let $(X,f) \in \md(\LL)$.
Define $\Hom_{KQ}(X,\tau^a(I_i))_b$ as the subspace
of $\Hom_{KQ}(X,\tau^a(I_i))$ consisting of all morphisms
$h$ such that
$$
0 = \tau^{b-a+1}(h) \circ \tau^{b-a}(f) \circ \cdots \circ \tau(f) 
\circ f\df X \to \tau^{b+1}(I_i).
$$

\begin{Lem}\label{liftXI1}
For $1 \le i \le n$ and $a \le b$ there is an isomorphism of vector spaces
$$
\Hom_{KQ}(X,\tau^a(I_i))_b \to \Hom_\LL((X,f),(I_{i,[a,b]},e_{i,[a,b]}))
$$
$$
h_a \mapsto \widetilde{h_a} 
:= \left[\bsm h_a\\h_{a+1}\\\vdots\\h_b \esm\right]
$$
where
$$
h_{a+j} := \tau^j(h_a) \circ \tau^{j-1}(f) \circ \cdots 
\tau^2(f) \circ \tau(f) \circ f\df X \to \tau^{a+j}(I_i)
$$
for $1 \le j \le b-a$.
\end{Lem}

\begin{proof}
Let
$h \in \Hom_\LL((X,f),(I_{i,[a,b]},e_{i,[a,b]}))$.
Thus 
$$
h = \left[\bsm h_a\\h_{a+1}\\\vdots\\h_b \esm\right]\df
X \to \bigoplus_{j=a}^b \tau^j(I_i)
$$
is a $KQ$-module homomorphism such that the diagram
$$
\xymatrix{
X \ar[r]^h \ar[d]^f& I_{i,[a,b]} \ar[d]^{e_{i,[a,b]}}\\
\tau(X) \ar[r]^{\tau(h)} & \tau(I_{i,[a,b]})
}
$$
commutes, in other words
$$
\left[\bsm h_{a+1}\\\vdots\\h_b\\0 \esm\right]
=
\left(\bsm 
0&1&&\\
 & \ddots & \ddots  &\\
 &   &       0&1\\
 &   &        &0
\esm\right)
\cdot 
\left[\bsm h_a\\h_{a+1}\\\vdots\\h_b \esm\right]
=
\left[\bsm \tau(h_a)f\\\tau(h_{a+1})f\\\vdots\\\tau(h_b)f \esm\right].
$$
Thus $h_a\df X \to \tau^a(I_i)$ determines $h_{a+1},\ldots,h_b$.
So the homomorphisms space 
$$
\Hom_\LL((X,f),(I_{i,[a,b]},e_{i,[a,b]}))
$$ 
can be identified
with the space of all homomorphisms $h_a\df X \to \tau^a(I_i)$ such that
$$
0 = \tau(h_b) \circ f = 
\tau^{b-a+1}(h_a) \circ \tau^{b-a}(f) \circ \cdots \circ \tau(f) 
\circ f\df X \to \tau^{b+1}(I_i).
$$
In this case, set
$$
\widetilde{h_a} = \left[\bsm h_a\\h_{a+1}\\\vdots\\h_b \esm\right]\df 
(X,f) \to (I_{i,[a,b]},e_{i,[a,b]})
$$
where
$
h_{a+j} = \tau^j(h_a) \circ \tau^{j-1}(f) \circ \cdots \tau^2(f) 
\circ \tau(f) \circ f
$ 
for $1 \le j \le b-a$.
\end{proof}

In the above lemma we call $\widetilde{h_a}$ the {\it lift} of 
$h_a\df X \to \tau^a(I_i)$.
The lift of a $KQ$-module homomorphism 
$$
h = (h_j)_j\df X \to \bigoplus_{j \in J} Y_j 
$$
with $Y_j$ indecomposable preinjective for all $j$
is defined by lifting every component $h_j$ of this homomorphism.
(Each indecomposable direct summand $Y_j$ is of the
form $\tau^a(I_i)$ for some $1 \le i \le n$ and $a \ge 0$.
Of course, we also have to specify 
with respect to which $b \ge a$ we want to lift 
$h_j\df X \to \tau^a(I_i)$.)

\begin{Cor}\label{liftXI2}
Let $(X,f) \in \CC_M$.
Then for $1 \le i \le n$ and $0 \le a \le t_i(M)$ 
we get isomorphisms of vector spaces
$$
\Hom_{KQ}(X,\tau^a(I_i)) \to 
\Hom_\LL((X,f),(I_{i,[a,t_i(M)]},e_{i,[a,t_i(M)]}))
$$
$$
h_a \mapsto \widetilde{h_a} 
= \left[\bsm h_a\\h_{a+1}\\\vdots\\h_{t_i(M)} \esm\right]
$$
where
$h_{a+j} = \tau^j(h_a) \circ \tau^{j-1}(f) \circ \cdots \tau^2(f) 
\circ \tau(f) 
\circ f$ for $1 \le j \le t_i(M)-a$.
\end{Cor}

\begin{proof}
In Lemma~\ref{liftXI1} take $b = t_i(M)$.
We have $X \in \add(M)$, but $\tau^{t_i(M)+1}(I_i)$ is not in
$\add(M)$.
Thus 
$$
\Hom_{KQ}(X,\tau^{t_i(M)+1}(I_i)) = 0.
$$
So there is no condition on the choice of $h_a$.
\end{proof}

Let again $(X,f) \in \md(\LL)$.
For $1 \le i \le n$ and $a \le b$, define $\Hom_{KQ}(\tau^b(I_i),X)_a$
as the subspace of $\Hom_{KQ}(\tau^b(I_i),X)$ consisting of all
morphisms $h$ such that 
$$
0 = f \circ \tau^{-1}(f) \circ \tau^{-2}(f) \circ \cdots \circ \tau^{-(b-a)}(f)
\circ \tau^{-(b-a)}(h)\df \tau^a(I_i) \to \tau(X).
$$

\begin{Lem}\label{liftIX1} 
There is an isomorphism of vector spaces 
$$
\Hom_{KQ}(\tau^b(I_i),X)_a \to \Hom_\LL((I_{i,[a,b]},e_{i,[a,b]}),(X,f))
$$
$$
h_b \mapsto \widetilde{h_b} = (h_a,h_{a+1},\ldots,h_b)
$$
where
$$
h_{b-j} := \tau^{-1}(f) \circ \tau^{-2}(f) \circ \cdots \circ \tau^{-j}(f)
\circ \tau^{-j}(h_b)\df \tau^{b-j}(I_i) \to X
$$ 
for $1 \le j \le b-a$.
\end{Lem}

\begin{proof}
Let
$
h \in \Hom_\LL((I_{i,[a,b]},e_{i,[a,b]}),(X,f))
$
for some $a \le b$.
Thus 
$$
h = (h_a,h_{a+1},\ldots,h_b)\df
\bigoplus_{j=a}^b \tau^j(I_i) \to X
$$
is a $KQ$-module homomorphism such that the diagram
$$
\xymatrix{
I_{i,[a,b]} \ar[r]^h \ar[d]^{e_{i,[a,b]}}& X \ar[d]^f\\
\tau(I_{i,[a,b]}) \ar[r]^{\tau(h)} & \tau(X)
}
$$
commutes. In other words
\begin{align*}
(0,\tau(h_a),\ldots,\tau(h_{b-1}))
&=
(\tau(h_a),\tau(h_{a+1}),\ldots,\tau(h_b)) \cdot
\left(\bsm 
0&1&&\\
 & \ddots & \ddots  &\\
 &   &       0&1\\
 &   &        &0
\esm\right) \\
&=
(fh_a,fh_{a+1},\ldots,fh_b).
\end{align*}
Thus $h_b\df \tau^b(I_i) \to X$ determines $h_a,\ldots,h_{b-1}$.
So the homomorphism space
$$
\Hom_\LL((X,f),(I_{i,[a,b]},e_{i,[a,b]}))
$$ 
can be identified
with the space of all homomorphisms $h_b\df \tau^b(I_i) \to X$ such that
$$
0 = f \circ h_a = f \circ \tau^{-1}(f) 
\circ \tau^{-2}(f) \circ \cdots \circ \tau^{-(b-a)}(f)
\circ \tau^{-(b-a)}(h_b) \df \tau^a(I_i) \to \tau(X).
$$
In this case, set
$$
\widetilde{h_b} = (h_a,h_{a+1},\ldots,h_b)\df 
(I_{i,[a,b]},e_{i,[a,b]}) \to (X,f)
$$
where
$
h_{b-j} = \tau^{-1}(f) \circ \tau^{-2}(f) \circ \cdots \circ \tau^{-j}(f)
\circ \tau^{-j}(h_b)
$ 
for $1 \le j \le b-a$.
\end{proof}

\begin{Cor}\label{liftIX2}
Let $(X,f) \in \md(\LL)$.
Then 
for $1 \le i \le n$ and $b \ge 0$ we get an isomorphism of vector spaces
$$
\Hom_{KQ}(\tau^b(I_i),X) \to \Hom_\LL((I_{i,[0,b]},e_{i,[0,b]}),(X,f))
$$
$$
h_b \mapsto \widetilde{h_b} = (h_0,h_1,\ldots,h_b)
$$
where
$h_{b-j} = \tau^{-1}(f) \circ \tau^{-2}(f) \circ \cdots \circ \tau^{-j}(f)
\circ \tau^{-j}(h_b)$ for $1 \le j \le b$.
\end{Cor}

\begin{proof}
In Lemma~\ref{liftIX1} take $a= 0$.
We have identified 
$
\Hom_\LL((I_{i,[0,b]},e_{i,[0,b]}),(X,f))
$ 
with the space of all homomorphisms $h_b\df \tau^b(I_i) \to X$ such that
$$
0 = fh_0\df I_i \to \tau(X),
$$
where $h_0$ is obtained from $h_b$ as described in Lemma~\ref{liftIX1}.
But for every $X \in \md(KQ)$ we have
$
\Hom_{KQ}(I_i,\tau(X)) = 0.
$
Thus there is no condition on the choice of $h_b$.
\end{proof}


{\Large\section{Construction of some $\CC_M$-complete rigid modules}
\label{completerigid}}


\subsection{The modules $T_M$ and $T_M^\vee$}\label{tia}
In this section,
let $M = M_1 \oplus \cdots \oplus M_r$ be a terminal $KQ$-module, and 
for $1 \le i \le n$ and $a \le b$ let 
$(I_{i,[a,b]},e_{i,[a,b]})$ be the $\LL$-module defined in Section
\ref{liftsection}.
For brevity, let $t_i := t_i(M)$.
Define$$
T_M := \bigoplus_{i=1}^n \bigoplus_{a=0}^{t_i} 
(I_{i,[a,t_i]},e_{i,[a,t_i]})
\text{\;\;\; and \;\;\;}
T_M^\vee := \bigoplus_{i=1}^n \bigoplus_{b=0}^{t_i} 
(I_{i,[0,b]},e_{i,[0,b]}).
$$
\index{$T_M$}\noindent
\index{$T_M^\vee$}\noindent

\begin{Lem}\label{constructionrigid}
For $1 \le i,j \le n$, $0 \le a \le t_i$ and $0 \le c \le t_j$ 
we have
$$
\Ext_\LL^1((I_{i,[a,t_i]},e_{i,[a,t_i]}),
(I_{j,[c,t_j]},e_{j,[c,t_j]})) = 0.
$$
\end{Lem}

\begin{proof}
By Lemma~\ref{extandhom} we know that
\begin{align*}
|\Ext_\LL^1((I_{i,[a,t_i]},e_{i,[a,t_i]}),(I_{j,[c,t_j]},e_{j,[c,t_j]}))| &=
|\Hom_\LL((I_{i,[a,t_i]},e_{i,[a,t_i]}),(I_{j,[c,t_j]},e_{j,[c,t_j]}))|\\
&+ |\Hom_\LL((I_{j,[c,t_j]},e_{j,[c,t_j]}),(I_{i,[a,t_i]},e_{i,[a,t_i]}))|\\
&- |\Hom_{KQ}(I_{i,[a,t_i]},I_{j,[c,t_j]})|\\
&- |\Hom_{KQ}(I_{j,[c,t_j]},I_{i,[a,t_i]})|\\  
&+ |\Ext_{KQ}^1(I_{i,[a,t_i]},I_{j,[c,t_j]})|\\ 
&+ |\Ext_{KQ}^1(I_{j,[c,t_j]},I_{i,[a,t_i]})|. 
\end{align*}
From Corollary~\ref{liftXI2} we get
$$
\Hom_\LL((I_{i,[a,t_i]},e_{i,[a,t_i]}),(I_{j,[c,t_j]},e_{j,[c,t_j]})) \cong
\Hom_{KQ}(I_{i,[a,t_i]},\tau^c(I_j)).
$$
Furthermore, the Auslander-Reiten formula yields
\begin{align*}
\Ext_{KQ}^1(I_{j,[c,t_j]},I_{i,[a,t_i]}) &\cong 
{\rm D}\Hom_{KQ}(I_{i,[a,t_i]},\tau(I_{j,[c,t_j]}))\\
&= {\rm D}\Hom_{KQ}\left(I_{i,[a,t_i]}, 
\tau\left(\bigoplus_{l=c}^{t_j} \tau^l(I_j)\right)\right)\\
&= {\rm D}\Hom_{KQ}\left(I_{i,[a,t_i]}, 
\bigoplus_{l=c+1}^{t_j+1} \tau^l(I_j)\right).
\end{align*}
Note that, since $\tau^{t_j+1}(I_j)\not\in \add(M)$, we have
$\Hom_{KQ}(I_{i,[a,t_i]},\tau^{t_j+1}(I_j)) = 0$.
This implies
\begin{multline*}
|\Ext_{KQ}^1(I_{j,[c,t_j]},I_{i,[a,t_i]})|=\\
|\Hom_{KQ}(I_{i,[a,t_i]},I_{j,[c,t_j]})| -
|\Hom_\LL((I_{i,[a,t_i]},e_{i,[a,t_i]}),(I_{j,[c,t_j]},e_{j,[c,t_j]}))|.
\end{multline*}
We also have a similar equality where $i$ and $j$ are exchanged, as
well as $a$ and $c$.
Summing up these two equalities, we get
$\Ext_\LL^1((I_{i,[a,t_i]},e_{i,[a,t_i]}),(I_{j,[c,t_j]},e_{j,[c,t_j]})) = 0$.
\end{proof}

By Corollary~\ref{liftXI2} there is an isomorphism
$$
\Phi_{[a,t_i],[c,t_j]}\df 
\bigoplus_{l=a}^{t_i} \Hom_{KQ}(\tau^l(I_i),\tau^c(I_j))
\to
\Hom_\LL((I_{i,[a,t_i]},e_{i,[a,t_i]}),(I_{j,[c,t_j]},e_{j,[c,t_j]}))
$$
defined by 
$$
h_{l,c} := (0,\ldots,0,h_c,0,\ldots,0) 
\mapsto
\widetilde{h_{l,c}} =
\left[\bsm h_{l,c}\\h_{l,c+1}\\\vdots\\h_{l,t_j} \esm\right]
$$
where
$$
h_c \in \Hom_{KQ}(\tau^l(I_i),\tau^c(I_j))
$$
and
$$
h_{l,c+k} := \tau^k(h_{l,c}) \circ \tau^{k-1}(e_{i,[a,t_i]}) 
\circ \cdots \circ \tau(e_{i,[a,t_i]}) \circ e_{i,[a,t_i]}
\df I_{i,[a,t_i]} \to \bigoplus_{u=c}^{t_j} \tau^u(I_j)
$$ 
for $1 \le k \le t_j-c$.
An easy calculation shows that
$$
\widetilde{h_{l,c}} = 
\left(\bbm 
&& &h_c&&&\\
&&&&\tau(h_c) &&\\
&&&&& \ddots &\\
&&&&&& \tau^{m}(h_c)\\
&&&&&&\\
&&&&&&\\
&&&&&&
\ebm\right)
\df (I_{i,[a,t_i]},e_{i,[a,t_i]}) \to (I_{j,[c,t_j]},e_{j,[c,t_j]})
$$
where $m = \max\{t_i-l,t_j-c\}$.
Here, the entries of the 
$(t_j-c+1) \times (t_i-a+1)$-matrix $\widetilde{h_{l,c}}$ 
are homomorphisms between
the indecomposable direct summands of the $KQ$-modules  
$I_{i,[a,t_i]}$ and $I_{j,[c,t_j]}$.
The only non-zero entries are the maps $\tau^k(h_c)$, 
$0 \le k \le m$.
(To be more precise, these are non-zero if and only if $h_c \not= 0$.)

\begin{Lem}
For $1 \le i \le n$ and $a \le b$ the endomorphism ring
$$
\End_\LL((I_{i,[a,b]},e_{i,[a,b]}))
$$ 
is local.
\end{Lem}

\begin{proof}
In the above situation, assume $a=c$, $i=j$ and $t_i=t_j$.
Let $a \le l \le t_i$.
If $a < l$, then
it follows easily that $\widetilde{h_{l,c}}$ is nilpotent
for all $h_c \in \Hom_{KQ}(\tau^l(I_i),\tau^c(I_j))$.
It is also clear that these homomorphisms form an ideal $I$ in
$\End_\LL((I_{i,[a,b]},e_{i,[a,b]}))$.
(If we write $\widetilde{h_{l,c}}$ again as a 
$(t_j-c+1) \times (t_i-a+1)$-matrix, then this matrix
is upper triangular with zero entries on the diagonal.
If $a=l$, then we obtain a diagonal matrix.)
Every ideal consisting only of nilpotent elements is contained 
in the radical of $\End_\LL((I_{i,[a,b]},e_{i,[a,b]}))$.
Now the factor algebra
$\End_\LL((I_{i,[a,b]},e_{i,[a,b]}))/I$ is 1-dimensional
with basis the residue class of
$$
\left(\bbm
1_{\tau^a(I_i)} &&&\\
&\tau(1_{\tau^a(I_i)})&&\\
&&\ddots&\\
&&& \tau^{t_i-a}(1_{\tau^a(I_i)})
\ebm\right).
$$
Here we use that $\Hom_{KQ}(X,X) \cong K$ for all indecomposable
preinjective $KQ$-modules $X$.
This finishes the proof.
\end{proof}

\begin{Cor}
For $1 \le i \le n$ and $a \le b$,
the $\LL$-module $(I_{i,[a,b]},e_{i,[a,b]})$ is indecomposable.
\end{Cor}

\begin{Cor}\label{completerigidT_M}
$T_M$ and $T_M^\vee$ are basic $\CC_M$-complete rigid $\LL$-modules.
\end{Cor}

\begin{proof}
Clearly, the modules $T_M$ and $T_M^\vee$ are contained in $\CC_M$.
By Lemma~\ref{constructionrigid} we know that $T_M$ is rigid.
Similarly one shows that $T_M^\vee$ is rigid.
Each $\LL$-module of the form $(I_{i,[a,b]},e_{i,[a,b]})$ is indecomposable,
and we have  $(I_{i,[a,b]},e_{i,[a,b]}) \cong (I_{j,[c,d]},e_{j,[c,d]})$ if and
only if $i=j$ and $[a,b] = [c,d]$.
Thus we get 
$$
\Sigma(T_M) = \Sigma(T_M^\vee) = r.
$$
Imitating the proof of \cite[Theorem 1.1]{GSExt} it is easy to show that 
$\rk(\CC_M) \le r$.
Thus we get $\rk(\CC_M) = \Sigma(M) = r$. 
This finishes the proof.
\end{proof}

\begin{Cor}
$\rk(\CC_M) = r$.
\end{Cor}

For later use, let us introduce the following abbreviations:
For $1 \le i \le n$ and $0 \le a \le b \le t_i$ set
\begin{align*}
T_{i,[a,b]}   &:=  (I_{i,[a,b]},e_{i,[a,b]}),\\
T_{i,a}      &:= (I_{i,[a,t_i]},e_{i,[a,t_i]}),\\
T_{i,b}^\vee &:= (I_{i,[0,b]},e_{i,[0,b]}).
\end{align*}
\index{$T_{i,[a,b]}$}\noindent
\index{$T_{i,a}$}\noindent
\index{$T_{i,b}^\vee$}\noindent

\subsection{The quivers of $\End_\LL(T_M)$ and $\End_\LL(T_M^\vee)$}
Let $\GG_M^*$ be defined as in Section~\ref{existrigid}.
As before, let $\GG_{T_M}$ be the quiver of $\End_\LL(T_M)$.

\begin{Lem}\label{quiverT_M}
We have $\GG_{T_M} = \GG_M^*$, where for $1 \le i \le n$ and
$0 \le a \le t_i$
the vertex $\tau^a(I_i)$ of $\GG_M^*$ corresponds to the vertex
$(I_{i,[a,t_i]},e_{i,[a,t_i]})$ of $\GG_{T_M}$.
\end{Lem}

\begin{proof}
Let $(I_{i,[a,t_i]},e_{i,[a,t_i]})$ and 
$(I_{j,[c,t_j]},e_{j,[c,t_j]})$ be indecomposable direct 
summands of $T_M$.

We want to construct a well behaved basis $B_{(i,a),(j,c)}$ of
$$
\Hom_\LL((I_{i,[a,t_i]},e_{i,[a,t_i]}),(I_{j,[c,t_j]},e_{j,[c,t_j]})).
$$
We write $B_{(i,a),(j,c)}$ as a disjoint union
$$
B_{(i,a),(j,c)} = \bigcup_{l=a}^{t_i} B_{(i,l),(j,c)}.
$$
where $B_{(i,l),(j,c)}$ are the images 
(under the map $\Phi_{[a,t_i],[c,t_j]}$)
of residue classes
of paths (in the path category of $\II_Q$) 
from $\tau^l(I_i)$ to $\tau^c(I_j)$.

Here we use that the mesh category of $\II_Q$ is obtained from
the path category by factoring out the mesh relations, and that
the full subcategory of indecomposable preinjective $KQ$-modules
is equivalent to the mesh category of $\II_Q$.
For details on mesh categories we refer to 
\cite[Chapter 10]{GR} and \cite[Lecture 1]{Ri3}.

Now it is easy to check that the homomorphisms 
$\widetilde{h_{l,c}}$ we constructed above
are irreducible in $\add(T_M)$ if and only if 
$h_c \in \Hom_{KQ}(\tau^l(I_i),\tau^c(I_j))$ is irreducible 
in $\II_Q$, 
or $l=a+1$, $i=j$, $c = a+1$ and 
$\widetilde{h_{l,c}}$ is a non-zero multiple of
$$
\left(\bbm 
0&1_{\tau^{a+1}(I_i)}&&&\\
0&&\tau(1_{\tau^{a+1}(I_i)}) &&\\
\vdots&&& \ddots &\\
0&&&& \tau^{m}(1_{\tau^{a+1}(I_i)})
\ebm\right)
\df (I_{i,[a,t_i]},e_{i,[a,t_i]}) \to (I_{i,[a+1,t_i]},e_{i,[a+1,t_i]})
$$
where $m = t_i-a-1$.
In other words, $h_c \in \Hom_{KQ}(\tau^{a+1}(I_i),\tau^{a+1}(I_i))$ 
is a non-zero multiple of $1_{\tau^{a+1}(I_i)}$. 
This implies $\GG_{T_M} = \GG_M^*$.
\end{proof}

The following Lemma is proved similarly as Lemma~\ref{quiverT_M}.

\begin{Lem}\label{quiverT_Mvee}
We have $\GG_{T_M^\vee} = \GG_M^*$, where for $1 \le i \le n$
and $0 \le b \le t_i$
the vertex $\tau^b(I_i)$ of $\GG_M^*$ corresponds to the vertex
$(I_{i,[0,b]},e_{i,[0,b]})$ of $\GG_{T_M^\vee}$.
\end{Lem}

Note that the $\CC_M$-projective direct summands of $T_M$ correspond to the 
rightmost vertices of $\GG_M^*$, whereas the $\CC_M$-projective summands
of $T_M^\vee$ correspond to the leftmost vertices of $\GG_M^*$.

One could also use covering methods to prove Lemma~\ref{quiverT_M}
and Lemma~\ref{quiverT_Mvee}, compare \cite{GLSAus}.
But note that in \cite{GLSAus} we only deal with
$Q$ being a Dynkin quiver and for $M$ we take the direct
sum of all indecomposable $KQ$-modules.
In this case, we have
$T_M = P_{Q^\op}$ and 
$T_M^\vee = I_{Q^\op}$,
where $P_Q$ and $I_Q$ are defined in \cite[Sections 1.7 and 1.2]{GLSAus}.

\subsection{Dimension vectors of some $\End_\LL(T_M)$-modules}
\label{dimvectors2}
As before, let $M = M_1 \oplus \cdots \oplus M_r$ be a terminal
$KQ$-module, and set $B:= \End_\LL(T_M)$.
For a $\LL$-module $(X,f) \in \CC_M$ 
we want to compute the dimension vector of the $B$-module
$
\Hom_\LL((X,f),T_M).
$
Since the indecomposable projective $B$-modules are
just the modules $\Hom_\LL(T_{i,a},T_M)$, $1 \le i \le n$, 
$0 \le a \le t_i$, we know that the entries of the dimension vector
$
\dimv(\Hom_\LL((X,f),T_M))
$
are 
$$
\dm \Hom_B(\Hom_\LL(T_{i,a},T_M),\Hom_\LL((X,f),T_M))
$$
where $1 \le i \le n$, $0 \le a \le t_i$.
We have
\begin{align*}
\Hom_B(\Hom_\LL(T_{i,a},T_M),\Hom_\LL((X,f),T_M)) &\cong
\Hom_\LL((X,f),T_{i,a})\\
&\cong \Hom_{KQ}(X,\tau^a(I_i)).
\end{align*}
The first isomorphism follows from
Corollary~\ref{fullyfaithful} and Lemma~\ref{fullyfaithful3}.
For the second isomorphism
we use Corollary~\ref{liftXI2}.

In other words, the entries of 
$
\dimv(\Hom_\LL((X,f),T_M))
$
are 
$\dm \Hom_{KQ}(X,M_s)$
where $1 \le s \le r$.
We can easily calculate $\dm \Hom_{KQ}(X,M_s)$
using the mesh category of $\II_Q$, see
\cite[Chapter 10]{GR}, \cite[Lecture 1]{Ri3}.

\subsection{An example of type $\A_3$}\label{example}
Let $Q$ be the quiver 
{\small
$$
\xymatrix{
1 & 2 \ar[l]\ar[r] & 3
}
$$
}\noindent
and let $M$ be the direct sum of all six indecomposable $KQ$-modules.
Thus $\GG_M = \GG_Q$ looks as follows:
{\small
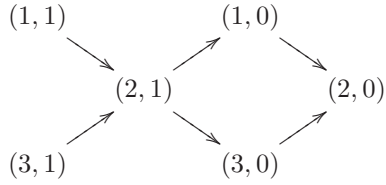
$$
\xymatrix@-1pc{
(1,1) \ar[dr] && (1,0) \ar[dr]
\\
& (2,1) \ar[dr]\ar[ur] && (2,0)
\\
(3,1) \ar[ur] && (3,0) \ar[ur]
}
$$
}\noindent
The following picture shows the quiver of 
$\End_\LL(T_M)$ where the vertices corresponding to the $T_{i,a}$
are labelled by the dimension vectors $\dimv(\Hom_\LL(T_{i,a},T_M))$. 
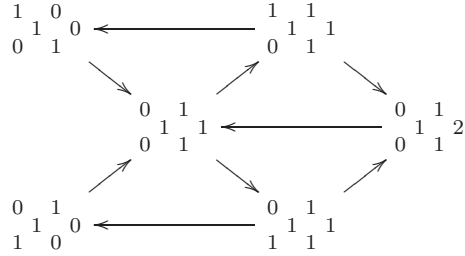
$$
\xymatrix@-1pc{
{\bsm1&&0\\&1&&0\\0&&1\esm} \ar[dr]
&& {\bsm1&&1\\&1&&1\\0&&1\esm} \ar[dr]\ar[ll]\\
&{\bsm0&&1\\&1&&1\\0&&1\esm} \ar[ur]\ar[dr]&& 
{\bsm0&&1\\&1&&2\\0&&1\esm} \ar[ll]\\
{\bsm0&&1\\&1&&0\\1&&0\esm}\ar[ur] &&
{\bsm0&&1\\&1&&1\\1&&1\esm}\ar[ll]\ar[ur]
}
$$
Similarly, the quiver of $\End_\LL(T_M^\vee)$ looks as follows:
$$
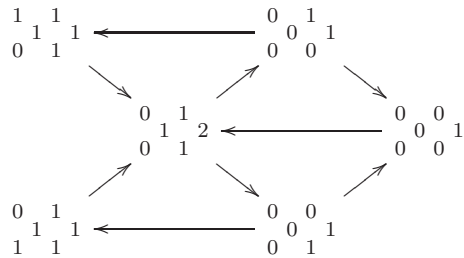
\xymatrix@-1pc{
{\bsm1&&1\\&1&&1\\0&&1\esm} \ar[dr]
&& {\bsm0&&1\\&0&&1\\0&&0\esm} \ar[dr]\ar[ll]\\
&{\bsm0&&1\\&1&&2\\0&&1\esm} \ar[ur]\ar[dr]&& 
{\bsm0&&0\\&0&&1\\0&&0\esm} \ar[ll]\\
{\bsm0&&1\\&1&&1\\1&&1\esm} \ar[ur] &&
{\bsm0&&0\\&0&&1\\0&&1\esm}\ar[ll]\ar[ur]
}
$$
The vertices corresponding to the $T_{i,b}^\vee$ are labelled by
the vectors $\dimv(\Hom_\LL(T_{i,b}^\vee,T_M))$.

\subsection{An example of type $\widetilde{\A}_2$}\label{A2tilde}
Let $Q$ be the quiver 
$$
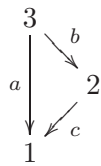
\xymatrix@-1pc{
3 \ar[dr]^b\ar[dd]_a\\
&2\ar[dl]^c\\
1
}
$$
and let $M$ be the terminal $KQ$-module with $t_i(M) = 1$
for all $i$.
Then the quiver of $\End_\LL(T_M)$ looks as follows:
$$
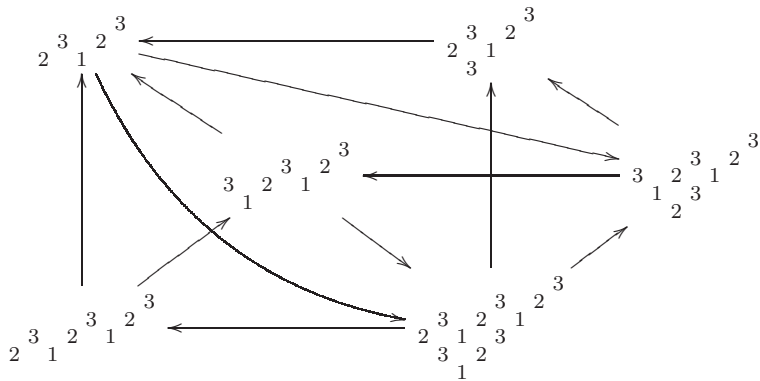
\xymatrix@-0.7pc{
{\bsm
&&&&3\\
&3&&2\\
2&&1
\esm} \ar@/_1.1cm/[ddrr]\ar[drrr]
&&
{\bsm
&&&&3\\
&3&&2\\
2&&1\\
&3
\esm
} \ar[ll] 
\\
&
{\bsm
&&&&&&3\\
&&&3&&2\\
3&&2&&1\\
&1
\esm} \ar[ul]\ar[dr]
&&
{\bsm
&&&&&&3\\
&&&3&&2\\
3&&2&&1\\
&1&&3\\
&&2
\esm
} \ar[ll]\ar[ul]
\\
{\bsm
&&&&&&&3\\
&&&&3&&2\\
&3&&2&&1\\
2&&1
\esm} \ar[uu]\ar[ur] 
&&
{\bsm
&&&&&&&3\\
&&&&3&&2\\
&3&&2&&1\\
2&&1&&3\\
&3&&2\\
&&1
\esm} \ar[ll]\ar[uu]\ar[ur]
}
$$
As vertices we display the indecomposable direct
summands of the $\LL$-module $T_M$.
The numbers can be interpreted as basis vectors or
as composition factors.
For example, 
$$
X := {\bsm
&&&&&&3\\
&&&3&&2\\
3&&2&&1\\
&1&&3\\
&&2
\esm
}
$$
stands for a 9-dimensional $\LL$-module
with dimension vector $(2,3,4)$.
More precisely, one could display the module $X$ as follows:
$$
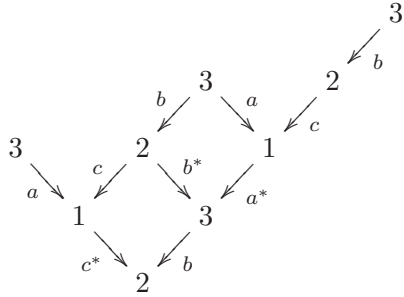
\xymatrix@-1pc{
&&&&&&3\ar[dl]^b\\
&&&3\ar[dl]_b\ar[dr]^a&&2\ar[dl]^c\\
3\ar[dr]_a && 2\ar[dr]^{b^*}\ar[dl]_c&&1\ar[dl]^{a^*}\\
&1\ar[dr]_{c^*}&&3\ar[dl]^b\\
&&2
}
$$
This picture shows how the different arrows
of the quiver $\overline{Q}$ of $\LL$ act
on the 9 basis vectors of the module.
For example, one can see immediately that
the socle of $X$ is isomorphic to $S_2$, and the
top is isomorphic to $S_3 \oplus S_3 \oplus S_3$.
One can also clearly see how the 2-dimensional indecomposable
injective $KQ$-module $I_2$ is ``glued from below'' to the
indecomposable $KQ$-module $\tau(I_2)$ using the arrows
$a^*,b^*,c^*$.
This gives a short exact sequence
$$
0 \to (I_2,0) \to X \to (\tau(I_2),0) \to 0
$$
of $\LL$-modules.

Using the same notation, the quiver of $\End_\LL(T_M^\vee)$ looks
as follows:
$$
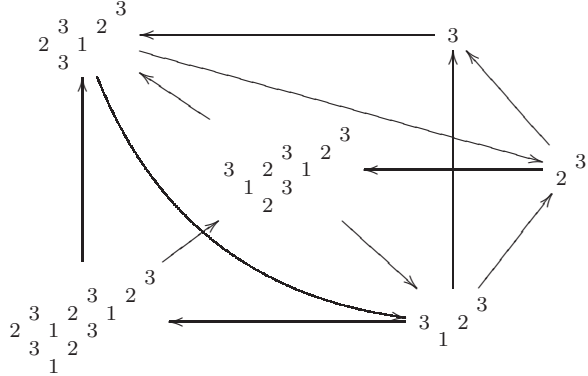
\xymatrix@-0.7pc{
{\bsm
&&&&3\\
&3&&2\\
2&&1\\
&3
\esm} \ar@/_1.1cm/[ddrr]\ar[drrr]
&&
{\bsm
3
\esm} \ar[ll] 
\\
&
{\bsm
&&&&&&3\\
&&&3&&2\\
3&&2&&1\\
&1&&3\\
&&2
\esm} \ar[ul]\ar[dr]
&&
{\bsm
&3\\
2
\esm} \ar[ll]\ar[ul]
\\
{\bsm
&&&&&&&3\\
&&&&3&&2\\
&3&&2&&1\\
2&&1&&3\\
&3&&2\\
&&1
\esm} \ar[uu]\ar[ur] 
&&
{\bsm
&&&3\\
3&&2\\
&1
\esm} \ar[ll]\ar[uu]\ar[ur]
}
$$


{\Large\section{$\CC_M$ is a cluster torsion class and is stably 2-Calabi-Yau}
\label{section7}}


\subsection{$\CC_M$ is a cluster torsion class}
As before, let $M = M_1 \oplus \cdots \oplus M_r$ be a terminal
$KQ$-module.
Set
$$
I_M = \bigoplus_{i=1}^n (I_{i,[0,t_i(M)]},e_{i,[0,t_i(M)]}).
$$
The following two lemmas are a direct consequence of 
Lemma~\ref{liftXI1} and Corollary~\ref{liftXI2}.

\begin{Lem}\label{homdim1}
For all $1 \le i \le n$ and $b \ge 0$, 
the $\LL$-module $(I_{i,[0,b]},e_{i,[0,b]})$ has a simple socle which
is isomorphic to $(S_i,0)$.
\end{Lem}

\begin{Lem}\label{homdim2}
For $X \in \CC_M$ we have
$
\dm \Hom_\LL(X,I_M) = \dm X.
$
\end{Lem}

\begin{Lem}\label{rightapprox1}
Let $(X,f) \in \CC_M$.
Then there exists a short exact sequence
$$
0 \to (X,f) \to (I,e) \to (Y,g) \to 0
$$
of $\LL$-modules with $(I,e) \in \add(I_M)$ and $(Y,g) \in \CC_M$.
\end{Lem}

\begin{proof}
Let
$$
h\df X \to \bigoplus_{i=1}^n I_i^{m_i} 
$$
be a monomorphism of $KQ$-modules.
Such a monomorphism exists, since $I_1,\cdots,I_n$ are the indecomposable
injective $KQ$-modules.
It follows from Corollary~\ref{liftXI2} that the lift 
$$
\widetilde{h}\df
(X,f) \to (I,e) := \bigoplus_{i=1}^n (I_{i,[0,t_i(M)]},e_{i,[0,t_i(M)]})^{m_i}
$$
of $h$ is a monomorphism of $\LL$-modules.
We denote its cokernel by $(Y,g)$.
Since $\CC_M$ is closed under factor modules, $(Y,g)$ 
is contained in $\CC_M$.
\end{proof}

\begin{Cor}\label{cogenerator}
$I_M$ is a cogenerator of $\CC_M$.
\end{Cor}

\begin{Lem}\label{Extinjective}
$I_M$ is $\CC_M$-injective.
\end{Lem}

\begin{proof}
It is enough to show that for $1 \le i \le n$ the module
$(I_{i,[0,t_i(M)]},e_{i,[0,t_i(M)]})$ is $\CC_M$-injective. 
Suppose $h'\df (X,f) \to (Y,g)$ is a monomorphism in $\CC_M$, and let 
$$
h\df (X,f) \to (I_{i,[0,t_i(M)]},e_{i,[0,t_i(M)]})
$$ 
be an arbitrary homomorphism.
By Corollary~\ref{liftXI2} we know that 
$$
h = \widetilde{h_0} = \left[\bsm h_0\\h_1\\ \vdots\\h_{t_i(M)} \esm\right]
$$
for some $KQ$-module homomorphism $h_0\df X \to I_i$.
Since $I_i$ is injective (as a $KQ$-module), 
and since $h'\df X \to Y$ is a monomorphism,
there exists some $KQ$-module homomorphism $h_0''\df Y \to I_i$ such that
$h_0'' \circ h' = h_0$.
$$
\xymatrix{
X \ar[d]_{h_0} \ar[r]^{h'} & Y \ar[dl]^{h_0''}\\
I_i
}
$$
We want to show that there exists a homomorphism 
$h''\df (Y,g) \to (I_{i,[0,t_i(M)]},e_{i,[0,t_i(M)]})$
such that the diagram
$$
\xymatrix{
(X,f) \ar[r]^{h'} \ar[d]^h& (Y,g) \ar[dl]^{h''}\\
(I_{i,[0,t_i(M)]},e_{i,[0,t_i(M)]})
}
$$
commutes.

Recall that a homomorphism $(X,f) \to (I_{i,[0,t_i(M)]},e_{i,[0,t_i(M)]})$ is
already determined by its component $X \to I_i$.

Let
$$
h'' := \widetilde{h_0''} = 
\left[\bsm h_0''\\h_1''\\ \vdots \\h_{t_i(M)}''\esm\right]
$$
be the lift of $h_0''$. 
It follows that the component $X \to I_i$ of the homomorphisms
$h'' \circ h'$ and $h$ is equal, namely $h_0''\circ h' = h_0$,
thus 
$h'' \circ h' = h$.

For brevity, let $I := (I_{i,[0,t_i(M)]},e_{i,[0,t_i(M)]})$.
Assume $Z \in \CC_M$.
We have to show that $\Ext_\LL^1(Z,I) = 0$.
Let 
$$
0 \to I \xrightarrow{f} E \to Z \to 0
$$
be a short exact sequence of $\LL$-modules.
By the above considerations, we know that $\Hom_\LL(f,I)$ is
surjective.
In particular, there exists a homomorphism $f'\df E \to I$ such that
$f'f = \id_I$.
Thus $f$ is a split monomorphism and the above sequence splits.
This finishes the proof.
\end{proof}

\begin{Lem}
If $C$ is a cogenerator of $\CC_M$, then $\add(C)$ contains all 
modules which are $\CC_M$-injective.
\end{Lem}

\begin{proof}
Let $I$ be $\CC_M$-injective.
Then there exists a short exact sequence
$$
0 \to I \xrightarrow{f} C' \to \Coker(f) \to 0
$$
of $\LL$-modules
with $C' \in \add(C)$.
We know that $\Coker(f) \in \CC_M$, because $\CC_M$ is closed under
factor modules.
Since $I$ is $\CC_M$-injective, the above sequence splits.
Therefore,
$I \in \add(C') \subseteq \add(C)$.
\end{proof}

Summarizing, we obtain the following:

\begin{Prop}\label{extprojinj}
If $M$ is a terminal $KQ$-module, then
$$
\add(I_M) = \{ \text{$\CC_M$-projectives} \} = 
\{ \text{$\CC_M$-injectives} \}.
$$
\end{Prop}

Now, let 
$$
T := \bigoplus_{i=1}^n \tau^{t_i(M)}(I_i).
$$
Recall that $T$ is a tilting module over $KQ$, and that
$$
\add(M) = \Gen(T) = \left\{ N \in \md(KQ) \mid \Ext_{KQ}^1(T,N) = 0 \right\}.
$$

\begin{Lem}\label{rightapprox}
Let $X$ be a $KQ$-module in $\add(M)$.
Then there exists a short exact sequence 
$$
0 \to T'' \to T' \xrightarrow{h} X \to 0
$$
of $KQ$-modules with $T',T'' \in \add(T)$
and $h$ a right $\add(T)$-approximation.
\end{Lem}

\begin{proof}
We deduce the result from the proof of \cite[Prop. 1.4 (b)]{Bo1}.
Let $h\df T' \to X$ be a right $\add(T)$-approximation of $X$.
Since $X \in \add(M) = \Gen(T)$, we know that $h$ is an
epimorphism.
Let $T'' = \Ker(h)$.
We obtain a short exact sequence
$$
0 \to T'' \to T' \xrightarrow{h} X \to 0.
$$
Applying $\Hom_{KQ}(T,-)$ to this sequence yields an exact sequence
$$
\Hom_{KQ}(T,T') \xrightarrow{\Hom_{KQ}(T,h)} \Hom_{KQ}(T,X)
\to \Ext_{KQ}^1(T,T'') \to \Ext_{KQ}^1(T,T') = 0. 
$$
Since $h$ is a right $\add(T)$-approximation, $\Hom_{KQ}(T,h)$
is surjective.
It follows that $\Ext_{KQ}^1(T,T'') = 0$.
Thus $T'' \in \add(M)$.
Every indecomposable direct summand of $T''$ maps non-trivially to
a module in $\add(T)$.
But the only modules in $\add(M)$ with this property lie
in $\add(T)$.
Thus $T'' \in \add(T)$.
This finishes the proof.
\end{proof}

\begin{Lem}\label{rightapprox2}
Let $(X,f) \in \CC_M$.
Then there exists a short exact sequence
$$
0 \to (Y,g) \to (I,e) \to (X,f) \to 0
$$
of $\LL$-modules with $(I,e) \in \add(I_M)$ and $(Y,g) \in \CC_M$.
\end{Lem}

\begin{proof}
Let 
$$
0 \to T'' \to T' \xrightarrow{h} X \to 0
$$
be the short exact sequence appearing in Lemma~\ref{rightapprox}.
It follows that
$$
T' = \bigoplus_{i=1}^n (\tau^{t_i(M)}(I_i))^{m_i}
$$
for some $m_i \ge 0$.
Set 
$$
(I,e) 
= \bigoplus_{i=1}^n(I_{i,[0,t_i(M)]},e_{i,[0,t_i(M)]})^{m_i}.
$$
Note that $(I,e) \in \add(I_M)$.
By Corollary~\ref{liftIX2} we can lift $h$ to a $\LL$-module homomorphism
$$
\widetilde{h}\df (I,e) \to (X,f).
$$
We denote the kernel of $\widetilde{h}$ by $(Y,g)$.
Thus we obtain a short exact sequence of $KQ$-modules
$$
0 \to Y \to I \xrightarrow{\widetilde{h}} X \to 0.
$$
Since $h$ occurs as a component of the homomorphism $\widetilde{h}$ and
since $h$ is a right $\add(T)$-approximation of $X$, we know that the
map 
$$
\Hom_{KQ}(T,\widetilde{h})\df \Hom_{KQ}(T,I) \to
\Hom_{KQ}(T,X)
$$ 
is surjective.
The module $I$ lies in $\add(M)$, thus 
$\Ext_{KQ}^1(T,I) = 0$.
So we get 
$$
\Ext_{KQ}^1(T,Y) = 0.
$$
This implies $Y \in \add(M)$, and therefore $(Y,g) \in \CC_M$.
\end{proof}

\begin{Cor}\label{generator}
$I_M$ is a generator of $\CC_M$.
\end{Cor}

\begin{Prop}\label{CMtorsion}
Let $M$ be a terminal $KQ$-module.
Then
$\CC_M$ is a cluster torsion class of $\nil(\LL)$ with
$
\rk(\CC_M) = \Sigma(M).
$
\end{Prop}

\begin{proof}
Combine Lemma~\ref{extensionclosed}, Lemma~\ref{factorclosed}, 
Proposition~\ref{extprojinj}, Corollary~\ref{cogenerator},
Corollary~\ref{generator}, Corollary~\ref{completerigidT_M} 
and Lemma~\ref{quiverT_M}.
\end{proof}

\begin{Cor}\label{rightapprox3}
Let $T$ be a $\CC_M$-maximal rigid $\LL$-module, and let $X \in \CC_M$.
Then there exists an exact sequence 
$
T'' \to T' \to X \to 0
$
with $T',T'' \in \add(T)$.
\end{Cor}

\begin{proof}
Every $\CC_M$-maximal rigid $\LL$-module contains $I_M$ as a direct
summand. Then use Lemma~\ref{rightapprox2} to get a surjective
map $h \df T' \to X$ with $T'\in \add(T)$ and $\Ker(h)\in\CC_M$, and 
a second time to get a surjective map $T''\to \Ker(h)$ with 
$T''\in\add(T)$. 
\end{proof}

\begin{Cor}\label{Frobenius2}
Let $M$ be a terminal $KQ$-module.
Then $\CC_M$ is a Frobenius category.
 \end{Cor}

\begin{proof}
Combine Proposition~\ref{extprojinj}, Lemma~\ref{rightapprox1} and 
Lemma~\ref{rightapprox2}.
\end{proof}

\subsection{The stable category ${\stCC}_M$ is 2-Calabi-Yau}
Let ${\stCC}_M$ be the {\it stable category} of $\CC_M$.
\index{${\stCC}_M$}
\index{stable category of $\CC_M$}
By definition the objects in ${\stCC}_M$ are the same as the
objects in $\CC_M$, and the morphisms spaces are the morphism
spaces in $\CC_M$ modulo morphisms factoring through 
$\CC_M$-projective-injective objects.
The category ${\stCC}_M$ is a triangulated 
category in a natural way \cite{H2}.
The shift is given by the relative syzygy functor
$$
\Omega_M^{-1}\df {\stCC}_M \to {\stCC}_M.
$$
We know that $\CC_M$ is closed under
extensions.
This implies
$$
\Ext_{\CC_M}^1(X,Y) = \Ext_\LL^1(X,Y)
$$ for all
objects $X$ and $Y$ in $\CC_M$.

Since $\CC_M$ is a Frobenius category, there is a functorial isomorphism
\begin{equation}\label{functorial1}
\underline{\CC}_M(X,\Omega_M^{-1}(Y)) \cong \Ext_{\CC_M}^1(X,Y)
\end{equation}
for all $X$ and $Y$ in $\CC_M$:
Every $f \in \CC_M(X,\Omega_M^{-1}(Y))$ gives rise to a commutative diagram
$$
\xymatrix{
\eta_f: & 0 \ar[r] & Y \ar@{=}[d] \ar[r] & E \ar[r] \ar[d] & X \ar[r] 
\ar[d]^{f} & 0\\
&0 \ar[r] & Y \ar[r] & I_Y \ar[r] & \Omega_M^{-1}(Y) \ar[r] & 0.
}
$$
The lower sequence is obtained from the embedding of $Y$ into
its injective hull $I_Y$ in $\CC_M$, and the
upper short exact sequence is just the pull-back of $f$.
Then $f \mapsto \eta_f$ yields the isomorphism (\ref{functorial1}).

Furthermore, using the canonical projective bimodule 
resolution of $\LL$, it is not
difficult to show that for all $\LL$-modules $X$ and $Y$ there exists
a functorial isomorphism
\begin{equation}\label{functorial2}
\Ext_\LL^1(X,Y) \cong {\rm D}\Ext_\LL^1(Y,X),
\end{equation}
see \cite[\S8]{GLSSemi2}.

Let $\T$ be a $K$-linear Hom-finite triangulated category
with shift functor $[1]$.
Then $\T$ is a {\it 2-Calabi-Yau category} 
\index{2-Calabi-Yau category}
if for all $X,Y \in \T$ there
is a functorial isomorphism
$$
\T(X,Y) \cong {\rm D}\T(Y,X[2]).
$$
If additionally $\T = \underline{\CC}$ 
for some Frobenius category $\CC$,
then $\T$ is called {\it algebraic}.
\index{algebraic 2-Calabi-Yau category}

\begin{Prop}
${\stCC}_M$ is an algebraic 2-Calabi-Yau category.
\end{Prop}

\begin{proof}
We have
\begin{align*}
\underline{\CC}_M(X,Y) &\cong \Ext_{\CC_M}^1(X,\Omega_M(Y))\\
                    &\cong {\rm D}\Ext_{\CC_M}^1(\Omega_M(Y),X)\\
                    &\cong {\rm D}\Ext_{\CC_M}^1(Y,\Omega_M^{-1}(X))\\
                    &\cong {\rm D}\underline{\CC}_M(Y,\Omega_M^{-2}(X)),   
\end{align*}
and all these isomorphisms are functorial.
\end{proof}


{\Large\section{Relative homology and $\CC$-maximal rigid modules}\label{section8}}


In this section, we recall some notions from relative homology theory
which, for Artin algebras, was developed by Auslander and Solberg
\cite{AS1,AS2}.

\subsection{Relative homology theory}\label{relhom1}
Let $A$ be a $K$-algebra, and let $X,Y,Z,T \in \md(A)$.
Set 
$$ 
F^T := \Hom_A(-,T)\df \md(A) \to \md(\End_A(T)).
$$
A short exact sequence
$$
0 \to X \to Y \to Z \to 0
$$
is $F^T$-{\it exact} if
\index{$F^T$-exact short exact sequence}
$0 \to F^T(Z) \to F^T(Y) \to F^T(X) \to 0$
is exact.
By $F^T(Z,X)$ we denote the set of equivalence classes of
$F^T$-exact sequences.

Let $\X_T$ be the subcategory of all
$X \in \md(A)$ such that
there exists an exact sequence
\begin{equation}\label{seq7}
0 \to X \to T_0 \xrightarrow{f_0} T_1 
\xrightarrow{f_1} T_2 \xrightarrow{f_2} \cdots
\end{equation}
where $T_i \in \add(T)$ for all $i$ and
$$
0 \to \Ker(f_i) \to T_i \to \Ima(f_i) \to 0
$$
are $F^T$-exact for all $i \ge 0$.
Sequence (\ref{seq7}) is an $F^T$-injective coresolution of $X$ in the
sense of \cite{AS2}.
Note that
$$
\add(T) \subseteq \X_T.
$$

For $X \in \X_T$ and $Z \in \md(A)$ let
$\Ext_{F^T}^i(Z,X)$, $i \ge 0$ 
be the cohomology groups obtained from the complex
\begin{equation}\label{seq8}
0 \to T_0 \xrightarrow{f_0} T_1 
\xrightarrow{f_1} T_2 \xrightarrow{f_2} \cdots
\end{equation}
by applying the functor $\Hom_A(Z,-)$.

\begin{Lem}[{\cite{AS1}}]
For $X \in \X_T$ and $Z \in \md(A)$ there is a functorial
isomorphism 
$$
\Ext_{F^T}^1(Z,X) = F^T(Z,X).
$$
\end{Lem}

\begin{Prop}[{\cite[Prop 3.7]{AS2}}]\label{seslift}
For $X \in \X_T$ and $Z \in \md(A)$ 
there is a functorial isomorphism
$$
\Ext_{F^T}^i(Z,X) \to \Ext_{\End_A(T)}^i(\Hom_A(X,T),\Hom_A(Z,T)) 
$$
for all $i \ge 0$.
\end{Prop}

\begin{Cor}\label{yoneda1}
For $X \in \X_T$ and $Z \in \md(A)$ 
there is a functorial isomorphism
$$
i_{Z,X,T}\df \Hom_A(Z,X) \to \Hom_{\End_A(T)}(\Hom_A(X,T),\Hom_A(Z,T))
$$
$$
h \mapsto (h' \mapsto h'h).
$$
\end{Cor}

If $X = Z$, we define $i_{X,T} := i_{X,X,T}$.

\begin{Cor}\label{yoneda2}
For $X \in \X_T$ and $Z \in \md(A)$
the map
$$
i_{X,T}\df \End_A(X) \to 
\Hom_{\End_A(T)}(\Hom_A(X,T),\Hom_A(X,T))
$$
is an anti-isomorphism of rings.
In other words, we get a ring isomorphism
$$
\End_A(X) \to \End_{\End_A(T)}(\Hom_A(X,T))^\op.
$$
\end{Cor}

\begin{proof}
It follows from the definitions that 
$i_{X,T}(h_1 \circ h_2) = i_{X,T}(h_2) \circ i_{X,T}(h_1)$.
\end{proof}

\begin{Cor}\label{fullyfaithful}
The functor 
$$
\Hom_A(-,T)\df \X_T \to \md(\End_A(T))
$$
is fully faithful.
In particular, $\Hom_A(-,T)$ has the following properties:
\begin{itemize}

\item[(i)]
If $X \in \X_T$ is indecomposable, then $\Hom_A(X,T)$ is indecomposable;

\item[(ii)]
If $\Hom_A(X,T) \cong \Hom_A(Y,T)$ for some $X,Y \in \X_T$, 
then $X \cong Y$.

\end{itemize}
\end{Cor}

Note that Corollary~\ref{fullyfaithful} follows already from
\cite[Section 3]{A}, see also 
\cite[Lemma 1.3 (b)]{APR}.

\begin{Cor}\label{fullyfaithful2}
Let $T \in \md(A)$, and let $\CC$ be an extension closed subcategory
of $\X_T$.
If 
$$
\psi\df 0 \to \Hom_A(Z,T) \xrightarrow{\Hom_A(g,T)} \Hom_A(Y,T)
\xrightarrow{\Hom_A(f,T)} \Hom_A(X,T) \to 0
$$
is a short exact sequence of $\End_A(T)$-modules with $X,Y,Z \in \CC$,
then 
$$
\eta\df 0 \to X \xrightarrow{f} Y \xrightarrow{g} Z \to 0 
$$
is a short exact sequence in $\md(A)$.
\end{Cor}

\begin{proof}
By Proposition~\ref{seslift} there exists an $F^T$-exact sequence
$$
\eta'\df 0 \to X \xrightarrow{f'} E \xrightarrow{g'} Z \to 0
$$
with $F^T(\eta') = \psi$.
Since $\CC$ is closed under extensions, we know that $E \in \CC$.
Now Corollary~\ref{fullyfaithful} implies that $E \cong Y$.
So there is a short exact sequence
$$
\eta''\df 0 \to X \xrightarrow{f''} Y \xrightarrow{g''} Z \to 0
$$
with $F^T(\eta'') = \psi$.
Again by Corollary~\ref{fullyfaithful} we get $f'' = f$ and $g'' = g$.
\end{proof}

\subsection{Relative homology for selfinjective torsion classes}
\label{relhom2}
The following lemma is stated in \cite[Lemma 5.1]{GLSRigid} for
preprojective algebras of Dynkin type.
But the same proof works for arbitrary preprojective algebras.

\begin{Lem}\label{mutation1}
Let $T$ and $X$ be rigid $\LL$-modules.
If
$$
0 \to X \xrightarrow{f} T' \to Y \to 0
$$
is a short exact sequence with $f$ a left $\add(T)$-approximation, 
then
$T \oplus Y$ is rigid.
\end{Lem}

\begin{Cor}\label{cormutation1}
Let $T$ and $X$ be rigid $\LL$-modules in a selfinjective torsion class 
$\CC$ of $\nil(\LL)$.
If $T$ is $\CC$-maximal rigid, then there exists a short exact sequence
$$
0 \to X \to T' \to T'' \to 0
$$
with $T',T'' \in \add(T)$.
\end{Cor}

\begin{proof}
In the situation of Lemma~\ref{mutation1}, if $T$ is $\CC$-maximal rigid,
we get $Y \in \add(T)$.
\end{proof}

\begin{Cor}\label{pdimone}
Let $T$ and $X$ be rigid $\LL$-modules in a selfinjective torsion class  
$\CC$ of $\nil(\LL)$.
If $T$ is $\CC$-maximal rigid, then $\Hom_\LL(X,T)$ is an $\End_\LL(T)$-module
with projective dimension at most one.
\end{Cor}

\begin{proof}
Applying $\Hom_\LL(-,T)$ to the short exact sequence
in Corollary~\ref{cormutation1} yields a projective resolution
$$
0 \to \Hom_\LL(T'',T) \to \Hom_\LL(T',T) \to \Hom_\LL(X,T) \to 0
$$
of the $\End_\LL(T)$-module $\Hom_\LL(X,T)$.
\end{proof}

\begin{Lem}\label{pdimone2}
Let $\CC$ be a selfinjective torsion class of $\nil(\LL)$.
If $T$ is a $\CC$-maximal 1-orthogonal $\LL$-module,
then for all $X \in \CC$ the $\End_\LL(T)$-module
$\Hom_\LL(X,T)$ has projective dimension at most one.
\end{Lem}

\begin{proof}
Let $X \in \CC$, and let $f\df X \to T'$ be a left $\add(T)$-approximation
of $X$.
Clearly, $f$ is injective, since $T$ is a cogenerator of $\CC$.
We obtain a short exact sequence
$$
0 \to X \xrightarrow{f} T' \to T'' \to 0
$$
where $T'' = \Coker(f)$.
Applying $\Hom_\LL(-,T)$ yields an exact sequence
$$
\eta\df 0 \to \Hom_\LL(T'',T) \to \Hom_\LL(T',T) \to \Hom_\LL(X,T) \to 0
$$
of $\End_\LL(T)$-modules.
It also follows that $\Ext_\LL^1(T'',T) = 0$.
Since $\CC$ is closed under factor modules, we know that $T'' \in \CC$.
This implies $T'' \in \add(T)$, because $T$ is $\CC$-maximal 1-orthogonal.
Thus $\eta$ is a projective resolution of $\Hom_\LL(X,T)$.
\end{proof}

\begin{Lem}\label{fullyfaithful3}
Let $\CC$ be a selfinjective torsion class of $\nil(\LL)$.
If $T$ is a $\CC$-maximal rigid $\LL$-module, then
$
\CC \subseteq \X_T.
$
\end{Lem}

\begin{proof}
For $X \in \CC$, let $f\df X \to T'$ be a 
left $\add(T)$-approximation, and let $Y$ be the cokernel 
of $f$.
Since $T$ is a cogenerator of $\CC$, we know that $f$ is injective.
The selfinjective torsion class $\CC$ is closed under factor modules,
thus $Y \in \CC$.
This yields the required $F^T$-injective coresolution of $X$.
\end{proof}


{\Large\section{Tilting and $\CC$-maximal rigid modules}\label{section9}}


In this section, we adapt some results due to Iyama \cite{I1,I2} to 
our situation of selfinjective torsion classes.

\begin{Thm}\label{maxorthogonal}
Let $M$ be a terminal $KQ$-module, and let
$T$ be a $\LL$-module in $\CC_M$
such that the following hold:
\begin{itemize}

\item[(i)]
$T$ is rigid;

\item[(ii)]
$T$ is a $\CC_M$-generator-cogenerator;

\item[(iii)]
$\gldim(\End_\LL(T)) \le 3$.

\end{itemize}
Then $T$ is $\CC_M$-maximal 1-orthogonal.
\end{Thm}

\begin{proof}
Let $X \in \CC_M$ with $\Ext_\LL^1(T,X) = 0$.
We have to show that $X \in \add(T)$.
By Corollary~\ref{rightapprox3}, there exists an exact sequence 
$
T'' \to T' \to X \to 0
$
with $T',T'' \in \add(T)$.
(For an arbitrary selfinjective torsion class of $\nil(\LL)$, the existence
of such an exact sequence is not known.)
Applying $\Hom_\LL(-,T)$ yields an exact sequence
$$
0 \to \Hom_\LL(X,T) \to \Hom_\LL(T',T) \to \Hom_\LL(T'',T) \to Z \to 0
$$
of $\End_\LL(T)$-modules.
Since $\gldim(\End_\LL(T)) \le 3$ we get
$\pdim(Z) \le 3$ and therefore 
$
\pdim(\Hom_\LL(X,T)) \le 1.
$
Let
$$
0 \to \Hom_\LL(T_2,T) \xrightarrow{G} \Hom_\LL(T_1,T) \xrightarrow{F}
\Hom_\LL(X,T) \to 0
$$
be a projective resolution of $\Hom_\LL(X,T)$.
Thus $T_1,T_2 \in \add(T)$.
Furthermore, we know by Corollary~\ref{yoneda1} that $F = \Hom_\LL(f,T)$ and
$G = \Hom_\LL(g,T)$ for some homomorphisms $f$ and $g$.
By Corollary~\ref{fullyfaithful2},
$$
0 \to X \xrightarrow{f} T_1 \xrightarrow{g} T_2 \to 0
$$
is a short exact sequence.
Since we assumed $\Ext_\LL^1(T,X) = 0$, we know that this sequence
splits.
Thus $X$ is isomorphic to a direct summand of $T_1$, and therefore
$X \in \add(T)$.
This finishes the proof.
\end{proof}

\begin{Thm}\label{Mcompletetilt}
Let $\CC$ be a selfinjective torsion class of $\nil(\LL)$.
If $T_1$ and $T_2$ are $\CC$-maximal rigid modules in $\CC$,
then $\Hom_\LL(T_2,T_1)$ is a classical tilting module over
$\End_\LL(T_1)$, and we have
$$
\End_{\End_\LL(T_1)}(\Hom_\LL(T_2,T_1)) \cong \End_\LL(T_2)^\op.
$$
\end{Thm}

\begin{proof}
Without loss of generality we assume 
$\Sigma(T_2) \ge \Sigma(T_1)$.
Set 
$$
T := \Hom_\LL(T_2,T_1)\ \text{ and }\ B:= \End_\LL(T_1).
$$
Let $f\df T_2 \to T_1'$ be a left $\add(T_1)$-approximation
of $T_2$.
Since $T_1$ is a cogenerator of $\CC$, we know that $f$ is a monomorphism.
Since $T_2$ is rigid we can use Lemma~\ref{mutation1} and get 
a short exact sequence
\begin{equation}\label{eq3}
0 \to T_2 \xrightarrow{f} T_1' \xrightarrow{g} T_1'' \to 0
\end{equation}
with $T_1',T_1'' \in \add(T_1)$.
This yields a projective resolution
\begin{equation}\label{eq5}
0 \to \Hom_\LL(T_1'',T_1) \xrightarrow{\Hom_\LL(g,T_1)} 
\Hom_\LL(T_1',T_1) \xrightarrow{\Hom_\LL(f,T_1)} \Hom_\LL(T_2,T_1)
\to 0.
\end{equation}
So the $B$-module $\Hom_\LL(T_2,T_1)$ has projective dimension 
at most one, which is the
first defining property of a classical tilting module.

Next, we show that $\Ext_B^1(T,T) = 0$.
Applying $\Hom_B(-,T)$ to Sequence (\ref{eq5}) yields
an exact sequence
\begin{multline}\label{eq4}
0 \to \End_B(T) \xrightarrow{F} 
\Hom_B(\Hom_\LL(T_1',T_1),T)
\xrightarrow{G} \Hom_B(\Hom_\LL(T_1'',T_1),T)\\
\to \Ext_B^1(T,T) \to 
\Ext_B^1(\Hom_\LL(T_1',T_1),T)
\end{multline}
where
$$
F = \Hom_B(\Hom_\LL(f,T_1),T) \text{ and }
G = \Hom_B(\Hom_\LL(g,T_1),T).
$$

\begin{Lem}\label{rigidtorigid1}
$\Ext_B^1(\Hom_\LL(T_1',T_1),T) = 0$.
\end{Lem}

\begin{proof}
This is clear, since $\Hom_\LL(T_1',T_1)$ is a projective $B$-module.
\end{proof}

\begin{Lem}\label{rigidtorigid2}
The map $G$ is surjective and 
$\Ext_B^1(T,T) = 0$.
\end{Lem}

\begin{proof}
One easily checks that the diagram
$$
\xymatrix{
\Hom_\LL(T_2,T_1') \ar[d]^{i_{T_2,T_1',T_1}} \ar[r]^{\Hom_\LL(T_2,g)} & 
\Hom_\LL(T_2,T_1'') \ar[d]^{i_{T_2,T_1'',T_1}}\\
\Hom_B(\Hom_\LL(T_1',T_1),T) \ar[r]^G & 
\Hom_B(\Hom_\LL(T_1'',T_1),T)
}
$$
is commutative.
The morphism $\Hom_\LL(T_2,g)$ is surjective, since $T_2$ is rigid.
Thus
$$
i_{T_2,T_1'',T_1} \circ \Hom_\LL(T_2,g) = G \circ i_{T_2,T_1',T_1}
$$
is surjective, since $i_{T_2,T_1'',T_1}$ is an isomorphism.
This implies that $G$ is surjective.
Now the result follows from Lemma~\ref{rigidtorigid1}.
\end{proof}

The number $\Sigma(T)$ of isomorphism classes of indecomposable direct summands
of $T$ is equal to $\Sigma(T_2)$.
We proved that $T$ is a partial tilting module over $B$.
This implies $\Sigma(T) \le \Sigma(T_1)$. 
By our assumption, $\Sigma(T_2) \ge \Sigma(T_1)$.
It follows that $\Sigma(T_1) = \Sigma(T_2)$.

Thus we proved that $T$ is a classical tilting module over $B$.
Now apply $\Hom_\LL(T_2,-)$ to Sequence (\ref{eq3}).
This yields a short exact sequence
$$
0 \to \End_\LL(T_2) \xrightarrow{\Hom_\LL(T_2,f)} \Hom_\LL(T_2,T_1')
\xrightarrow{\Hom_\LL(T_2,g)} \Hom_\LL(T_2,T_1'') \to 0.
$$

\begin{Lem}\label{rigidtorigid3}
There exists an anti-isomorphism of ring 
$$
\xi\df \End_\LL(T_2) \to \End_{\End_\LL(T_1)}(T)
$$
such that the diagram
$$
\xymatrix{
0 \ar[r] & \End_\LL(T_2) \ar[r]^-{\Hom_\LL(T_2,f)} \ar[d]^\xi &
\Hom_\LL(T_2,T_1') \ar[r]^{\Hom_\LL(T_2,g)} \ar[d]^{i_{T_2,T_1',T_1}} &
\Hom_\LL(T_2,T_1'') \ar[r] \ar[d]^{i_{T_2,T_1'',T_1}} & 0\\
0 \ar[r] & \End_B(T) \ar[r]^-{F} &
\Hom_B(\Hom_\LL(T_1',T_1),T) \ar[r]^G &
\Hom_B(\Hom_\LL(T_1'',T_1),T) \ar[r] & 0
}
$$
commutes and has exact rows.
\end{Lem}

\begin{proof}
Set $\xi(h)(h'') = h''h$ for all $h \in \End_\LL(T_2)$ and 
$h'' \in \Hom_\LL(T_2,T_1)$.
Now one easily checks that
$$
(F \circ \xi)(h) = (i_{T_2,T_1',T_1} \circ\Hom_\LL(T_2,f))(h)\df
h' \mapsto h'fh.
$$
\end{proof}

Lemma~\ref{rigidtorigid3} implies that
$
\End_{\End_\LL(T_1)}(T) \cong \End_\LL(T_2)^\op.
$
This finishes the proof of Theorem~\ref{Mcompletetilt}.
\end{proof}

\begin{Cor}
Let $\CC$ be a selfinjective torsion class of $\nil(\LL)$.
If $T_1$ and $T_2$ are $\CC$-maximal rigid $\LL$-modules, 
then $\Sigma(T_1) = \Sigma(T_2)$.
\end{Cor}

\begin{Cor}\label{maxcompl}
Let $\CC$ be a selfinjective torsion class of $\nil(\LL)$.
For a $\LL$-module $T$ the following are equivalent:
\begin{itemize}

\item
$T$ is $\CC$-maximal rigid;

\item
$T$ is $\CC$-complete rigid.

\end{itemize}
\end{Cor}


{\Large\section{A functor from $\CC_M$ to the cluster category $\CC_Q$}
\label{clustercat}}


\subsection{A triangle equivalence}
Assume in this section that $M$ is a terminal $KQ$-module
with $t_i(M) = 1$ for all $i$.
(Note that this assumption excludes the linearly oriented quiver
of Dynkin type $\A_n$.)
Thus
$$
M  = \bigoplus_{i=1}^n (I_i \oplus \tau(I_i))
$$
where $I_1,\ldots,I_n$ are the indecomposable injective $KQ$-modules.
By $\CC_Q$ 
\index{$\CC_Q$}
we denote the {\it cluster category} 
\index{cluster category}
associated to $Q$.
Cluster  categories were invented by
Buan, Marsh, Reineke, Reiten and Todorov \cite{BMRRT}.
Keller \cite{Ke}
proved that they are triangulated categories in a natural way.

\begin{Thm}
Under the assumptions above,
the categories ${\stCC}_{M}$ and $\CC_Q$ are triangle equivalent.
\end{Thm}

\begin{proof}
We proved already that $\stCC_M$ is an algebraic 2-Calabi-Yau category.
According to an important result by Keller and Reiten 
\cite{KR}, it is enough to
construct a $\CC_M$-maximal 1-orthogonal module $T$ in $\CC_M$ such that
the quiver of the stable endomorphism algebra $\End_{{\stCC}_M}(T)$ 
is isomorphic to $Q^\op$. Using Lemma~\ref{quiverT_M},
it is easy to check that the module $T_M$ we constructed in Section 
\ref{completerigid} has this property.
\end{proof}

The proof of Keller and Reiten's theorem is quite
involved and it does not seem to provide an explicit functor. 
Here we present 
an elementary construction of a $K$-linear functor 
$G\df \CC_M \to \CC_Q$
such that the kernel of $G$ consists precisely of the morphisms which 
factor through $\CC_M$-projective-injective modules.
Thus we obtain a $K$-linear equivalence $\underline{G}\df \stCC_M \to \CC_Q$. 
Note however that we do not discuss the possible triangulated structures of
$\stCC_M$ and $\CC_Q$.

\subsection{Derived categories of path algebras}\label{2sl:pre}

Let us review a few facts about the derived category of a path algebra
which we will use without further reference. 
This material can be found in Happel's book \cite{H2}. 
Write $\DD := \DD^b(\md(KQ))$ for the bounded derived category of
$\md(KQ)$.
Recall that 
$$
\DD = \bigvee_{i\in\Z}(\md(KQ))[i]
$$ 
since $KQ$ is 
hereditary, see also Figure~\ref{fig:D-F}. 
As usual, we identify $\md(KQ)$ with the full subcategory 
$\md(KQ)[0]$ of $\DD$.

If $I \in \md(KQ)$ is injective, then 
$\tau_\DD^{-1}(I) = (\nu^{-1}(I))[1]$ may be considered as a complex of
projective $KQ$-modules which is concentrated in degree $-1$. 
Here $\nu\df \md(KQ) \to \md(KQ)$ is the Nakayama functor,
see for example \cite{ASS,Ri1}.
More generally, if
$$
0 \to L \to I' \xrightarrow{\iota} I'' \to 0
$$
is an injective resolution of a $KQ$-module $L$, then 
$$
\tau_\DD^{-1}(L) = (\nu^{-1}(I') \xrightarrow{\nu^{-1}(\iota)} \nu^{-1}(I''))[1] 
$$
may be viewed as a complex of projectives concentrated in degree $-1$ and $0$.
This is essentially the same as saying that $\tau_\DD^{-1}$ is 
the right derived functor 
$$
\mathbf{R}\Hom_{KQ}({\rm D}KQ,-)[1] \cong 
\mathbf{R}\Hom_{KQ}({\rm D}KQ[-1],-).
$$
Here, we consider the injective cogenerator 
${\rm D}KQ := \Hom_K(KQ,K)$ of $\md(KQ)$
as a bimodule. 

In particular, if $L \in \md(KQ)$ 
has no projective direct summand, then
$\tau_\DD(L) = \tau_Q(L)$. 
This follows from the usual construction of the
Auslander-Reiten translation $\tau = \tau_Q$ in $\md(KQ)$.

\subsection{Cluster categories}
Let us review the construction of the cluster
category $\CC_Q$ as a $K$-linear category. 
It is by definition the orbit
category of $\DD$ by the action of the group 
$\ebrace{F}$ generated by the self-equivalence 
$F := \tau_\DD^{-1} \circ [1]$ of $\DD$. 
Keller \cite{Ke} proved that this is in fact a triangulated category.

Now, let $\F$ be the full subcategory of $\DD$ which consists of all objects
which are isomorphic to a complex 
$0 \to I' \to I'' \to 0$ 
of injective $KQ$-modules 
concentrated in degree $0$ and $1$. 
So each object $C$ in $\F$ is naturally
of the form $C_{{\rm inj},1} \oplus C_{{\rm mod},0}$ where $C_{{\rm mod},0}$ 
is isomorphic to a $KQ$-module
concentrated in degree $0$, and $C_{{\rm inj},1}$ 
is isomorphic to an injective $KQ$-module concentrated in degree $1$.

Thus the indecomposable objects in $\F$ are just the indecomposable
$KQ$-modules $L$ and the shifts $I_i[-1]$ of the indecomposable
injective $KQ$-modules $I_1,\ldots,I_n$.
(Recall that we identify $\md(KQ)$ and $\md(KQ)[0]$.)

Note that $F(\F)$ consists of those objects in $\DD$ which are isomorphic
to a complex
$0 \to P' \to P'' \to 0$ 
of projective $KQ$-modules 
concentrated in degree
$-2$ and $-1$. 
In fact, $F(I' \to I'') = (\nu^{-1}(I') \to \nu^{-1}(I''))[2]$ which
is a complex of projectives concentrated in degree $-2$ and $-1$, where
$\nu$ is again the Nakayama functor of $\md(KQ)$.

We conclude that we may consider $\F$ as a fundamental domain of the
action of the group $\ebrace{F}$ on $\DD$. 
Thus we can identify 
the objects of $\CC_Q$ and $\F$. 
Note that with this identification we
have
$$
\CC_Q(X,Y) = \DD(X,Y)\oplus \DD(X,F(Y))
$$ 
for $X,Y \in \F$. 
The composition is given by
$$
(\phi_0,\phi_1)\circ (\psi_0,\psi_1) = (\phi_0\psi_0, 
(F\phi_0)\psi_1+\phi_1\psi_0)
$$
Recall that for $M,N \in \md(KQ)$ we have $\Hom_\DD(M,N[i]) = 0$ 
unless $i \in \{ 0,1 \}$, see also Figure~\ref{fig:D-F}.
\begin{figure}
\setlength{\unitlength}{.8mm}
\begin{picture}(165,45)(-10,-10)
\linethickness{1pt}
\path(-1,0)(9,0)(14,20)(11,25)(-1,25) 
\put(2,12){\makebox(0,0)[b]{$\mathcal{I}[-1]$}}
\path(30,25)(14,25)(17,20)(12,0)(30,0)
\put(20,13){\makebox(0,0)[b]{$\cP$}}
\spline(32,32)(32,5)(35,-1)(46,-1)(48,5)(48,32)
\put(39,13){\makebox(0,0)[b]{$\RR$}}
\path(51,0)(61,0)(66,20)(63,25)(51,25) 
\put(55,13){\makebox(0,0)[b]{$\mathcal{I}$}}
\path(82,25)(66,25)(69,20)(64,0)(82,0)
\put(77,12){\makebox(0,0)[b]{$\cP[1]$}}
\spline(85,32)(85,5)(88,-1)(99,-1)(101,5)(101,32)
\put(93,12){\makebox(0,0)[b]{$\RR[1]$}}
\path(104,0)(114,0)(120,20)(116,25)(104,25) 
\put(108,12){\makebox(0,0)[b]{$\mathcal{I}[1]$}}
\path(139,25)(119,25)(123,20)(117,0)(139,0)
\put(129,12){\makebox(0,0)[b]{$\cP[2]$}}
\dottedline{2}(64,25)(67,20)(62,0)(58,-5)(11,-5)(8,0)(13,20)(10,25)(15,35)%
(58,35)(64,25)
\put(37,-10){\makebox(0,0)[b]{$\F$}}
\dottedline{2}(64,25)(70,35)(114,35)(121,25)(124,20)(118,0)(113,-5)(65,-5)%
(62,0)
\put(94,-11){\makebox(0,0)[b]{$F(\F)$}}
\end{picture}
\caption{$\DD^b(\md(KQ))$ and $\F$}
\label{fig:D-F}
\end{figure}

\subsection{Description of $\CC_M$}
Recall that the objects in $\CC_M$ are of the form
$X = (I'' \oplus \tau(I'),f)$, where $I'$ and $I''$ are 
injective $KQ$-modules 
and 
$$
f \in \Hom_{KQ}(I''\oplus \tau(I'),\tau(I'')\oplus \tau^2(I')).
$$
For obvious reasons we can and will identify $f$ with a homomorphism
$f\df \tau(I') \to \tau(I'')$. 
If $Y = (J''\oplus \tau(J'), g)$ is another object
in $\CC_M$, then we have
$$
\CC_M(X,Y) = \left\{ \left(\bsm \vph''& \widetilde{\vph}\\0 &\vph' \esm\right) 
\in
\Hom_{KQ}(I''\oplus \tau(I') ,J''\oplus \tau(J') ) \mid 
g \circ \vph' = \tau(\vph'') \circ f \right\}.
$$
Thus the diagram
$$
\xymatrix{
\tau(I') \ar[d]^f\ar[r]^{\vph'} & \tau(J') \ar[d]^g\\
\tau(I'')\ar[r]^{\tau(\vph'')} & \tau(J'')
}
$$
commutes.
Note that there is no condition on 
$\widetilde{\vph} \in \Hom_{KQ}(\tau(I'),J'')$.

\subsection{Description of the functor $G$}
Using the above conventions and notations we define $G\df \CC_M \to \CC_Q$
on  objects as
$$
G(X) := (0 \to I' \xrightarrow{\tau^{-1}(f)} I'' \to 0)\in\F.
$$
For a morphism $\vph\in\CC_M(X,Y)$ we define
$$
G(\vph) := \left( (\tau^{-1}(\vph'),\vph''), 
\tau^{-1}_\DD(\widetilde{\vph})\right)
$$
The first component, $(\tau^{-1}(\vph'),\vph'')$ is by the definition
of the homomorphisms in $\CC_M$ a homomorphism between complexes of 
injective modules.
So this is well defined. 
As for the second component consider the following diagram 
for morphisms $G(X)\to FG(Y)$ in $\DD$:
$$
\xymatrix@+1.0pc{
0 \ar[r]\ar[d] & \tau_D^{-1}(J')\ar[d]^{\tau_\DD^{-2}(g)}\\
I' \ar[r]^>>>>>>{\tau_\DD^{-1}(\widetilde{\vph})}\ar[d]_{\tau^{-1}(f)} &
\tau_\DD^{-1}(J'') \ar[d]\\
I'' \ar[r] & 0
}
$$

\begin{Thm}
The functor $G$ is an epivalence. 
For $\vph \in \CC_M(X,Y)$ we have
$G(\vph)=0$
if and only if there exists 
$$
\eta = 
\left(\bbm 
\eta_2 & 0\\ 
\widetilde{\eta} &\eta_1
\ebm\right)
\in 
\DD(I''\oplus\tau(I'),  \tau^{-1}_\DD(J'')\oplus J')
$$
such that
\begin{equation}\label{eq:phifact}
\vph = 
\left(\bbm\vph'' & \widetilde{\vph}\\ 0 & \vph' \ebm\right) =
\left(\bbm 
\tau^{-1}(g)\circ\widetilde{\eta}&\tau^{-1}(g)\circ\eta_1 +\tau(\eta_1)\circ f\\ 
0&\tau(\widetilde{\eta})\circ f
\ebm\right).
\end{equation}
Moreover, the condition \eqref{eq:phifact} is equivalent to the condition
that $\vph$ factors through a $\CC_M$-projective-injective module.
This implies that
$\underline{G}\df \stCC_M \to \CC_Q$ is an equivalence.
\end{Thm}

\begin{proof}
(a) 
Obviously, $G$ is dense. 
On morphisms, $G$ is surjective in the first 
component because $(\tau^{-1}(\vph'),\vph'')$ can be any morphism between
the two complexes of injectives which are concentrated in degree $0$
and $1$. 
Moreover, in the derived category $\DD$ homomorphisms between
bounded complexes of injectives are just given by morphisms of
complexes modulo homotopy.

(b) 
In order to see that $G$ is also surjective in the second component, we 
consider the standard triangles
$$
I' \xrightarrow{\tau^{-1}(f)} I'' \to G(X)[1] \to I'[1]
$$ 
and
$$
\tau_\DD^{-1}(J') \xrightarrow{\tau^{-2}(g)} \tau_\DD^{-1}(J'') \to 
FG(Y) \to \tau_\DD^{-1}(J')[1]
$$
in $\DD$.
For $i \not= j$ and $a,b \in \{',''\}$ we get
$\DD(I^{a}[i],\tau_\DD^{-1}(J^{b})[j]) = 0$.
From the corresponding long exact sequences we obtain the following commutative
diagram with exact rows and columns:
$$
\def\objectstyle{\scriptstyle}
\xymatrix{
\DD(I'',\tau_\DD^{-1}(J')) \ar[r]\ar[d] & 
\DD(I'',\tau_\DD^{-1}(J'')) \ar[r]\ar[d]  
& \DD(I'',FG(Y))\ar[r]\ar[d]   & 0 \ar[d]\\
\DD(I',\tau_\DD^{-1}(J'))\ar[r]\ar[d] & \DD(I',\tau_\DD^{-1}(J'')) 
\ar[r]\ar[d]  &
\DD(I',FG(Y))\ar[r]\ar[d]   & 0\ar[d]\\
\DD(G(X),\tau_\DD^{-1}(J')) \ar[r]\ar[d] & \DD(G(X),\tau_\DD^{-1}(J'')) 
\ar[r]\ar[d] &
\DD(G(X),FG(Y)) \ar[r]\ar[d] & \DD(G(X),\tau_\DD^{-1}(J')[1]\ar[d])\\
0 \ar[r] & 0 \ar[r] &\DD(I''[-1],FG(Y))\ar[r]  & 0
}
$$
Thus, $\DD(I''[1],FG(Y)) = 0 = \DD(G(X),\tau_\DD^{-1}(J')[1])$, 
and we conclude that
$$
\DD(G(X),FG(Y)) \cong \frac{\DD(I',\tau_\DD^{-1}(J''))}%
{(\tau_\DD^{-2}(g))\DD(I',\tau_\DD^{-1}(J')) + 
\DD(I'',\tau_\DD^{-1}(J''))(\tau^{-1}(f))}.
$$

(c) 
Our claim on the kernel of $G$ follows
from the end of steps (a) and (b), respectively.
Now one can use our results in Section~\ref{liftsection}
to describe the morphisms in $\CC_M$ which factor through 
$\CC_M$-projective-injectives.
It follows
that this is equivalent to
the description of the kernel of $G$ in \eqref{eq:phifact}.
\end{proof}

In practice, the functor $G\df \CC_M \to \CC_Q$ is (at least on
objects) easy to handle:
Take an indecomposable $KQ$-module $L$, and let
$$
0 \to L \to I' \xrightarrow{f} I'' \to 0
$$
be a minimal injective resolution of $L$.
Define a $\LL$-module
$$
\widetilde{L} := \left( I'' \oplus \tau(I'),
\left(\bsm 0 & \tau(f)\\0&0\esm\right) 
\right).
$$
Then we have $G(\widetilde{L}) = L$.
In particular, if $L = I_i$ is injective, then
$\widetilde{L} = (\tau(I_i),0)$.
Note that there is a short exact sequence
$$
0 \to (I'',0) \to \widetilde{L} \to (\tau(I'),0) \to 0.
$$
of $\LL$-modules.

Next, let $L := I_i[-1]$ be the $[-1]$-shift of an indecomposable
injective $KQ$-module $I_i$.
Set 
$\widetilde{L} := (I_i,0)$.
Again, we have $G(\widetilde{L}) = L$.

This describes the preimages of the indecomposable 
objects in $\CC_Q$ under the equivalence
$\underline{G}\df \stCC_M \to \CC_Q$.

\noindent
{\bf Example}:
Let $Q$ be the quiver
{\small
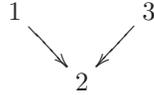
$$
\xymatrix@-0.8pc{
1 \ar[dr]&&3\ar[dl]\\
&2
}
$$
}\noindent
and let
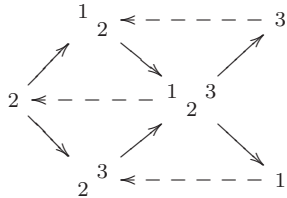
$$
\xymatrix@-1.0pc{
&{\bsm1\\&2\esm} \ar[dr] && 
{\bsm 3\esm}\ar@{-->}[ll]\\
{\bsm 2\esm} \ar[ur]\ar[dr] && 
{\bsm1&&3\\&2&\esm} \ar[ur]\ar[dr]\ar@{-->}[ll]\\
&{\bsm&3\\2\esm} \ar[ur] && 
{\bsm 1\esm} \ar@{-->}[ll]
}
$$
be the Auslander-Reiten quiver of $KQ$.
(The dotted arrows show how the Auslander-Reiten translation
$\tau$ acts on the non-projective indecomposable $KQ$-modules.)
Then
$$
0 \to \bsm 1\\&2\esm \to \bsm 1&&3\\&2\esm \xrightarrow{f} \bsm 3\esm \to 0
$$
is a minimal injective resolution of the $KQ$-module
$L := \bsm 1\\&2\esm$, where $f$ is just the obvious projection map.
It follows that
$$
\widetilde{L} = \left(\bsm 3\esm \oplus \bsm 2\esm,
\left( \bsm 0 & \tau(f) \\ 0 & 0 \esm\right) 
\right) = \bsm 2\\&3\esm.
$$
Note that $\tau(f)\df \bsm 2\esm \to \bsm 1\\&2\esm$ is
the obvious inclusion map.
(Here we are using the same notation as in Section~\ref{A2tilde}:
The numbers $1,2,3$ correspond to composition factors
of $KQ$-modules and $\LL$-modules, respectively.
For example $\bsm 2\\&3\esm$ is the 2-dimensional indecomposable
$\LL$-module with top $S_2$ and socle $S_3$.)

\newpage

{\Large\part{Mutations}
\label{part3}}



{\Large\section{Mutation of $\CC$-maximal rigid modules}\label{section11}}


\begin{Prop}\label{mutation2}
Let $T \oplus X$ be a basic rigid $\LL$-module 
such that the following
hold:
\begin{itemize}

\item $X$ is indecomposable;

\item $X \in \Cogen(T)$.

\end{itemize}
Then there exists a short exact sequence
$$
0 \to X \xrightarrow{f} T' \xrightarrow{g} Y \to 0
$$
such that the following hold:
\begin{itemize}

\item
$f$ is a minimal left $\add(T)$-approximation;

\item
$g$ is a minimal right $\add(T)$-approximation;

\item $T \oplus Y$ is basic rigid;

\item $Y$ is indecomposable and $X \not\cong Y$.

\end{itemize}
\end{Prop}

\begin{proof}
Let $f\df X \to T'$ be a minimal left $\add(T)$-approximation of $X$.
Since $X \in \Cogen(T)$, it follows that $f$ is a monomorphism.
Now copy the proof of \cite[Proposition 5.6]{GLSRigid}.
\end{proof}

In the situation of the above proposition, we call $\{ X,Y \}$ an
{\it exchange pair associated to} $T$, and we write
$$
\mu_X(T \oplus X) = T \oplus Y.
$$
We say that $T \oplus Y$ is the mutation of $T \oplus X$ in direction $X$.
The short exact sequence
\[
0 \to X \xrightarrow{f} T' \xrightarrow{g} Y \to 0
\]
is the {\it exchange sequence} starting in $X$ and ending in $Y$.

\begin{Prop}\label{mutation3}
Let $X$ and $Y$ be indecomposable rigid $\LL$-modules with 
$$
\dm \Ext_\LL^1(Y,X) = 1,
$$
and let
$$
0 \to X \xrightarrow{f} E \xrightarrow{g} Y \to 0
$$
be a non-split short exact sequence.
Then the following hold:
\begin{itemize}

\item[(i)]
$E \oplus X$ and $E \oplus Y$ are rigid
and $X,Y \notin \add(E)$.

\item[(ii)]
If we assume additionally that $T \oplus X$ and
$T \oplus Y$ are basic $\CC$-maximal rigid $\LL$-modules for some 
selfinjective torsion class $\CC$ of $\nil(\LL)$,
then 
$f$ is a minimal left $\add(T)$-approximation
and
$g$ is a minimal right $\add(T)$-approximation.

\end{itemize}
\end{Prop}

\begin{proof}
If $X$ and $Y$ are in some selfinjective torsion class $\CC$, then
$E \in \CC$, since $\CC$ is closed under extensions.
Now copy the proof of \cite[Proposition 5.7]{GLSRigid}.
\end{proof}

\begin{Cor}\label{cormutation3}
Let $\CC$ be a selfinjective torsion class of $\nil(\LL)$.
Let $\{X,Y\}$ be an exchange pair associated to some basic rigid $\LL$-module
$T$ such that $T \oplus X$ and $T \oplus Y$ are $\CC$-maximal rigid, 
and assume $\dm \Ext_\LL^1(Y,X) = 1$.
Then
\[
\mu_Y(\mu_X(T \oplus X)) = T \oplus X.
\]
\end{Cor}

\begin{proof}
Copy the proof of \cite[Corollary 5.8]{GLSRigid}.
\end{proof}


{\Large\section{Endomorphism algebras of $\CC$-maximal
rigid modules}\label{endosection}}


In this section,
let $\CC$ be a selfinjective torsion class of $\nil(\LL)$.
We denote by $I_\CC$ its $\CC$-projective generator-cogenerator.
We work mainly with 
basic rigid $\LL$-modules in $\CC$.
However, all our results on their endomorphism algebras are Morita
invariant, thus they 
hold for endomorphism algebras of arbitrary rigid $\LL$-modules in $\CC$.

Let $A$ be a $K$-algebra, and
let $M = M_1^{n_1} \oplus \cdots \oplus M_t^{n_t}$ 
be a finite-dimensional $A$-module, where the $M_i$
are pairwise non-isomorphic indecomposable
modules and $n_i \ge 1$.
As before let $S_i = S_{M_i}$ be the simple $\End_A(M)$-module
corresponding to $M_i$.
Then
$\Hom_A(M_i,M)$ is the indecomposable projective $\End_A(M)$-module with
top $S_i$.
The basic facts on the quiver $\GG_M$ of the endomorphism algebra
$\End_A(M)$ are collected in \cite[Section 3.2]{GLSRigid}.

\begin{Thm}[\cite{Ig}]\label{Igusa}
Let $A$ be a finite-dimensional $K$-algebra.
If $\gldim(A) < \infty$,
then the quiver of $A$ has no loops.
\end{Thm}

\begin{Prop}[{\cite[Proposition 3.11]{GLSRigid}}]\label{2cycleresult}
Let $A$ be a finite-dimensional $K$-algebra.
If $\gldim(A) < \infty$ and if the quiver of $A$ has a 2-cycle, 
then $\Ext_A^2(S,S) \not= 0$ 
for some simple $A$-module $S$.
\end{Prop}

\begin{Lem}[{\cite[Lemma 6.1]{GLSRigid}}]\label{noloopcrit2}
Let $\{ X,Y \}$ be an exchange pair associated to a
basic rigid $\LL$-module $T$.
Then the following are equivalent:
\begin{itemize}

\item
The quiver of $\End_\LL(T \oplus X)$ has no loop at $X$;

\item
Every non-isomorphism $X \to X$ factors through $\add(T)$;

\item
$\dm \Ext_\LL^1(Y,X) = 1$.

\end{itemize}
\end{Lem}

\begin{Lem}\label{simplesocle}
Let $T$ be a basic $\CC$-maximal rigid $\LL$-module.
If the quiver of $\End_\LL(T)$ has no loops, then
every indecomposable $\CC$-projective module has a simple
socle.
\end{Lem}

\begin{proof}
Let $P$ be an indecomposable $\CC$-projective module.
Let $h\df P \to T'$ be a minimal left $\add(T/P)$-approximation,
and set $X := P/\Ker(h)$.
Since $P$ is $\CC$-projective, $\Ker(h) \not= 0$.

Let $U$ be a non-zero submodule of $P$, and set $X' := P/U$.
Since $\CC$ is closed under factor modules, we get $X' \in \CC$.
By $p\df P \to X'$ we denote the canonical projection morphism.

We know that $T$ is a cogenerator of $\CC$.
Thus there exists a monomorphism
$$
\phi = \left[\bsm\phi_1\\\vdots\\\phi_m\\\theta\esm\right]\df
X' \to P^m \oplus T''
$$
with $T'' \in \add(T/P)$.
Since $U \not= 0$, none of the homomorphisms $\phi_i\df X' \to P$
is invertible.
In particular, none of the $\phi_i$ is an epimorphism.
The image of 
$$
\phi \circ p = 
\left[\bsm\phi_1 \circ p\\\vdots\\\phi_m \circ p\\\theta 
\circ p\esm\right]\df
P \to P^m \oplus T''
$$
is isomorphic to $X'$, and $\phi_i \circ p\df P \to P$ is not invertible
for all $i$.
Since the quiver of $\End_\LL(T)$ has no loops,
there exist homomorphisms 
$\phi_i'\df P \to T_i'$ and $\phi_i''\df T_i' \to P$ with
$T_i' \in \add(T/P)$ such that
$$
\phi_i \circ p = \phi_i'' \circ \phi_i'
$$
for all $i$.
Set
$$
\phi' = 
\left[\bsm\phi_1'\\\vdots\\\phi_m'\\\theta 
\circ p\esm\right]\df
P \to \left(\bigoplus_{i=1}^m T_i'\right) \oplus T''.
$$
It follows that $\phi \circ p = \phi'' \circ \phi'$ where
$$
\phi'' = \left[\bsm\phi_1''&&&\\&\ddots&&\\&&\phi_m''&\\&&&1_{T''}\esm\right].
$$
Thus
the image of $\phi'$ has at least the dimension of $X'$,
and we have 
$$
\Ker(\phi') \subseteq \Ker(\phi \circ p) = \Ker(p) = U.
$$
Now $h$ is a left $\add(T/P)$-approximation, thus
$\phi'$ factors through $h$.
Therefore $\dm X' \le \dm \Ima(h) = \dm X$.
It follows that $\Ker(h)$ must be simple.

Next, assume that $U_1$ and $U_2$ are simple submodules of $P$
with $U_1 \not= U_2$.
Thus there exists a monomorphism $P \to P/U_1 \oplus P/U_2$. 
From the above considerations we know that $P/U_1$ and $P/U_2$
are both in $\Cogen(T/P)$.
This implies that $P$ is in $\Cogen(T/P)$, a contradiction.
We conclude that $P$ has a simple socle.
\end{proof}

\begin{Prop}\label{gldim}
Let $T$ be a basic $\CC$-maximal rigid $\LL$-module.
If the quiver of $\End_\LL(T)$ has no loops, then
\[
\gldim(\End_\LL(T)) = 3.
\]
\end{Prop}

\begin{proof}
Set $B = \End_\LL(T)$.
By assumption, the quiver of $B$ has no loops.
It follows that $\Ext_B^1(S,S) = 0$ for all simple
$B$-modules $S$.
Let
$
T = T_1 \oplus \cdots \oplus T_r
$
with $T_i$ indecomposable for all $i$.
As before,
denote the simple $B$-module corresponding to $T_i$ by $S_{T_i}$.

Assume that $X = T_i$ is not $\CC$-projective.
Let $\{ X,Y \}$ be the exchange pair associated to $T/X$.
Note that $I_\CC \in \add(T/X)$.
This implies $X \in \Gen(T/X)$ and $X \in \Cogen(T/X)$.

By Lemma~\ref{noloopcrit2} we have $\dm \Ext_\LL^1(Y,X) = 1$.
Let
\[
0 \to X \to T' \to Y \to 0
\]
and
\[
0 \to Y \to T'' \to X \to 0
\]
be the corresponding non-split short exact sequences.
As in the proof of \cite[Proposition 6.2]{GLSRigid}
we obtain a minimal projective resolution
$$
0 \to \Hom_\LL(X,T) \to \Hom_\LL(T'',T) \to \Hom_\LL(T',T) \to \Hom_\LL(X,T) 
\to S_X \to 0. 
$$
In particular, $\pdim_B(S_X) = 3$.

Next, assume that $P = T_i$ is $\CC$-projective.
By Lemma~\ref{simplesocle} we know that $P$ has a simple socle, say $S$.
As in \cite[Proposition 9.4]{GLSRigid} one shows that $X := P/S$ is rigid.
Note also that $X \in \CC$.
Let $f\df X \to T'$ be a minimal left $\add(T/P)$-approximation.
It is easy to show that $f$ is injective. 
We get a short exact sequence
\[
0 \to X \xrightarrow{f} T' \to Y \to 0.
\]
It follows that
$Y \in \add(T)$.
The projection $\pi\df P \to X$ yields an exact sequence
\[
P \xrightarrow{h} T' \to Y \to 0
\]
where $h = f\pi$.
One can easily check that
$h$ is a minimal left $\add(T/P)$-approximation.
Applying $\Hom_\LL(-,T)$ to this sequence gives a projective
resolution
\[
0 \to \Hom_\LL(Y,T) \to \Hom_\LL(T',T) \xrightarrow{\Hom_\LL(h,T)} 
\Hom_\LL(P,T) \to S_P \to 0.
\]
This implies $\pdim(S_P) \le 2$.
For details we refer to 
the proof of \cite[Proposition 6.2]{GLSRigid}.
This finishes the proof.
\end{proof}

Recall the definition of a cluster torsion class of $\nil(\LL)$ (see
Section~\ref{defclustertorsionclass}).
The statements in the following theorem are presented in the order
in which we prove them.

\begin{Thm}\label{quivershape}
Let $\CC$ be a cluster torsion class of $\nil(\LL)$.
Let $T$ be a basic $\CC$-maximal rigid $\LL$-module, and set $B = \End_\LL(T)$.
Then the following hold:
\begin{itemize}

\item[(1)]
The quiver of $B$ has no loops;

\item[(2)]
$\gldim(B) = 3$;

\item[(3)]
For all simple $B$-modules $S$ we have
$\Ext_B^1(S,S) = 0$ and
$\Ext_B^2(S,S) = 0$;

\item[(4)]
The quiver of $B$ has no 2-cycles.

\end{itemize}
\end{Thm}

\begin{proof}
By Theorem~\ref{Mcompletetilt} we know that $\End_\LL(T_\CC)$ and
$\End_\LL(T)$ are derived equivalent, since every $\CC$-complete
rigid module is obviously $\CC$-maximal rigid.
Since the quiver of $\End_\LL(T_\CC)$ has no loops, 
Proposition~\ref{gldim} implies that $\gldim(\End_\LL(T_\CC)) = 3 < \infty$.
This implies
$\gldim(\End_\LL(T)) < \infty$.
Thus by Theorem~\ref{Igusa} the quiver of $\End_\LL(T)$ has no loops.
Then again Proposition~\ref{gldim} yields $\gldim(\End_\LL(T)) = 3$.
This proves (1) and (2).

Since the quiver of $B$ has no loops, we have $\Ext_B^1(S,S) = 0$
for all simple $B$-modules $S$.
Let $X$ be a direct summand of $T$ such that
$X$ is not $\CC$-projective.
In the proof of Proposition~\ref{gldim}, we constructed a projective resolution
\[
0 \to \Hom_\LL(X,T) \to \Hom_\LL(T'',T) \to \Hom_\LL(T',T) 
\to \Hom_\LL(X,T) \to S_X \to 0, 
\]
and we also know that $X \notin \add(T'')$.
Thus applying $\Hom_B(-,S_X)$ to this resolution yields
$\Ext_B^2(S_X,S_X) = 0$.
Next, assume
$P$ is an indecomposable $\CC$-projective direct summand of $T$. 
As in the proof of Proposition~\ref{gldim} we have a projective resolution
\[
0 \to \Hom_\LL(Y,T) \to \Hom_\LL(T',T) \xrightarrow{\Hom_\LL(h,T)} 
\Hom_\LL(P,T) \to S_P \to 0
\]
where $P \notin \add(T')$.
Since the module $T'$ projects onto $Y$, we conclude that
$P \notin \add(Y)$.
Applying $\Hom_B(-,S_P)$ to the above resolution of $S_P$ yields
$\Ext_B^2(S_P,S_P) = 0$.
This finishes the proof of (3).

We proved that for all simple $B$-modules $S$ we have
$\Ext_B^2(S,S) = 0$.
We also know that $\gldim(B) = 3 < \infty$.
Then it follows from Proposition~\ref{2cycleresult} that 
the quiver of $B$ cannot have 2-cycles.
Thus (4) holds. 
This finishes the proof.
\end{proof}

\begin{Cor}\label{corquivershape}
Let $\CC$ be a cluster torsion class of $\nil(\LL)$.
Let $T$ be a basic $\CC$-maximal rigid $\LL$-module, and let $X$ be
an indecomposable direct summand of $T$ which is not $\CC$-projective.
Let
$$
0 \to X \to T' \to Y \to 0
$$
be the corresponding exchange sequence starting in $X$.
Then the following hold:
\begin{itemize}

\item
We have $\dm \Ext_\LL^1(Y,X) = \dm \Ext_\LL^1(X,Y) = 1$, and 
the exchange sequence ending in $X$ is of the form
\[
0 \to Y \to T'' \to X \to 0
\]
for some $T'' \in \add(T/X)$;

\item
The simple $\End_\LL(T)$-module $S_X$ has a minimal projective resolution
of the form
\[
0 \to \Hom_\LL(X,T) \to \Hom_\LL(T'',T) \to \Hom_\LL(T',T) \to \Hom_\LL(X,T) 
\to S_X \to 0;
\]

\item
We have $\add(T') \cap \add(T'') = 0$.

\end{itemize}
\end{Cor}

\begin{proof}
Copy the proof of \cite[Corollary 6.5]{GLSRigid}.
\end{proof}

\begin{Thm}\label{threeequiv}
Let $M$ be a terminal $KQ$-module.
For a $\LL$-module $T$ in $\CC_M$ the following are equivalent:
\begin{itemize}

\item[(i)]
$T$ is $\CC_M$-maximal rigid;

\item[(ii)]
$T$ is $\CC_M$-complete rigid;

\item[(iii)]
$T$ is $\CC_M$-maximal 1-orthogonal.

\end{itemize}
\end{Thm}

\begin{proof}
Since $\CC_M$ is a selfinjective torsion class of $\nil(\LL)$, we know 
from Corollary~\ref{maxcompl} that
(i) and (ii) are equivalent.
Every $\CC_M$-maximal 1-orthogonal module is obviously
$\CC_M$-maximal rigid.
Vice versa, assume that $T$ is $\CC_M$-maximal rigid.
We know that there exists some $\CC_M$-complete rigid module
$T_M$ such that the quiver of $\End_\LL(T_M)$ has no loops.
By Theorem~\ref{quivershape} we get that $\gldim(\End_\LL(T)) = 3$.
Thus we can use Theorem~\ref{maxorthogonal} and 
get that $T$ is $\CC_M$-maximal 1-orthogonal.
\end{proof}

We conjecture that Theorem~\ref{threeequiv} can be generalized 
to the case where $\CC$ is a cluster torsion class of $\nil(\LL)$.
Note however that in this article (and also in \cite{GLSRigid}) 
we do not make any use of the fact 
that every $\CC_M$-maximal rigid module is $\CC_M$-maximal 1-orthogonal.

\begin{Prop}\label{CalabiYau}
Let $\CC$ be a cluster torsion class of $\nil(\LL)$. 
Let $T$ be a basic $\CC$-maximal rigid $\LL$-module, and let $X$ be
an indecomposable direct summand of $T$ which is not $\CC$-projective.
Set $B = \End_\LL(T)$.
Then for any simple $B$-module $S$ we have
$$
\dm \Ext_B^{3-i}(S_X,S) = \dm \Ext_B^i(S,S_X)
$$
where $0 \le i \le 3$.
\end{Prop}

\begin{proof}
Copy the proof of \cite[Proposition 6.6]{GLSRigid}.
\end{proof}


{\Large\section{From mutation of modules to 
mutation of matrices}\label{graphmutation}}


In this section, let $\CC$ be a cluster torsion class
of $\nil(\LL)$.

Let $T = T_1 \oplus \cdots \oplus T_r$ be a basic $\CC$-maximal rigid
$\LL$-module with $T_i$ indecomposable for all $i$.
Without loss of generality we assume that 
$T_{r-n+1}, \ldots, T_r$ are $\CC$-projective.
For $1 \le i \le r$ let $S_i = S_{T_i}$ be the simple $\End_\LL(T)$-module
corresponding to $T_i$.
The matrix
$$
C_T = (c_{ij})_{1 \le i,j \le r}
$$
\index{$C_T$ (where $T$ is maximal rigid)}\noindent
where 
$$
c_{ij} = \dm \Hom_{\End_\LL(T)}(\Hom_\LL(T_i,T),\Hom_\LL(T_j,T)) 
= \dm \Hom_\LL(T_j,T_i)
$$
is the Cartan matrix of the algebra $\End_\LL(T)$.

By Theorem~\ref{quivershape} we know that $\gldim(\End_\LL(T)) = 3$.
As in \cite[Section 7]{GLSRigid} this implies that
\begin{equation}\label{RTCT}
R_T =(r_{ij})_{1 \le i,j \le r} = C_T^{-t}
\end{equation}
\index{$R_T$ (where $T$ is maximal rigid)}\noindent
is the matrix of the Ringel form of $\End_\LL(T)$, where 
$$
r_{ij} = \bil{S_i,S_j} = \sum_{i=0}^3 (-1)^i \, 
\dm \Ext_{\End_\LL(T)}^i(S_i,S_j). 
$$

\begin{Lem}\label{rijequation}
Assume that $i \le r-n$ or $j \le r-n$.
Then the following hold:
\begin{itemize} 

\item
$r_{ij} = \dm \Ext_{\End_\LL(T)}^1(S_j,S_i) - 
\dm \Ext_{\End_\LL(T)}^1(S_i,S_j)$;

\item
$r_{ij} = - r_{ji}$;

\item
$
r_{ij} =
\begin{cases}
\text{number of arrows $j \to i$ in $\GG_T$}    & 
\text{if $r_{ij} > 0$},\\
-(\text{number of arrows $i \to j$ in $\GG_T$}) & 
\text{if $r_{ij} < 0$},\\
0                                   & \text{otherwise}.
\end{cases}
$

\end{itemize}
\end{Lem}

\begin{proof}
Copy the proof of \cite[Lemma 7.3]{GLSRigid}.
\end{proof}

Recall that $B(T) = B(\GG_T) = (t_{ij})_{1 \le i,j \le r}$ is the $r \times r$-matrix 
defined by
\[
t_{ij} = 
(\text{number of arrows $j \to i$ in $\GG_T$}) 
-(\text{number of arrows $i \to j$ in $\GG_T$}).
\]
Let 
$B(T)^\circ = (t_{ij})$ and $R_T^\circ = (r_{ij})$
be the $r \times (r-n)$-matrices obtained from
$B(T)$ and $R_T$, respectively, by deleting the last $n$ columns.
As a consequence of Lemma~\ref{rijequation} we get the following:

\begin{Cor}\label{corrijequation}
$R_T^\circ = B(T)^\circ$.
\end{Cor}

The dimension vector of the indecomposable projective 
$\End_\LL(T)$-module
$\Hom_\LL(T_i,T)$ is the $i$th column of the matrix $C_T$.

For $1 \le k \le r-n$ let
\begin{equation}\label{seq1}
0 \to T_k \to T' \to T_k^* \to 0
\end{equation}
and
\begin{equation}\label{seq2}
0 \to T_k^* \to T'' \to T_k \to 0
\end{equation}
be exchange sequences associated to the direct summand $T_k$
of $T$.
Keeping in mind the remarks in \cite[Section 3.2]{GLSRigid},
it follows from Lemma~\ref{rijequation} that
$$
T'  = \bigoplus_{r_{ik} > 0} T_i^{r_{ik}}
\;\;\;\; \text{ and }  \;\;\;\;
T'' = \bigoplus_{r_{ik} < 0} T_i^{-r_{ik}}.
$$
Set 
$$
T^* = \mu_{T_k}(T) = T_k^* \oplus T/T_k.
$$

For an $m \times m$-matrix $B$ and some
$k \in [1,m]$ we define an $m \times m$-matrix $S = S(B,k) = (s_{ij})$
by
\[
s_{ij} =
\begin{cases}
-\delta_{ij} + \dfrac{|b_{ij}|-b_{ij}}{2} & \text{if $i = k$},\\
\delta_{ij} & \text{otherwise}.
\end{cases}
\]
By $S^t$ we denote the transpose of the matrix $S = S(B,k)$.

Now let $S = S(R_T,k)$.
The proofs of the following proposition and its corollaries are 
identical to the ones in \cite{GLSRigid}.

\begin{Prop}\label{quivermutation3}
With the above notation we have
\[
C_{T^*} = SC_TS^t.
\]
\end{Prop}

\begin{Cor}\label{quivermutation4}
$R_{T^*} = S^tR_TS$.
\end{Cor}

\begin{Cor}\label{quivermutation5}
$R_{T^*}^\circ = \mu_k(R_T^\circ)$.
\end{Cor}

Now we combine Corollary~\ref{corrijequation} and Corollary 
\ref{quivermutation5} and obtain the following theorem:

\begin{Thm}\label{thmquivermutation}
$B(\mu_{T_k}(T))^\circ = \mu_k(B(T)^\circ)$.
\end{Thm}

In particular, applying Theorem~\ref{thmquivermutation} to the 
cluster torsion class $\CC_M$ we have proved Theorem~\ref{main4}.


{\Large\section{Mutations of clusters via dimension vectors}
\label{section15}}


Let $\CC$ be a cluster torsion class of $\nil(\LL)$ of rank $r$.
Let $T_\CC$ be a fixed basic $\CC$-maximal rigid module and
set $B=\End_{\LL}(T_{\CC})$. 
In this section we prove that every indecomposable rigid module
$X$ in $\CC$ is determined by the dimension vector ${\mathbf d}_X$ of the 
$B$-module $\Hom_\LL(X,T_{\CC})$.
If $\{X,Y\}$  is an exchange pair associated to
$U=U_1\oplus\cdots\oplus U_{r-1}$, 
we also give an easy combinatorial rule to calculate ${\mathbf d}_Y$ in terms
of ${\mathbf d}_X$ and the vectors ${\mathbf d}_{U_i}$.

\subsection{Dimension vectors of rigid modules}
Let $A$ be a finite-dimensional $K$-algebra, and assume that
$K$ is algebraically closed.
For $d \ge 1$ let $A^d$ be the free $A$-module of rank $d$.
Let $U$ be an $A$-module which is isomorphic to a submodule of $A^d$,
and set
$$
{\mathbf e} = \dimv(A^d) - \dimv(U).
$$
By $\md(A,{\mathbf e})$ we denote the affine variety of $A$-modules 
with dimension vector ${\mathbf e}$.

The {\it Richmond stratum} $\cS(U,A^d)$ is the 
subset of $\md(A,{\mathbf e})$ consisting of the modules $M$ such that
there exists a short exact sequence
$$
0 \to U \to A^d \to M \to 0.
$$

\begin{Thm}[{\cite[Theorem 1]{R}}]
The Richmond stratum
$\cS(U,A^d)$ is a smooth, irreducible, locally closed subset 
of $\md(A,{\mathbf e})$, and
$$
\dm \cS(U,A^d) = \dm \Hom_A(U,A^d) - \dm \End_A(U).
$$
\end{Thm}

\begin{Cor}\label{rigiddim}
Assume that $\gldim(A) < \infty$.
Let $M$ and $N$ be $A$-modules such
that the following hold:
\begin{itemize}

\item
$\dimv(M) = \dimv(N)$;

\item
$M$ and $N$ are rigid;

\item
$\pdim(M) \le 1$ and $\pdim(N) \le 1$;
\end{itemize}
Then $M \cong N$.
\end{Cor}

\begin{proof}
Let ${\mathbf d} = (d_1,\ldots,d_n)$ be the dimension vector of 
the modules $M$ and $N$,
and set $d = d_1 + \cdots + d_n$.
So there are epimorphisms $f\df A^d \to M$ and $g\df A^d \to N$.
Since the projective dimensions of $M$ and $N$ are at most one,
we get two short exact sequences
$$
0 \to P' \to A^d \to M \to 0 
\text{\;\;\; and \;\;\;}
0 \to P'' \to A^d \to N \to 0
$$
with $P'$ and $P''$ projective modules which have the same dimension
vector.

Since $\gldim(A) < \infty$, the Cartan matrix of $A$ is invertible.
Thus the dimension vectors of the indecomposable projective $A$-modules
are linearly independent.
These two facts yield that $P'$ and $P''$ are isomorphic.
Since $M$ and $N$ are rigid, the orbits $\orb_M$ and $\orb_N$ are
dense in the Richmond stratum
$\cS(P',A^d)$,  thus $\orb_M = \orb_N$ and therefore
$M \cong N$.
Here we use the fact that Richmond strata are irreducible.
\end{proof}

Corollary~\ref{rigiddim} is in some sense optimal as the following 
two examples show.
Let $Q$ be the quiver with two vertices $1$ and $2$, and two arrows
$a\df 1 \to 2$ and $b\df 2 \to 1$.

Let $A_1 = KQ/I_1$ where the ideal $I_1$ of the path algebra $KQ$ is
generated by the path $ba$.
Let $M = \bsm 1\\2 \esm$ and $N = \bsm 2\\1 \esm$.
Obviously, $\dimv(M) = \dimv(N) = (1,1)$ and $M \not\cong N$.
The following hold:
\begin{itemize}

\item
$M$ and $N$ are rigid;

\item
$\pdim(M) = 1$ and $\pdim(N) = 2$;

\item
$\gldim(A) = 2$.

\end{itemize}

Next, let $A_2 = KQ/I_2$ where the ideal $I_2$ of the path algebra $KQ$ is
generated by the paths $ab$ and $ba$.
Define the modules $M$ and $N$ as above.
Then the following hold:
\begin{itemize}

\item
$M$ and $N$ are rigid;

\item
$\pdim(M) = \pdim(N) = 0$;

\item
$\gldim(A) = \infty$.

\end{itemize}

\begin{Cor}
Let $X$ and $Y$ be indecomposable rigid modules in $\CC$.
If ${\mathbf d}_X={\mathbf d}_Y$ then $X \cong Y$.
\end{Cor}
\begin{proof}
Let $M=\Hom_\LL(X,T_\CC)$ and $N=\Hom_\LL(Y,T_\CC)$.
Since $M$ and $N$ are direct summands of tilting modules over $B$
they are rigid, and
by Corollary~\ref{pdimone} they have projective
dimension at most one. 
Hence by  Corollary~\ref{rigiddim} we have $M\cong N$.
Now applying Lemma~\ref{fullyfaithful3} and Corollary~\ref{fullyfaithful} we get
that $X \cong Y$.
\end{proof}

\subsection{Mutations via dimension vectors} 
We now explain how to calculate mutations of clusters via dimension vectors.
We start with some notation:
For 
${\mathbf d} = (d_1,\ldots,d_r)$ and ${\mathbf e} = (e_1,\ldots,e_r)$
in $\Z^r$ define
$$
\max\{{\mathbf d},{\mathbf e}\} := (f_1,\ldots,f_r)
$$
where $f_i = \max\{d_i,e_i\}$ for $1 \le i \le r$.
Set ${\rm Max}\{{\mathbf d},{\mathbf e}\} := {\mathbf d}$
if $d_i \ge e_i$ for all $i$.
In this case, we write
${\mathbf d} \ge {\mathbf e}$.
Of course, ${\rm Max}\{{\mathbf d},{\mathbf e}\} = {\mathbf d}$
implies
$\max\{{\mathbf d},{\mathbf e}\} = {\mathbf d}$.
By $|{\mathbf d}|$ we denote the sum of the entries of ${\mathbf d}$.

Let $\GG^*$ be a quiver 
as constructed in Section~\ref{preinjective}.
We assume that $\GG^*$ has $r$ vertices.
Now replace each vertex $i$ of $\GG^*$ by some
${\mathbf d}_i \in \Z^r$.
Thus we obtain a new quiver $(\GG^*)'$ whose 
vertices are elements in $\Z^r$.

For $k \not= (i,0)$ define
the mutation $\mu_{{\mathbf d}_k}((\GG^*)')$ of $(\GG^*)'$ at
the vertex ${\mathbf d}_k$ in two steps:
\begin{itemize}

\item[(1)]
Replace the vertex ${\mathbf d}_k$ of $(\GG^*)'$ by
$$
{\mathbf d}_k^* := -{\mathbf d}_k + 
\max\left\{\sum_{{\mathbf d}_k \to {\mathbf d}_i} {\mathbf d}_i,
\sum_{{\mathbf d}_j \to {\mathbf d}_k}{\mathbf d}_j \right\}
$$
where the sums are taken over all arrows in $(\GG^*)'$ which start,
respectively end in the vertex ${\mathbf d}_k$;

\item[(2)]
Change the arrows of $(\GG^*)'$ following Fomin and Zelevinsky's
quiver mutation rule for the vertex ${\mathbf d}_k$.

\end{itemize}

Thus starting with $(\GG^*,({\mathbf d}_i)_i)$ we can use iterated mutation
and obtain quivers whose vertices
are elements in $\Z^r$.
It is an important and interesting question, if these quivers
parametrize the seeds or clusters of the cluster algebra
$\cA(B(\GG^*)^\circ)$ associated to $\GG^*$, and if the elements
in $\Z^r$ appearing as vertices are in bijection with the cluster
variables of $\cA(B(\GG^*)^\circ)$.

For example,
if for each $i$ we choose ${\mathbf d}_i = -{\mathbf e}_i$, where ${\mathbf e}_i$ is
the $i$th canonical basis vector of $\Z^r$, then
the resulting vertices (i.e. elements in $\Z^r$) are the denominator
vectors of the cluster variables of
$\cA(B(\GG^*)^\circ)$,
compare \cite[Section 7, Equation (7.7)]{FZ3}.
(The variables attached to the vertices $(i,0)$ serve as 
(non-invertible) coefficients.
To obtain the denominator vectors as defined in \cite{FZ3} one
has to ignore the entries corresponding to these $n$ coefficients.)
It is an open problem, if these denominator vectors actually parametrize
the cluster variables of $\cA(B(\GG^*)^\circ)$.

We will show that for an appropriate choice of the initial vectors ${\mathbf d}_i$,
the quivers obtained by iterated mutation of $(\GG^*)'$ are
in bijection with the seeds and clusters of $\cA(B(\GG^*)^\circ)$.
All resulting vertices (including the ${\mathbf d}_i$) will be elements
in $\N^r$, and we will show that for our particular choice of
initial vectors, we can
use ``{\rm Max}'' instead of ``{\rm max}'' in the formula above.
(This holds for all iterated mutations.)

\begin{Prop}\label{dimcount}
Let $T$ and $R$ be $\CC$-maximal rigid $\LL$-modules, and assume
that $R$ is basic.
Let 
$$
\eta'\df
0 \to R_k \xrightarrow{f'} R' \xrightarrow{g'} R_k^* \to 0 
\text{\;\;\;and \;\;\;}
\eta''\df
0 \to R_k^* \xrightarrow{f''} R'' \xrightarrow{g''} R_k \to 0
$$
be the two exchange sequences associated to an
indecomposable direct summand $R_k$ of $R$ which is not 
$\CC$-projective.
Then 
$
\dm \Hom_\LL(R',T) \not= \dm \Hom_\LL(R'',T),
$
and we have
$$
\dimv(\Hom_\LL(R_k,T)) +  \dimv(\Hom_\LL(R_k^*,T)) = 
\max\{ \dimv(\Hom_\LL(R',T)), \dimv(\Hom_\LL(R'',T)) \}.
$$
Furthermore, the following are equivalent:
\begin{itemize}

\item[(i)]
$\eta'$ is $F^T$-exact;

\item[(ii)] 
$\dm \Hom_\LL(R',T) > \dm \Hom_\LL(R'',T)$;
\item[(iii)]
$\dimv(\Hom_\LL(R',T)) \ge \dimv(\Hom_\LL(R'',T))$.

\end{itemize}
\end{Prop}

\begin{proof}
Set $B = \End_\LL(T)$. 
By \cite[Lemma 2.2]{H3} we may assume without loss of 
generality that
$
\Ext_B^1(\Hom_\LL(R_k^*,T),\Hom_\LL(R_k,T)) = 0.
$
By Proposition~\ref{seslift},
\begin{align*}
1 = \dm \Ext_\LL^1(R_k^*,R_k) &\ge \dm \Ext_{F^T}^1(R_k^*,R_k)\\
&= \dm \Ext_B^1(\Hom_\LL(R_k,T),\Hom_\LL(R_k^*,T)) > 0.
\end{align*}
This implies $\Ext_\LL^1(R_k^*,R_k) = \Ext_{F^T}^1(R_k^*,R_k)$. 
Thus $\eta'$ is $F^T$-exact, and
$$
\eta\df 
0 \to \Hom_\LL(R_k^*,T) \xrightarrow{\Hom_\LL(g',T)} \Hom_\LL(R',T) 
\xrightarrow{\Hom_\LL(f',T)} \Hom_\LL(R_k,T) \to 0
$$
is a (non-split) short exact sequence.
If we apply $\Hom_\LL(-,T)$ to $\eta''$,
we obtain an exact sequence
$$
0 \to \Hom_\LL(R_k,T) \xrightarrow{\Hom_\LL(g'',T)} \Hom_\LL(R'',T) 
\xrightarrow{\Hom_\LL(f'',T)} \Hom_\LL(R_k^*,T).
$$
Now $\Hom_\LL(f'',T)$ cannot be an epimorphism, since that would 
yield a non-split extension and we know that 
$\Ext_B^1(\Hom_\LL(R_k^*,T),\Hom_\LL(R_k,T)) = 0$.
Thus for dimension reasons we get
$\dm \Hom_\LL(R',T) > \dm \Hom_\LL(R'',T)$.
Using the functors $\Hom_B(P,-)$ where $P$ runs through the
indecomposable projective $B$-modules, it also follows that
$\dimv(\Hom_\LL(R',T)) > \dimv(\Hom_\LL(R'',T))$.
Finally, the formula for dimension vectors follows from the
exactness of $\eta$.
\end{proof}

Proposition~\ref{dimcount} yields
an easy combinatorial rule for the mutation of $\CC$-maximal rigid
modules.
Let 
$
T = T_1 \oplus \cdots \oplus T_r
$ 
be a $\CC$-maximal rigid $\LL$-module.
We assume that $T_{r-n+1},\ldots,T_r$
are $\CC$-projective.
For $1 \le i \le r$ let
$
{\mathbf d}_i := \dimv(\Hom_\LL(T_i,T_\CC)).
$

As before, let $\GG_T$ be the quiver of $\End_\LL(T)$.
The vertices of $\GG_T$ are labelled by the modules $T_i$.
For each $i$ we replace the vertex labelled by $T_i$ by the
dimension vector ${\mathbf d}_i$.
The resulting quiver is denoted by $\GG_T'$.

For $k \in [1,r-n]$ let
$$
0 \to T_k \to T' \to T_k^* \to 0
\text{\;\;\; and \;\;\;}
0 \to T_k^* \to T'' \to T_k \to 0
$$
be the two resulting exchange sequences.
We can now easily compute the dimension vector of the
$\End_\LL(T_\CC)$-module $\Hom_\LL(T_k^*,T_\CC)$, namely 
Proposition~\ref{dimcount} yields that
$$
{\mathbf d}_k^* := \dimv(\Hom_\LL(T_k^*,T_\CC)) = 
\begin{cases}
- {\mathbf d}_k + 
\sum_{{\mathbf d}_k \to {\mathbf d}_i}{\mathbf d}_i & \text{if
$\sum_{{\mathbf d}_k \to {\mathbf d}_i} |{\mathbf d}_i| >
\sum_{{\mathbf d}_j \to {\mathbf d}_k} |{\mathbf d}_j|$},\\
- {\mathbf d}_k +  \sum_{{\mathbf d}_j \to {\mathbf d}_k}{\mathbf d}_j & 
\text{otherwise},
\end{cases}
$$
where the sums are taken over all arrows in $\GG_T'$ which start,
respectively end in the vertex ${\mathbf d}_k$.
More precisely, we have 
\begin{equation}\label{denform}
{\mathbf d}_k^* = -{\mathbf d}_k + 
\max\left\{\sum_{{\mathbf d}_k \to {\mathbf d}_i} {\mathbf d}_i,
\sum_{{\mathbf d}_j \to {\mathbf d}_k}{\mathbf d}_j \right\}
\end{equation}
and we know that
\begin{equation}\label{denform2}
\max\left\{\sum_{{\mathbf d}_k \to {\mathbf d}_i} {\mathbf d}_i,
\sum_{{\mathbf d}_j \to {\mathbf d}_k}{\mathbf d}_j \right\}
=
{\rm Max}\left\{\sum_{{\mathbf d}_k \to {\mathbf d}_i}{\mathbf d}_i,
\sum_{{\mathbf d}_j \to {\mathbf d}_k}{\mathbf d}_j \right\}.
\end{equation}
%

\subsection{Example}\label{ex14.3}
Let $M$ be a terminal $KQ$-module such that 
$\GG_M$ is the quiver  
{\small
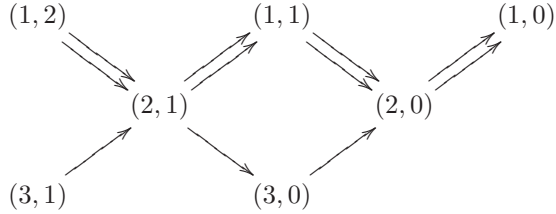
$$
\xymatrix@-0.5pc{ 
(1,2) \ar@<0.5ex>[dr]\ar@<-0.5ex>[dr] &&
(1,1) \ar@<0.5ex>[dr]\ar@<-0.5ex>[dr] && (1,0)
\\
&(2,1) \ar@<0.5ex>[ur]\ar@<-0.5ex>[ur]\ar[dr] && 
(2,0) \ar@<0.5ex>[ur]\ar@<-0.5ex>[ur]
\\
(3,1) \ar[ur] && 
(3,0) \ar[ur]
}
$$
}\noindent
which appeared already in Section~\ref{preinjective}.
Let 
$
T_M = T_1 \oplus \cdots \oplus T_7.
$ 
As always we assume without loss of generality that
$T_5,T_6,T_7$ are $\CC_M$-projective.
The following picture shows the quiver $\GG_{T_M}'$.
Its vertices are the dimension vectors of the
$\End_\LL(T_M)$-modules $\Hom_\LL(T_i,T_M)$.
These dimension vectors can be constructed easily 
by standard calculations inside the mesh category
of $\N Q^\op$, see also Section~\ref{dimvectors2}.
The dimension vectors associated to the indecomposable $\CC_M$-projectives
are labelled in red colour.
$$
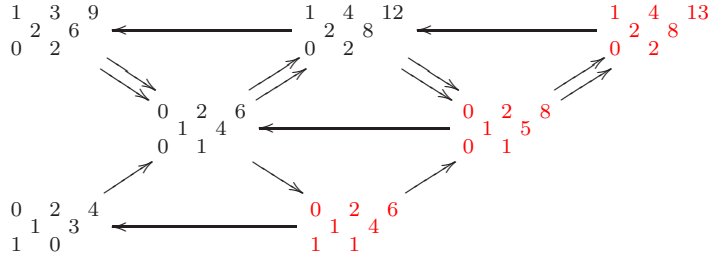
\xymatrix@-1pc{
{\bsm1&&3&&9\\&2&&6\\0&&2\esm} \ar@<0.5ex>[dr] \ar@<-0.5ex>[dr]
&& {\bsm1&&4&&12\\&2&&8\\0&&2\esm} \ar@<0.5ex>[dr] \ar@<-0.5ex>[dr]\ar[ll]&&
{\red\bsm1&&4&&13\\&2&&8\\0&&2\esm} \ar[ll]\\
&{\bsm0&&2&&6\\&1&&4\\0&&1\esm} \ar@<0.5ex>[ur]\ar@<-0.5ex>[ur]\ar[dr]&& 
{\red\bsm0&&2&&8\\&1&&5\\0&&1\esm} \ar[ll]\ar@<0.5ex>[ur]\ar@<-0.5ex>[ur]\\
{\bsm0&&2&&4\\&1&&3\\1&&0\esm}\ar[ur] &&
{\red\bsm0&&2&&6\\&1&&4\\1&&1\esm}\ar[ll]\ar[ur]
}
$$
Compare this also to the example in Section~\ref{example}.

Now let us mutate the $\LL$-module $T_k$ where
$$
\dimv(\Hom_\LL(T_k,T_M)) = {\bsm1&&4&&12\\&2&&8\\0&&2\esm}.
$$
We have to look at all arrows starting and ending in the corresponding
vertex of $\GG_{T_M}'$, and add up the entries of the 
attached dimension vectors, as explained in the previous section.
Since 
$$
\left|{\red\bsm1&&4&&13\\&2&&8\\0&&2\esm}\right| 
+ 2 \cdot \left|{\bsm0&&2&&6\\&1&&4\\0&&1\esm}\right| =  58 > 57 =
\left|{\bsm1&&3&&9\\&2&&6\\0&&2\esm}\right| + 
2 \cdot \left|{\red\bsm0&&2&&8\\&1&&5\\0&&1\esm}\right|,
$$ 
we get
$$
\dimv(\Hom_\LL(T_k^*,T_M)) = 
{\red\bsm1&&4&&13\\&2&&8\\0&&2\esm} 
+ 2 \cdot {\bsm0&&2&&6\\&1&&4\\0&&1\esm} 
-  {\bsm1&&4&&12\\&2&&8\\0&&2\esm}
= {\bsm0&&4&&13\\&2&&8\\0&&2\esm}
$$
and the quiver $\GG_{\mu_{T_k}(T_M)}'$ looks as follows:
$$
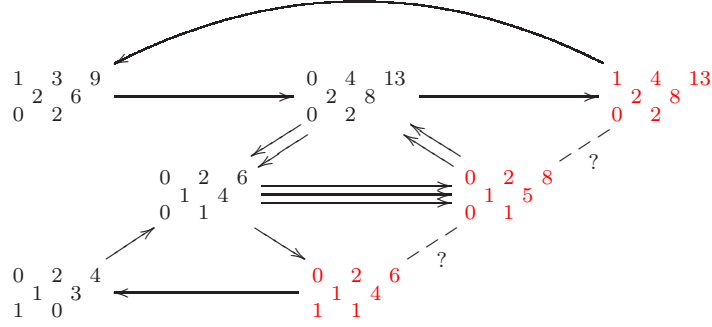
\xymatrix@-1pc{
&&&&\\
{\bsm1&&3&&9\\&2&&6\\0&&2\esm} \ar[rr]
&& {\bsm0&&4&&13\\&2&&8\\0&&2\esm} \ar@<0.5ex>[dl] \ar@<-0.5ex>[dl]\ar[rr]&&
{\red\bsm1&&4&&13\\&2&&8\\0&&2\esm} \ar@/_3pc/[llll]\\
&{\bsm0&&2&&6\\&1&&4\\0&&1\esm} \ar[dr]
\ar@<0.7ex>[rr]\ar@<-0.7ex>[rr]\ar[rr] && 
{\red\bsm0&&2&&8\\&1&&5\\0&&1\esm} \ar@{--}_{?}[ur]
\ar@<0.5ex>[ul]\ar@<-0.5ex>[ul]\\
{\bsm0&&2&&4\\&1&&3\\1&&0\esm}\ar[ur] &&
{\red\bsm0&&2&&6\\&1&&4\\1&&1\esm}\ar[ll]\ar@{--}_{?}[ur]
}
$$
Note that we cannot control how the arrows between vertices corresponding
to the three indecomposable
$\CC_M$-projectives
behave under mutation.
But this does not matter, because these arrows are not needed for the
mutation of seeds and clusters.
In the picture, we indicate the missing information by lines
of the form $\xymatrix{\ar@{--}[r]&}$.
This process can be iterated, and our theory says that each of the resulting
dimension vectors determines uniquely a cluster variable.

\subsection{Characterization of $Q$-split exact sequences}
\label{Qsplit}
In this section, let
$M = M_1 \oplus \cdots \oplus M_r$ be a terminal $KQ$-module.
We need the following result.

\begin{Lem}\label{lemma15.2}
Let $N_1,N_2 \in \add(M)$.
If
$\dm \Hom_{KQ}(N_1,M_i) = \dm \Hom_{KQ}(N_2,M_i)$ for all
$1 \le i \le r$, then $N_1 \cong N_2$.
\end{Lem}

\begin{proof}
For $i = 1,2$
we have  $\Hom_{KQ}(N_i,N) = 0$
for all indecomposable $KQ$-modules $N$ with $N \not\in \add(M)$.
It is a well known result by Auslander that
for any finite-dimensional algebra $A$ the numbers
$\dm \Hom_A(X,Z)$, where $Z$ runs through all finite-dimensional
indecomposable $A$-modules, determine a finite-dimensional
$A$-module $X$ uniquely up to isomorphism.
Applying this to $X = N_i$ yields the result.
\end{proof}

As before let 
$
\pi_Q\df \md(\LL) \to \md(KQ)
$
be the restriction functor, which is obviously exact.
Let 
$$
T_M = \bigoplus_{i=1}^n \bigoplus_{a=0}^{t_i} T_{i,a}
$$
be the $\CC_M$-complete rigid module we constructed before.
Set $B := \End_\LL(T_M)$.

We know that the contravariant 
functor $\Hom_\LL(-,T_M)$ yields an embedding
$$
\CC_M \to \md(B).
$$
If $X \in \CC_M$,
then the entries of the dimension vector $\dimv_B(\Hom_\LL(X,T_M))$ are
$$
\dm \Hom_B(\Hom_\LL(T_{i,a},T_M),\Hom_\LL(X,T_M)) =
\dm \Hom_{KQ}(\pi_Q(X),\tau^a(I_i))
$$
where $1 \le i \le n$ and $0 \le a \le t_i$, compare 
Section~\ref{dimvectors2}.
This together with Lemma~\ref{lemma15.2} yields 
the following result:

\begin{Lem}\label{lemma15.1}
For $\LL$-modules $X,Y \in \CC_M$ the following
are equivalent:
\begin{itemize}

\item
$\pi_Q(X) \cong \pi_Q(Y)$;

\item
$\dimv_B(\Hom_\LL(X,T_M)) = \dimv_B(\Hom_\LL(Y,T_M))$.

\end{itemize}
\end{Lem}

A short exact sequence 
$
\eta\df 0 \to X \to Y \to Z \to 0
$
of $\LL$-modules is called $Q$-{\it split} 
\index{$Q$-split short exact sequence}
if the short exact sequence
$
0 \to \pi_Q(X) \to \pi_Q(Y) \to \pi_Q(Z) \to 0
$
splits.

\begin{Prop}\label{QTexact}
For a short exact sequence 
$\eta\df 0 \to X \to Y \to Z \to 0$ of $\LL$-modules in $\CC_M$
the following are equivalent:
\begin{itemize}

\item
$\eta$ is $F^{T_M}$-exact;

\item
$\eta$ is $Q$-split.

\end{itemize}
\end{Prop}

\begin{proof}
Clearly, $\eta$ is $F^{T_M}$-exact if and only if 
$$
\dimv_B(\Hom_\LL(X,T_M)) + \dimv_B(\Hom_\LL(Z,T_M)) =
\dimv_B(\Hom_\LL(Y,T_M)).
$$
By Lemma~\ref{lemma15.1} this is equivalent to 
$\pi_Q(X) \oplus \pi_Q(Z) \cong \pi_Q(Y)$.
This is the case if and only if $\eta$ is $Q$-split.
\end{proof}

Proposition~\ref{QTexact} singles out a distinguished class of short exact
sequences of modules over preprojective algebras.
We believe that these sequences are important and deserve much attention.

For each indecomposable direct summand $(I_{i,[a,t_i]},e_{i,[a,t_i]})$
of $T_M$ we know its projection to $\md(KQ)$, namely
$$
\pi_Q(I_{i,[a,t_i]},e_{i,[a,t_i]}) = I_{i,[a,t_i]} = 
\bigoplus_{j=a}^{t_i} \tau^j(I_i).
$$
Using the mesh category of $\II_Q$ we can compute 
$\dimv_B(\Hom_\LL((I_{i,[a,t_i]},e_{i,[a,t_i]}),T_M))$.
Thus, 
combining Propositions~\ref{dimcount} and \ref{QTexact} 
we can inductively determine the $KQ$-module $\pi_Q(R)$ for each
cluster monomial $\delta_R$ in $\RR(\CC_M,T_M)$.

\subsection{Example}
Let $\LL$ be of Dynkin type $\A_3$.
Then the short exact sequences
$$
\eta'\df
0 \to \bsm1&&3\\&2&\esm \to \bsm&2&\\1&&3\\&2&\esm \to 
\bsm2\esm \to 0
\text{\;\;\; and \;\;\;}
\eta''\df
0 \to \bsm2\esm \to \bsm1&\\&2\esm \oplus \bsm&3\\2&\esm \to 
\bsm1&&3\\&2&\esm \to 0
$$
are exchange sequences in $\md(\LL)$.
There are four Dynkin quivers of type $\A_3$.
In each case, we determine if $\eta'$ or $\eta''$ is
$Q$-split:

\begin{center}
\begin{tabular}{|c|c|c|}
\hline
$Q$ & $\eta'$ & $\eta''$\\ \hline
$1 \xleftarrow{} 2 \xleftarrow{} 3$   & - & $Q$-split \\ \hline
$1 \xrightarrow{} 2 \xrightarrow{} 3$ & - & $Q$-split \\ \hline
$1 \xleftarrow{} 2 \xrightarrow{} 3$  & - & $Q$-split \\ \hline
$1 \xrightarrow{} 2 \xleftarrow{} 3$  & $Q$-split & - \\ \hline
\end{tabular}
\end{center}


{\Large\section{The algebra $\End_\LL(T_M)$ is quasi-hereditary}
\label{quasihereditary}}


\subsection{The partial ordering of tilting modules}\label{tiltingord}
Let $\T_A$ (resp. $\T_A^{\rm cl}$)
be a set of representatives of
isomorphism classes of all basic tilting modules (resp. basic classical
tilting modules) 
over $A$.

For a tilting module $T \in \T_A$ let
\begin{align*}
T^\perp &:= \left\{ Y \in \md(A) \mid \Ext_A^i(T,Y) = 0
\text{ for all } i \ge 1 \right\},\\
{^\perp}T &:= \left\{ X \in \md(A) \mid \Ext_A^i(X,T) = 0
\text{ for all } i \ge 1 \right\}.
\end{align*}
From now on we use $\perp$ only in this sense, so there
is no danger of confusing it with the notation in Section~\ref{perp}.

For $R,T \in \T_A$ define $R \le T$ if
$R^\perp \subseteq T^\perp$.
Thus $(\T_A,\le)$ and also $(\T_A^{\rm cl},\le)$ 
become partially ordered sets.
It follows that there is a unique maximal element, namely we have
$T \le {_A}A$
for all $T \in \T_A$.
Minimal elements need not exist in general.

Riedtmann and Schofield \cite{RS2} define a quiver ${\mathcal K}_A$ as follows:
The vertices of ${\mathcal K}_A$ are the elements in $\T_A$, and
there is an arrow $T \to R$ if and only if the following hold:
\begin{itemize}

\item[(i)]
$T = N \oplus X$ and $R = N \oplus Y$ with $X$ and $Y$ indecomposable 
and $X \not\cong Y$;

\item[(ii)]
There exists a short exact sequence
$$
0 \to X \to N' \to Y \to 0
$$
with $N' \in \add(N)$.

\end{itemize}
In this case, we have $\Ext_A^1(X,Y) = 0$.
Here are some known facts:
\begin{itemize}

\item[(a)]
If $T$ and $R$ are basic tilting modules satisfying (i), then
there is an arrow $T \to R$ or $R \to T$ in ${\mathcal K}_A$;

\item[(b)]
The quiver
${\mathcal K}_A$ is the Hasse quiver of the partially order set
$(\T_A,\le)$, see \cite{HU3};

\item[(c)]
If there is an arrow $T \to R$ in ${\mathcal K}_A$, then $R \le T$ and
$\pdim(R) \ge \pdim(T)$;

\item[(d)]
If $R \in T^\perp$, then $R \le T$, see
\cite[Lemma 2.1, (a)]{HU4} and  
\cite[Proof of Theorem 2.1]{HU3}.

\end{itemize}
%

\subsection{Quasi-hereditary algebras}\label{qhalgebras}
Let $A$ be a finite-dimensional algebra.
By $P_1,\ldots,P_r$ and $Q_1,\ldots,Q_r$ and
$S_1,\ldots,S_r$ we denote the indecomposable 
projective, indecomposable injective and simple $A$-modules,
respectively,
where $S_i = {\rm top}(P_i)  = \soc(Q_i)$. 

For a class $\U$ of $A$-modules let $\F(\U)$ 
\index{$\F(\U)$ (class of $\U$-filtered modules)}
be the class of all $A$-modules 
$X$ which
have a filtration 
$$
X = X_0 \supseteq X_1 \supseteq \cdots \supseteq X_t = 0
$$
of submodules 
such that all factors $X_{j-1}/X_j$ belong to $\U$ for all $1 \le j \le t$.
Such a filtration is called a $\U$-{\it filtration} of $X$.
We call these modules the $\U$-{\it filtered modules}.

Let $\Delta_i$ be the largest factor module of $P_i$ in 
$\F(S_1,\ldots,S_i)$,
and set 
$$
\Delta = \{ \Delta_1,\ldots,\Delta_r \}.
$$
\index{$\Delta$ (set of standard modules)}\noindent
The modules $\Delta_i$ are called {\it standard modules}.
\index{standard module}
The algebra $A$ is called {\it quasi-hereditary} 
\index{quasi-hereditary algebra}
if 
$\End_A(\Delta_i) \cong K$ for all $i$, and if
${_A}A$ belongs to $\F(\Delta)$.
Quasi-hereditary algebras first occured in Cline, Parshall
and Scott's \cite{CPS} study of highest weight categories. 

Note that the definition of a quasi-hereditary algebra 
depends on the chosen ordering of the simple modules.
If we reorder them, it could happen that our algebra 
is no longer quasi-hereditary.

Now assume $A$ is a quasi-hereditary algebra, and let $\F(\Delta)$ be the
subcategory of $\Delta$-filtered $A$-modules.
For $X \in \F(\Delta)$ let $[X:\Delta_i]$ be the number of times
that $\Delta_i$ occurs as a factor in a $\Delta$-filtration
of $X$.
Then 
$$
\dimv_\Delta(X) = ([X:\Delta_1],\ldots,[X:\Delta_r]) \in \N^r
$$ 
is the $\Delta$-{\it dimension vector}
\index{$\Delta$-dimension vector}
of $X$.

Let $\nabla_i$ be the largest submodule of $Q_i$ in $\F(S_1,\ldots,S_i)$,
and let 
$$
\nabla = \{ \nabla_1,\ldots,\nabla_n\}.
$$
The modules $\nabla_i$ are called {\it costandard modules}.

Let $A$ be a quasi-hereditary algebra.
The following results (and the missing definitions) 
can be found in \cite{Ri5,Ri4}:
\begin{itemize}

\item[(i)]
There is a unique (up to isomorphism) basic tilting module
$T^{\Delta\cap\nabla}$ 
\index{$T^{\Delta\cap\nabla}$}
over $A$ such that 
$$
\add(T^{\Delta\cap\nabla}) = \F(\Delta) \cap \F(\nabla).
$$

\item[(ii)]
We have $\F(\Delta) = {^\perp}(T^{\Delta\cap\nabla})$ and
$\F(\nabla) = (T^{\Delta\cap\nabla})^\perp$.

\item[(iii)]
$\F(\Delta)$ is closed under extensions and
under direct summands.

\item[(iv)]
$\F(\Delta)$ is a resolving and functorially finite subcategory
of $\md(A)$.

\item[(v)]
$\F(\Delta)$ has Auslander-Reiten sequences.

\item[(vi)]
$[P_i:\Delta_j] = [\nabla_j:S_i]$ for all $1 \le i,j \le r$,
where $[\nabla_j:S_i]$ is the Jordan-H\"older multiplicity of
$S_i$ in $\nabla_j$.

\item[(vii)]
If $X \in \F(\Delta)$, then
$[X:\Delta_i] = \dm \Hom_A(X,\nabla_i)$ for all $i$.

\item[(viii)]
$\Hom_A(\Delta_i,\Delta_j) = 0$ for all $i > j$.

\item[(ix)]
$\Ext_A^1(\Delta_i,\Delta_j) = 0$ for all $i \ge j$.

\item[(x)]
The $\F(\Delta)$-projective modules are the projective $A$-modules.
The $\F(\nabla)$-injective modules are the injective $A$-modules.

\item[(xi)]
The $\F(\Delta)$-injective modules are the modules
in $\add(T^{\Delta\cap\nabla})$.
The $\F(\nabla)$-projective modules are the modules
in $\add(T^{\Delta\cap\nabla})$.

\item[(xii)]
If $\Ext_A^1(X,\nabla_i) = 0$ for all $i$, then
$X \in \F(\Delta)$.
Similarly, if $\Ext_A^1(\Delta_i,Y) = 0$ for all $i$, then
$Y \in \F(\nabla)$.

\end{itemize}
The module $T^{\Delta\cap\nabla}$ 
is called the {\it characteristic tilting module} of $A$.
\index{characteristic tilting module}
In general, $T^{\Delta\cap\nabla}$ is not a classical tilting module.
The endomorphism algebra $\End_A(T^{\Delta \cap \nabla})$ is
called the {\it Ringel dual} of $A$.
It is again a quasi-hereditary algebra in a natural way, see \cite{Ri5}.

\subsection{$\GG_M$-adapted orderings}\label{Madapted}
Let $M = M_1 \oplus \cdots \oplus M_r$ be a terminal $KQ$-module.
An ordering
$x(1) < x(2) < \cdots < x(r)$ 
of the vertices of the quiver $\GG_M$ is called 
$\GG_M$-{\it adapted}
\index{$\GG_M$-adapted ordering}
if the following hold:
If there exists an oriented path from $x(j)$ to 
$x(i)$ in $\GG_M$ we must have $j>i$.
Such orderings always exist since $\GG_M$ is a quiver without oriented
cycles.

\subsection{The algebra $\End_\LL(T_M)$ is quasi-hereditary}
As before, let $M = M_1 \oplus \cdots \oplus M_r$ be a
terminal $KQ$-module, and let $t_i = t_i(M)$ for all $1 \le i \le r$.
For brevity set
$$
B := \End_\LL(T_M).
$$

Recall that there is an inclusion functor
$
\iota_Q\df \md(KQ) \to \md(\LL)
$
defined by $\iota_Q(X) = (X,0)$.
For every $\LL$-module $(X,f)$ we have
$
\iota_Q(\pi_Q(X,f)) = (X,0).
$

Assume that $x(1) < x(2) < \cdots < x(r)$ is a $\GG_M$-adapted ordering
of the vertices of $\GG_M$, compare Section~\ref{adaptedordering}.
Thus we get a bijection 
$$
{\rm v}\df \{x(1),\ldots,x(r)\} \to \{M_1,\ldots,M_r\}.
$$
Set $M(x(i)) := {\rm v}(x(i))$.
It follows that
$\Hom_{KQ}(M(x(i)),M(x(j))) = 0$
if $i < j$.

For each $x(j)$ we have 
$$
M(x(j)) = \tau^a(I_i)
$$
for some uniquely determined $1 \le i \le n$ and $0 \le a \le t_i$.
Define
\begin{align*}
M_{i,a}       &:= M(x(j)),\\
P_{x(j)}      &:= P_{i,a} := \Hom_\LL(T_{i,a},T_M),\\
\Delta_{x(j)} &:= \Delta_{i,a} := \Hom_\LL((M(x(j)),0),T_M),\\
S_{x(j)}      &:= S_{i,a} := {\rm top}(P_{x(j)}).
\end{align*}
\index{$M_{i,a}$}\noindent
\index{$M(x(j))$}\noindent
\index{$\Delta_{x(j)}$}\noindent
\index{$\Delta_{i,a}$}\noindent
For the definition of $T_{i,a}$ we refer to Section~\ref{tia}.
Recall that $P_{x(j)}$ is an indecomposable projective
$B$-module, so its top $S_{x(j)}$ is simple.
We get an ordering
$$
S_{x(1)} < S_{x(2)} < \cdots < S_{x(r)}
$$
of the simple $B$-modules.

Set
$$
\Delta := \left\{ \Delta_{x(1)},\ldots,\Delta_{x(r)} \right\}.
$$

\begin{Lem}\label{qhlemma1}
For each $x(j)$ the following hold:
\begin{itemize}

\item[(i)]
$\pdim(\Delta_{x(j)}) \le 1$;

\item[(ii)]
${\rm top}(\Delta_{x(j)}) \cong S_{x(j)}$;

\item[(iii)]
$\End_B(\Delta_{x(j)}) \cong K$.

\end{itemize}
\end{Lem}

\begin{proof}
For $1 \le i \le n$ and $0 \le a \le t_i$ we get a
$Q$-split short exact sequence
$$
0 \to (\tau^a(I_i),0) \to T_{i,a} \to T_{i,a+1} \to 0
$$
of $\LL$-modules.
This follows easily from the construction of $T_{i,a}$.
Applying $\Hom_\LL(-,T_M)$ yields a short exact sequence
\begin{equation}\label{deltares}
0 \to P_{i,a+1} \to P_{i,a} \to \Delta_{i,a} \to 0
\end{equation}
of $B$-modules.
Here we set $T_{i,t_i+1} = 0$ and $P_{i,t_i+1} = 0$.
Clearly, the exact sequence (\ref{deltares}) is a projective resolution
of $\Delta_{i,a}$.
Thus $\pdim(\Delta_{i,a}) \le 1$.
Furthermore, $P_{i,a}$ is an indecomposable projective module and
therefore has a simple top.
It follows that
${\rm top}(\Delta_{i,a}) \cong S_{i,a}$.
Finally, we have
$$
\End_B(\Delta_{x(j)}) \cong \End_\LL((M(x(j)),0)) \cong \End_{KQ}(M(x(j))) 
\cong K.
$$
(The $KQ$-module $M(x(j))$ is indecomposable preinjective, and therefore its
endomorphism ring is $K$.)
\end{proof}

\begin{Lem}\label{qhlemma2}
We have ${_B}B \in \F(\Delta)$.
\end{Lem}

\begin{proof}
The short exact sequence (\ref{deltares}) in the
proof of Lemma~\ref{qhlemma1} yields a filtration
$$
0 = P_{i,t_i+1} \subset P_{i,t_i} \subset \cdots \subset
P_{i,a+1} \subset P_{i,a}
$$
such that
$
P_{i,k}/P_{i,k+1} \cong \Delta_{i,k}
$
for all $a \le k \le t_i$.
Since each indecomposable projective $B$-module is of the form
$P_{i,a}$ for some $i$ and $a$, this implies
${_B}B \in \F(\Delta)$.
\end{proof}

\begin{Lem}\label{qhlemma3}
A simple $B$-module $S_{x(i)}$ occurs with
multiplicity 
$$
[\Delta_{x(j)}:S_{x(i)}] = \dm \Hom_{KQ}(M(x(j)),M(x(i)))
$$
in every
composition series of $\Delta_{x(j)}$.
\end{Lem}

\begin{proof}
Clearly, for each $i$ and $a$ we have
$$
[\Delta_{x(j)}:S_{i,a}] = \dm \Hom_B(\Hom_\LL(T_{i,a},T_M),\Delta_{x(j)}).
$$
Then by the considerations in Section~\ref{dimvectors2} we know that
$S_{i,a}$ occurs 
$$
\dm \Hom_{KQ}(M(x(j)),\tau^a(I_i))
$$ 
times
in every composition series of $\Delta_{x(j)}$.
\end{proof}

Let $(X,f) \in \CC_M$, and let
$$
X = M(x(1))^{m_{x(1)}} \oplus \cdots \oplus M(x(r))^{m_{x(r)}}
$$ 
be a direct sum decomposition of $X$ into indecomposables.
We assume that $X \not= 0$.
Let $k$ be minimal such that $m_{x(k)} > 0$.
It follows that
$$
\Hom_{KQ}(M(x(k)),\tau M(x(j))) = 0
$$
for all $j$ with $m_{x(j)} > 0$.
For some direct summand $M(x(k))$ of $X$ let
$\iota\df M(x(k)) \to X$ be the canonical inclusion map,
and let $\pi\df X \to X/M(x(k))$ be the corresponding projection.
Furthermore, let $i\df X/M(x(k)) \to X$ be the obvious inclusion.

We obtain a short exact sequence 
\begin{equation}\label{Qsplitsequence}
\xymatrix{
0 \ar[d]^0 \ar[r] & M(x(k)) \ar[d]^0 \ar[r]^>>>>>>>>{\iota} & X \ar[d]^f 
\ar[r]^>>>>>>>>{\pi} & X/M(x(k)) \ar[r] \ar[d]^{f'} & 0 \ar[d]^0 \\
\tau(0) \ar[r]        & \tau(M(x(k))) \ar[r]^>>>>>{\tau(\iota)}& 
\tau(X) \ar[r]^>>>>>{\tau(\pi)} &
\tau(X/M(x(k))) \ar[r] & \tau(0) 
}
\end{equation}
of $\LL$-modules,
where 
$f' = \tau(\pi) \circ f \circ i$.
Clearly, the short exact sequence (\ref{Qsplitsequence})
is $Q$-split.
It follows that it stays exact if we apply $\Hom_\LL(-,T_M)$.

For an algebra $A$ let $\cP_{\le 1}(A)$ be the subcategory
of all modules $X \in \md(A)$ with $\pdim(X) \le 1$.

The main result of this section is the following:

\begin{Thm}\label{qh1}
Let $M$ be a terminal $KQ$-module.
The following hold:
\begin{itemize}

\item[(i)]
The algebra $B:= \End_\LL(T_M)$ is quasi-hereditary with standard modules
$$
\Delta = \left\{ \Delta_{x(1)},\ldots,\Delta_{x(r)} \right\};
$$

\item[(ii)]
$\F(\Delta) = \Hom_\LL(\CC_M,T_M)$;

\item[(iii)]
$T^{\Delta\cap\nabla} = \Hom_\LL(T_M^\vee,T_M)$;

\item[(iv)]
$\F(\Delta) \subseteq \cP_{\le 1}(B)$.

\end{itemize}
\end{Thm}

\begin{proof}
(i): 
By Lemma~\ref{qhlemma3} we know that
$[\Delta_{x(j)}:S_{x(i)}] \not= 0$ implies
$j \ge i$.
Furthermore, we have $S_{i,a+1} > S_{i,a}$.
Using this and the short exact sequence
$$
0 \to P_{i,a+1} \to P_{i,a} \to \Delta_{i,a} \to 0
$$
and the fact ${\rm top}(P_{i,a+1}) \cong S_{i,a+1}$, we get
that $\Delta_{i,a}$ is the largest factor module of $P_{i,a}$
in $\F(\{ S \mid S \le S_{i,a} \})$
where $S$ runs through the simple $B$-modules.
By Lemma~\ref{qhlemma2} we know that ${_B}B \in \F(\Delta)$.
Now Lemma~\ref{qhlemma1}, (iii) yields that $B$ is quasi-hereditary.

(ii):
For $X,Z \in \CC_M$ and we have a functorial isomorphism
$$
\Ext_{F^{T_M}}^1(Z,X) \to \Ext_B^1(\Hom_\LL(X,T_M),\Hom_\LL(Z,T_M)).
$$
Thus the image of the functor
$$
F := \Hom_\LL(-,T_M)\df \CC_M \to \md(B)
$$
is extension closed.
Clearly, for all $x(j)$ the standard module 
$\Delta_{x(j)}$ is in the image of
$F$.
It follows that $\F(\Delta) \subseteq \Ima(F)$.

Using the short exact sequence (\ref{Qsplitsequence}) and induction
on the number of indecomposable direct summands of $X$ one shows
that 
$$
\Hom_\LL((X,f),T_M) \in \F(\Delta)
$$ 
for all $(X,f) \in \CC_M$.
Here one uses that $\F(\Delta)$ is closed under extensions.
Thus $\Ima(F) \subseteq \F(\Delta)$.

(iii):
Let $1 \le i \le n$ and $0 \le b \le t_i$.
Recall that $T_{i,b}^\vee = (I_{i,[0,b]},e_{i,[0,b]})$
and
$$
T_M^\vee = \bigoplus_{i=1}^n \bigoplus_{b=0}^{t_i} T_{i,b}^\vee.
$$
Thus $\Hom_\LL(T_{i,b}^\vee,T_M) \in \F(\Delta)$.
To prove that $\Hom_\LL(T_{i,b}^\vee,T_M) \in \F(\nabla)$
we have to show that
$$
\Ext_B^1(\Delta_{x(j)},\Hom_\LL(T_{i,b}^\vee,T_M)) = 0
$$
for all $x(j)$, compare Section~\ref{qhalgebras}, (xii).
This is equivalent to showing that
$$
\Ext_{F^{T_M}}^1(T_{i,b}^\vee,(M(x(j)),0)) = 0
$$
for all $x(j)$.
In other words, we have to show that every $Q$-split short
exact sequence
\begin{equation}\label{Qsplitsequence2}
0 \to (M(x(j)),0) \xrightarrow{f} E \to T_{i,b}^\vee \to 0
\end{equation}
splits.
Without loss of generality we can assume that 
$$
f = \left[\bbm1\\0\ebm\right]\df M(x(j)) \to M(x(j)) \oplus I_{i,[0,b]}.
$$

Let $N$ be the terminal $KQ$-module defined by $t_k(N) = b$ for
all  $1 \le k \le n$.
(If we are in the Dynkin case, it can happen that
$\tau^b(I_j) = 0$ for some $j$.
In this case, let $t_j(N)$ be the minimal $l$ such that
$\tau^l(I_j) \not= 0$.)

Case 1: If $M(x(j)) \in \add(N)$, then $T_{i,b}^\vee = T_{i,0}$
is $\CC_M$-projective-injective.
Thus the sequence (\ref{Qsplitsequence2}) splits.

Case 2: Assume that $M(x(j)) \notin \add(N)$.
This implies $\Hom_{KQ}(I_{i,[0,b]},\tau(M(x(j)))) = 0$.
Since the sequence (\ref{Qsplitsequence2}) is $Q$-split we know that
$E$ is isomorphic to $(M(x(j)) \oplus I_{i,[0,b]},h)$
where $h$ is of the form
$$
h = \left(\bbm 0 & h'\\ h''&e_{i,[0,b]}\ebm\right)\df
M(x(j)) \oplus I_{i,[0,b]} \to \tau(M(x(j))) \oplus \tau(I_{i,[0,b]}).
$$
It follows that $h'' = 0$, otherwise $f$ would not yield a homomorphism
in $\CC_M$.
We also have $h' = 0$, since
$\Hom_{KQ}(I_{i,[0,b]},\tau(M(x(j)))) = 0$.
This implies that the short exact sequence (\ref{Qsplitsequence2}) splits.

So we proved that 
$$
\Hom_\LL(T_M^\vee,T_M) \in \F(\Delta) \cap \F(\nabla).
$$
Since $T_M^\vee$ has the correct number of isomorphism classes
of indecomposable direct summands, namely $r$, we know
that $\Hom_\LL(T_M^\vee,T_M) = T^{\Delta\cap\nabla}$.
This finishes the proof of (iii).

(iv): This follows directly from Lemma~\ref{qhlemma1}, (i).
\end{proof}

\begin{Cor}
Let $M$ be a terminal $KQ$-module.
Then $\Hom_\LL(-,T_M)$ yields an anti-equivalence
$$
\CC_M \to \F(\Delta).
$$
\end{Cor}

\begin{proof}
Combine Corollary~\ref{fullyfaithful}, Lemma~\ref{fullyfaithful3} and
Theorem~\ref{qh1}, (ii).
\end{proof}

One can easily construct examples of the form $B = \End_\LL(T_M)$
such that $\F(\Delta)$ is a proper subcategory
of $\cP_{\le 1}(B)$, see Example 2 in Section~\ref{D4ex2}.

Next, we describe the $B$-modules $\nabla_{i,a}$.
By Section~\ref{qhalgebras}, (vi) and the proof of Lemma~\ref{qhlemma2} 
we know that
$$
[\nabla_{i,a}:S_{j,b}] = [P_{j,b}:\Delta_{i,a}] =
\begin{cases}
1 & \text{if $i = j$ and $b \le a \le t_j$},\\
0 & \text{otherwise}.
\end{cases}
$$
Thus the $B$-modules $\nabla_{i,a}$ are very easy to describe,
namely $\nabla_{i,a}$ is a serial $B$-module which has
a unique composition series
$$
0 = U_{i,a+1} \subset U_{i,a} \subset U_{i,a-1} \subset \cdots \subset
U_{i,1} \subset U_{i,0} = \nabla_{i,a}
$$
such that
$U_{i,k}/U_{i,k+1} \cong S_{i,k}$ for all $0 \le k \le a$.
In particular, we have
$$
{\rm top}(\nabla_{i,a}) \cong S_{i,0}.
$$
There is no other indecomposable $B$-module $U$ with
$\dimv(U) = \dimv(\nabla_{i,a})$, since the support of
$\nabla_{i,a}$ is a quiver of type ${\mathbb A}_{a+1}$.

It follows from the proof of Lemma~\ref{qhlemma2} and
Section~\ref{qhalgebras}, (ix) that
each $\Delta$-filtration of the indecomposable projective
$B$-module $P_{i,a}$ is structured as follows:
$$
\begin{tabular}{c}
$\Delta_{i,a}$\\
\hline
$\Delta_{i,a+1}$\\
\hline
$\cdots$\\
\hline
$\Delta_{i,t_i}$
\end{tabular}
$$
(We just display the factors of the $\Delta$-filtration of
$P_{i,a}$.)

Next, let us analyse the structure of the characteristic tilting module
$$
T^{\Delta\cap\nabla} = \bigoplus_{i=1}^n \bigoplus_{b=0}^{t_i}
\Hom_\LL(T_{i,b}^\vee,T_M)
$$
in more detail:
It follows easily from the definitions that for all $1 \le i \le n$
and $0 \le b \le t_i$ there is a $Q$-split short exact sequence
$$
0 \to T_{i,b}^\vee \to T_{i,0} \to T_{i,b+1} \to 0
$$
of $\LL$-modules.
Applying $\Hom_\LL(-,T_M)$ yields a short exact sequence
$$
0 \to P_{i,b+1} \to P_{i,0} \to \Hom_\LL(T_{i,b}^\vee,T_M) \to 0
$$
of $B$-modules.
(Again, we set $T_{i,t_i+1} = 0$ and $P_{i,t_i+1} = 0$.)
It follows that each $\Delta$-filtration of 
the $B$-module $\Hom_\LL(T_{i,b}^\vee,T_M)$
has the following structure:
$$
\begin{tabular}{c}
$\Delta_{i,0}$\\
\hline
$\Delta_{i,1}$\\
\hline
$\cdots$\\
\hline
$\Delta_{i,b}$
\end{tabular}
$$
Thus, it is enough to know the structure of the 
$\F(\Delta)$-projective-injective $B$-modules
$P_{i,0} = \Hom_\LL(T_{i,0},T_M)$, in order to describe all
indecomposable direct summands of ${_B}B$ and $T^{\Delta\cap\nabla}$.

Next, let
$Q_M^\op$ 
\index{$Q_M^\op$}
be the full subquiver of $\GG_{T_M}$
with vertices $T_{i,[t_i,t_i]}$ where $1 \le i \le n$.
(For example, if $t_i = t_j$ for all $i$ and $j$, then
$Q_M^\op \cong Q^\op$.)
Let $M'$ be the terminal $KQ_M^\op$-module defined by
$t_i(M') = t_i(M)$ for all $1 \le i \le n$.
Then one can check that
$$
\End_\LL(T_M^\vee)^\op \cong \End_\LL(T_{M'}).
$$
Note that $\End_\LL(T_M^\vee)^\op$ is the endomorphism
algebra of our characteristic tilting module $T^{\Delta \cap \nabla}$.
In other words,
$\End_\LL(T_M^\vee)^\op$ is the Ringel dual of $\End_\LL(T_M)$.
It follows that $\End_\LL(T_M^\vee)^\op$ (and therefore also
$\End_\LL(T_M^\vee)$) is again a quasi-hereditary algebra.

We conclude that $\End_\LL(T_M)$ belongs to a rather special and
interesting class of quasi-hereditary algebras:
\begin{itemize}

\item[(a)]
The characteristic tilting module $T^{\Delta\cap\nabla}$ has
projective dimension one.

\item[(b)]
Each indecomposable projective $\End_\LL(T_M)$-module
and each indecomposable direct summand
of $T^{\Delta\cap\nabla}$ is ``$\Delta$-serial'', i.e. it
has a unique $\Delta$-filtration.
In particular, its $\Delta$-dimension vector
has only entries 0 or 1.

\item[(c)]
All modules in $\nabla$ are serial modules.

\end{itemize}

\begin{Lem}
The characteristic tilting module
$T^{\Delta\cap\nabla}$ is the unique minimal element in the poset 
$(\T_B^{\rm cl}\cap\F(\Delta),\le)$.
\end{Lem}

\begin{proof}
Let $T \in \T_B^{\rm cl} \cap \F(\Delta)$ with $T \le T^{\Delta\cap\nabla}$.
This implies
$\Ext_B^1(T^{\Delta\cap\nabla},T) = 0$. 
By Section~\ref{qhalgebras}, (ii) we get
$$
T \in \F(\Delta) \cap \F(\nabla) = \add(T^{\Delta\cap\nabla}).
$$
Therefore $T = T^{\Delta\cap\nabla}$.
This shows that $T^{\Delta\cap\nabla}$ is a minimal element
in $(\T_B^{\rm cl}\cap\F(\Delta),\le)$.

To show uniqueness,
assume $T$ is a minimal element in $(\T_B^{\rm cl}\cap\F(\Delta),\le)$.
Since $T \in \F(\Delta)$ we know again by Section~\ref{qhalgebras}, (ii)
that $T^{\Delta\cap\nabla} \in T^\perp$.
Now Section~\ref{tiltingord}, (d) yields that 
$T^{\Delta\cap\nabla} \le T$.
This implies $T =  T^{\Delta\cap\nabla}$.
\end{proof}

\begin{Lem}
The modules $\Hom_\LL(T_{i,0},T_M)$, $1 \le i \le n$ are
the indecomposable $\F(\Delta)$-projective-injectives modules.
\end{Lem}

\begin{proof}
The modules $\Hom_\LL(T_{i,0},T_M)$ are
the only indecomposable projective $B$-modules, which are 
direct summands of
$T^{\Delta\cap\nabla}$.
Therefore, by
Section~\ref{qhalgebras}, (xi) we know that
$\Hom_\LL(T_{i,0},T_M)$ is $\F(\Delta)$-injective.
Furthermore, any $\F(\Delta)$-projective module is
projective by Section~\ref{qhalgebras}, (x).
This finishes the proof.
\end{proof}

\begin{Cor}
If $T$ is a tilting module in $\F(\Delta)$, then
$T$ has a direct summand isomorphic
to
$$
\bigoplus_{i=1}^n \Hom_\LL(T_{i,0},T_M).
$$
\end{Cor}

It is straightforward to construct examples where
the $B$-module $\Hom_\LL(T_{i,0},T_M)$ 
is not injective in $\md(B)$, see Example 2 in Section~\ref{D4ex2}.

\begin{Conj}
The Hasse diagram of $(\T_B^{\rm cl} \cap \F(\Delta),\le)$ is
connected.
\end{Conj}

\begin{Lem}\label{comparedimv}
For a $\LL$-module $(X,f) \in \CC_M$ the following hold:
\begin{itemize}

\item[(i)]
$[\Hom_\LL((X,f),T_M):\Delta_{x(j)}] = [\Hom_\LL((X,0),T_M):\Delta_{x(j)}]$;

\item[(ii)]
$\dimv_\Delta(\Hom_\LL((X,f),T_M)) = 
\dimv_\Delta(\Hom_\LL((X,0),T_M))$;

\item[(iii)]
$\dimv(\Hom_\LL((X,f),T_M)) = \dimv(\Hom_\LL((X,0),T_M)).
$

\end{itemize}
\end{Lem}

\begin{proof}
Using the short exact sequence (\ref{Qsplitsequence})
one shows by induction on the number of indecomposable
direct summands of $X$ that
(i) holds.
It follows from the definitions
that (i) and (ii) are equivalent.
Modules having the same $\Delta$-dimension vector, also have
the same dimension vectors, i.e. (ii) implies (iii).
\end{proof}

\subsection{Examples of type $\A_3$}
Every finite-dimensional algebra $A$ with global dimension at most two
is quasi-hereditary, i.e. one can find an ordering
of the simple $A$-modules such that $A$ becomes quasi-hereditary with respect
to that ordering.
If $T$ is a $\CC_M$-maximal rigid $\LL$-module, then 
$\gldim(\End_\LL(T)) = 3$.
It seems to be difficult to determine when $\End_\LL(T)$ is
quasi-hereditary and when not.

Even if $Q$ is a quiver of type $\A_3$ 
there are maximal rigid modules whose endomorphism
algebra is not quasi-hereditary:
Let $T$
be the maximal rigid $\LL$-module
$$
{\bsm2&\\&3\esm} \oplus {\bsm&2&\\1&&3\esm} \oplus 
{\bsm&2\\1&\esm} \oplus
{\bsm1&&\\&2&\\&&3\esm} \oplus {\bsm&2&\\1&&3\\&2&\esm} \oplus 
{\bsm&&3\\&2&\\1&&\esm}.
$$
The quiver of $\End_\LL(T)$ looks as follows:
$$
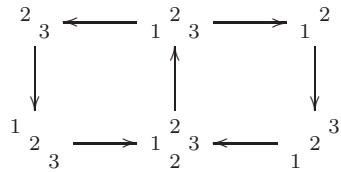
\xymatrix{
{\bsm2&\\&3\esm} \ar[d] & {\bsm&2&\\1&&3\esm} \ar[l]\ar[r] & 
{\bsm&2\\1&\esm} \ar[d]\\
{\bsm1&&\\&2&\\&&3\esm} \ar[r] & {\bsm&2&\\1&&3\\&2&\esm} \ar[u] & 
{\bsm&&3\\&2&\\1&&\esm} \ar[l]
}
$$
It is not difficult to show that $\End_\LL(T)$ is not a quasi-hereditary
algebra.

As another example (also in the type $\A_3$ case),
let $T$
be the maximal rigid $\LL$-module
$$
{\bsm2&\\&3\esm} \oplus {\bsm&2&\\1&&3\esm} \oplus 
{\bsm3\esm} \oplus
{\bsm1&&\\&2&\\&&3\esm} \oplus {\bsm&2&\\1&&3\\&2&\esm} \oplus 
{\bsm&&3\\&2&\\1&&\esm}.
$$
The quiver of $\End_\LL(T)$ looks as follows:
$$
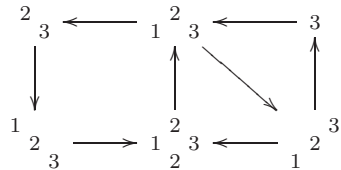
\xymatrix{
{\bsm2&\\&3\esm} \ar[d] & {\bsm&2&\\1&&3\esm} \ar[dr]\ar[l] & 
{\bsm3\esm} \ar[l]\\
{\bsm1&&\\&2&\\&&3\esm} \ar[r] & {\bsm&2&\\1&&3\\&2&\esm} \ar[u] & 
{\bsm&&3\\&2&\\1&&\esm} \ar[l]\ar[u]
}
$$
The algebra $\End_\LL(T)$ can be equipped with the structure
of a quasi-hereditary algebra.
However, it turns out that for each of these structures, there
are standard modules $\Delta_i$ with $\pdim(\Delta_i) = 2$.
Thus the image of the functor 
$\Hom_\LL(-,T)$ does not coincide with the category $\F(\Delta)$ of
$\Delta$-filtered $\End_\LL(T)$-modules.

\subsection{An example of type $\A_4$}

Let $Q$ be the quiver
{\small
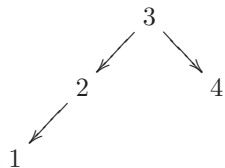
$$
\xymatrix@-0.8pc{
&& 3 \ar[dl]\ar[dr]\\
& 2 \ar[dl] && 4\\
1
}
$$
}\noindent
of Dynkin type $\A_4$.

The Auslander-Reiten quiver ${\rm AR}(KQ)$ of $KQ$ looks as follows:
{\small
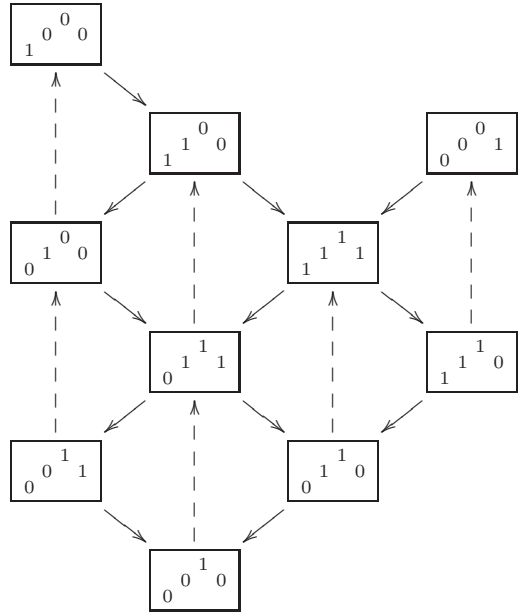
$$
\xymatrix@-1pc{
\fbox{$\bsm&&0&\\&0&&0\\1&&&\esm$} \ar[dr]\\
& \fbox{$\bsm&&0&\\&1&&0\\1&&&\esm$} \ar[dr]\ar[dl] && 
\fbox{$\bsm&&0&\\&0&&1\\0&&&\esm$} \ar[dl]\\
\fbox{$\bsm&&0&\\&1&&0\\0&&&\esm$} \ar@{-->}[uu]\ar[dr] && 
\fbox{$\bsm&&1&\\&1&&1\\1&&&\esm$}  \ar[dl]\ar[dr]\\
& \fbox{$\bsm&&1&\\&1&&1\\0&&&\esm$} \ar@{-->}[uu] \ar[dl]\ar[dr] && 
\fbox{$\bsm&&1&\\&1&&0\\1&&&\esm$}  \ar@{-->}[uu]\ar[dl] \\
\fbox{$\bsm&&1&\\&0&&1\\0&&&\esm$} \ar@{-->}[uu] \ar[dr] && 
\fbox{$\bsm&&1&\\&1&&0\\0&&&\esm$} \ar@{-->}[uu] \ar[dl]\\
& \fbox{$\bsm&&1&\\&0&&0\\0&&&\esm$} \ar@{-->}[uu] 
}
$$
}\noindent
As usual, the solid arrows correspond to the irreducible
maps between the indecomposable $KQ$-modules, and the
dotted arrows describe the Auslander-Reiten translation in
$\md(KQ)$.

Next, let $M(x(1)),\ldots,M(x(10))$ be the 10 indecomposable
$KQ$-modules labelled in such a way that 
$x(1) < x(2) < \cdots < x(10)$
is a $\GG_M$-adapted ordering of the vertices of $\GG_M$, where $M$
is just the direct sum of all indecomposable $KQ$-modules.
For example, 
{\small
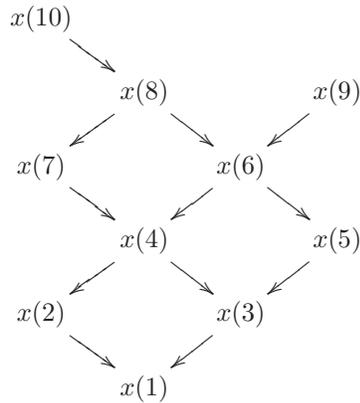
$$
\xymatrix@-1pc{
x(10) \ar[dr]\\
& x(8) \ar[dr]\ar[dl] && 
x(9) \ar[dl]\\
x(7) \ar[dr] && 
x(6) \ar[dl]\ar[dr]\\
& x(4) \ar[dl]\ar[dr] && 
x(5) \ar[dl] \\
x(2) \ar[dr] && 
x(3) \ar[dl]\\
& x(1)
}
$$
}\noindent
is such a labelling.
This labelling of the indecomposable $KQ$-modules will be fixed for
the rest of this section.

Now we replace the dotted arrows in ${\rm AR}(KQ)$ by solid arrows,
and for simplicity we label the vertices by $i$ instead of $x(i)$.
In this way we obtain the quiver $\GG_{T_M} = \GG_M^*$ 
of the algebra $\End_\LL(T_M)$,
where $M := M(x(1)) \oplus \cdots \oplus M(x(10))$: 
{\small
$$
\xymatrix@-0.5pc{
10 \ar[dr]\\
& 8 \ar[dr]\ar[dl] && 
9 \ar[dl]\\
7 \ar[uu]\ar[dr] && 
6 \ar[dl]\ar[dr]\\
& 4 \ar[uu] \ar[dl]\ar[dr] && 
5 \ar[uu]\ar[dl] \\
2 \ar[uu] \ar[dr] && 
3 \ar[uu] \ar[dl]\\
& 1 \ar[uu] 
}
$$
}\noindent
Note that $\End_\LL(T_M)$ contains as a subalgebra the
Auslander algebra of $KQ$
(just delete the solid arrows which came from the dotted arrows
of ${\rm AR}(KQ)$).

For each vertex $x(i)$ of the quiver of $\End_\LL(T_M)$ 
let $\Delta_{x(i)}$ be the associated standard module.
It turns out that $\Delta_{x(i)}$ is just the indecomposable
projective module over the Auslander algebra of $KQ$ considered
as an $\End_\LL(T_M)$-module.

For example, as a quiver representation the $\End_\LL(T_M)$-module 
$\Delta_{x(8)}$ looks
as follows:
{\small
$$
\xymatrix@-0.5pc{
0 \ar[dr]\\
& K \ar[dr]^\id \ar[dl]_\id && 0 \ar[dl]\\
K \ar[uu]\ar[dr]^\id && K \ar[dl]_\id \ar[dr]^\id\\
& K \ar[uu]_{\rm 0} \ar[dl]\ar[dr]^\id && K \ar[uu]\ar[dl]_\id \\
0 \ar[uu] \ar[dr] && K \ar[uu]_{\rm 0} \ar[dl] \\
& 0 \ar[uu] 
}
$$
}\noindent
The arrows without a label are just zero maps.
Another way of displaying this module would be
$$
\bsm&8\\7&&6\\&4&&5\\&&3\esm
$$
The numbers correspond to composition factors.
Here is a table of all the modules $\Delta_{x(i)}$ and 
$\nabla_{x(i)}$:
\begin{center}
\begin{tabular}{|c|c|c|c|c|c|c|c|c|c|c|}\hline
$\Delta_{x(i)}$ &
$\bsm1\esm$ & 
$\bsm2&\\&1\esm$ &
$\bsm&3\\1\esm$ &
$\bsm&4\\2&&3\\&1\esm$ &
$\bsm&&5\\&3\\1\esm$ &
$\bsm&&6\\&4&&5\\2&&3\\&1\esm$ &
$\bsm7\\&4\\&&3\esm$ &
$\bsm&8\\7&&6\\&4&&5\\&&3\esm$ &
$\bsm&&&9\\&&6\\&4\\2\esm$ &
$\bsm&\\10\\&8\\&&6\\&&&5\\&\esm$\\ \hline
$\nabla_{x(i)}$ &
$\bsm1\esm$ & 
$\bsm2\esm$ &
$\bsm&3\esm$ &
$\bsm1\\4\esm$ &
$\bsm5\esm$ &
$\bsm3\\6\esm$ &
$\bsm2\\7\esm$ &
$\bsm&\\1\\4\\8\\&\esm$ &
$\bsm5\\9\esm$ &
$\bsm2\\7\\10\esm$\\ \hline 
\end{tabular}
\end{center}
Here are pictures describing the structure of the indecomposable
projective $\End_\LL(T_M)$-modules.
The factors of a $\Delta$-filtration of the modules $P_{x(i)}$ 
are marked by different colours.

The top $S_{x(i)}$
of a standard module $\Delta_{x(i)}$ is displayed as
$\fbox{i}$.
Note that
$P_{x(8)} \subset P_{x(4)} \subset P_{x(1)}$,
$P_{x(10)} \subset P_{x(7)} \subset P_{x(2)}$,
$P_{x(6)} \subset P_{x(3)}$
and $P_{x(9)} \subset P_{x(5)}$.
{\small
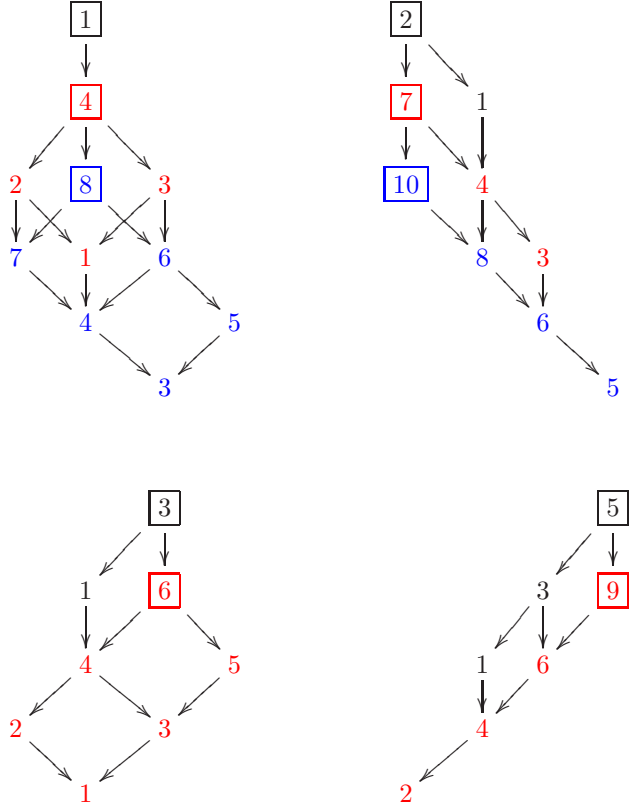
$$
\xymatrix@-1pc{
&\fbox{$1$} \ar[d] && &&
& \fbox{$2$}\ar[d]\ar[dr]&&&
\\
&\red{\fbox{$4$}} \ar[dl]\ar[d]\ar[dr]&& &&
& \red{\fbox{$7$}}\ar[d]\ar[dr]&1\ar[d]&&
\\
\red{2}\ar[d]\ar[dr] & \blue{\fbox{$8$}} \ar[dl]\ar[dr] & 
\red{3} \ar[dl]\ar[d]& &&
& \blue{\fbox{$10$}}\ar[dr]&\red{4}\ar[d]\ar[dr]&&
\\
\blue{7}\ar[dr] & \red{1} \ar[d] & \blue{6} \ar[dl]\ar[dr]& &&
&&\blue{8}\ar[dr]&\red{3}\ar[d]&
\\
& \blue{4} \ar[dr] &&\blue{5}\ar[dl] &&
&&&\blue{6} \ar[dr]&
\\
&&\blue{3} && &&&&&\blue{5}\\
&&&\\
&&\fbox{$3$}\ar[dl]\ar[d]& &&
&&&&\fbox{$5$}\ar[dl]\ar[d]
\\
&1\ar[d]&\red{\fbox{$6$}}\ar[dl]\ar[dr]& &&
&&&3\ar[dl]\ar[d]&\red{\fbox{$9$}}\ar[dl]
\\
&\red{4}\ar[dl]\ar[dr]&&\red{5}\ar[dl] &&
&&1\ar[d]&\red{6}\ar[dl]&
\\
\red{2}\ar[dr]&&\red{3}\ar[dl]& &&
&&\red{4}\ar[dl]&&
\\
&\red{1}&& &&
&\red{2}&&&
}
$$
}\noindent

Next, we would like to
construction $T_M$ explicitly.
First, observe that the
Auslander-Reiten quiver ${\rm AR}(KQ^\op)$ of $KQ^\op$
looks as follows:
{\small
$$
\xymatrix@-1pc{
{\rm AR}(KQ^\op)\df &&&x(1)\ar[dl]\ar[dr]&&\\
&&x(2)\ar[dr]&&x(3)\ar[dl]\ar[dr]&\\
&&&x(4)\ar@{-->}[uu]\ar[dl]\ar[dr]&&x(5)\ar[dl]\\
&&x(7)\ar@{-->}[uu]\ar[dr]&&x(6)\ar@{-->}[uu]\ar[dl]\ar[dr]&\\
&&&x(8)\ar@{-->}[uu]\ar[dl]&&x(9)\ar@{-->}[uu]\\
&&x(10)\ar@{-->}[uu]&&&\\
Q^\op\df&&&&3\ar[dl]\ar[dr]&\\
&&&2\ar[dl]&&4&\\
&&1&&&
}
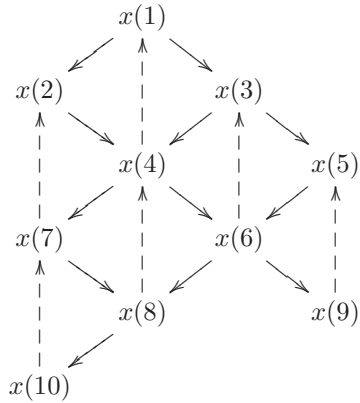
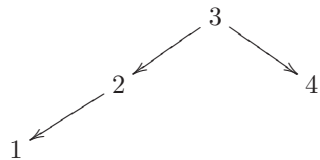
$$
}\noindent
Again, we will just write $i$ instead of $x(i)$.

Denote the Auslander algebra of $KQ^\op$ by
$C$, and let
$P(i)$, $1 \le i \le r$ 
be the indecomposable projective $C$-modules.
Set 
$$
T_i :=  F_\lambda(P(i))
\text{\;\;\; and \;\;\;}
T_M = T_1 \oplus \cdots \oplus T_r.
$$ 
Here $F_\lambda$ is the push-down functor $\md(\widetilde{\LL}) \to \md(\LL)$
where $\widetilde{\LL}$ is the obvious covering (with Galois group
$\Z$) of $\LL$.
We can consider $C$ as a subalgebra of $\widetilde{\LL}$, so the
expression $F_\lambda(P(i))$ makes sense.
The following table illustrates how the modules $P(i)$ and $T_i$ look like:
\begin{center}
\begin{tabular}{|c|c|c|c|c|c|c|c|c|c|c|}\hline
$P(i)$ &
$\bsm&\\&1\\2&&3\\&4&&5\\&&6\\&\esm$ & 
$\bsm2\\&4\\&&6\\&&&9\esm$ &
$\bsm&&3\\&4&&5\\7&&6\\&8\esm$ &
$\bsm&4\\7&&6\\&8&&9\esm$ &
$\bsm&&&5\\&&6\\&8\\10\esm$ &
$\bsm&&6\\&8&&9\\10\esm$ &
$\bsm7\\&8\esm$ &
$\bsm&8\\10\esm$ &
$\bsm9\esm$ &
$\bsm10\esm$\\ \hline 
$T_i$ &
$\bsm&\\&2\\1&&3\\&2&&4\\&&3\\&\esm$ & 
$\bsm1\\&2\\&&3\\&&&4\esm$ &
$\bsm&&3\\&2&&4\\1&&3\\&2\esm$ &
$\bsm&2\\1&&3\\&2&&4\esm$ &
$\bsm&&&4\\&&3\\&2\\1\esm$ &
$\bsm&&3\\&2&&4\\1\esm$ &
$\bsm1\\&2\esm$ &
$\bsm&2\\1\esm$ &
$\bsm4\esm$ &
$\bsm1\esm$\\ \hline 
\end{tabular}
\end{center}
Note that $Q$ could be identified with the full subquiver given by the
vertices $\{10,8,6,9\}$.
Then
one easily checks that the restriction of $T_i$ to the full subquiver
given by $\{10,8,6,9\}$ is just the module $M(x(i))$.

\subsection{Examples of type $\D_4$}\label{D4ex2}
Let $Q$ be the quiver
{\small
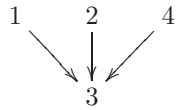
$$
\xymatrix@-0.5pc{
1\ar[dr]&2\ar[d]&4\ar[dl]\\
&3
}
$$
}
of Dynkin type $\D_4$.
The Auslander-Reiten quiver of $KQ$ looks as follows:
{\small
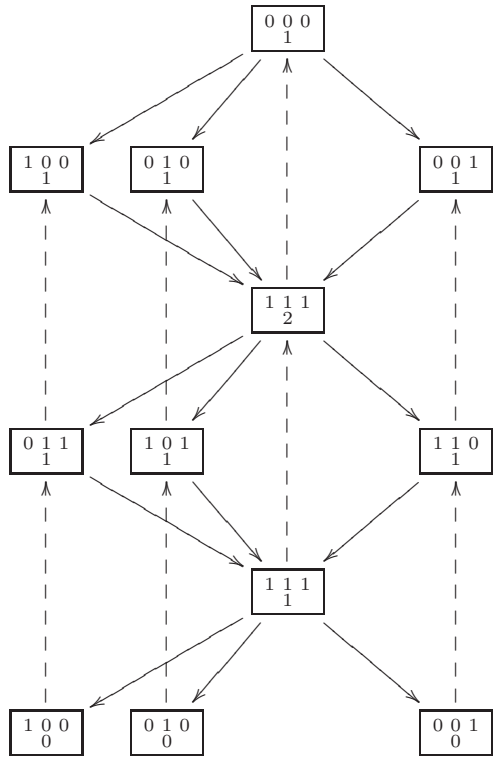
$$
\xymatrix@-1pc{
&&\fbox{$\bsm0&0&0\\&1&\esm$}\ar[ddll]\ar[ddl]\ar[ddrr]\\
&&&&&\\
\fbox{$\bsm1&0&0\\&1&\esm$}\ar[ddrr] &
\fbox{$\bsm0&1&0\\&1&\esm$}\ar[ddr]&&&
\fbox{$\bsm0&0&1\\&1&\esm$}\ar[ddll]\\
&&&&&\\
&&\fbox{$\bsm1&1&1\\&2&\esm$}\ar[ddll]\ar[ddl]\ar[ddrr]\ar@{-->}[uuuu]\\
&&&&&\\
\fbox{$\bsm0&1&1\\&1&\esm$}\ar[ddrr]\ar@{-->}[uuuu] &
\fbox{$\bsm1&0&1\\&1&\esm$}\ar[ddr]\ar@{-->}[uuuu]&&&
\fbox{$\bsm1&1&0\\&1&\esm$}\ar[ddll]\ar@{-->}[uuuu]\\
&&&&&\\
&&\fbox{$\bsm1&1&1\\&1&\esm$}\ar[ddll]\ar[ddl]\ar[ddrr]\ar@{-->}[uuuu]\\
&&&&&\\
\fbox{$\bsm1&0&0\\&0&\esm$}\ar@{-->}[uuuu]&
\fbox{$\bsm0&1&0\\&0&\esm$}\ar@{-->}[uuuu]&&&
\fbox{$\bsm0&0&1\\&0&\esm$}\ar@{-->}[uuuu]
}
$$
}\noindent
As before,
we label the indecomposable $KQ$-modules 
$M(x(1)),\ldots,M(x(12))$
such that 
$x(1) < x(2) < \cdots < x(12)$ is a $\GG_M$-adapted ordering:
{\small
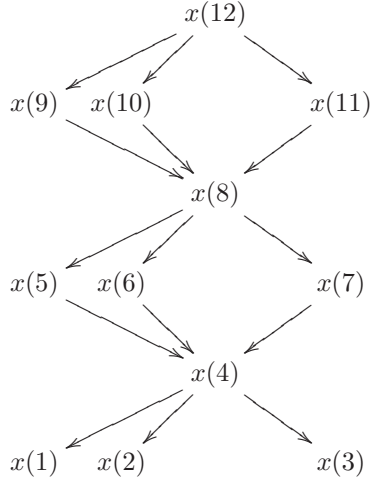
$$
\xymatrix@-1.5pc{
&&x(12)\ar[ddll]\ar[ddl]\ar[ddrr]\\
&&&&&\\
x(9)\ar[ddrr] & x(10)\ar[ddr]&&&
x(11)\ar[ddll]\\
&&&&&\\
&&x(8)\ar[ddll]\ar[ddl]\ar[ddrr]\\
&&&&&\\
x(5)\ar[ddrr] &
x(6)\ar[ddr] &&&
x(7)\ar[ddll]\\
&&&&&\\
&&x(4)\ar[ddll]\ar[ddl]\ar[ddrr]\\
&&&&&\\
x(1) &
x(2) &&&
x(3)
}
$$
}\noindent
For simplicity we write again just $i$ instead of $x(i)$.

\noindent{\bf Example 1}: 
Let 
$$
M := \bigoplus_{i=1}^{12} M(x(i)).
$$
The following picture describes the indecomposable projective
$\End_\LL(T_M)$-modules $P_{x(1)}$ and $P_{x(4)}$.
Up to symmetry, the projectives $P_{x(2)}$ and $P_{x(3)}$ look
like $P_{x(1)}$.
The factors of a $\Delta$-filtration of each of these modules
is highlighted by different colours.
{\small
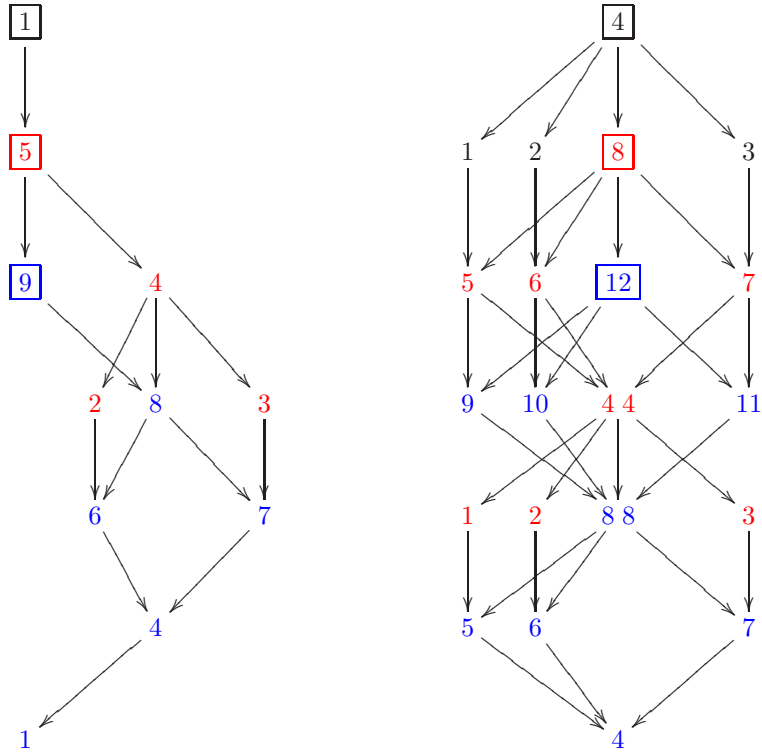
$$
\xymatrix@-1pc{
\fbox{$1$} \ar[dd]&&&&&&&&
&&\fbox{$4$} \ar[ddll]\ar[ddl]\ar[dd]\ar[ddrr]\\
&&&&&\\
\red{\fbox{$5$}}\ar[dd]\ar[ddrr] &&&&&&&&
1\ar[dd]& 2 \ar[dd]& \red{\fbox{$8$}}\ar[ddll]\ar[ddl]\ar[dd]\ar[ddrr]
&&3\ar[dd]\\
&&&&&\\
\blue{\fbox{$9$}}\ar[ddrr] && \red{4}\ar[ddrr]\ar[dd]\ar[ddl] &&&&&&
\red{5}\ar[dd]\ar[ddrr]&\red{6}\ar[dd]\ar[ddr]&
\blue{\fbox{$12$}}\ar[ddll]\ar[ddl]\ar[ddrr]&&
\red{7}\ar[ddll]\ar[dd]\\
&&&&&\\
&\red{2}\ar[dd] & \blue{8}\ar[ddl]\ar[ddrr] && \red{3}\ar[dd]&&&&
\blue{9}\ar[ddrr]&\blue{10}\ar[ddr]&\red{4\;4}
\ar[ddll]\ar[ddl]\ar[dd]\ar[ddrr]&&
\color{blue}{11}\ar[ddll]\\
&&&&&&\\
& \blue{6}\ar[ddr] &&& \blue{7}\ar[ddll]&&&&
\color{red}{1}\ar[dd]&\color{red}{2}\ar[dd]&\color{blue}{8\;8}
\ar[ddll]\ar[ddl]\ar[ddrr]&&
\color{red}{3}\ar[dd]\\
&&&&&\\
&&\blue{4}\ar[ddll]&&&&&&
\color{blue}{5}\ar[ddrr]&\color{blue}{6}\ar[ddr]&&&
\color{blue}{7}\ar[ddll]\\
&&&&&\\
\blue{1}&&&&&&&&
&&\color{blue}{4}
}
$$
}\noindent
Note that
$$
\dimv(P_{x(4)}) = (2,2,2,4,2,2,2,3,1,1,1,1).
$$
Thus $\dm P_{x(4)} = 23$.
We know that $P_{x(4)} \in \F(\Delta) \cap \F(\nabla)$.
We leave it as an easy exercise to work out that
$P_{x(4)}$ has a $\nabla$-filtration of the following structure:
$$
\begin{tabular}{c}
$\nabla_{12}$\\
\hline
$\nabla_9 \oplus \nabla_{10} \oplus \nabla_{11}$\\
\hline
$\nabla_8 \oplus \nabla_8$\\
\hline
$\nabla_5 \oplus \nabla_6 \oplus \nabla_7$\\
\hline
$\nabla_4$
\end{tabular}
$$

\noindent{\bf Example 2}: 
We keep the notation from above, but now we define
$$
M := \bigoplus_{i=1}^{8} M(x(i)).
$$
Again, set $B := \End_\LL(T_M)$.
The following picture describes the indecomposable projective
$B$-modules $P_{x(i)}$.
The factors of a $\Delta$-filtration are marked by different colours.
Note that
$P_{x(8)} \subset P_{x(4)}$,
$P_{x(5)} \subset P_{x(1)}$,
$P_{x(6)} \subset P_{x(2)}$
and $P_{x(7)} \subset P_{x(3)}$.
{\small
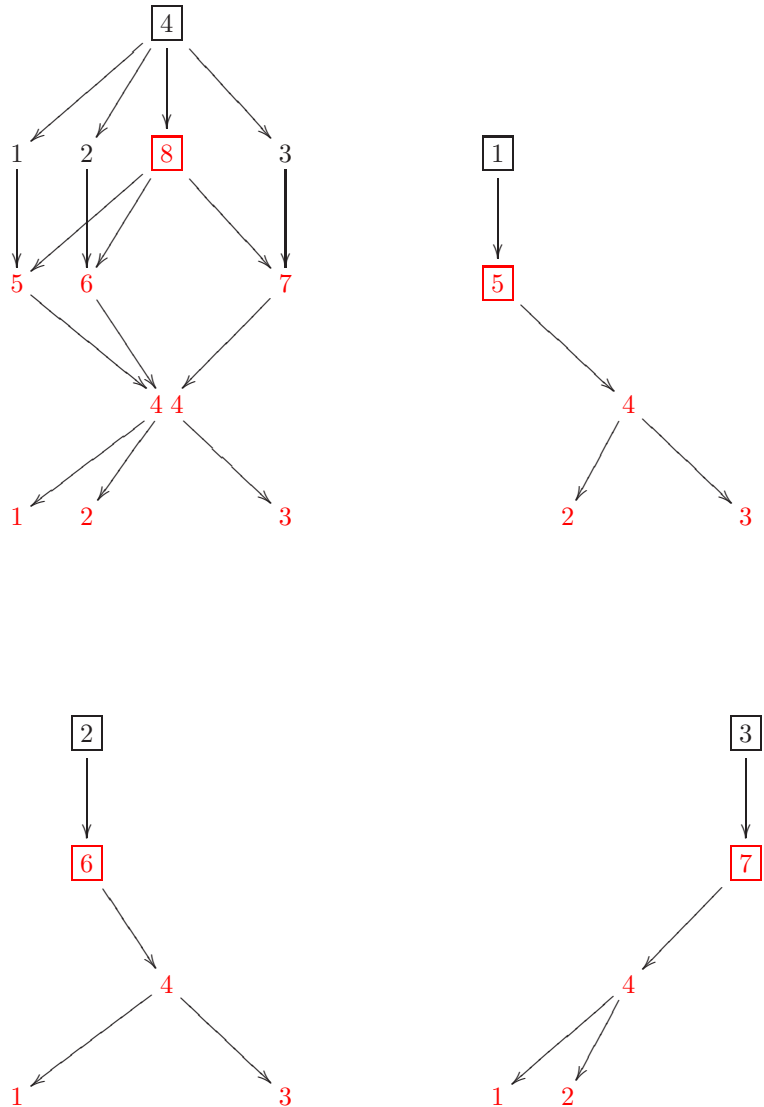
$$
\xymatrix@-1pc{
&&\fbox{$4$} \ar[ddll]\ar[ddl]\ar[dd]\ar[ddrr]\\
&&&&&\\
1\ar[dd]& 2 \ar[dd]& \red{\fbox{$8$}}\ar[ddll]\ar[ddl]\ar[ddrr]
&&3\ar[dd] 
&&&&\fbox{$1$}\ar[dd]
\\
&&&&&\\
\red{5}\ar[ddrr]&\red{6}\ar[ddr]&&&
\red{7}\ar[ddll]
&&&&\red{\fbox{$5$}}\ar[ddrr]
\\
&&&&&\\
&&\red{4\;4}
\ar[ddll]\ar[ddl]\ar[ddrr]&&
&&&&&&\red{4}\ar[ddl]\ar[ddrr]
\\
&&&&&\\
\color{red}{1}&\color{red}{2}&&&
\color{red}{3}
&&&&&\red{2}&&&\red{3}\\
&&&&&&&&&&\\
&&&&&&&&&&\\
&&&&&&&&&&\\
&\fbox{$2$}\ar[dd]
&&&&&&&&&&&\fbox{$3$}\ar[dd]
\\
&&&&&&&&&&\\
&\red{\fbox{$6$}}\ar[ddr]&&&&&&&&&&&\red{\fbox{$7$}}\ar[ddll]\\
&&&&&&&&&&\\
&& \red{4}\ar[ddll]\ar[ddrr]&&&&&&&&\red{4}\ar[ddll]\ar[ddl]\\
&&&&&&&&&&\\
\red{1} &&&&\red{3}&&&&\red{1}& \red{2}
}
$$
}\noindent
Note that none of the modules $P_{x(i)}$ has a simple socle.
Thus there are no non-zero projective-injective $B$-modules.

This time there is a $\nabla$-filtration of $P_{x(4)}$ which is structured as
follows:
$$
\begin{tabular}{c}
$\nabla_8$\\
\hline
$\nabla_5 \oplus \nabla_6 \oplus \nabla_7$\\
\hline
$\nabla_4 \oplus \nabla_4$\\
\hline
$\nabla_1 \oplus \nabla_2 \oplus \nabla_3$\\
\end{tabular}
$$

There is an obvious embedding 
$\iota\df P_{x(1)} \to P_{x(4)}$ such that
the cokernel $Z$ of $\iota$ looks as follows:
{\small
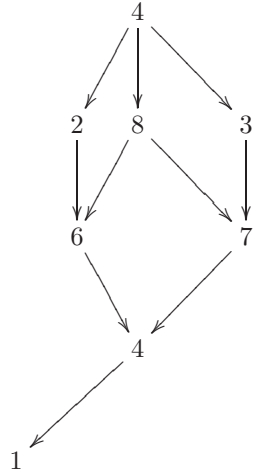
$$
\xymatrix@-1pc{
&&4\ar[ddl]\ar[dd]\ar[ddrr]\\
&&&&&&&&&&\\
&2 \ar[dd] & 8 \ar[ddl]\ar[ddrr]&&3\ar[dd]\\
&&&&&&&&&&\\
&6\ar[ddr]&&&7\ar[ddll]\\
&&&&&&&&&&\\
&&4\ar[ddll]\\
&&&&&&&&&&\\
1
}
$$
}\noindent
Thus we obtain a short exact sequence
$$
0 \to P_{x(1)} \xrightarrow{\iota} P_{x(4)} \to Z \to 0.
$$
This shows that $\pdim(Z) \le 1$.

Next, we apply $\Hom_B(-,Z)$ to the short exact sequence above.
Up to scalars there is only one homomorphism 
$f\df P_{x(1)} \to Z$. 
(The image of $f$ is the socle $S_{x(1)}$ of $Z$.)
It is also obvious that $f$ factors through $\iota$, i.e.
there exists a homomorphism $g\df P_{x(4)} \to Z$ such that 
$g \circ \iota = f$.
(The image of $g$ is the unique 2-dimensional submodule of 
$Z$.)
It follows that 
$$
\Hom_B(P_{x(4)},Z) \xrightarrow{\Hom_B(\iota,Z)} \Hom_B(P_{x(1)},Z)
$$
is surjective.
Since $P_{x(4)}$ is projective, we have
$\Ext_B^1(P_{x(4)},Z) = 0$.
It follows that $\Ext_B^1(Z,Z) = 0$.

It is also clear that $Z$ does not lie in $\F(\Delta)$.
(The top of $Z$ is isomorphic to $S_{x(4)}$, so if $Z \in \F(\Delta)$,
then $\Delta_{x(4)}$ has to be a factor module of $Z$, which is
not the case.)

Taking the Bongartz completion of $Z$ we obtain
a classical tilting module in $\T_B^{\rm cl}$ which is not
contained in $\F(\Delta)$.
In particular, $\F(\Delta)$ and $\cP_{\le 1}(B)$ do
not coincide.


{\Large\section{Mutations of clusters via $\Delta$-dimension vectors}
\label{section16}}


For ${\mathbf d} = (d_1,\ldots,d_r)$ and ${\mathbf e} = (e_1,\ldots,e_r)$
in $\Z^r$ define
$$
{\mathbf d} \cdot {\mathbf e} := \sum_{i=1}^r d_ie_i.
$$

Let $M = M_1 \oplus \cdots \oplus M_r$ be a terminal
$KQ$-module, and
let $x(1) < x(2) < \cdots < x(r)$ be a $\GG_M$-adapted ordering
of the vertices of $\GG_M$.
Set $B := \End_\LL(T_M)$.

We know that the $K$-dimension of the standard module $\Delta_{x(i)}$ is
\begin{align*}
\dm \Delta_{x(i)} &= \dm \Hom_\LL((M(x(i)),0),T_M) \\
&= \sum_{j=1}^r \dm \Hom_{KQ}(M(x(i)),M(x(j))) \\ 
&= \sum_{j=1}^i \dm \Hom_{KQ}(M(x(i)),M(x(j))) 
\end{align*}
for all $i$.
So $\dm \Delta_{x(i)}$ can be calculated in the mesh category.
Define
$$
d_\Delta := (\dm \Delta_{x(1)},\ldots,\dm \Delta_{x(r)}).
$$
As before, for a $\LL$-module $X \in \CC_M$ let 
$\dimv_\Delta(\Hom_\LL(X,T_M))$ be
the $\Delta$-dimension vector of the $B$-module
$\Hom_\LL(X,T_M)$.

Now let $T = T_1 \oplus \cdots \oplus T_r$ 
be a basic $\CC_M$-maximal rigid $\LL$-module, and
suppose that
$T_k$ is not $\CC_M$-projective-injective. 
Then we can mutate $T$ in direction $T_k$.
We obtain two exchange sequences
$$
0 \to T_k \to T' \to T_k^* \to 0
\text{\;\;\; and \;\;\;}
0 \to T_k^* \to T'' \to T_k \to 0
$$
with $T',T'' \in \add(T/T_k)$.

For brevity, set 
$$
{\mathbf d}_i := \dimv_\Delta(\Hom_\LL(T_i,T_M))
$$
for all $1 \le i \le r$.
Similarly to the definition of $\GG_T'$ in Section~\ref{section15}
let $\GG_T''$ be the quiver which is obtained from the quiver
of $\End_\LL(T)$ be replacing the vertex corresponding to $T_i$
by the $\Delta$-dimension vector $\dimv_\Delta(\Hom_\LL(T_i,T_M))$.

\begin{Prop}
The $\Delta$-dimension vector of the $B$-module $\Hom_\LL(T_k^*,T_M)$
is
$$
{\mathbf d}_k^* :=
\begin{cases}
- {\mathbf d}_k + \sum_{{\mathbf d}_k \to {\mathbf d}_i} {\mathbf d}_i &
\text{if
$
\sum_{{\mathbf d}_k \to {\mathbf d}_i} {\mathbf d}_i \cdot d_\Delta >
\sum_{{\mathbf d}_j \to {\mathbf d}_k} {\mathbf d}_j \cdot d_\Delta
$},\\
- {\mathbf d}_k +  
\sum_{{\mathbf d}_j \to {\mathbf d}_k} {\mathbf d}_j &
\text{otherwise}.
\end{cases}
$$
Here
the sums are taken over all arrows of the quiver
of $\GG_T''$ which start,
respectively end in the vertex ${\mathbf d}_k$.
\end{Prop}

\begin{proof}
This follows immediately from our results in Section~\ref{section15}.
\end{proof}

In the example in Section~\ref{ex14.3}, the 
graph $\GG_{T_M}''$ looks as follows:
$$
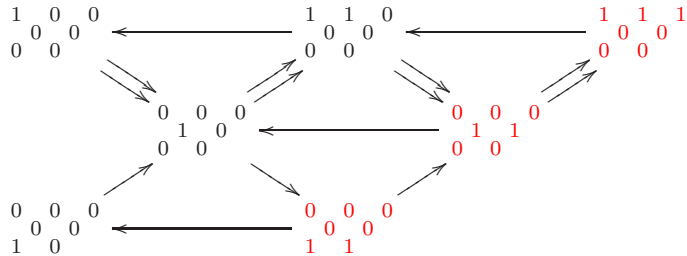
\xymatrix@-1pc{
{\bsm1&&0&&0\\&0&&0\\0&&0\esm} \ar@<0.5ex>[dr] \ar@<-0.5ex>[dr]
&& {\bsm1&&1&&0\\&0&&0\\0&&0\esm} \ar@<0.5ex>[dr] \ar@<-0.5ex>[dr]\ar[ll]&&
{\red\bsm1&&1&&1\\&0&&0\\0&&0\esm} \ar[ll]\\
&{\bsm0&&0&&0\\&1&&0\\0&&0\esm} \ar@<0.5ex>[ur]\ar@<-0.5ex>[ur]\ar[dr]&& 
{\red\bsm0&&0&&0\\&1&&1\\0&&0\esm} \ar[ll]\ar@<0.5ex>[ur]\ar@<-0.5ex>[ur]\\
{\bsm0&&0&&0\\&0&&0\\1&&0\esm}\ar[ur] &&
{\red\bsm0&&0&&0\\&0&&0\\1&&1\esm}\ar[ll]\ar[ur]
}
$$
The $\Delta$-dimension vectors 
associated to the indecomposable $\CC_M$-projectives
are labelled in red colour.

An easy calculation in the mesh category shows that
$$
d_\Delta = \bsm 23&&6&&1\\&14&&3\\11&&4\esm.
$$
Again, we mutate the $\LL$-module $T_k$ where
$$
\dimv_\Delta(\Hom_\LL(T_k,T_M)) = {\bsm1&&1&&0\\&0&&0\\0&&0\esm}.
$$
We get
$$
\dimv_\Delta(\Hom_\LL(T_k^*,T_M)) = 
- {\bsm1&&1&&0\\&0&&0\\0&&0\esm} +
{\bsm1&&1&&1\\&2&&0\\0&&0\esm} =
{\bsm0&&0&&1\\&2&&0\\0&&0\esm}
$$
since
$$
{\bsm1&&1&&1\\&0&&0\\0&&0\esm}\cdot d_\Delta + 
2 \cdot {\red\bsm0&&0&&0\\&1&&0\\0&&0\esm}\cdot d_\Delta = 58 >
57 =
{\red\bsm1&&0&&0\\&0&&0\\0&&0\esm} \cdot d_\Delta
+ 2 \cdot {\bsm0&&0&&0\\&1&&1\\0&&0\esm} \cdot d_\Delta.
$$ 

In this easy example, the calculation of the above inequality 
was not necessary.
Using $\Delta$-dimension vectors one can see immediately
that
$$
- {\bsm1&&1&&0\\&0&&0\\0&&0\esm} +
{\bsm1&&0&&0\\&0&&0\\0&&0\esm} + 
2 \cdot {\red\bsm0&&0&&0\\&1&&1\\0&&0\esm}
$$
contains a negative entry, so the other option, namely taking the arrows 
in $\GG_T''$ ending in 
$$
{\bsm1&&1&&0\\&0&&0\\0&&0\esm},
$$ 
is the correct one.

This shows that in some situations 
it can be more convenient to calculate mutations
using $\Delta$-dimension vectors in contrast to using ordinary
dimension vectors as in Section~\ref{section15}.


{\Large\section{A sequence of mutations from $T_M$ to $T_M^\vee$}
\label{section17}}


\subsection{The algorithm}\label{algorithm}
Let $M = M_1 \oplus \cdots \oplus M_r$ be a terminal $KQ$-module,
and (as before) let $Q_M^\op$ 
\index{$Q_M^\op$}
be the full subquiver of $\GG_{T_M}$
with vertices $T_{i,[t_i,t_i]}$ where $1 \le i \le n$.

Without loss of generality we assume that the vertices $1,\ldots,n$ of 
$Q_M$ 
are numbered in such a way that $1$ is a sink of $Q_M$
and $k+1$
is a sink of the quiver
$$
\sigma_k \cdots \sigma_2\sigma_1(Q_M)
$$
for all $1 \le k \le n-1$.
We call $1 < 2 < \cdots < n$ a $Q_M$-{\it adapted ordering} of the vertices
of $Q_M$.
(If $k$ is a vertex of a quiver $\GG$, then $\sigma_k(\GG)$
is the quiver obtained from $\GG$ by reversing all arrows ending or
starting at $k$.
For brevity 
we just wrote $i$ for the vertex $T_{i,[t_i,t_i]}$ of $Q_M$ and
$Q_M^\op$.)

Now we describe an algorithm which will yield a directed path 
in the Hasse quiver of the partial ordering $\T_B^{\rm cl}$ from
the (unique) maximal element $T_M$ to the (unique) minimal element $T_M^\vee$.
The proof is done by induction on the number $r-n$ of indecomposable
non-injective direct summands of $M$.
This is left as an exercise to the reader.

In the following, we just ignore the symbols of the form 
$T_{i,[a,b]}$ in case $a<0$ or $b<0$.

\noindent{\bf Step 1}:
We mutate the following 
$$
r_1 := \sum_{i=1}^n t_i 
$$ 
vertices of $\GG_{T_M}$ in the given order:
\begin{align*}
T_{1,[t_1,t_1]},
T_{1,[t_1-1,t_1]},
T_{1,[t_1-2,t_1]},
\ldots,
T_{1,[1,t_1]},\\
T_{2,[t_2,t_2]},
T_{2,[t_2-1,t_2]},
T_{2,[t_2-2,t_2]},
\ldots,
T_{2,[1,t_2]},\\
\cdots\\
T_{n,[t_n,t_n]},
T_{n,[t_n-1,t_n]},
T_{n,[t_n-2,t_n]},
\ldots,
T_{n,[1,t_n]}
\end{align*}
We obtain a new quiver $\GG_{T_M}^1$ with $r_1$ new vertices
\begin{align*}
T_{1,[t_1-1,t_1-1]},
T_{1,[t_1-2,t_1-1]},
T_{1,[t_1-3,t_1-1]},
\ldots,
T_{1,[0,t_1-1]},\\
T_{2,[t_2-1,t_2-1]},
T_{2,[t_2-2,t_2-1]},
T_{2,[t_2-3,t_2-1]},
\ldots,
T_{2,[0,t_2-1]},\\
\cdots\\
T_{n,[t_n-1,t_n-1]},
T_{n,[t_n-2,t_n-1]},
T_{n,[t_n-3,t_n-1]},
\ldots,
T_{n,[0,t_n-1]}
\end{align*}

\noindent{\bf Step 2}:
We mutate the following 
$$
r_2 := \sum_{i=1}^n \max\{0,t_i-1\}
$$ 
vertices of $\GG_{T_M}^1$ in the following order:
\begin{align*}
T_{1,[t_1-1,t_1-1]},
T_{1,[t_1-2,t_1-1]},
T_{1,[t_1-3,t_1-1]},
\ldots,
T_{1,[1,t_1-1]},\\
T_{2,[t_2-1,t_2-1]},
T_{2,[t_2-2,t_2-1]},
T_{2,[t_2-3,t_2-1]},
\ldots,
T_{2,[1,t_2-1]},\\
\cdots\\
T_{n,[t_n-1,t_n-1]},
T_{n,[t_n-2,t_n-1]},
T_{n,[t_n-3,t_n-1]},
\ldots,
T_{n,[1,t_n-1]}
\end{align*}
We obtain a new quiver $\GG_{T_M}^2$ with $r_2$ 
new vertices
\begin{align*}
T_{1,[t_1-2,t_1-2]},
T_{1,[t_1-3,t_1-2]},
T_{1,[t_1-4,t_1-2]},
\ldots,
T_{1,[0,t_1-2]},\\
T_{2,[t_2-2,t_2-2]},
T_{2,[t_2-3,t_2-2]},
T_{2,[t_2-4,t_2-2]},
\ldots,
T_{2,[0,t_2-2]},\\
\cdots\\
T_{n,[t_n-2,t_n-2]},
T_{n,[t_n-3,t_n-2]},
T_{n,[t_n-4,t_n-2]},
\ldots,
T_{n,[0,t_n-2]}
\end{align*}

\noindent{\bf Step k}:
We mutate the following 
$$
r_k := \sum_{i=1}^n \max\{0,t_i-(k-1)\}
$$ 
vertices of $\GG_{T_M}^{k-1}$ in the following order:
\begin{align*}
T_{1,[t_1-(k-1),t_1-(k-1)]},
T_{1,[t_1-k,t_1-(k-1)]},
T_{1,[t_1-(k+1),t_1-(k-1)]},
\ldots,
T_{1,[1,t_1-(k-1)]},\\
T_{2,[t_2-(k-1),t_2-(k-1)]},
T_{2,[t_2-k,t_2-(k-1)]},
T_{2,[t_2-(k+1),t_2-(k-)1]},
\ldots,
T_{2,[1,t_2-(k-1)]},\\
\cdots\\
T_{n,[t_n-(k-1),t_n-(k-1)]},
T_{n,[t_n-k,t_n-(k-1)]},
T_{n,[t_n-(k+1),t_n-(k-1)]},
\ldots,
T_{n,[1,t_n-(k-1)]}
\end{align*}
We obtain a new quiver $\GG_{T_M}^k$ with $r_k$ 
new vertices
\begin{align*}
T_{1,[t_1-k,t_1-k]},
T_{1,[t_1-(k+1),t_1-k]},
T_{1,[t_1-(k+2),t_1-k]},
\ldots,
T_{1,[0,t_1-k]},\\
T_{2,[t_2-k,t_2-k]},
T_{2,[t_2-(k+1),t_2-k]},
T_{2,[t_2-(k+2),t_2-k]},
\ldots,
T_{2,[0,t_2-2]},\\
\cdots\\
T_{n,[t_n-k,t_n-k]},
T_{n,[t_n-(k+2),t_n-k]},
T_{n,[t_n-(k+3),t_n-k]},
\ldots,
T_{n,[0,t_n-k]}
\end{align*}
The algorithm stops when
all vertices are of the form
$T_{i,[0,b]}$ where $1 \le i \le n$ and $0 \le b \le t_i$.
This will happen after 
$$
r(M) := \sum_{i=1}^n \frac{t_i(t_i+1)}{2}
$$ 
mutations.

As an example, assume
$Q$ is a Dynkin quiver of type ${\mathbb E}_8$,
and let $M$ be the direct sum of all 120 indecomposable
$KQ$-modules.
In this case, we get $t_i = 14$ for all 8 vertices $i$ of $Q$.
Then our algorithm says that starting with $T_M$
we can reach $T_M^\vee$ after
$r(M) = 8 \cdot 105 = 840$ 
mutations.

Note that if we start with our initial 
maximal rigid module $T_M$, and if we only perform the $r(M)$ 
mutations described
in the algorithm, then we obtain the subset
$$
\left\{ T_{i,[a,b]} \mid 1 \le i \le n, 0 \le a \le b \le t_i  
\right\}
$$
of the set of indecomposable rigid modules of $\CC_M$.
In particular, this subset contains all 
modules
$(M_i,0)$ where $1 \le i \le r$, namely
$$
\left\{ (M_i,0) \mid 1 \le i \le r \right\} = 
\left\{ T_{l,[c,c]} \mid 1 \le l \le n, 0 \le c \le t_l \right\}.
$$
It follows that a given rigid module of $\CC_M$ has 
at most $n$ indecomposable direct summands of 
the form
$(M_i,0)$.
(Otherwise we would get a rigid $\C Q$-module with more
than $n$ isomorphism classes of indecomposable direct
summands, and this is not possible.)

\subsection{Generalized determinantal identities}
Recall from Section~\ref{sect_dualsemi}
the definition of $\delta_X$ for $X\in \nil(\LL)$.
The known multiplicative properties of these functions are reviewed
below in Theorem~\ref{Th:ClusChar}.
In view of these properties,
the previous explicit sequence of mutations from $T_M$
to $T_M^\vee$ yields an interesting family of 
identities satisfied by the $\delta_{T_{i,[c,d]}}$.
To avoid cumbersome notation,
in the following theorem, for $\LL$-modules $X$ and $Y$ we 
write $X \cdot Y$ instead of $\delta_X \cdot \delta_Y$.
If $c > d$, then set $T_{i,[c,d]} := 1$.

\begin{Thm}[Generalized determinantal identities]\label{detthm}
Let $M = M_1 \oplus \cdots \oplus M_r$ be a terminal $KQ$-module.
Then for $1 \le i \le n$ and $1 \le a \le b \le t_i$ we have
\begin{multline}\label{detform1}
T_{i,[a-1,b]} \cdot T_{i,[a,b-1]} = T_{i,[a,b]} \cdot T_{i,[a-1,b-1]}\\
- \prod_{i \to j} T_{j,[a+(t_j-t_i),b+(t_j-t_i)]}
\prod_{k \to i} T_{k,[a-1+(t_k-t_i),b-1+(t_k-t_i)]}
\end{multline}
where the products are taken over all arrows of the quiver $Q_M^\op$
which start and end in $i$, respectively.
\end{Thm}

\begin{proof}
This follows immediately from Theorem~\ref{Th:ClusChar} and
the algorithm described in
Section~\ref{algorithm}.
Formula~(\ref{detform1}) is just an exchange relation
corresponding to the mutation of $T_{i,[a,b]}$
with $T_{i,[a,b]}^* = T_{i,[a-1,b-1]}$.
\end{proof}

If $M$ is a terminal $KQ$-module such that all $t_i$'s are equal
to a fixed $t$, then the formula in Theorem~\ref{detthm}
simplifies as follows: 
\begin{equation}\label{Tiformula1}
T_{i,[a-1,b]} \cdot T_{i,[a,b-1]} = T_{i,[a,b]} \cdot T_{i,[a-1,b-1]}
- \prod_{i \to j} T_{j,[a,b]}
\prod_{k \to i} T_{k,[a-1,b-1]}
\end{equation}
where the products are taken over all arrows of the quiver $Q_M^\op
= Q^\op$,
which start and end in $i$, respectively.

\begin{Rem}
{\rm
Fomin and Zelevinsky \cite[Theorem 1.17]{FZ5} prove generalized
determinantal identities associated to pairs of Weyl group elements
for all Dynkin cases (including the non-simply laced cases).
Using the material of Sections~\ref{coordinateunipotent}, \ref{adaptedordering} 
below, the formula~(\ref{detform1}) can be seen as a generalization of 
some of their identities to the symmetric Kac-Moody case.
}
\end{Rem}

Note that the intervals appearing on the right hand side of 
Formula~(\ref{detform1}) all have length $b-a+1$.
On the left hand side we have the interval $[a-1,b]$ of length $b-a+2$
and the interval $[a,b-1]$ of length $b-a$.
Thus we obtain a recursive description of any 
$\delta_{T_{i,[a,b]}}$ 
in terms of the $\delta_{T_{l,[c,c]}}$.
This shows that every $\delta_{T_{i,[a,b]}}$ is a rational function
of the $\delta_{T_{l,[c,c]}}$.
Recall that each $\delta_{T_{l,[c,c]}}$ is of the form $\delta_{(M_i,0)}$ for some $i$.

In fact, 
we will show that for any $\LL$-module
$X \in \CC_M$ we have
$\delta_X \in \C[\delta_{(M_1,0)},\ldots,\delta_{(M_r,0)}]$,
see Theorem~\ref{proofmain1}.
In particular,
for all
$1 \le i \le n$ and $0 \le a \le b \le t_i$ the
rational function $\delta_{T_{i,[a,b]}}$ is
a polynomial
in $\delta_{(M_1,0)},\ldots,\delta_{(M_r,0)}$.

Another proof of the polynomiality of the $\delta_{T_{i,[a,b]}}$ 
was found by Kedem and Di~Francesco \cite{DFK},
using ideas of Fomin and Zelevinsky (in particular \cite[Lemma 4.2]{BFZ}).
We thank these four mathematicians
for communicating their insights to us at MSRI in March 2008.

\subsection{Examples}
Let $Q$ be the quiver
{\footnotesize
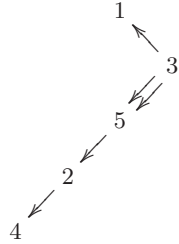
$$
\xymatrix@-1.2pc{
&&1 \\
&&& 3 \ar[ul]\ar@<0.5ex>[dl]\ar@<-0.5ex>[dl]\\
&& 5 \ar[dl]\\
& 2 \ar[dl]\\
4
}
$$
}\noindent
Let $M$ be the terminal $KQ$-module with $t_1 = t_3 = 3$,
$t_2 = t_5 = 2$ and $t_4 = 1$.
In the following we just write $i,\![a,\!b]$ instead of $T_{i,[a,b]}$.
Then the quiver $\GG_{T_M}$ of $\End_\LL(T_M)$ looks as follows:
{\footnotesize
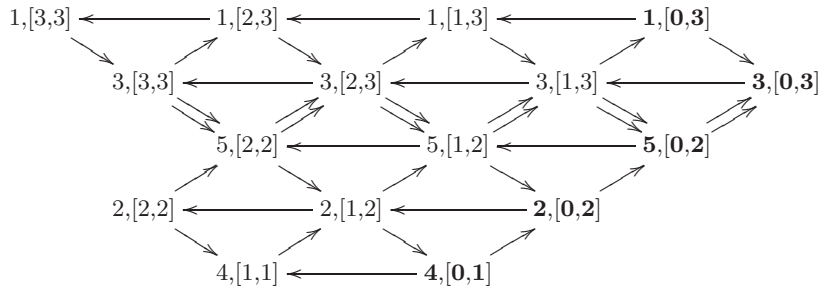
$$
\xymatrix@-1.2pc{
{1,\![3,\!3]}\ar[dr] &&
{1,\![2,\!3]} \ar[ll]\ar[dr] &&
{1,\![1,\!3]} \ar[dr]\ar[ll] &&
\mathbf{1,\![0,\!3]} \ar[ll]\ar[dr] \\
& {3,\![3,\!3]} \ar@<0.5ex>[dr]\ar@<-0.5ex>[dr]\ar[ur]&&
{3,\![2,\!3]} \ar[ll] \ar[ur]\ar@<0.5ex>[dr]\ar@<-0.5ex>[dr] &&
{3,\![1,\!3]} \ar[ll]\ar[ur]\ar@<0.5ex>[dr]\ar@<-0.5ex>[dr] &&
\mathbf{3,\![0,\!3]} \ar[ll] \\
&& {5,\![2,\!2]} \ar[dr]\ar@<0.5ex>[ur]\ar@<-0.5ex>[ur] &&
{5,\![1,\!2]} \ar[ll]\ar[dr]\ar@<0.5ex>[ur]\ar@<-0.5ex>[ur]  &&
\mathbf{5,\![0,\!2]} \ar[ll]\ar@<0.5ex>[ur]\ar@<-0.5ex>[ur] \\
& {2,\![2,\!2]}\ar[ur] \ar[dr]&& 
{2,\![1,\!2]} \ar[ll]\ar[ur]\ar[dr] && 
\mathbf{2,\![0,\!2]} \ar[ll]\ar[ur] \\
&& {4,\![1,\!1]} \ar[ur]&& 
\mathbf{4,\![0,\!1]} \ar[ll] \ar[ur] 
}
$$
}\noindent
Clearly, $Q_M^\op$ looks like this:
{\footnotesize
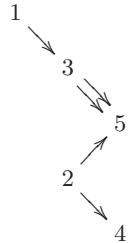
$$
\xymatrix@-1.2pc{
1 \ar[dr]\\
& 3 \ar@<0.5ex>[dr]\ar@<-0.5ex>[dr]\\
&& 5\\
& 2 \ar[dr]\ar[ur]\\
&&4
}
$$
}\noindent
(Note that we have chosen the numbering of the vertices
of $Q$ in such a way that the ordering of the vertices of 
$Q_M$ is $Q_M$-adapted.)

Starting at $\GG_{T_M}$, 
it takes 19 mutations to reach the quiver 
of $\End_\LL(T_M^\vee)$.
In the following pictures we perform sometimes two mutations
at the same time, in case these mutations do not affect each other.
The vertices we are going to mutate are marked in different colours.

{\bf Step 1}:

{\footnotesize
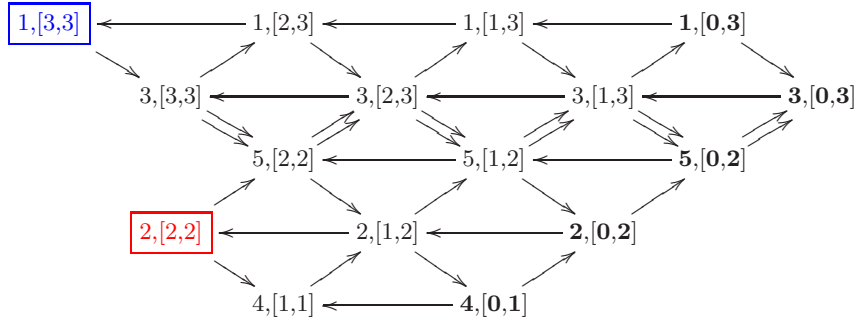
$$
\xymatrix@-1.2pc{
\blue{\fbox{1,[3,3]}}\ar[dr] &&
{1,\![2,\!3]} \ar[ll]\ar[dr] &&
{1,\![1,\!3]} \ar[dr]\ar[ll] &&
\mathbf{1,\![0,\!3]} \ar[ll]\ar[dr] \\
& {3,\![3,\!3]} \ar@<0.5ex>[dr]\ar@<-0.5ex>[dr]\ar[ur]&&
{3,\![2,\!3]} \ar[ll] \ar[ur]\ar@<0.5ex>[dr]\ar@<-0.5ex>[dr] &&
{3,\![1,\!3]} \ar[ll]\ar[ur]\ar@<0.5ex>[dr]\ar@<-0.5ex>[dr] &&
\mathbf{3,\![0,\!3]} \ar[ll] \\
&& {5,\![2,\!2]} \ar[dr]\ar@<0.5ex>[ur]\ar@<-0.5ex>[ur] &&
{5,\![1,\!2]} \ar[ll]\ar[dr]\ar@<0.5ex>[ur]\ar@<-0.5ex>[ur]  &&
\mathbf{5,\![0,\!2]} \ar[ll]\ar@<0.5ex>[ur]\ar@<-0.5ex>[ur] \\
& \red{\fbox{2,[2,2]}}\ar[ur] \ar[dr]&& 
{2,\![1,\!2]} \ar[ll]\ar[ur]\ar[dr] && 
\mathbf{2,\![0,\!2]} \ar[ll]\ar[ur] \\
&& {4,\![1,\!1]} \ar[ur]&& 
\mathbf{4,\![0,\!1]} \ar[ll] \ar[ur] 
}
$$
}\noindent

{\footnotesize
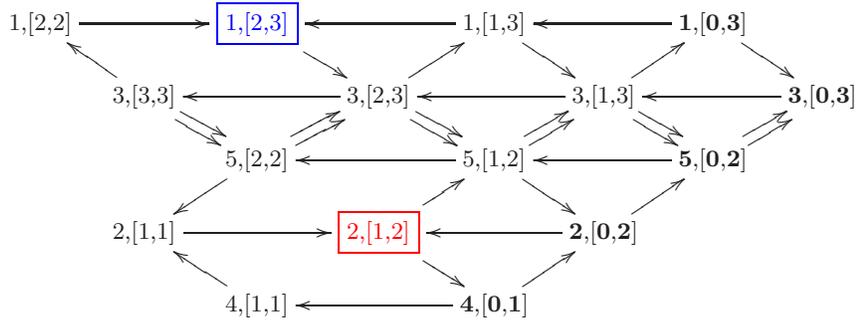
$$
\xymatrix@-1.2pc{
{1,\![2,\!2]}\ar[rr] &&
\blue{\fbox{1,[2,3]}} \ar[dr] &&
{1,\![1,\!3]} \ar[dr]\ar[ll] &&
\mathbf{1,\![0,\!3]} \ar[ll]\ar[dr] \\
& {3,\![3,\!3]} \ar@<0.5ex>[dr]\ar@<-0.5ex>[dr]\ar[ul]&&
{3,\![2,\!3]} \ar[ll] \ar[ur]\ar@<0.5ex>[dr]\ar@<-0.5ex>[dr] &&
{3,\![1,\!3]} \ar[ll]\ar[ur]\ar@<0.5ex>[dr]\ar@<-0.5ex>[dr] &&
\mathbf{3,\![0,\!3]} \ar[ll] \\
&& {5,\![2,\!2]} \ar[dl]\ar@<0.5ex>[ur]\ar@<-0.5ex>[ur] &&
{5,\![1,\!2]} \ar[ll]\ar[dr]\ar@<0.5ex>[ur]\ar@<-0.5ex>[ur]  &&
\mathbf{5,\![0,\!2]} \ar[ll]\ar@<0.5ex>[ur]\ar@<-0.5ex>[ur]\\
& {2,\![1,\!1]}\ar[rr] && 
\red{\fbox{2,[1,2]}} \ar[ur]\ar[dr] && 
\mathbf{2,\![0,\!2]} \ar[ll]\ar[ur]\\
&& {4,\![1,\!1]} \ar[ul]&& 
\mathbf{4,\![0,\!1]} \ar[ll]\ar[ur]
}
$$
}\noindent

{\footnotesize
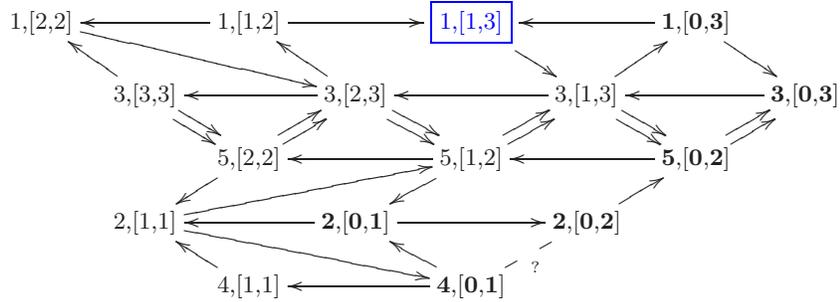
$$
\xymatrix@-1.2pc{
{1,\![2,\!2]}\ar[drrr] &&
{1,\![1,\!2]} \ar[rr]\ar[ll] &&
\blue{\fbox{1,[1,3]}} \ar[dr] &&
\mathbf{1,\![0,\!3]} \ar[ll]\ar[dr]\\
& {3,\![3,\!3]} \ar@<0.5ex>[dr]\ar@<-0.5ex>[dr]\ar[ul]&&
{3,\![2,\!3]} \ar[ll] \ar[ul]\ar@<0.5ex>[dr]\ar@<-0.5ex>[dr] &&
{3,\![1,\!3]} \ar[ll]\ar[ur]\ar@<0.5ex>[dr]\ar@<-0.5ex>[dr] &&
\mathbf{3,\![0,\!3]} \ar[ll]\\
&& {5,\![2,\!2]} \ar[dl]\ar@<0.5ex>[ur]\ar@<-0.5ex>[ur] &&
{5,\![1,\!2]} \ar[ll]\ar[dl]\ar@<0.5ex>[ur]\ar@<-0.5ex>[ur]  &&
\mathbf{5,\![0,\!2]} \ar[ll]\ar@<0.5ex>[ur]\ar@<-0.5ex>[ur] \\
& {2,\![1,\!1]}\ar[drrr]\ar[urrr] && 
\mathbf{2,\![0,\!1]} \ar[rr]\ar[ll] && 
\mathbf{2,\![0,\!2]} \ar[ur]\\
&& {4,\![1,\!1]} \ar[ul]&& 
\mathbf{4,\![0,\!1]} \ar[ll]\ar@{--}_?[ur]\ar[ul]
}
$$
}\noindent

{\footnotesize
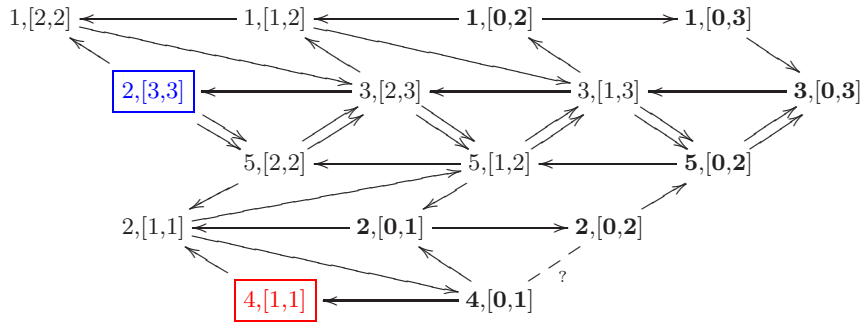
$$
\xymatrix@-1.2pc{
{1,\![2,\!2]}\ar[drrr] &&
{1,\![1,\!2]} \ar[ll]\ar[drrr] &&
\mathbf{1,\![0,\!2]} \ar[rr]\ar[ll] &&
\mathbf{1,\![0,\!3]} \ar[dr]\\
& \blue{\fbox{2,[3,3]}} \ar@<0.5ex>[dr]\ar@<-0.5ex>[dr]\ar[ul]&&
{3,\![2,\!3]} \ar[ll] \ar[ul]\ar@<0.5ex>[dr]\ar@<-0.5ex>[dr] &&
{3,\![1,\!3]} \ar[ll]\ar[ul]\ar@<0.5ex>[dr]\ar@<-0.5ex>[dr] &&
\mathbf{3,\![0,\!3]} \ar[ll]\\
&& {5,\![2,\!2]} \ar[dl]\ar@<0.5ex>[ur]\ar@<-0.5ex>[ur] &&
{5,\![1,\!2]} \ar[ll]\ar[dl]\ar@<0.5ex>[ur]\ar@<-0.5ex>[ur]  &&
\mathbf{5,\![0,\!2]} \ar[ll] \ar@<0.5ex>[ur]\ar@<-0.5ex>[ur] \\
& {2,\![1,\!1]}\ar[drrr]\ar[urrr] && 
\mathbf{2,\![0,\!1]} \ar[rr]\ar[ll] && 
\mathbf{2,\![0,\!2]} \ar[ur]\\
&& \red{\fbox{4,[1,1]}} \ar[ul]&& 
\mathbf{4,\![0,\!1]} \ar[ll]\ar[ul]\ar@{--}_?[ur] 
}
$$
}\noindent

{\footnotesize
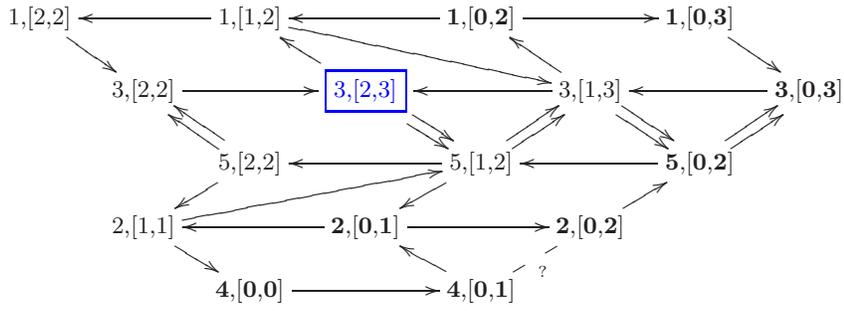
$$
\xymatrix@-1.2pc{
{1,\![2,\!2]}\ar[dr] &&
{1,\![1,\!2]} \ar[ll]\ar[drrr] &&
\mathbf{1,\![0,\!2]} \ar[rr]\ar[ll] &&
\mathbf{1,\![0,\!3]} \ar[dr] \\
& {3,\![2,\!2]} \ar[rr] &&
\blue{\fbox{3,[2,3]}} \ar[ul]\ar@<0.5ex>[dr]\ar@<-0.5ex>[dr] &&
{3,\![1,\!3]} \ar[ll]\ar[ul]\ar@<0.5ex>[dr]\ar@<-0.5ex>[dr] &&
\mathbf{3,\![0,\!3]} \ar[ll]\\
&& {5,\![2,\!2]} \ar[dl]\ar@<0.5ex>[ul]\ar@<-0.5ex>[ul] &&
{5,\![1,\!2]} \ar[ll]\ar[dl]\ar@<0.5ex>[ur]\ar@<-0.5ex>[ur]  &&
\mathbf{5,\![0,\!2]} \ar[ll]\ar@<0.5ex>[ur]\ar@<-0.5ex>[ur] \\
& {2,\![1,\!1]}\ar[dr]\ar[urrr] && 
\mathbf{2,\![0,\!1]} \ar[rr]\ar[ll] && 
\mathbf{2,\![0,\!2]} \ar[ur]\\
&& \mathbf{4,\![0,\!0]} \ar[rr]&& 
\mathbf{4,\![0,\!1]} \ar[ul] \ar@{--}_?[ur] 
}
$$
}\noindent

{\footnotesize
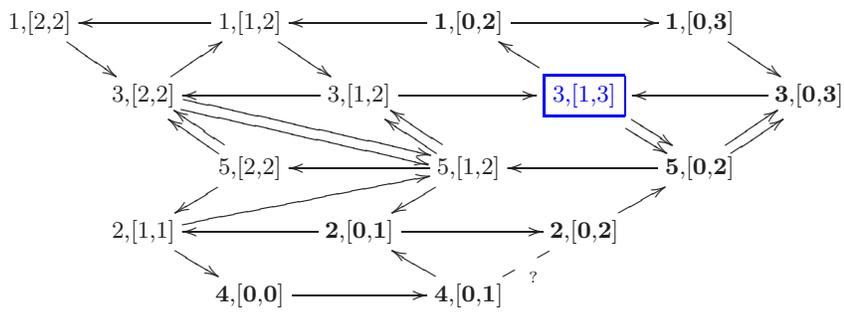
$$
\xymatrix@-1.2pc{
{1,\![2,\!2]}\ar[dr] &&
{1,\![1,\!2]} \ar[ll]\ar[dr] &&
\mathbf{1,\![0,\!2]} \ar[rr]\ar[ll] &&
\mathbf{1,\![0,\!3]} \ar[dr] \\
& {3,\![2,\!2]}\ar[ur]\ar@<0.5ex>[drrr]\ar@<-0.5ex>[drrr]  &&
{3,\![1,\!2]} \ar[ll]\ar[rr] &&
\blue{\fbox{3,[1,3]}} \ar[ul]\ar@<0.5ex>[dr]\ar@<-0.5ex>[dr] &&
\mathbf{3,\![0,\!3]} \ar[ll]\\
&& {5,\![2,\!2]} \ar[dl]\ar@<0.5ex>[ul]\ar@<-0.5ex>[ul] &&
{5,\![1,\!2]} \ar[ll]\ar[dl]\ar@<0.5ex>[ul]\ar@<-0.5ex>[ul]  &&
\mathbf{5,\![0,\!2]} \ar[ll] \ar@<0.5ex>[ur]\ar@<-0.5ex>[ur] \\
& {2,\![1,\!1]}\ar[dr]\ar[urrr] && 
\mathbf{2,\![0,\!1]} \ar[rr]\ar[ll] && 
\mathbf{2,\![0,\!2]} \ar[ur] \\
&& \mathbf{4,\![0,\!0]} \ar[rr]&& 
\mathbf{4,\![0,\!1]} \ar[ul] \ar@{--}_?[ur]
}
$$
}\noindent

{\footnotesize
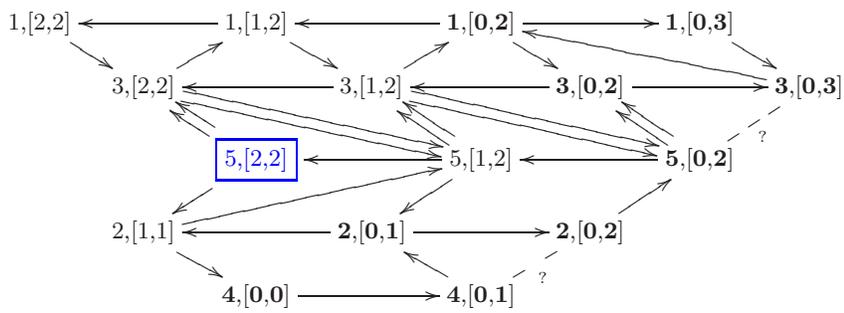
$$
\xymatrix@-1.2pc{
{1,\![2,\!2]}\ar[dr] &&
{1,\![1,\!2]} \ar[ll]\ar[dr] &&
\mathbf{1,\![0,\!2]} \ar[rr]\ar[ll]\ar[dr] &&
\mathbf{1,\![0,\!3]} \ar[dr]\\
& {3,\![2,\!2]}\ar[ur]\ar@<0.5ex>[drrr]\ar@<-0.5ex>[drrr]  &&
{3,\![1,\!2]} \ar[ll]\ar[ur]\ar@<0.5ex>[drrr]\ar@<-0.5ex>[drrr] &&
\mathbf{3,\![0,\!2]} \ar[ll]\ar[rr] &&
\mathbf{3,\![0,\!3]} \ar[ulll]\\
&& \blue{\fbox{5,[2,2]}} \ar[dl]\ar@<0.5ex>[ul]\ar@<-0.5ex>[ul] &&
{5,\![1,\!2]} \ar[ll]\ar[dl]\ar@<0.5ex>[ul]\ar@<-0.5ex>[ul]  &&
\mathbf{5,\![0,\!2]} \ar[ll]\ar@<0.5ex>[ul]\ar@<-0.5ex>[ul] \ar@{--}_?[ur]\\
& {2,\![1,\!1]}\ar[dr]\ar[urrr] && 
\mathbf{2,\![0,\!1]} \ar[rr]\ar[ll] && 
\mathbf{2,\![0,\!2]} \ar[ur]\\
&& \mathbf{4,\![0,\!0]} \ar[rr]&& 
\mathbf{4,\![0,\!1]} \ar[ul]\ar@{--}_?[ur] 
}
$$
}\noindent

{\footnotesize
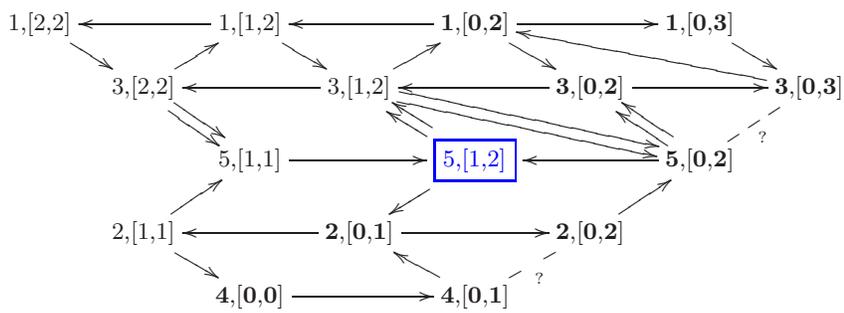
$$
\xymatrix@-1.2pc{
{1,\![2,\!2]}\ar[dr] &&
{1,\![1,\!2]} \ar[ll]\ar[dr] &&
\mathbf{1,\![0,\!2]} \ar[rr]\ar[ll]\ar[dr] &&
\mathbf{1,\![0,\!3]} \ar[dr]\\
& {3,\![2,\!2]}\ar[ur]\ar@<0.5ex>[dr]\ar@<-0.5ex>[dr]  &&
{3,\![1,\!2]} \ar[ll]\ar[ur]\ar@<0.5ex>[drrr]\ar@<-0.5ex>[drrr] &&
\mathbf{3,\![0,\!2]} \ar[ll]\ar[rr] &&
\mathbf{3,\![0,\!3]} \ar[ulll]\\
&& {5,\![1,\!1]} \ar[rr] &&
\blue{\fbox{5,[1,2]}} \ar[dl]\ar@<0.5ex>[ul]\ar@<-0.5ex>[ul]  &&
\mathbf{5,\![0,\!2]} \ar[ll]\ar@<0.5ex>[ul]\ar@<-0.5ex>[ul] \ar@{--}_?[ur]\\
& {2,\![1,\!1]}\ar[dr]\ar[ur] && 
\mathbf{2,\![0,\!1]} \ar[rr]\ar[ll] && 
\mathbf{2,\![0,\!2]} \ar[ur]\\
&& \mathbf{4,\![0,\!0]} \ar[rr]&& 
\mathbf{4,\![0,\!1]} \ar[ul] \ar@{--}_?[ur] 
}
$$
}\noindent

{\bf Step 2}:

{\footnotesize
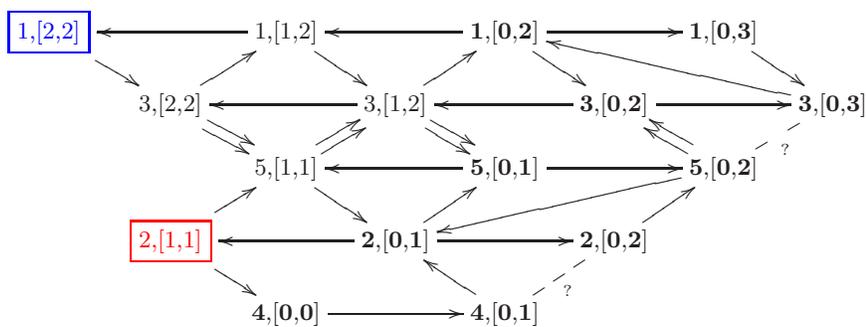
$$
\xymatrix@-1.2pc{
\blue{\fbox{1,[2,2]}}\ar[dr] &&
{1,\![1,\!2]} \ar[ll]\ar[dr] &&
\mathbf{1,\![0,\!2]} \ar[rr]\ar[ll]\ar[dr] &&
\mathbf{1,\![0,\!3]} \ar[dr]\\
& {3,\![2,\!2]}\ar[ur]\ar@<0.5ex>[dr]\ar@<-0.5ex>[dr]  &&
{3,\![1,\!2]} \ar[ll]\ar[ur]\ar@<0.5ex>[dr]\ar@<-0.5ex>[dr] &&
\mathbf{3,\![0,\!2]} \ar[ll]\ar[rr] &&
\mathbf{3,\![0,\!3]} \ar[ulll]\\
&& {5,\![1,\!1]} \ar@<0.5ex>[ur]\ar@<-0.5ex>[ur]\ar[dr] &&
\mathbf{5,\![0,\!1]} \ar[ll]\ar[rr]&&
\mathbf{5,\![0,\!2]} \ar@<0.5ex>[ul]\ar@<-0.5ex>[ul]\ar[dlll]\ar@{--}_?[ur]\\
& \red{\fbox{2,[1,1]}}\ar[dr]\ar[ur] && 
\mathbf{2,\![0,\!1]} \ar[rr]\ar[ll]\ar[ur] && 
\mathbf{2,\![0,\!2]} \ar[ur]\\
&& \mathbf{4,\![0,\!0]} \ar[rr]&& 
\mathbf{4,\![0,\!1]} \ar[ul]\ar@{--}_?[ur] 
}
$$
}\noindent

{\footnotesize
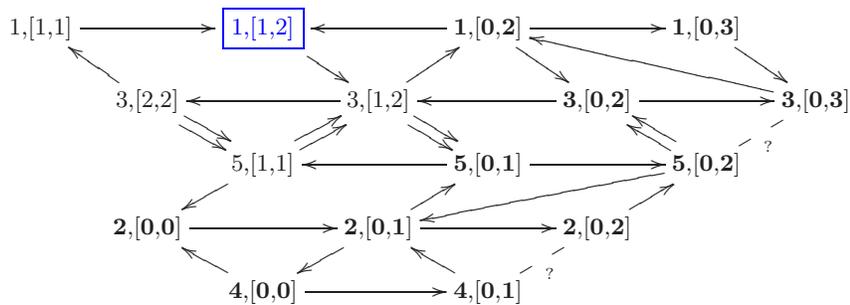
$$
\xymatrix@-1.2pc{
{1,\![1,\!1]}\ar[rr] &&
\blue{\fbox{1,[1,2]}} \ar[dr] &&
\mathbf{1,\![0,\!2]} \ar[rr]\ar[ll]\ar[dr] &&
\mathbf{1,\![0,\!3]} \ar[dr]\\
& {3,\![2,\!2]}\ar[ul]\ar@<0.5ex>[dr]\ar@<-0.5ex>[dr]  &&
{3,\![1,\!2]} \ar[ll]\ar[ur]\ar@<0.5ex>[dr]\ar@<-0.5ex>[dr] &&
\mathbf{3,\![0,\!2]} \ar[ll]\ar[rr] &&
\mathbf{3,\![0,\!3]} \ar[ulll]\\
&& {5,\![1,\!1]} \ar@<0.5ex>[ur]\ar@<-0.5ex>[ur]\ar[dl] &&
\mathbf{5,\![0,\!1]} \ar[ll]\ar[rr]&&
\mathbf{5,\![0,\!2]} \ar@<0.5ex>[ul]\ar@<-0.5ex>[ul]\ar[dlll]\ar@{--}_?[ur]\\
& \mathbf{2,\![0,\!0]} \ar[rr] && 
\mathbf{2,\![0,\!1]} \ar[rr]\ar[ur]\ar[dl] && 
\mathbf{2,\![0,\!2]} \ar[ur]\\
&& \mathbf{4,\![0,\!0]} \ar[rr]\ar[ul]&& 
\mathbf{4,\![0,\!1]} \ar[ul] \ar@{--}_?[ur] 
}
$$
}\noindent

{\footnotesize
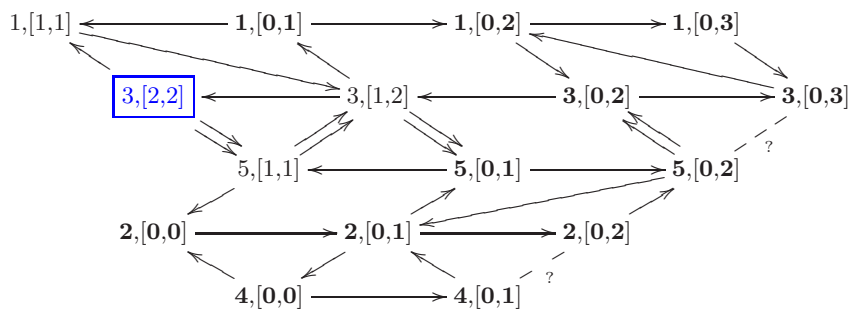
$$
\xymatrix@-1.2pc{
{1,\![1,\!1]}\ar[drrr] &&
\mathbf{1,\![0,\!1]} \ar[rr]\ar[ll] &&
\mathbf{1,\![0,\!2]} \ar[rr]\ar[dr] &&
\mathbf{1,\![0,\!3]} \ar[dr]\\
& \blue{\fbox{3,[2,2]}}\ar[ul]\ar@<0.5ex>[dr]\ar@<-0.5ex>[dr]  &&
{3,\![1,\!2]} \ar[ll]\ar[ul]\ar@<0.5ex>[dr]\ar@<-0.5ex>[dr] &&
\mathbf{3,\![0,\!2]} \ar[ll]\ar[rr] &&
\mathbf{3,\![0,\!3]} \ar[ulll]\\
&& {5,\![1,\!1]} \ar@<0.5ex>[ur]\ar@<-0.5ex>[ur]\ar[dl] &&
\mathbf{5,\![0,\!1]} \ar[ll]\ar[rr]&&
\mathbf{5,\![0,\!2]} \ar@<0.5ex>[ul]\ar@<-0.5ex>[ul]\ar[dlll]\ar@{--}_?[ur]\\
& \mathbf{2,\![0,\!0]} \ar[rr] && 
\mathbf{2,\![0,\!1]} \ar[rr]\ar[ur]\ar[dl] && 
\mathbf{2,\![0,\!2]} \ar[ur]\\
&& \mathbf{4,\![0,\!0]} \ar[rr]\ar[ul]&& 
\mathbf{4,\![0,\!1]} \ar[ul]\ar@{--}_?[ur] 
}
$$
}\noindent

{\footnotesize
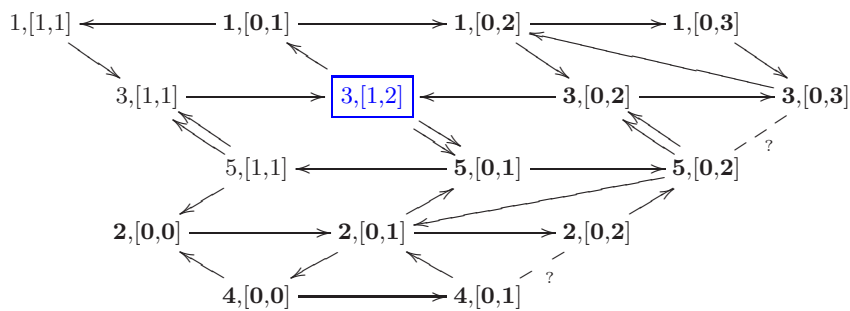
$$
\xymatrix@-1.2pc{
{1,\![1,\!1]}\ar[dr] &&
\mathbf{1,\![0,\!1]} \ar[rr]\ar[ll] &&
\mathbf{1,\![0,\!2]} \ar[rr]\ar[dr] &&
\mathbf{1,\![0,\!3]} \ar[dr]\\
& {3,\![1,\!1]} \ar[rr] &&
\blue{\fbox{3,[1,2]}} \ar[ul]\ar@<0.5ex>[dr]\ar@<-0.5ex>[dr] &&
\mathbf{3,\![0,\!2]} \ar[ll]\ar[rr] &&
\mathbf{3,\![0,\!3]} \ar[ulll]\\
&& {5,\![1,\!1]} \ar@<0.5ex>[ul]\ar@<-0.5ex>[ul]\ar[dl] &&
\mathbf{5,\![0,\!1]} \ar[ll]\ar[rr]&&
\mathbf{5,\![0,\!2]} \ar@<0.5ex>[ul]\ar@<-0.5ex>[ul]\ar[dlll]\ar@{--}_?[ur]\\
& \mathbf{2,\![0,\!0]} \ar[rr] && 
\mathbf{2,\![0,\!1]} \ar[rr]\ar[ur]\ar[dl] && 
\mathbf{2,\![0,\!2]} \ar[ur]\\
&& \mathbf{4,\![0,\!0]} \ar[rr]\ar[ul]&& 
\mathbf{4,\![0,\!1]} \ar[ul]\ar@{--}_?[ur] 
}
$$
}\noindent

{\footnotesize
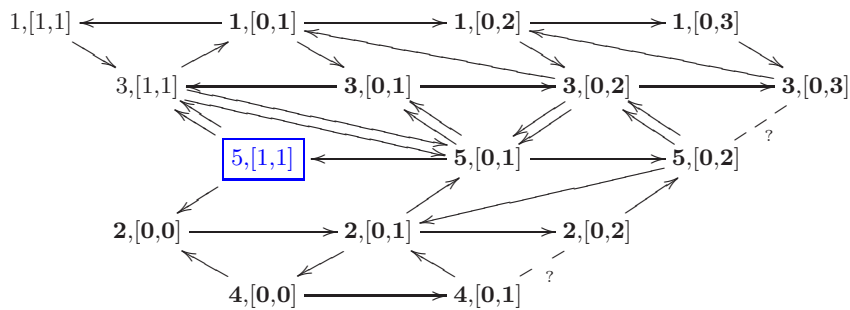
$$
\xymatrix@-1.2pc{
{1,\![1,\!1]}\ar[dr] &&
\mathbf{1,\![0,\!1]} \ar[rr]\ar[ll]\ar[dr] &&
\mathbf{1,\![0,\!2]} \ar[rr]\ar[dr] &&
\mathbf{1,\![0,\!3]} \ar[dr]\\
& {3,\![1,\!1]} \ar[ur] \ar@<0.5ex>[drrr]\ar@<-0.5ex>[drrr] &&
\mathbf{3,\![0,\!1]} \ar[ll]\ar[rr]&&
\mathbf{3,\![0,\!2]} \ar[rr]\ar[ulll] \ar@<0.5ex>[dl]\ar@<-0.5ex>[dl] &&
\mathbf{3,\![0,\!3]} \ar[ulll]\\
&& \blue{\fbox{5,[1,1]}} \ar@<0.5ex>[ul]\ar@<-0.5ex>[ul]\ar[dl] &&
\mathbf{5,\![0,\!1]} \ar[ll]\ar[rr]\ar@<0.5ex>[ul]\ar@<-0.5ex>[ul] &&
\mathbf{5,\![0,\!2]} \ar@<0.5ex>[ul]\ar@<-0.5ex>[ul]\ar[dlll]\ar@{--}_?[ur]\\
& \mathbf{2,\![0,\!0]} \ar[rr] && 
\mathbf{2,\![0,\!1]} \ar[rr]\ar[ur]\ar[dl] && 
\mathbf{2,\![0,\!2]} \ar[ur]\\
&& \mathbf{4,\![0,\!0]} \ar[rr]\ar[ul]&& 
\mathbf{4,\![0,\!1]} \ar[ul]\ar@{--}_?[ur] 
}
$$
}\noindent

{\bf Step 3}:

{\footnotesize
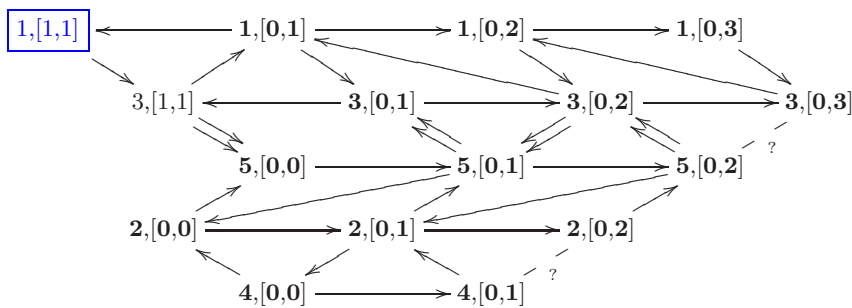
$$
\xymatrix@-1.2pc{
\blue{\fbox{1,[1,1]}}\ar[dr] &&
\mathbf{1,\![0,\!1]} \ar[rr]\ar[ll]\ar[dr] &&
\mathbf{1,\![0,\!2]} \ar[rr]\ar[dr] &&
\mathbf{1,\![0,\!3]} \ar[dr]\\
& {3,\![1,\!1]} \ar[ur] \ar@<0.5ex>[dr]\ar@<-0.5ex>[dr] &&
\mathbf{3,\![0,\!1]} \ar[ll]\ar[rr]&&
\mathbf{3,\![0,\!2]} \ar[rr]\ar[ulll] \ar@<0.5ex>[dl]\ar@<-0.5ex>[dl] &&
\mathbf{3,\![0,\!3]} \ar[ulll]\\
&& \mathbf{5,\![0,\!0]} \ar[rr]&&
\mathbf{5,\![0,\!1]} \ar[rr]\ar@<0.5ex>[ul]\ar@<-0.5ex>[ul]\ar[dlll] &&
\mathbf{5,\![0,\!2]} \ar@<0.5ex>[ul]\ar@<-0.5ex>[ul]\ar[dlll]\ar@{--}_?[ur]\\
& \mathbf{2,\![0,\!0]} \ar[rr]\ar[ur] && 
\mathbf{2,\![0,\!1]} \ar[rr]\ar[ur]\ar[dl] && 
\mathbf{2,\![0,\!2]} \ar[ur]\\
&& \mathbf{4,\![0,\!0]} \ar[rr]\ar[ul]&& 
\mathbf{4,\![0,\!1]} \ar[ul]\ar@{--}_?[ur] 
}
$$
}\noindent

{\footnotesize
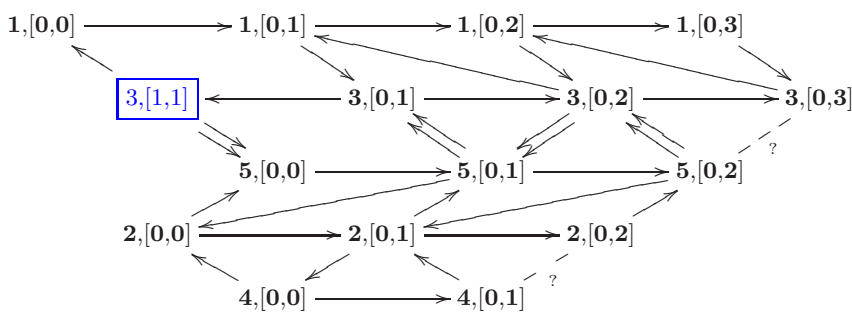
$$
\xymatrix@-1.2pc{
\mathbf{1,\![0,\!0]}\ar[rr] &&
\mathbf{1,\![0,\!1]} \ar[rr]\ar[dr] &&
\mathbf{1,\![0,\!2]} \ar[rr]\ar[dr] &&
\mathbf{1,\![0,\!3]} \ar[dr]\\
& \blue{\fbox{3,[1,1]}} \ar[ul] \ar@<0.5ex>[dr]\ar@<-0.5ex>[dr] &&
\mathbf{3,\![0,\!1]} \ar[ll]\ar[rr]&&
\mathbf{3,\![0,\!2]} \ar[rr]\ar[ulll] \ar@<0.5ex>[dl]\ar@<-0.5ex>[dl] &&
\mathbf{3,\![0,\!3]} \ar[ulll]\\
&& \mathbf{5,\![0,\!0]} \ar[rr]&&
\mathbf{5,\![0,\!1]} \ar[rr]\ar@<0.5ex>[ul]\ar@<-0.5ex>[ul]\ar[dlll] &&
\mathbf{5,\![0,\!2]} \ar@<0.5ex>[ul]\ar@<-0.5ex>[ul]\ar[dlll]\ar@{--}_?[ur]\\
& \mathbf{2,\![0,\!0]} \ar[rr]\ar[ur] && 
\mathbf{2,\![0,\!1]} \ar[rr]\ar[ur]\ar[dl] && 
\mathbf{2,\![0,\!2]} \ar[ur]\\
&& \mathbf{4,\![0,\!0]} \ar[rr]\ar[ul]&& 
\mathbf{4,\![0,\!1]} \ar[ul]\ar@{--}_?[ur] 
}
$$
}\noindent

{\footnotesize
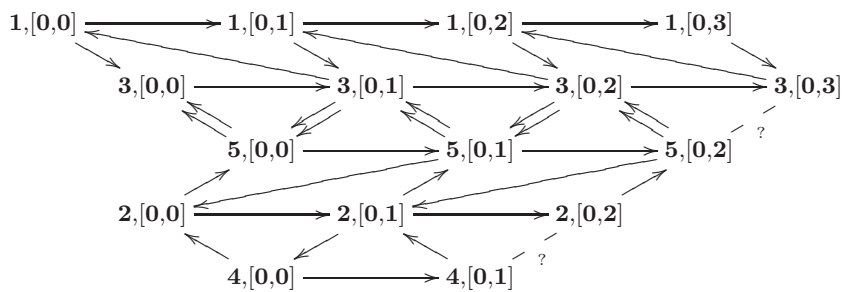
$$
\xymatrix@-1.2pc{
\mathbf{1,\![0,\!0]}\ar[rr]\ar[dr] &&
\mathbf{1,\![0,\!1]} \ar[rr]\ar[dr] &&
\mathbf{1,\![0,\!2]} \ar[rr]\ar[dr] &&
\mathbf{1,\![0,\!3]} \ar[dr]\\
& \mathbf{3,\![0,\!0]} \ar[rr]&&
\mathbf{3,\![0,\!1]} \ar[rr]\ar@<0.5ex>[dl]\ar@<-0.5ex>[dl]\ar[ulll]&&
\mathbf{3,\![0,\!2]} \ar[rr]\ar[ulll] \ar@<0.5ex>[dl]\ar@<-0.5ex>[dl] &&
\mathbf{3,\![0,\!3]} \ar[ulll]\\
&& \mathbf{5,\![0,\!0]} \ar[rr]\ar@<0.5ex>[ul]\ar@<-0.5ex>[ul]&&
\mathbf{5,\![0,\!1]} \ar[rr]\ar@<0.5ex>[ul]\ar@<-0.5ex>[ul]\ar[dlll] &&
\mathbf{5,\![0,\!2]} \ar@<0.5ex>[ul]\ar@<-0.5ex>[ul]\ar[dlll]\ar@{--}_?[ur]\\
& \mathbf{2,\![0,\!0]} \ar[rr]\ar[ur] && 
\mathbf{2,\![0,\!1]} \ar[rr]\ar[ur]\ar[dl] && 
\mathbf{2,\![0,\!2]} \ar[ur]\\
&& \mathbf{4,\![0,\!0]} \ar[rr]\ar[ul]&& 
\mathbf{4,\![0,\!1]} \ar[ul]\ar@{--}_?[ur] 
}
$$
}

We finally arrived at 
a quiver whose vertices correspond to the indecomposable direct summands of
$T_M^\vee$.
Since we know how the quiver of $\End_\LL(T_M^\vee)$ looks like,
we also obtain the missing arrows marked by $\xymatrix{\ar@{--}[r] &}$.
Namely we have an arrow
{\footnotesize
$\xymatrix{\mathbf{3,\![0,\!3]} \ar@<0.5ex>[r]\ar@<-0.5ex>[r] & \mathbf{5,\![0,\!2]}}$
} 
and two arrows
{\footnotesize
$\xymatrix{\mathbf{2,\![0,\!2]} \ar[r] & \mathbf{4,\![0,\!1]}}$
}.

Next, we discuss another example, which illustrates
why we call the formula in Theorem~\ref{detthm} a
generalized determinantal identity.
This time we will display explicitly the $\LL$-modules
occuring in our sequence of mutations from $T_M$ to $T_M^\vee$.
Let $Q$ be the quiver
{\small
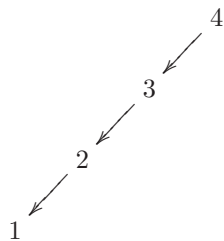
$$
\xymatrix@-0.8pc{
&&&4 \ar[dl] \\
&&3 \ar[dl] \\
&2 \ar[dl] \\
1
}
$$
}\noindent
The Auslander-Reiten quiver of $\C Q$ is:
{\small
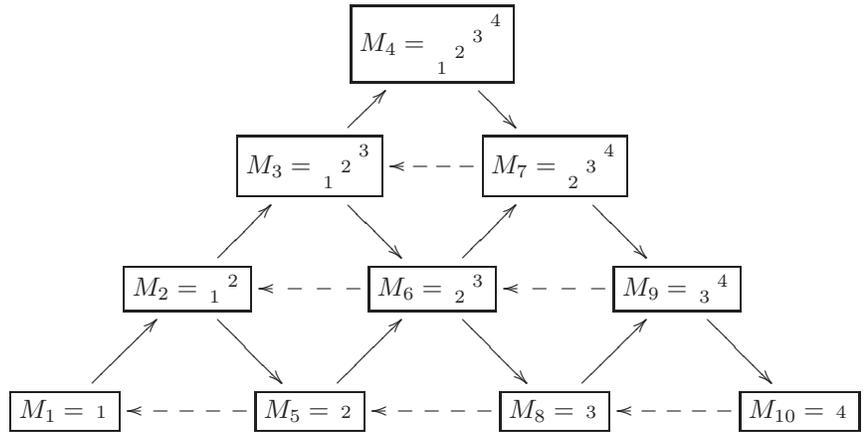
$$
\xymatrix@!@-3.8pc{
&&&{\fbox{$M_4 = \bsm &&&4\\&&3\\&2\\1 \esm$}} \ar[dr]\\
&&{\fbox{$M_3 = \bsm &&3\\&2\\1 \esm$}} \ar[ur]\ar[dr] && 
{\fbox{$M_7 = \bsm &&4\\&3\\2 \esm$}}\ar[dr] \ar@{-->}[ll]\\
& {\fbox{$M_2 = \bsm &2\\1 \esm$}} \ar[ur]\ar[dr] && 
{\fbox{$M_6 = \bsm &3\\2 \esm$}} \ar[ur]\ar[dr] \ar@{-->}[ll] && 
{\fbox{$M_9 = \bsm &4\\3 \esm$}}\ar[dr] \ar@{-->}[ll] \\
{\fbox{$M_1 = \bsm 1 \esm$}} \ar[ur] && 
{\fbox{$M_5 = \bsm 2 \esm$}} \ar[ur] \ar@{-->}[ll] && 
{\fbox{$M_8 = \bsm 3 \esm$}} \ar[ur] \ar@{-->}[ll] && 
{\fbox{$M_{10}= \bsm 4 \esm$}}  \ar@{-->}[ll]
}
$$
}

Let us display a sequence of mutations from $T_M$ to $T_M^\vee$:
{\small
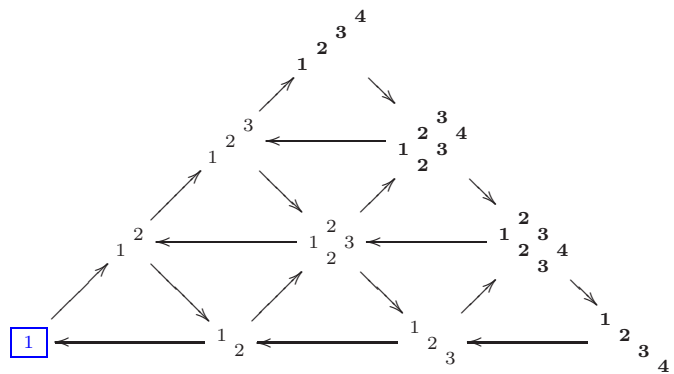
$$
\xymatrix@!@-1.8pc{
&&&{{\bsm &&&\mathbf 4\\&&\mathbf 3\\&\mathbf 2\\\mathbf 1 \esm}} \ar[dr]
\\
&&{{\bsm &&3\\&2\\1 \esm}} \ar[ur]\ar[dr] && 
{{\bsm &&\mathbf 3\\&\mathbf 2&&\mathbf 4\\\mathbf 1&&\mathbf 3\\&\mathbf 2 \esm}}\ar[dr]\ar[ll]
\\
& {{\bsm &2\\1 \esm}} \ar[ur]\ar[dr] && 
{{\bsm &2\\1&&3\\&2 \esm}} \ar[ur]\ar[dr]\ar[ll] && 
{{\bsm &\mathbf 2\\\mathbf 1&&\mathbf 3\\&\mathbf 2&&\mathbf 4\\&&\mathbf 3 \esm}}\ar[dr]\ar[ll]
\\
\blue{\fbox{$\bsm 1 \esm$}} \ar[ur] && 
{{\bsm 1\\&2 \esm}} \ar[ur]\ar[ll] && 
{{\bsm 1\\&2\\&&3 \esm}} \ar[ur]\ar[ll] && 
{{\bsm \mathbf 1\\&\mathbf 2\\&&\mathbf 3\\&&& \mathbf 4 \esm}} \ar[ll]
}
$$
}

{\small
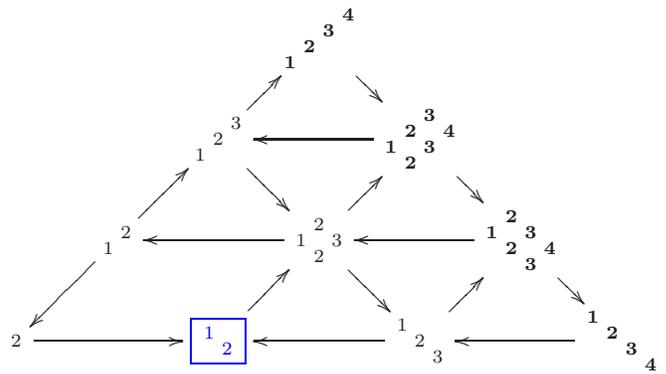
$$
\xymatrix@!@-1.8pc{
&&&{{\bsm &&&\mathbf 4\\&&\mathbf 3\\&\mathbf 2\\\mathbf 1 \esm}} \ar[dr]
\\
&&{{\bsm &&3\\&2\\1 \esm}} \ar[ur]\ar[dr] && 
{{\bsm &&\mathbf 3\\&\mathbf 2&&\mathbf 4\\\mathbf 1&&\mathbf 3\\&\mathbf 2 \esm}}\ar[dr]\ar[ll]
\\
& {{\bsm &2\\1 \esm}} \ar[ur]\ar[dl] && 
{{\bsm &2\\1&&3\\&2 \esm}} \ar[ur]\ar[dr]\ar[ll] && 
{{\bsm &\mathbf 2\\\mathbf 1&&\mathbf 3\\&\mathbf 2&&\mathbf 4\\&&\mathbf 3 \esm}}\ar[dr]\ar[ll]
\\
{{\bsm 2 \esm}} \ar[rr] && 
\blue{\fbox{$\bsm 1\\&2 \esm$}} \ar[ur] && 
{{\bsm 1\\&2\\&&3 \esm}} \ar[ur]\ar[ll] && 
{{\bsm \mathbf 1\\&\mathbf 2\\&&\mathbf 3\\&&&\mathbf 4 \esm}} \ar[ll]
}
$$
}

{\small
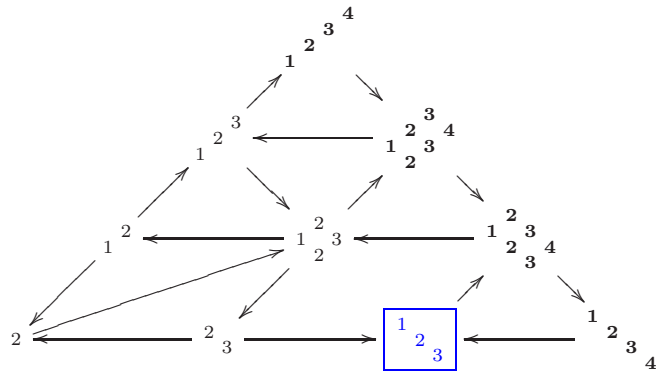
$$
\xymatrix@!@-1.8pc{
&&&{{\bsm &&&\mathbf 4\\&&\mathbf 3\\&\mathbf 2\\\mathbf 1 \esm}} \ar[dr]
\\
&&{{\bsm &&3\\&2\\1 \esm}} \ar[ur]\ar[dr] && 
{{\bsm &&\mathbf 3\\&\mathbf 2&&\mathbf 4\\\mathbf 1&&\mathbf 3\\&\mathbf 2 \esm}}\ar[dr]\ar[ll]
\\
& {{\bsm &2\\1 \esm}} \ar[ur]\ar[dl] && 
{{\bsm &2\\1&&3\\&2 \esm}} \ar[ur]\ar[ll]\ar[dl] && 
{{\bsm &\mathbf 2\\\mathbf 1&&\mathbf 3\\&\mathbf 2&&\mathbf 4\\&&\mathbf 3 \esm}}\ar[dr]\ar[ll]
\\
{{\bsm 2 \esm}} \ar[rrru]&& 
{{\bsm 2\\&3 \esm}} \ar[rr]\ar[ll] && 
\blue{\fbox{$\bsm 1\\&2\\&&3 \esm$}} \ar[ur] && 
{{\bsm \mathbf 1\\&\mathbf 2\\&&\mathbf 3\\&&&\mathbf 4 \esm}} \ar[ll]
}
$$
}

{\small
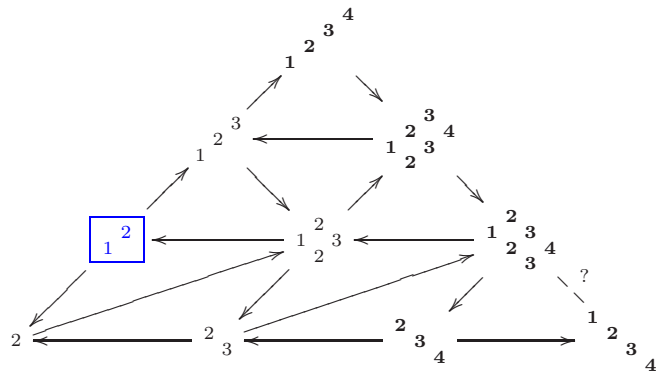
$$
\xymatrix@!@-1.8pc{
&&&{{\bsm &&&\mathbf 4\\&&\mathbf 3\\&\mathbf 2\\\mathbf 1 \esm}} \ar[dr]
\\
&&{{\bsm &&3\\&2\\1 \esm}} \ar[ur]\ar[dr] && 
{{\bsm &&\mathbf 3\\&\mathbf 2&&\mathbf 4\\\mathbf 1&&\mathbf 3\\&\mathbf 2 \esm}}\ar[dr]\ar[ll]
\\
& \blue{\fbox{$\bsm &2\\1 \esm$}} \ar[ur]\ar[dl] && 
{{\bsm &2\\1&&3\\&2 \esm}} \ar[ur]\ar[ll]\ar[dl] && 
{{\bsm &\mathbf 2\\\mathbf 1&&\mathbf 3\\&\mathbf 2&&\mathbf 4\\&&\mathbf 3 \esm}}\ar@{--}[dr]^?
\ar[ll]\ar[dl]
\\
{{\bsm 2 \esm}} \ar[rrru]&& 
{{\bsm 2\\&3 \esm}} \ar[ll]\ar[rrru] && 
{{\bsm \mathbf 2\\&\mathbf 3\\&&\mathbf 4 \esm}} \ar[rr]\ar[ll] && 
{{\bsm \mathbf 1\\&\mathbf 2\\&&\mathbf 3\\&&&\mathbf 4 \esm}}
}
$$
}

{\small
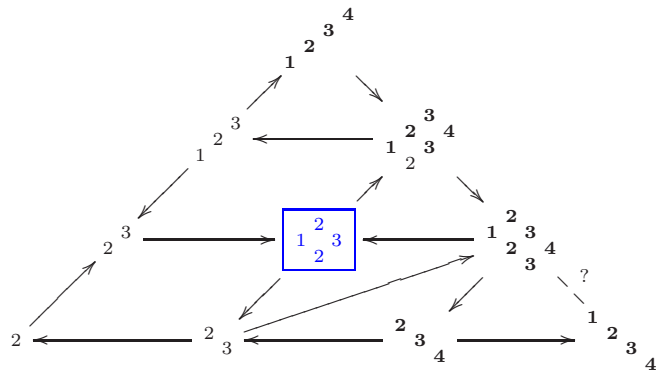
$$
\xymatrix@!@-1.8pc{
&&&{{\bsm &&&\mathbf 4\\&&\mathbf 3\\&\mathbf 2\\\mathbf 1 \esm}} \ar[dr]
\\
&&{{\bsm &&3\\&2\\1 \esm}} \ar[dl]\ar[ur] && 
{{\bsm &&\mathbf 3\\&\mathbf 2&&\mathbf 4\\\mathbf 1&&\mathbf 3\\&2 \esm}}\ar[dr]\ar[ll]
\\
& {{\bsm &3\\2 \esm}} \ar[rr] && 
\blue{\fbox{$\bsm &2\\1&&3\\&2 \esm$}} \ar[ur]\ar[dl] && 
{{\bsm &\mathbf 2\\\mathbf 1&&\mathbf 3\\&\mathbf 2&&\mathbf 4\\&&\mathbf 3 \esm}}\ar@{--}[dr]^?
\ar[ll]\ar[dl]
\\
{{\bsm 2 \esm}} \ar[ur]&& 
{{\bsm 2\\&3 \esm}} \ar[ll]\ar[rrru] && 
{{\bsm \mathbf 2\\&\mathbf 3\\&&\mathbf 4 \esm}} \ar[rr]\ar[ll] && 
{{\bsm \mathbf 1\\&\mathbf 2\\&&\mathbf 3\\&&&\mathbf 4 \esm}}
}
$$
}

{\small
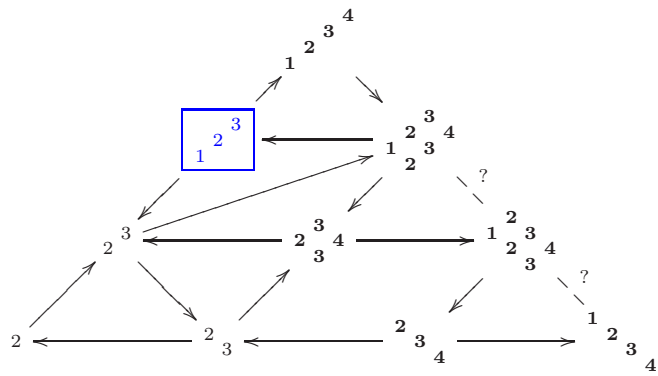
$$
\xymatrix@!@-1.8pc{
&&&{{\bsm &&&\mathbf 4\\&&\mathbf 3\\&\mathbf 2\\\mathbf 1 \esm}} \ar[dr]
\\
&&\blue{\fbox{$\bsm &&3\\&2\\1 \esm$}} \ar[dl]\ar[ur] && 
{{\bsm &&\mathbf 3\\&\mathbf 2&&\mathbf 4\\\mathbf 1&&\mathbf 3\\&\mathbf 2 \esm}}\ar@{--}[dr]^?
\ar[ll]\ar[dl]
\\
& {{\bsm &3\\2 \esm}}\ar[dr]\ar[rrru] && 
{{\bsm &\mathbf 3\\\mathbf 2&&\mathbf 4\\&\mathbf 3 \esm}} \ar[rr]\ar[ll]&& 
{{\bsm &\mathbf 2\\\mathbf 1&&\mathbf 3\\&\mathbf 2&&\mathbf 4\\&&\mathbf 3 \esm}}\ar@{--}[dr]^?
\ar[dl]
\\
{{\bsm 2 \esm}} \ar[ur]&& 
{{\bsm 2\\&3 \esm}} \ar[ll]\ar[ur] && 
{{\bsm \mathbf 2\\&\mathbf 3\\&&\mathbf 4 \esm}} \ar[rr]\ar[ll] && 
{{\bsm \mathbf 1\\&\mathbf 2\\&&\mathbf 3\\&&&\mathbf 4 \esm}}
}
$$
}

{\small
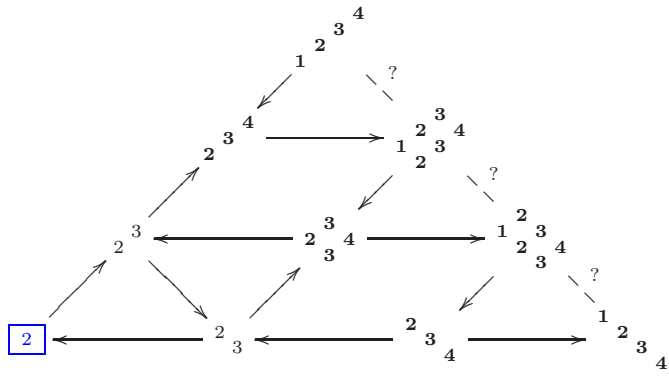
$$
\xymatrix@!@-1.8pc{
&&&{{\bsm &&&\mathbf 4\\&&\mathbf 3\\&\mathbf 2\\\mathbf 1 \esm}} \ar[dl]\ar@{--}[dr]^?
\\
&&{{\bsm &&\mathbf 4\\&\mathbf 3\\\mathbf 2 \esm}} \ar[rr] && 
{{\bsm &&\mathbf 3\\&\mathbf 2&&\mathbf 4\\\mathbf 1&&\mathbf 3\\&\mathbf 2 \esm}}\ar@{--}[dr]^?
\ar[dl]
\\
& {{\bsm &3\\2 \esm}}\ar[dr]\ar[ur] && 
{{\bsm &\mathbf 3\\\mathbf 2&&\mathbf 4\\&\mathbf 3 \esm}} \ar[rr]\ar[ll]&& 
{{\bsm &\mathbf 2\\\mathbf 1&&\mathbf 3\\&\mathbf 2&&\mathbf 4\\&&\mathbf 3 \esm}}\ar@{--}[dr]^?
\ar[dl]
\\
\blue{\fbox{$\bsm 2 \esm$}} \ar[ur]&& 
{{\bsm 2\\&3 \esm}} \ar[ll]\ar[ur] && 
{{\bsm \mathbf 2\\&\mathbf 3\\&&\mathbf 4 \esm}} \ar[rr]\ar[ll] && 
{{\bsm \mathbf 1\\&\mathbf 2\\&&\mathbf 3\\&&&\mathbf 4 \esm}}
}
$$
}

{\small
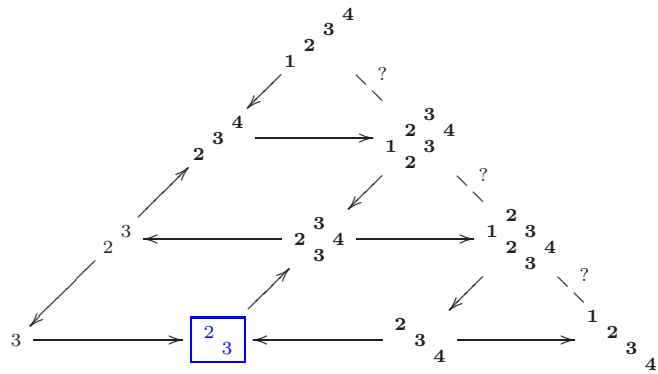
$$
\xymatrix@!@-1.8pc{
&&&{{\bsm &&&\mathbf 4\\&&\mathbf 3\\&\mathbf 2\\\mathbf 1 \esm}} \ar[dl]\ar@{--}[dr]^?
\\
&&{{\bsm &&\mathbf 4\\&\mathbf 3\\\mathbf 2 \esm}} \ar[rr] && 
{{\bsm &&\mathbf 3\\&\mathbf 2&&\mathbf 4\\\mathbf 1&&\mathbf 3\\&\mathbf 2 \esm}}\ar@{--}[dr]^?
\ar[dl]
\\
& {{\bsm &3\\2 \esm}} \ar[ur]\ar[dl] && 
{{\bsm &\mathbf 3\\\mathbf 2&&\mathbf 4\\&\mathbf 3 \esm}} \ar[rr]\ar[ll]&& 
{{\bsm &\mathbf 2\\\mathbf 1&&\mathbf 3\\&\mathbf 2&&\mathbf 4\\&&\mathbf 3 \esm}}\ar@{--}[dr]^?
\ar[dl]
\\
{{\bsm 3 \esm}} \ar[rr]&& 
\blue{\fbox{$\bsm 2\\&3 \esm$}} \ar[ur] && 
{{\bsm \mathbf 2\\&\mathbf 3\\&&\mathbf 4 \esm}} \ar[rr]\ar[ll] && 
{{\bsm \mathbf 1\\&\mathbf 2\\&&\mathbf 3\\&&&\mathbf 4 \esm}}
}
$$
}

{\small
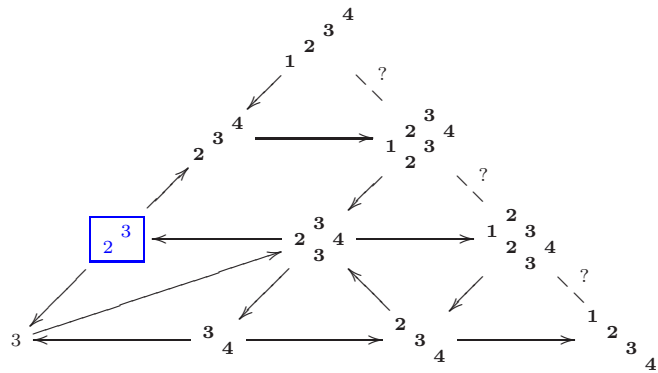
$$
\xymatrix@!@-1.8pc{
&&&{{\bsm &&&\mathbf 4\\&&\mathbf 3\\&\mathbf 2\\\mathbf 1 \esm}} \ar[dl]\ar@{--}[dr]^?
\\
&&{{\bsm &&\mathbf 4\\&\mathbf 3\\\mathbf 2 \esm}} \ar[rr] && 
{{\bsm &&\mathbf 3\\&\mathbf 2&&\mathbf 4\\\mathbf 1&&\mathbf 3\\&\mathbf 2 \esm}}\ar@{--}[dr]^?
\ar[dl]
\\
& \blue{\fbox{$\bsm &3\\2 \esm$}} \ar[ur]\ar[dl] && 
{{\bsm &\mathbf 3\\\mathbf 2&&\mathbf 4\\&\mathbf 3 \esm}} \ar[rr]\ar[ll]\ar[dl]&& 
{{\bsm &\mathbf 2\\\mathbf 1&&\mathbf 3\\&\mathbf 2&&\mathbf 4\\&&\mathbf 3 \esm}}\ar@{--}[dr]^?
\ar[dl]
\\
{{\bsm 3 \esm}}\ar[rrru] && 
{{\bsm \mathbf 3\\&\mathbf 4 \esm}} \ar[rr]\ar[ll] && 
{{\bsm \mathbf 2\\&\mathbf 3\\&&\mathbf 4 \esm}} \ar[rr]\ar[ul] && 
{{\bsm \mathbf 1\\&\mathbf 2\\&&\mathbf 3\\&&&\mathbf 4 \esm}}
}
$$
}

{\small
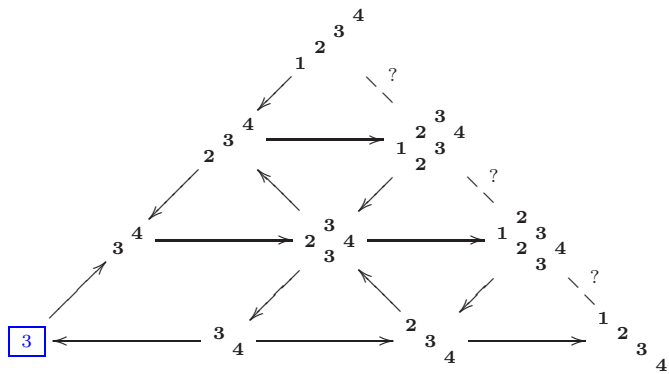
$$
\xymatrix@!@-1.8pc{
&&&{{\bsm &&&\mathbf 4\\&&\mathbf 3\\&\mathbf 2\\\mathbf 1 \esm}} \ar[dl]\ar@{--}[dr]^?
\\
&&{{\bsm &&\mathbf 4\\&\mathbf 3\\\mathbf 2 \esm}} \ar[rr]\ar[dl] && 
{{\bsm &&\mathbf 3\\&\mathbf 2&&\mathbf 4\\\mathbf 1&&\mathbf 3\\&\mathbf 2 \esm}}\ar@{--}[dr]^?
\ar[dl]
\\
& {{\bsm &\mathbf 4\\\mathbf 3 \esm}} \ar[rr] && 
{{\bsm &\mathbf 3\\\mathbf 2&&\mathbf 4\\&\mathbf 3 \esm}} \ar[rr]\ar[dl]\ar[ul]&& 
{{\bsm &\mathbf 2\\\mathbf 1&&\mathbf 3\\&\mathbf 2&&\mathbf 4\\&&\mathbf 3 \esm}}\ar@{--}[dr]^?
\ar[dl]
\\
\blue{\fbox{$\bsm 3 \esm$}} \ar[ur] && 
{{\bsm \mathbf 3\\&\mathbf 4 \esm}} \ar[rr]\ar[ll] && 
{{\bsm \mathbf 2\\&\mathbf 3\\&&\mathbf 4 \esm}} \ar[rr]\ar[ul] && 
{{\bsm \mathbf 1\\&\mathbf 2\\&&\mathbf 3\\&&&\mathbf 4 \esm}}
}
$$
}

{\small
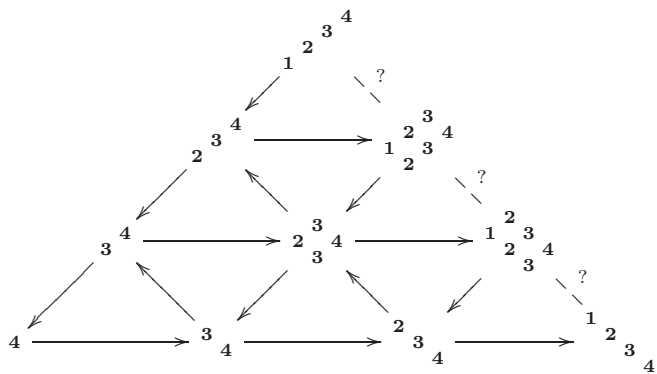
$$
\xymatrix@!@-1.8pc{
&&&{{\bsm &&&\mathbf 4\\&&\mathbf 3\\&\mathbf 2\\\mathbf 1 \esm}} \ar[dl]\ar@{--}[dr]^?
\\
&&{{\bsm &&\mathbf 4\\&\mathbf 3\\\mathbf 2 \esm}} \ar[rr]\ar[dl] && 
{{\bsm &&\mathbf 3\\&\mathbf 2&&\mathbf 4\\\mathbf 1&&\mathbf 3\\&\mathbf 2 \esm}}\ar@{--}[dr]^?
\ar[dl]
\\
& {{\bsm &\mathbf 4\\\mathbf 3 \esm}} \ar[rr]\ar[dl] && 
{{\bsm &\mathbf 3\\\mathbf 2&&\mathbf 4\\&\mathbf 3 \esm}} \ar[rr]\ar[dl]\ar[ul]&& 
{{\bsm &\mathbf 2\\\mathbf 1&&\mathbf 3\\&\mathbf 2&&\mathbf 4\\&&\mathbf 3 \esm}}\ar@{--}[dr]^?
\ar[dl]
\\
{{\bsm \mathbf 4 \esm}} \ar[rr] && 
{{\bsm \mathbf 3\\&\mathbf 4 \esm}} \ar[rr]\ar[ul] && 
{{\bsm \mathbf 2\\&\mathbf 3\\&&\mathbf 4 \esm}} \ar[rr]\ar[ul] && 
{{\bsm \mathbf 1\\&\mathbf 2\\&&\mathbf 3\\&&&\mathbf 4 \esm}}
}
$$
}

The first quiver on the list above is the quiver of
the endomorphism algebra $\End_A(T_M)$, and the
last quiver on the list is the quiver of $\End_\LL(T_M^\vee)$.

Define a matrix 
$$
X :=
\left(
\bbm 
1 & x_1 & x_2 & x_3 & x_4\\
0 &   1 & x_5 & x_6 & x_7\\
0 &   0 &   1 & x_8 & x_9\\
0 &   0 &   0 &   1 & x_{10}\\
0 &   0 &   0 &   0 &   1
\ebm
\right)
$$
The matrix $X$ has columns and rows indexed by $[1,5]$.
For two subsets $I,J \subseteq [1,5]$ with the same cardinality
$k$ let 
$$
\Delta_{I,J} \in \C[x_1,\ldots,x_{10}]
$$
be the minor of $X$ with respect to the rows $I$ and columns $J$.
For example, 
$$
\Delta_{23,35} = {\rm det}\left( \bbm x_5 & x_7\\1&x_9\ebm\right) 
= x_5x_9 - x_7. 
$$
If $I = [1,b] \cup I'$ and $J = [1,b] \cup J'$,
then we have $\Delta_{I,J} = \Delta_{I',J'}$.
Also, for some subsets $I$ and $J$ of $[1,5]$, the
minor $\Delta_{I,J}$ will be zero.
If $I = [1,k]$, we will just write $\Delta_J$ instead of $\Delta_{I,J}$.

We now take the quiver of $\End_\LL(T_M)$ as displayed in the
list above and replace the
vertices by polynomials in $\C[x_1,\ldots,x_{10}]$ as follows:
{\small
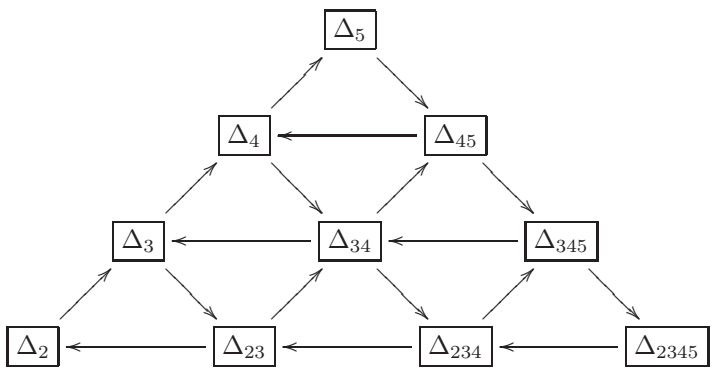
$$
\xymatrix@!@-1.8pc{
&&&{\fbox{$\Delta_{5}$}} \ar[dr]
\\
&&{\fbox{$\Delta_{4}$}} \ar[ur]\ar[dr] && 
{\fbox{$\Delta_{45}$}}\ar[dr]\ar[ll]
\\
& {\fbox{$\Delta_{3}$}} \ar[ur]\ar[dr] && 
{\fbox{$\Delta_{34}$}} \ar[ur]\ar[dr]\ar[ll] && 
{\fbox{$\Delta_{345}$}}\ar[dr]\ar[ll]
\\
{\fbox{$\Delta_{2}$}} \ar[ur] && 
{\fbox{$\Delta_{23}$}} \ar[ur]\ar[ll] && 
{\fbox{$\Delta_{234}$}} \ar[ur]\ar[ll] && 
{\fbox{$\Delta_{2345}$}} \ar[ll]
}
$$
}

There is an isomorphism
$$
\eta\df \C[x_1,\ldots,x_{10}] \to 
\C[\delta_{(M_1,0)},\ldots,\delta_{(M_{10},0)}] = \RR(\CC_M)
$$
which is defined by $x_i \mapsto \delta_{(M_i,0)}$.
Note that
\begin{center}
\begin{tabular}{llll}
$x_1 = \Delta_2$, & $x_2 = \Delta_3$, & $x_3 = \Delta_4$, & $x_4 = \Delta_5$,\\
$x_5 = \Delta_{13}$, & $x_6 = \Delta_{14}$, & $x_7 = \Delta_{15}$,\\
$x_8 = \Delta_{124}$, & $x_9 = \Delta_{125}$,\\ 
$x_{10} = \Delta_{1235}$
\end{tabular}
\end{center}

Now one can easily show that for $1 \le i \le 4$ and
$0 \le a \le b \le i-1$ we have
$$
\eta\df \Delta_{[1,i-a],[1,i-b-1] \cup [4-b+1,4-a+1]} 
\mapsto \delta_{T_{i,[a,b]}}.
$$
(We assume that $[1,0] := \varnothing$.)
Thus all cluster variables appearing in the above sequence of mutations
from $T_M$ to $T_M^\vee$ are images of minors of the matrix $X$.

On the other hand, there are cluster variables $\delta_R$
in $\RR(\CC_M)$ which are not images
of minors of the matrix $X$, compare
\cite[Section 13.1]{GLSSemi1}.

\newpage


{\Large\part{Cluster algebras}
\label{part4}}


In this part, let $K = \C$ be the field of complex numbers.


{\Large\section{Kac-Moody Lie algebras and semicanonical bases}
\label{semican}}


In this section we recall known results on Kac-Moody
Lie algebras and semicanonical bases.

\subsection{Kac-Moody Lie algebras}
Let $\GG=(\GG_0,\GG_1,\gamma)$ be a finite graph (without loops).  
It has as set of vertices $\GG_0$, edges $\GG_1$ and 
$\gamma\df\GG_1\to \cP_2(\GG_0)$ determining
the adjacency of the edges; 
here $\cP_2(\GG_0)$ denotes the set
of two-element subsets of $\GG_0$. 
If $\GG_0 = \{ 1,2,\ldots,n \}$ we can assign
to $\GG$ a {\em symmetric generalized Cartan matrix} 
$C_\GG = (c_{ij})_{i,j = 1,2,\ldots n}$, which is an $n \times n$-matrix
with integer entries
$$
c_{ij} :=
\begin{cases}
2 & \text{if $i = j$},\\
-\abs{\gamma^{-1}(\{i,j\})} & \text{if $i \neq j$}.
\end{cases}
$$ 
Obviously, the assignment
$\GG \mapsto C_\GG$ induces a bijection between isomorphism classes of
graphs with vertex set $\{ 1,2,\ldots,n \}$ and symmetric generalized 
Cartan matrices
in $\Z^{n\times n}$ up to simultaneous permutation of rows and columns.

If $Q=(Q_0,Q_1,s,t)$ is a quiver without oriented cycles (in particular
without loops) its underlying graph $\abs{Q}:=(Q_0,Q_1,q)$ is given
by $q(a)=\{s(a),t(a)\}$ for all $a \in Q_1$ i.e. it is obtained by 
``forgetting'' the orientation of the edges.
We write then also $C_Q := C_{\abs{Q}}$.

It will be convenient for us to consider $\g := \g_Q := \g(C_Q)$  the 
(symmetric) {\it Kac-Moody Lie algebra} 
\index{Kac-Moody Lie algebra}
\index{$\g$}
associated to $Q$, which is defined as follows:
Let $\h$ be a $\C$-vector space of dimension $2n-\operatorname{rank}(C_Q)$,
and let $\Pi := \{ \alpha_1,\ldots,\alpha_n \} \subset \h^*$ and 
$\Pi^\vee := \{ \alpha_1^\vee,\ldots,\alpha_n^\vee \} \subset \h$
be linearly independent subsets of the vector spaces $\h^*$ and $\h$,
respectively, such that
$$
\alpha_i(\alpha_j^\vee) = c_{ij}
$$
for all $i,j$.
Then $\g = (\g,[-,-])$ is the Lie algebra over $\C$
generated by $\h$ and the symbols $e_i$
and $f_i\ (1 \le i \le n)$ satisfying the following 
defining relations:
\begin{itemize}

\item[(L1)]
$[h,h'] = 0$ for all $h,h' \in \h$,

\item[(L2)]
$[h,e_i] = \alpha_i(h)e_i$, and
$[h,f_i] = -\alpha_i(h)f_i$,

\item[(L3)]
$[e_i,f_i] = \alpha_i^\vee$ and
$[e_i,f_j] = 0$ for all $i \not= j$,

\item[(L4)]
$(\ad(e_i)^{1-c_{ij}})(e_j) = 0$ for all $i \not= j$,

\item[(L5)]
$(\ad (f_i)^{1-c_{ij}})(f_j) = 0$ for all $i \not= j$.
\end{itemize}
(For $x,y \in \g$ and $m \ge 1$
we set $\ad(x)(y) := \ad(x)^1(y) := [x,y]$
and
$\ad(x)^{m+1}(y) := \ad(x)^m([x,y])$.)

The Lie algebra $\g$ is finite-dimensional if and only if
$Q$ is a Dynkin quiver.
In this case, $\g$ is the usual simple Lie algebra associated to $Q$.

Conversely, if $\g = \g(C)$ is a Kac-Moody Lie algebra defined by a symmetric
generalized Cartan matrix $C$, we say that $\g$ is of type $\GG$ if $C=C_\GG$.
This is well defined for symmetric Kac-Moody Lie algebras.
We call $\GG$ the {\it Dynkin graph} of $\g$.

For $\alpha \in \h^*$ let
$$
\g_\alpha := 
\{ x \in \g \mid [h,x] = \alpha(h)x \text{ for all $h \in \h$} \}.
$$
\index{$\g_\alpha$}\noindent
One can show that $\dm \g_\alpha < \infty$ for all $\alpha$.
By 
$
R := \Z\alpha_1 \oplus \cdots \oplus \Z\alpha_n
$
we denote the {\it root lattice} 
\index{root lattice}
of $\g$.
Define 
$
R^+ :=  \N\alpha_1 \oplus \cdots \oplus \N\alpha_n.
$
The {\it roots} of $\g$ 
\index{root of $\g$}
are defined as the elements in
$$
\Delta := \{  \alpha \in R \setminus \{0\} \mid \g_\alpha \not= 0 \}.
$$
\index{$\Delta$ (set of roots)}\noindent
Set $\Delta^+ := \Delta \cap R^+$ and $\Delta^- := \Delta \cap (-R^+)$.
\index{$\Delta^+$}
\index{$\Delta^-$}
One can show that $\Delta = \Delta^+ \cup \Delta^-$.
The elements in $\Delta^+$ and $\Delta^-$ are the {\it positive roots}
\index{positive root}
and the {\it negative roots}, 
\index{negative root}
respectively.
The elements in $\{ \alpha_1,\ldots,\alpha_n \}$ are positive roots
of $\g$ and are called {\it simple roots}.
\index{simple root}

One has the triangular decomposition
$
\g = \n_- \oplus \h \oplus \n
$
with 
$$
\n_- = \bigoplus_{\alpha \in \Delta^+} \g_{-\alpha}
\text{\;\;\; and \;\;\;}
\n = \bigoplus_{\alpha \in \Delta^+} \g_\alpha.
$$
\index{$\n$}\noindent
\index{$\n_-$}\noindent
The Lie algebra $\n$ is generated by $e_1,\ldots,e_n$
with defining relations (L4).
Set $\n_\alpha := \g_\alpha$ if $\alpha \in R^+ \setminus \{0\}$.

For $1 \le i \le n$ define an element $s_i$
in the automorphism group ${\rm Aut}(\h^*)$ of $\h^*$ by
$$
s_i(\alpha) := \alpha - \alpha(\alpha_i^\vee)\alpha_i
$$
for all $\alpha \in \h^*$.
The subgroup $W \subset {\rm Aut}(\h^*)$ generated by $s_1,\ldots,s_n$
is the {\it Weyl group} 
\index{Weyl group}
of $\g$.
The elements $s_i$ are called
{\it Coxeter generators} 
\index{Coxeter generator}
of $W$.
The identity element of $W$ is denoted by 1.
The {\it length} $l(w)$ 
\index{length of a Weyl group element}
of some $w \not= 1$ in $W$ is 
the smallest number $t \ge 1$ such that $w = s_{i_t}\cdots s_{i_2} s_{i_1}$
for some $1 \le i_j \le n$.
In this case $(i_t,\ldots,i_2, i_1)$ is a
{\it reduced expression} for $w$.
\index{reduced expression for a Weyl group element}
Let $R(w)$ be the set of all reduced expressions for $w$.
We set $l(1) = 0$.

Let us repeat the following definition from Section~\ref{wadapt}.
A Weyl group element $w$ is $Q$-{\it adaptable} 
\index{$Q$-adaptable Weyl group element}
if there
exists a reduced expression $(i_t,\ldots, i_2, i_1) \in R(w)$ such
that $i_1$ is a sink of $Q$, and $i_{k+1}$ is a sink of
$\si_{i_k} \cdots \si_{i_2}\si_{i_1}(Q)$ for all
$1 \le k \le t-1$.
If this is the case, we say that $(i_t,\ldots,i_2,i_1)$ is a
$Q$-{\it adapted} reduced expression.
A Weyl group element $w$ is {\it adaptable} 
\index{adaptable Weyl group element}
if there exists an orientation $Q$ of the Dynkin graph
of $\g$ such that $w$ is $Q$-adaptable.

A root $\alpha \in \Delta$ is a {\it real root} 
\index{real root}
if
$\alpha = w(\alpha_i)$
for some $w \in W$ and some $i$.
It is well known that $\dm \g_\alpha = 1$ if $\alpha$ is a real root. 
By $\Delta_{\rm re}$ 
\index{$\Delta_{\rm re}$}
we denote the set of real roots of $\g$.
Define $\Delta_{\rm re}^+ := \Delta_{\rm re} \cap \Delta^+$.
\index{$\Delta_{\rm re}^+$}

Finally, let us fix a $\Z$-basis 
$\{\vpi_j\mid 1\le j \le 2n - \operatorname{rank}(C_Q)\}$ of
$\mathfrak{h}^*_\Z$ such that
\[
\vpi_j(\alpha^\vee_i)=\delta_{ij},\qquad
(1\le i\le n,\ 1\le j\le 2n-\operatorname{rank}(C_Q)),  
\]
(see \cite[\S 6.1.6]{Ku}).
The $\vpi_j$ are the {\it fundamental weights}. 
\index{fundamental weight}
We denote by 
\[
P := \bigoplus_{i=1}^{2n-\operatorname{rank}(C_Q)} \Z\vpi_i
\] 
the integral weight lattice, and we set
\[
P^+ := \{ \nu\in P\mid \nu(\alpha_i^\vee)\ge 0\text{ for } 1\leq i\leq n \}.
\]

\subsection{The universal enveloping algebra $U(\n)$}\label{universal}
The universal enveloping algebra $U(\n)$ 
\index{$U(\n)$}
of the Lie algebra
$\n$ is the associative $\C$-algebra defined by generators
$E_1,\ldots,E_n$ and relations
$$
\sum_{k=0}^{1-c_{ij}} (-1)^k E_i^{(k)}E_jE_i^{(1-c_{ij}-k)} = 0
$$
for all $i \not= j$,
where the $c_{ij}$ are the entries of the generalized Cartan matrix
$C_Q$, and 
$$
E_i^{(k)} := E_i^k/k!.
$$
We have a canonical embedding 
$\iota\df \n \to U(\n)$ which maps $e_i$ to $E_i$ for all $1 \le i \le n$.
We consider $\n$ as a subspace of $U(\n)$, and we also identify
$e_i$ and $E_i$.

Let 
\[
J= \left\{
\begin{array}{ll}
\N_1&\mbox{if}\ \dm\n = \infty,\\
{[1,d]}&\mbox{if}\ \dm\n = d.
\end{array}
\right.
\]
Let ${\rm P} := \{ p_i \mid i \in J \}$ be a $\C$-basis of $\n$
such that ${\rm P} \cap \n_\alpha$ is a basis of $\n_\alpha$
for all positive roots $\alpha$.
We assume that $\{ e_1,\ldots,e_n \} \subset {\rm P}$.
Thus $e_i$ is a basis vector of the
(1-dimensional) space $\n_{\alpha_i}$.
 
For $k \ge 0$ define $p_i^{(k)} := p_i^k/k!$.
Let $\N^{(J)}_{\phantom{0}}$ be the set of
tuples $(m_i)_{i \in J}$ of natural numbers $m_i$ such that
$m_i = 0$ for all but finitely many $m_i$.
For ${\bf m} = (m_i)_{i \ge 1} \in \N^{(J)}_{\phantom{0}}$ define
$$
p_{\bf m} := p_1^{(m_1)} p_2^{(m_2)}\cdots p_s^{(m_s)}
$$
where $s$ is chosen such that $m_j = 0$ for all $j > s$.

\begin{Thm}[Poincar{\'e}-Birkhoff-Witt]\label{PBW1}
The set
$$
\cP := \left\{ p_{\bf m} \mid {\bf m} \in \N^{(J)}_{\phantom{0}} \right\}
$$
is a $\C$-basis of $U(\n)$.
\end{Thm}

The basis $\cP$ is called a {\it PBW-basis} of $U(\n)$.
\index{PBW-basis of $U(\n)$}
For $d = (d_1,\ldots,d_n) \in \N^n$ let
$U_d$ be the subspace of $U(\n)$ spanned by the 
elements of the form $e_{i_1}e_{i_2} \cdots e_{i_m}$,
where for each $1 \le i \le n$ the set
$\{ i_k \mid i_k = i, 1\le k \le m \}$ contains exactly
$d_i$ elements.

It follows that
$$
U(\n) = \bigoplus_{d \in \N^n} U_d.
$$
This turns $U(\n)$ into an $\N^n$-graded algebra.

Furthermore, $U(\n)$ is a cocommutative Hopf algebra
with comultiplication
$$
\Delta\df U(\n) \to U(\n) \otimes U(\n)
$$
defined by $\Delta(x) := 1 \otimes x + x \otimes 1$
for all $x \in \n$.
It is easy to check that 
\begin{equation}\label{eqPBWdelta}
\Delta(p_{\bf m}) = \sum_{\bf k} p_{\bf k} \otimes p_{{\bf m}-{\bf k}},
\end{equation}
where the sum is over all tuples ${\bf k} = (k_i)_{i\ge 1}$ with
$0\le k_i \le m_i$ for every $i$.

By $U_d^*$ we denote the vector space dual of $U_d$.
Define the {\it graded dual} of $U(\n)$ 
\index{graded dual of $U(\n)$}
by
$$
U(\n)^*_{\rm gr} := \bigoplus_{d \in \N^n} U_d^*.
$$
\index{$U(\n)^*_{\rm gr}$}\noindent
It follows that
$U(\n)^*_{\rm gr}$ is a commutative associative $\C$-algebra
with multiplication defined via the comultiplication
$\Delta$ of $U(\n)$:
For $f',f'' \in U(\n)^*_{\rm gr}$ and
$x \in U(\n)$, we have
$$
(f' \cdot f'')(x) = \sum_{(x)} f'(x_{(1)})f''(x_{(2)}),
$$
where (using the Sweedler notation) we write
$$
\Delta(x) = \sum_{(x)} x_{(1)} \otimes x_{(2)}.
$$

Let
$
\cP^* := 
\left\{ p_{\bf m}^* \mid {\bf m} \in \N^{(J)}_{\phantom{0}} \right\}
$
be the dual PBW-basis of $U(\n)^*_{\rm gr}$, where 
$$
p_{\bf m}^*(p_{\bf n}) :=
\begin{cases}
1 & \text{if ${\bf m} = {\bf n}$},\\
0 & \text{otherwise}.
\end{cases}
$$
The element in $\cP^*$ corresponding to $p_i \in {\rm P}$ is
denoted by $p_i^*$.
It follows from (\ref{eqPBWdelta}) that 
$$
p_{\bf m}^* \cdot p_{\bf n}^* = p_{{\bf m}+{\bf n}}^*,
$$
that is, each element $p_{\bf m}^*$ in $\cP^*$
is equal to a monomial in the $p_i^*$'s.
Hence, the graded dual $U(\n)^*_{\rm gr}$ can be identified with the
polynomial algebra $\C[p_1^*,p_2^*,\ldots]$ (with countably 
many variables $p_i^*$).

\subsection{The Lie algebra $\n(w)$}\label{nw}
Let 
$$
\widehat{\n} := \prod_{\alpha \in \Delta^+} \n_\alpha
$$
\index{$\widehat{\n}$}\noindent
be the completion of the Lie algebra $\n$.

A subset $\Theta \subseteq \Delta^+$ is {\it bracket closed}
\index{bracket closed subset of $\Delta^+$}
if for all $\alpha,\beta \in \Theta$ with $\alpha+\beta \in \Delta^+$
we have $\alpha+\beta \in \Theta$.
One calls $\Theta$
{\it bracket coclosed}
\index{bracket coclosed subset of $\Delta^+$}
if $\Delta^+ \setminus \Theta$ is bracket closed.

For $w \in W$ set
$$
\Delta_w^+ := \left\{ \alpha \in \Delta^+ \mid w(\alpha) < 0 \right\}.
$$
\index{$\Delta_w^+$}\noindent
It follows that for each $(i_t,\ldots,i_2, i_1) \in R(w)$ we have
$$
\Delta_w^+  = \{\alpha_{i_1}, s_{i_1}(\alpha_{i_2}), \ldots,
s_{i_1}s_{i_2}\cdots s_{i_{t-1}}(\alpha_{i_t}) \}.
$$
The set $\Delta_w^+$ contains $l(w)$ positive roots, all
of these are real roots. See for example~\cite[1.3.14]{Ku}.
The next lemma is well known.

\begin{Lem}\label{bracketclosed1}
For every $w \in W$, the set
$\Delta_w^+$ is bracket closed
and bracket coclosed.
\end{Lem}

Let 
$$
\n(w) := \bigoplus_{\alpha \in \Delta_w^+} \n_\alpha 
$$
\index{$\n(w)$}\noindent
be the {\it nilpotent Lie algebra associated to} $w$.
Thus $\dm \n(w) = l(w)$.

The following lemma is also well known (see \cite{Be}).

\begin{Lem}\label{bracketclosed2}
Let $M = M_1 \oplus \cdots \oplus M_r$ be a terminal $\C Q$-module,
and let $x(1) < x(2) < \cdots < x(r)$ be a $\GG_M$-adapted ordering 
of the vertices of $\GG_M$.
For $1 \le j \le r$ set $i_j = i$ where $i$ is
the vertex of $Q$ with $x(j) = (i,a)$ for some $a$.
Then 
$$
w:= w(M) := s_{i_r}\cdots s_{i_2} s_{i_1}
$$ 
is $Q^\op$-adaptable, and
we have
$$
\Delta_w^+ = \Delta_M^+ := \{ \dimv(M_1),\ldots,\dimv(M_r) \}.
$$
Moreover, $w$ does not depend on the choice of the adapted ordering.
\end{Lem}

Let us discuss an example.
Let $Q$ be the quiver
{\small
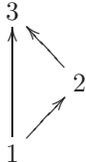
$$
\xymatrix@-0.8pc{
3 \\
& 2 \ar[ul]\\
1 \ar[ur]\ar[uu]
}
$$
}\noindent
We take a terminal $\C Q$-module $M$ defined by
$t_1 = 2$ and $t_2 = t_3 = 1$.
Thus the quiver $\GG_M$ looks as follows:
$$
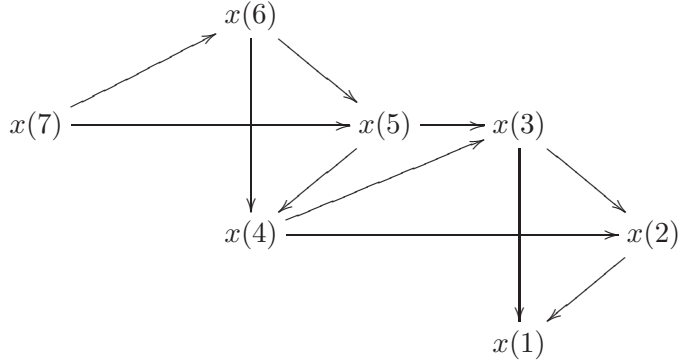
\xymatrix{
&& x(6) \ar[dr]\ar[dd] \\
x(7) \ar[rrr]\ar[urr] &&& x(5) \ar[r]\ar[dl] & x(3) \ar[dr]\ar[dd]
\\
&& x(4) \ar[rrr]\ar[rru] &&& x(2) \ar[dl]\\
&&&&x(1) 
}
$$
The vertices are labelled such that the ordering 
$x(1) < x(2) < \cdots < x(7)$ is $\GG_M$-adapted.
Thus
$w = w(M) := s_1s_3s_2s_1s_3s_2s_1$ 
is a $Q^\op$-adaptable Weyl group element.
We get
$$
\Delta_w^+ := \left\{ 
\bsm0\\&0\\1\esm,\,
\bsm0\\&1\\1\esm,\,
\bsm1\\&1\\2\esm,\,
\bsm1\\&2\\2\esm,\,
\bsm2\\&2\\3\esm,\,
\bsm2\\&3\\3\esm,\,
\bsm3\\&3\\4\esm
\right\}.
$$
These are just the dimension vector of the indecomposable
direct summands of $M$.
For example, we have
$$
s_1s_2s_3(\alpha_1)  = s_1s_2\left(\bsm1\\&0\\1\esm\right)
= s_1\left(\bsm1\\&2\\1\esm\right)
= \left(\bsm1\\&2\\2\esm\right).
$$
The part of the preinjective component of $\C Q$, where
the indecomposable direct summands of $M$ lie, looks like this:
$$
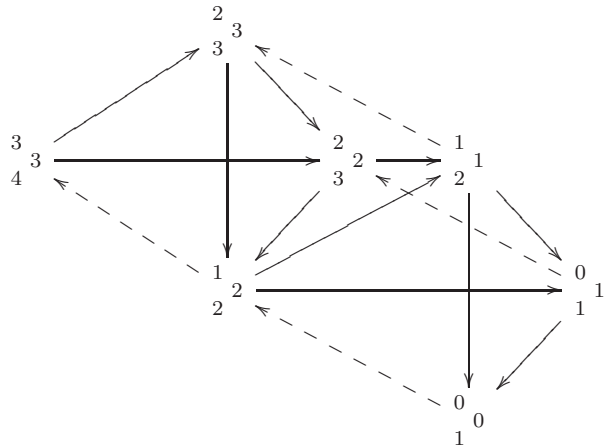
\xymatrix{
&& {\bsm2\\&3\\3\esm} \ar[dr]\ar[dd] \\
{\bsm3\\&3\\4\esm} \ar[rrr]\ar[urr] &&& {\bsm2\\&2\\3\esm} \ar[r]\ar[dl] & 
{\bsm1\\&1\\2\esm} \ar[dr]\ar[dd]
\ar@{-->}[ull]\\
&& {\bsm1\\&2\\2\esm} \ar@{-->}[ull]\ar[rrr]\ar[rru] &&& 
{\bsm0\\&1\\1\esm} \ar[dl]\ar@{-->}[ull]\\
&&&& {\bsm0\\&0\\1\esm} \ar@{-->}[ull]
}
$$

\subsection{Semicanonical bases}\label{section_semican}
Recall that for each dimension vector $d = (d_1,\ldots,d_n)$
we defined the affine variety $\LL_d$ of nilpotent $\LL$-modules
with dimension vector $d$.
A subset $C$ of $\LL_d$ is said to be constructible if it is
a finite union of locally closed subsets.
A function 
$$
f\df \LL_d \to \C
$$ 
is {\it constructible} 
\index{constructible function $\LL_d \to \C$}
if 
the image $f(\LL_d)$ is finite and 
$f^{-1}(m)$ is a constructible
subset of $\LL_d$ for all $m \in \C$.
The set of constructible functions on $\LL_d$ is denoted by 
$M(\LL_d)$.
\index{$M(\LL_d)$}
This is a $\C$-vector space. 

Recall that the group $\GL_d$ acts on $\LL_d$ by conjugation.
By $M(\LL_d)^{\GL_d}$ 
\index{$M(\LL_d)^{\GL_d}$}
we denote the subspace of
$M(\LL_d)$ consisting of the constructible functions which are
constant on the $\GL_d$-orbits in $\LL_d$.
Set
$$
\widetilde{\MM} := \bigoplus_{d \in \N^n} M(\LL_d)^{\GL_d}.
$$
\index{$\widetilde{\MM}$}\noindent
For $f' \in M(\LL_{d'})^{\GL_{d'}}$, $f'' \in M(\LL_{d''})^{\GL_{d''}}$
and $d = d'+d''$
we define a constructible function 
$$
f := f' \star f''\df \LL_d \to \C
$$ 
\index{$f' \star f''\df \LL_d \to \C$}\noindent
in $M(\LL_d)^{\GL_d}$ by
$$
f(X) := \sum_{m \in \C} m \, \chi_c\left(\left\{ 
U \subseteq X \mid f'(U)f''(X/U) = m
\right\} \right)
$$
for all $X \in \LL_d$, where $U$ runs over 
the points of the Grassmannian of all submodules
of $X$ with $\dimv(U) = d'$.
Here, for a constructible subset $V$ of a complex variety 
we denote by $\chi_c(V)$ its
(topological) Euler characteristic with respect to cohomology with compact 
support.
This turns $\widetilde{\MM}$ into an associative $\C$-algebra.

\begin{Rem}\label{newconv}
{\rm
Note that the product $\star$ defined here is opposite to 
the convolution product we have used in 
\cite{GLSSemi1,GLSVerma,GLSSemi2}.
This new convention turns out to be better adapted to our choice
of categorifying $\C[N(w)]$ and $\C[N^w]$ by categories
closed under factor modules.
It is also compatible with our choice in \cite{GLSFlag}
of categorifying coordinate rings of partial flag varieties 
by categories closed under submodules (see Remark~\ref{finrem}).
}
\end{Rem}

For the canonical basis vector $e_i := \dimv(S_i)$ we know
that $\LL_{e_i}$ is just a point, which (as a $\LL$-module)
is isomorphic to the simple module $S_i$.
Define ${\bf 1}_i\df \LL_{e_i} \to \C$ by ${\bf 1}_i(S_i) := 1$.
By $\MM$ we denote the subalgebra of $\widetilde{\MM}$ generated
by the functions ${\bf 1}_i$ where $1 \le i \le n$.
Set $\MM_d := \MM \cap M(\LL_d)^{\GL_d}$.
It follows that
$$
\MM := \bigoplus_{d \in \N^n} \MM_d
$$
\index{$\MM$}\noindent
is an $\N^n$-graded $\C$-algebra.

\begin{Thm}[{Lusztig \cite{Lu2}}]\label{lusztig1}
There is an isomorphism of $\N^n$-graded $\C$-algebras
$$
U(\n) \to \MM
$$
defined by $E_i \mapsto {\bf 1}_i$ for $1 \le i \le n$.
\end{Thm}

Let $\irr(\LL_d)$ be the set of irreducible components
of $\LL_d$.

\begin{Thm}[{Lusztig \cite{Lu2}}]\label{lusztig2}
For each $Z \in \irr(\LL_d)$ there
is a unique $f_Z\df \LL_d \to \C$ in $\MM_d$
such that $f_Z$ takes value 1 on some dense open subset of
$Z$ and value 0 on some dense open subset of any other irreducible
component $Z'$ of $\LL_d$.
Furthermore, the set 
$$
\cS := \left\{ f_Z \mid Z \in \irr(\LL_d), d \in \N^n \right\}
$$
\index{$\cS$}\noindent
is a $\C$-basis of $\MM$.
\end{Thm}

The basis $\cS$ is called the {\it semicanonical basis} 
\index{semicanonical basis}
of $\MM$.
By Theorem~\ref{lusztig1} we just identify $\MM$ and $U(\n)$ and
consider $\cS$ also as a basis of $U(\n)$.

Now we turn to the graded dual 
$$
U(\n)^*_{\rm gr} = \bigoplus_{d \in \N^n} U_d^*
$$
of $U(\n)$.
Let $\MM_d^*$ be the dual space of $\MM_d$, and set
$$
\MM^* := \bigoplus_{d \in \N^n} \MM_d^*.
$$
\index{$\MM^*$}\noindent

For $X \in \LL_d$ define a linear form
$$
\delta_X\df \MM_d \to \C
$$
\index{$\delta_X\df \MM_d \to \C$}\noindent
by $\delta_X(f) := f(X)$.
It is not difficult to show that
the map $X \mapsto \delta_X$ from $\LL_d$ to $\MM_d^*$ is constructible, i.e.
this map has a finite image and the preimage of each
element in $\MM_d^*$ is constructible in $\LL_d$.
So on every irreducible component 
$Z \in \irr(\LL_d)$ there is a Zariski open set on which
this map is constant.
Define $\rho_Z := \delta_X$ for any $X$ in this open set.
\index{$\rho_Z$}
The $\C$-vector space $\MM_d^*$ 
is spanned by the functions
$\delta_X$ with $X \in \LL_d$.
Then by construction
$$
\cS^* := \left\{ \rho_Z \mid Z \in \irr(\LL_d), d \in \N^n \right\}
$$
\index{$\cS^*$}\noindent
is the basis of $\MM^* \equiv U(\n)^*_{\rm gr}$ 
dual to Lusztig's semicanonical basis 
\index{dual semicanonical basis}
$\cS$
of $U(\n)$.

\subsection{A cluster character}\label{subsectclusterchar}
The map $X \mapsto \delta_X$ from $\nil(\LL)$ to $U(\n)^*_{\text{gr}}$ 
has the following multiplicative properties.

\begin{Thm}[{\cite[Lemma 7.2]{GLSSemi1} 
and~\cite[Theorem 2]{GLSSemi2}}]\label{Th:ClusChar}
Let $X,Y \in \nil(\LL)$.
\begin{itemize}

\item[(i)] 
We have
\[
\delta_X\delta_Y = \delta_{X\oplus Y}.
\]

\item[(ii)]
If $\dm \Ext_\LL^1(X,Y) = 1$ with
\[
0 \to X\to E'\to Y\to 0 \text{ and } 0\to Y\to E''\to X\to 0
\]
the corresponding non-split short exact sequences, then
\[
\delta_X\delta_Y = \delta_{E'} + \delta_{E''}.
\]
\end{itemize}
\end{Thm}

In fact, the main result of~\cite{GLSSemi2} is a more general multiplication
formula (without the restriction to the case
$\dm \Ext^1_\LL(X,Y) = 1$).
However, in this paper we do not need this.

\subsection{Dual Verma modules}\label{dualVerma}
We present some of the results of~\cite{GLSVerma} in a form
convenient for our present purpose. 
For $i=1,\ldots, n$, let
$\hI_i$ denote the injective $\LL$-module with simple socle
$S_i$. 
If $\LL$ is not of Dynkin type, $\hI_i$ is infinite-dimensional.
For $\nu \in P^+$ we write
\[
\hI_\nu := \bigoplus_{i=1}^n \hI_i^{\,\nu(\alpha_i^\vee)}.
\]
For $i=1,\ldots,n$ and a nilpotent $\LL$-module $X$ we denote by
$\G(i,X)$ the variety of submodules $Y$ of $X$ such that 
$X/Y \cong S_i$.
Similarly, if  
$$
\soc(X) = \bigoplus_{i=1}^n S_i^{m_i}
$$ 
and 
$\nu \in P^+$ is such that $\nu(\alpha_i^{\vee})\geq m_i$ for $1 \le i \le n$, 
then we have an embedding $X \hookrightarrow \hI_\nu$.
In this case, we denote by $\G(i,\nu,X)$ the variety of submodules $Y$ of 
$\hI_\nu$ such that $X\subset Y$ and $Y/X \cong S_i$.  
Hence, if
$\dimv(X) = \beta$ and $f \in \MM_{\beta-\alpha_i}$, 
we can form the following sum
\[
S = \sum_{m\in\Z} m\chi_c(\{ Y \in \G(i,X)\mid f(Y)=m \}).
\]
For convenience we shall denote such an expression by an
integral, for example, 
\[
S = \int_{Y\in\G(i,X)} f(Y).
\]
Similarly, there exists a partition 
$$
\G(i,X) = \bigsqcup_{j=1}^m A_j
$$ 
into constructible subsets such that
$\delta_Y=\delta_{Y'}$ for all $Y,Y'\in A_j$.  
Then, choosing  arbitrary $Y_j\in A_j$ for $j=1, \ldots, m$,
we can also denote by an integral the following element of
$\MM_{\beta-\alpha_i}^*$
\[
\int_{Y\in\G(i,X)}\delta_Y =\sum_{j=1}^m \chi_c(A_j)\delta_{Y_j}.
\]

\begin{Thm} 
Let $\lambda\in P$ be an integral weight,
and let $M_{\rm low}(\lambda)$ be the lowest  weight Verma 
{\em right} $U(\g)$-module,
with lowest weight $\lambda$.
Under the  identifications 
$$
M_{\rm low}(\lambda) \equiv U(\n) \equiv \MM
$$ 
we equip $\MM$ with the structure of a right $U(\g)$-module as follows.
The generators $e_i\in\n,\ f_i\in\n_-,\ h\in \h$ act
on
$g\in\MM(\beta)$ by
\begin{align*}
(g \cdot e_i)(X')  &= \int_{Y\in\G(i,X')} g(Y),\\
(g \cdot f_i)(X) &= \int_{Y\in\G(i,\nu,X)}g(Y) - 
(\nu-\lambda)(\alpha_i^\vee)g(X\oplus S_i),\\
g\cdot h &= (\lambda -\beta)(h) g,
\end{align*}
where $X'\in\LL_{\beta+\alpha_i}$, $X\in\LL_{\beta-\alpha_i}$ and 
$\nu\in P^+$ are as above.
\end{Thm}

Note that $g\cdot e_i=g*\mathbf{1}_i$ by our convention for the multiplication
in $\MM$. Moreover, 
the formula for $g\cdot f_i\in\MM_{\beta-\alpha_i}$ is 
in fact independent of the choice of $\nu$.

Recall, that for each $\h$-diagonalizable right $U(\g)$-module 
$$
M = \bigoplus_{\mu\in \h^*} M_\mu
$$ 
one can consider the {\em dual} representation
$$
M^* = \bigoplus_{\mu\in\h^*} M^*_\mu
$$ 
defined by
$M^*_\mu := \Hom_\C (M_{\mu},\C)$ and
\[
(x\cdot\phi)(m):=\phi( m\cdot x),
\qquad
(x\in\g,\ m\in M).
\]
Consider the canonical epimorphism from the Verma module
$M_{\rm low}(\lambda)$ to the irreducible lowest weight right 
$U(\n)$-module $L_{\rm low}(\lambda)$. For the
corresponding dual representations we obtain an inclusion
$$
L^*_{\rm low}(\lambda) \hookrightarrow M_{\rm low}^*(\lambda).
$$ 
It is well known that
$L_{\rm low}^*(\lambda)$ is isomorphic to the irreducible highest weight 
left module $L(\lambda)$. This yields the following realization of $L(\lambda)$
in terms of $\delta$-functions.

\begin{Thm} \label{thm:irred-l}
Let $\lambda\in P^+$ be a dominant weight.
The subspace of $U(\n)^*_{\rm gr}$ 
spanned
by all $\delta_X$ such that $X$ is isomorphic to a submodule of $\hI_\lambda$
carries the above-mentioned structure of an irreducible highest weight module
$L(\lambda)$. For such $X$ with $\dimv(X)=\beta$ 
the action of the Chevalley generators of $U(\g)$ is 
given by
\begin{align*}
e_i\cdot \delta_X &=\int_{Y\in\G(i,X)} \delta_Y,\\
f_i\cdot \delta_X &=\int_{Y'\in\G(i,\lambda,X)}\delta_{Y'},\\
h\cdot \delta_X &=(\lambda-\beta)(h)\delta_X.
\end{align*}
\end{Thm}

Note that $U(\n)_{\rm gr}^*$ carries also a {\em right} $U(\n)$-module
structure coming from the left regular representation of $U(\n)$. In
order to describe it, we introduce the following definition.
For $X \in \LL_\beta$ we denote by $\G'(i,X)$ the variety of submodules
$Y$ of $X$ such that $\dimv(Y) = \alpha_i$. Each element of this space
is isomorphic to $S_i$ and clearly $\G'(i,X)$ is a projective space.
It is easy to see that
\[
\delta_X \cdot e_i = \int_{S\in\G'(i,X)} \delta_{X/S}.
\]
Under the above identification 
$M^*_{\rm low}(\lambda) \equiv U(\n)_{\rm gr}^*$, the
subspace of $U(\n)^*_{\rm gr}$ carrying the
$U(\g)$-module $L(\lambda)$ can be described as follows.

\begin{Cor}\label{embedding1}
We have
\[
L(\lambda) = \left\{\phi\in U(\n)^*_{\rm gr}\mid
 \phi\cdot e_i^{\lambda(\alpha_i^\vee)+1}=0,\ (i=1,\ldots, n) \right\}.
\]
\end{Cor}

\begin{proof}
The nilpotent $\LL$-module $X$ is isomorphic to a submodule of $\hI_\lambda$
if and only if ${\delta_X \cdot e_i^{\lambda(\alpha_i^\vee)+1}=0}$ 
for every $i$.
The claim then follows from Theorem~\ref{thm:irred-l}.   
\end{proof}


{\Large\section{A dual PBW-basis and a dual semicanonical basis 
for $\cA(\CC_M)$}\label{PBWsection}}


In this section we prove Theorem~\ref{main7} and Theorem~\ref{basesthm}.
We also deduce from these results the existence of semicanonical bases
for the cluster algebras $\widetilde{\RR}(\CC_M)$ and
$\underline{\RR}(\CC_M)$ obtained by inverting and specializing
coefficients, respectively.

\subsection{Proof of Theorem~\ref{main7}}\label{clustercattoclusteralg}

By the definition of the cluster algebra $\cA(\CC_M,T)$, 
its initial seed
is $(\yy, B(T)^\circ)$. 
Let $\F=\C(y_1,\ldots,y_r)$.
Since $T$ is rigid, by Theorem~\ref{Th:ClusChar}~(i)
every monomial in the $\delta_{T_i}$ belongs to
the dual semicanonical basis $\cS^*$, hence the $\delta_{T_i}$ 
are algebraically
independent and $(\delta_{T_1}, \ldots , \delta_{T_r})$ is a transcendence
basis of the subfield $\G$ it generates inside the fraction field
of $U(\n)^*_{\text{gr}}$.
Let $\iota\df \F\to\G$ be the field isomorphism defined by
$\iota(y_i)=\delta_{T_i}\ (1\le i \le r)$.
Combining Theorem~\ref{main4} and Theorem~\ref{Th:ClusChar}~(ii) we see that
the cluster variable $z$ of $\cA(\CC_M,T)$ obtained from the
initial seed $(\yy, B(T)^\circ)$ through a sequence of seed mutations
in successive directions $k_1,\ldots,k_s$ will be mapped by $\iota$
to $\delta_X$, where $X \in \CC_M$ is the indecomposable rigid module obtained
by the same sequence of mutations of rigid modules.
It follows that $\iota$ restricts to an isomorphism from
$\cA(\CC_M,T)$ to $\RR(\CC_M,T)$.
This isomorphism is completely determined by the images of the 
elements $y_i$, hence the unicity. 
The cluster monomials are mapped to elements $\delta_R$ where 
$R$ is a (not necessarily maximal or basic) 
rigid module in $\CC_M$, hence an element of $\cS^*$.
This finishes the proof of Theorem~\ref{main7}.

\subsection{Proof of Theorem~\ref{basesthm}}
Let $M = M_1 \oplus \cdots \oplus M_r$ be a terminal $\C Q$-module.
We fix for this section a $\GG_M$-adapted ordering 
$x(1) < x(2) < \cdots < x(r)$
of the vertices of $\GG_M$ and we may assume that $M_i$ corresponds to the
vertex $x(i)$ for $1 \leq i\leq r$. 
Moreover, we write $\alpha(i) := \dimv(M_i)$.

We have
$$
\C[\delta_{(M_1,0)},\ldots, \delta_{(M_r,0)}] \subseteq 
\RR(\CC_M) \subseteq
\Span_\C\ebrace{\delta_X \mid X \in \CC_M},
$$
where the first inclusion follows from the observation that each of the
$\LL$-modules $(M_i,0)$ for $1 \le i \le r$ is the direct summand of a
maximal rigid module on the mutation path from $T_M$ to $T_M^\vee$, 
see Section~\ref{section17}. 
The second inclusion follows from the observation
that $\Span_\C\ebrace{\delta_X \mid X \in \CC_M}$ is an algebra. 
This follows
from the fact that $\CC_M$ is an additive category together with
Theorem~\ref{Th:ClusChar}~(i).

For each $M' \in \add(M)$ 
we can choose a $\C Q$-module homomorphism 
$$
f_{M'}\df M' \to \tau(M')
$$
such that $(M',f_{M'})$ is generic in $\pi^{-1}_{Q,d}(\orb_{M'})$,
where $d = \dimv(M')$.
It follows that $\delta_{(M',f_{M'})}$ belongs to the dual semicanonical basis
$\cS^*$ of $ U(\mathfrak n)^*_{\rm gr}$.
If $M' = M_i$ is an indecomposable direct summands of $M$,
then $\Hom_{\C Q}(M',\tau(M')) = 0$ and therefore
$f_{M'}=0$.

The following theorem is a slightly more explicit statement
of Theorem~\ref{basesthm}:

\begin{Thm}\label{proofmain1}
Let $M$ be a terminal $\C Q$-module.
Then the following hold:
\begin{itemize}

\item[(i)]
We have
$$
\RR(\CC_M) = \C[\delta_{(M_1,0)},\ldots,\delta_{(M_r,0)}]
= \Span_\C\ebrace{\delta_X \mid X \in \CC_M};
$$

\item[(ii)]
The set
$$
\cP_M^* := \left\{ \delta_{(M',0)} \mid M' \in \add(M) \right\}
$$
\index{$\cP_M^*$}\noindent
is a $\C$-basis of $\RR(\CC_M)$;

\item[(iii)]
The subset
$$
\cS_M^* := \left\{ \delta_{(M',f_{M'})} \mid M' \in \add(M) \right\}
$$
\index{$\cS_M^*$}\noindent
of the dual semicanonical basis
is a $\C$-basis of $\RR(\CC_M)$, and all
cluster monomials of $\RR(\CC_M)$ belong
to $\cS_M^*$.

\end{itemize}
\end{Thm}

The basis $\cP_M^*$ will be called 
{\it dual PBW-basis} 
\index{dual PBW-basis of $\RR(\CC_M)$}
of
$\RR(\CC_M)$,
and $\cS_M^*$ the 
{\it dual semicanonical basis} 
\index{dual semicanonical basis of $\RR(\CC_M)$}
of
$\RR(\CC_M)$.
The proof of this theorem will be given after a series of lemmas.

Let $\MM_{Q,d}$ 
\index{$\MM_{Q,d}$}
be the $\C$-vector space of
$\GL_d$-invariant constructible functions
$\rep(Q,d) \to \C$.
Set
$$
\MM_Q := \bigoplus_{d \in \N^n} \MM_{Q,d}.
$$
\index{$\MM_Q$}\noindent
Note that the affine space $\rep(Q,d)$ of representations of $Q$ with
dimension vector $d$ can be viewed as an irreducible component of $\LL_d$.
In fact, $\dm \rep(Q,d) = \dm \LL_d$, and we have
$$
\rep(Q,d) = \{ (f_a,f_{a^*})_{a \in Q_1} \in \LL_d \mid 
f_{a^*} = 0 \text{ for all } a \in Q_1 \},
$$
compare Section~\ref{quivers}.
Thus we have a natural restriction 
$
\Res_Q\df \MM \to \MM_Q,
$
\index{$\Res_Q\df \MM \to \MM_Q$}\noindent
which is an algebra homomorphism.

\begin{Lem}[\cite{Lu1,S}]\label{prooflemma1}
Let $f \in \MM_d$.
If $\Res_Q(f) = 0$, then $f=0$.
\end{Lem}

Let
$$
\n = \bigoplus_{d \in \Delta^+} \n_d
$$
be the root space decomposition of $\n$.
We consider $\n$ as a subspace of the universal enveloping
algebra $U(\n)$.
Since we identify $U(\n)$ and $\MM$,
we can think of an element $f$ in $\n_d$ as a constructible function
$f\df \LL_d \to \CC$ in $\MM_d$.

\begin{Lem}\label{prooflemma2}
Let $f \in \n_d$.
If 
$d \not\in \{ \alpha(i) \mid 1 \le i \le r \}$,
then 
$$
f(X) = 0 \text{ for all } X \in \CC_M.
$$
\end{Lem}

\begin{proof}
We know from Lemma~\ref{bracketclosed1} and Lemma~\ref{bracketclosed2} that
$\{ \alpha(i) \mid 1 \le i \le r \}$ 
is a bracket closed subset of $\Delta^+$.
Thus if $X \in \CC_M$ has a dimension vector $d'$ in $\Delta^+$ 
we must have
$d' \in  \{ \alpha(i) \mid 1 \le i \le r \}$.
Therefore, since $f \in \MM_d$ with $d \in \Delta^+$
and $d \not\in \{ \alpha(i) \mid 1 \le i \le r \}$, we have 
$f(X) = 0$ for every $X \in \CC_M$.
\end{proof}

Next, we construct a PBW-basis of $U(\n)$ as follows.
Choose a $\C$-basis 
$$
{\rm P} := \{ p_j \mid j \in J \}
$$
of $\n$ consisting of weight vectors for the adjoint representation.
We number it so that
for $1 \le j \le r$ the vector $p_j$ belongs to $\n_{\alpha(j)}$.
We denote by  
$\cP = \{ p_{\bf m} \mid {\bf m} \in \N^{(J)}_{\phantom{0}} \}$ 
the PBW-basis of $U(\n)$ with respect to ${\rm P}$ (see Theorem~\ref{PBW1}).

\begin{Lem}\label{prooflemma3}
Let $p_\mm \in \cP$ where $\mm = (m_j)_{j\in J}$.
If $m_j > 0$ for some $j > r$, then
$$
p_\mm(X) = 0 \text{ for all } X \in \CC_M.
$$
Equivalently, $\delta_X(p_\mm) = 0$ for all $X \in \CC_M$.
\end{Lem}

\begin{proof}
We regard $p_\mm$ as an element of $\MM$, hence as a convolution
product
\[
p_\mm = p_1^{(m_1)}\star p_2^{(m_2)} \star \cdots \star p_s^{(m_s)}.
\] 
Let us assume that $s > r$ and $m_s > 0$.
It follows that 
$p_\mm = p \star p_s$ where 
$$
p := \frac{1}{m_s} \left(p_1^{(m_1)} 
\star p_{2}^{(m_{2})} \star \cdots \star 
p_{s-1}^{(m_{s-1})} \star p_s^{(m_s-1)}\right).
$$
Now let $X \in \CC_M$.
Then 
$$
p_{\mathbf m}(X) = (p \star p_s)(X) = 
\sum_{m \in \C} m \, \chi_c(\{ U \subseteq X \mid p(U)p_s(X/U) = m \}).
$$
Since $\CC_M$ is closed under factor modules, we get $X/U \in \CC_M$
for all submodules $U$ of $X$.
Now Lemma~\ref{prooflemma2} yields $p_s(X/U) = 0$ for all such $U$.
Thus we proved that $p_\mm(X) = 0$ for all $X \in \CC_M$.
\end{proof}

We denote by $\cP^*$ the dual of the PBW-basis $\cP$ 
of $U(\n)$.
Define
$$
\cP_M^* := \left\{ (p_1^*)^{m_1} \cdots (p_r^*)^{m_r} \mid 
m_i \ge 0 \right\} \subseteq \cP^*.
$$

Recall that the dimension vectors of
the indecomposable representations of the quiver $Q$ are known (\cite{K2}):
\begin{Thm}[Kac]\label{KacThm}
There is an indecomposable representation in $\rep(Q,d)$
if and only if $d \in \Delta^+$.
Furthermore, there is (up to isomorphism) exactly 
one indecomposable representation
in $\rep(Q,d)$ if and only if $d \in \Delta_{\rm re}^+$.
\end{Thm}

Let $\Delta\df U(\n) \to U(\n) \otimes U(\n)$
be the comultiplication of $U(\n)$.
For $f \in U(\n)$ it is well known that
$$
\Delta(f) = 1 \otimes f + f \otimes 1
$$ 
if and only if $f \in \n$, see for example \cite{D}.
As before, we identify $U(\n)$ and $\MM$, and we denote
the comultiplication of $\MM$ also by $\Delta$.
If $f \in \MM_d$, $X' \in \LL_{d'}$ and $X'' \in \LL_{d''}$ with $d = d'+d''$,
then 
$$
f(X' \oplus X'') = \Delta(f)(X',X''),
$$
see \cite{GLSSemi1}.
Thus, if $f \in \n_d$ for some $d$, 
and $X \in \LL_d$ is not indecomposable, we can write
$X = X' \oplus X''$ with
$X' \not= 0 \not= X''$ and
$$
f(X)=
(1 \otimes f)(X',X'') + (f \otimes 1)(X',X'') = f(X') + f(X'')
=0,
$$
since $\dimv(X')\not = d \not = \dimv(X'')$.
Therefore we proved the following result:

\begin{Lem}\label{prooflemma5}
For a dimension vector $d$, let $\LL_d^{\rm ind}$ be the constructible
subset of $\LL_d$ consisting of the indecomposable $\LL$-modules in $\LL_d$.
If $f \in \n_d$, then ${\rm supp}(f) \subseteq \LL_d^{\rm ind}$.
\end{Lem}

We know by Lemma~\ref{prooflemma1} 
that the map $\Res_Q\df \MM_d \to \MM_{Q,d}$
is injective.
Thus, if $0 \not= f \in \MM_d$, then $\Res_Q(f) \not= 0$.
Let $d$ be a real root.
Thus ${\rm ind}(Q,d) = \orb_{M'}$ for some $\C Q$-module $M'$.
Therefore, if $f \in \n_d$, then 
$$
\Res_Q(f) = c_{M'}{\mathbf 1}_{\orb_{M'}}
$$
for some $c_{M'} \not= 0$.
In particular, for the functions $p_1,\ldots,p_r$
in ${\rm P}$ we get 
\[
\Res_Q(p_i) = c_i {\mathbf 1}_{\orb(M_i)}
\]
(for typographical reasons we write here $\orb(M_i)=\orb_{M_i}$).
We can assume (after possibly rescaling the $p_i$) that
$c_i = 1$ for all $1 \le i \le r$.

\begin{Lem}\label{prooflemma6}
We have 
$$
\Res_{Q}(p_1^{(m_1)} \star \cdots \star p_r^{(m_r)}) =
{\mathbf 1}_{\orb(M_1)}^{(m_1)} \star \cdots \star 
{\mathbf 1}_{\orb(M_r)}^{(m_r)} = {\mathbf 1}_{\orb_{M'}}
$$
where $M' := M_1^{m_1} \oplus \cdots \oplus M_r^{m_r}$.
\end{Lem}

\begin{proof}
Since $x(1) < x(2) < \cdots < x(r)$ is a $\GG_M$-adapted ordering
of the vertices of $\GG_M$, we get that
$\Ext_{\C Q}^1(M_i,M_j) = 0$ for all $i \ge j$.
Now the result follows from the definition of the multiplication
$\star$ of constructible functions.
\end{proof}

\begin{Lem}\label{prooflemma7}
We have 
$
\cP_M^* = \left\{ \delta_{(M',0)} \mid M' \in \add(M) \right\}.
$
\end{Lem}

\begin{proof}
For $M',M'' \in \add(M)$ with
$M' := M_1^{m_1} \oplus \cdots \oplus M_r^{m_r}$
we have
\begin{align*}
\delta_{(M'',0)}(p_1^{(m_1)} \star \cdots \star p_r^{(m_r)}) &=
(p_1^{(m_1)} \star \cdots \star p_r^{(m_r)})(M'',0) \\
&= \Res_Q(p_1^{(m_1)} \star \cdots \star p_r^{(m_r)})(M'') \\
&= {\mathbf 1}_{\orb_{M'}}(M'') 
= 
{\begin{cases}
1 & \text{if $M' \cong M''$},\\
0 & \text{otherwise}.
\end{cases}}
\end{align*}
This follows from Lemma~\ref{prooflemma6}, and the fact that
$(M'',0)$ is
the image of $M''$ under the natural inclusion
$\rep(Q) \to \nil(\LL)$ which sends a $\C Q$-module $L$ to
the $\LL$-module $(L,0)$.
Hence, using also Lemma~\ref{prooflemma3}, the claim is proved.
\end{proof}

\begin{proof}[Proof of Theorem~\ref{proofmain1}]
Let $X \in \CC_M$.
By Lemma~\ref{prooflemma3} and Lemma~\ref{prooflemma7}, 
$\delta_X$ is a linear combination of
dual PBW-basis vectors of the form $\delta_{(M',0)}$ with
$M' \in \add(M)$.
Hence
$\delta_X \in \C[\delta_{(M_1,0)},\ldots,\delta_{(M_r,0)}]$,
and
$$
\Span_\C\ebrace{\delta_X \mid X \in \CC_M}
\subseteq
\C[\delta_{(M_1,0)},\ldots,\delta_{(M_r,0)}]
\subseteq
\RR(\CC_M). 
$$
Using the known reverse inclusions we
get (i) and (ii) of Theorem~\ref{proofmain1}.

Next, let $M' \in \add(M)$.
Thus $\delta_{(M',f_{M'})}$ is an element in the dual
semicanonical basis~$\cS^*$.
Since $(M',f_{M'}) \in \CC_M$, we know that$\delta_{(M',f_{M'})} \in \RR(\CC_M)$.
For dimension reasons this implies that
$$
\cS_M^* := 
\left\{ \delta_{(M',f_{M'})} \mid M' \in \add(M) \right\}
= \cS^* \cap \RR(\CC_M)
$$
\index{$\cS_M^*$}\noindent
is a $\C$-basis of $\RR(\CC_M)$.
By what we proved before, the
cluster monomials of $\RR(\CC_M)$ are a subset
of $\cS_M^*$.
This proves (iii).
\end{proof}

\subsection{Example}
Let us discuss an example of base change between 
$\cP_M^*$ and $\cS_M^*$.
Let $Q$ be the quiver
{\small
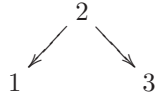
$$
\xymatrix@-0.8pc{
& 2\ar[dl]\ar[dr]\\
1 && 3
}
$$
}\noindent
The Auslander-Reiten quiver of $\md(\C Q)$ looks as follows:
$$
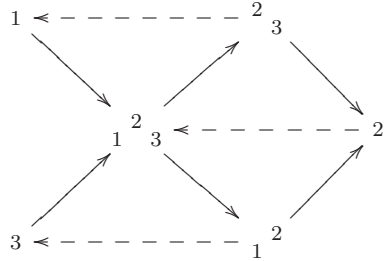
\xymatrix{
{\bsm 1\esm} \ar[dr] && {\bsm 2\\&3\esm} \ar[dr]\ar@{-->}[ll]\\
& {\bsm &2\\1&&3\esm} \ar[ur]\ar[dr] && {\bsm 2\esm}\ar@{-->}[ll]\\
{\bsm 3\esm} \ar[ur]&& {\bsm &2\\1\esm}\ar[ur]\ar@{-->}[ll] 
}
$$
The $\LL$-modules $T_{i,[a,b]}$ are the following:
\begin{align*}
T_{1,[0,0]} &= {\bsm &2\\1\esm}, & 
T_{1,[1,1]} &= {\bsm 3\esm}, &
T_{1,[0,1]} &= {\bsm &&3\\&2\\1\esm},\\
T_{2,[0,0]} &= {\bsm 2\esm}, & 
T_{2,[1,1]} &= {\bsm &2\\1&&3\esm}, &
T_{2,[0,1]} &= {\bsm &2\\1&&3\\&2\esm},\\
T_{3,[0,0]} &= {\bsm 2\\&3\esm}, & 
T_{3,[1,1]} &= {\bsm 1\esm}, &
T_{3,[0,1]} &= {\bsm 1\\&2\\&&3\esm}.
\end{align*}
Let $M = M_1 \oplus \cdots \oplus M_6$ be the direct sum
of all six indecomposable $\C Q$-modules.
For each indecomposable $\LL$-module of the form $(M_i,0)$ one
easily checks that 
$\F_{\ii,X}$ is either empty or
a single point, so $\chi_c(\F_{\ii,X})$ is either 0 or 1.

Thus the functions $\delta_{(M_i,0)}$ are known.
Using Theorem~\ref{detthm} we get
\begin{align*}
\delta_{T_{1,[0,1]}} &= \delta_{T_{1,[1,1]}}\cdot\delta_{T_{1,[0,0]}} 
- \delta_{T_{2,[1,1]}},\\
\delta_{T_{2,[0,1]}} &= \delta_{T_{2,[1,1]}}\cdot\delta_{T_{2,[0,0]}}
-\delta_{T_{1,[0,0]}}\cdot\delta_{T_{3,[0,0]}},\\
\delta_{T_{3,[0,1]}} &= \delta_{T_{3,[1,1]}}\cdot\delta_{T_{3,[0,0]}}
- \delta_{T_{2,[1,1]}}.
\end{align*}
The initial cluster of our cluster algebra $\RR(\CC_M)$ looks
as follows:
$$
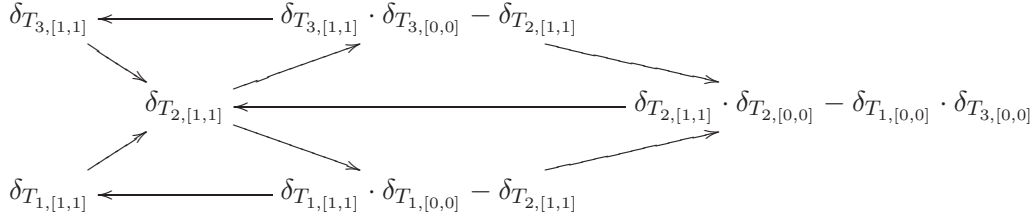
\xymatrix@=0.5cm{
{\delta_{T_{3,[1,1]}}} \ar[dr] && {\delta_{T_{3,[1,1]}}\cdot
\delta_{T_{3,[0,0]}}
- \delta_{T_{2,[1,1]}}}  \ar[dr]\ar[ll]\\
&{\delta_{T_{2,[1,1]}}} \ar[ur]\ar[dr] &&
{\delta_{T_{2,[1,1]}}\cdot\delta_{T_{2,[0,0]}}
-\delta_{T_{1,[0,0]}}\cdot\delta_{T_{3,[0,0]}}}\ar[ll] \\
{\delta_{T_{1,[1,1]}}} \ar[ur]&& 
{\delta_{T_{1,[1,1]}}\cdot\delta_{T_{1,[0,0]}} 
- \delta_{T_{2,[1,1]}}} \ar[ur] \ar[ll]
}
$$
The cluster variables in $\RR(\CC_M)$
are
$$
\left\{
\delta_{T_{i,[c,c]}} \;,\; \delta_{T_{i,[0,1]}} \mid
1 \le i \le 3 \text{ and } c = 0,1 \right\}
\cup
\left\{
\delta_{\bsm 1\\&2\esm},
\delta_{\bsm &3\\2\esm},
\delta_{\bsm 1&&3\\&2\esm}
\right\}.
$$
(Here we consider the three coefficients $\delta_{T_{i,[0,1]}}$
also as cluster variables.)
Beginning with our initial cluster we can mutate several times
and get
\begin{align*}
\delta_{\bsm 1\\&2\esm} &= \delta_{\bsm 1\esm}\cdot\delta_{\bsm 2\esm} -
\delta_{\bsm&2\\1\esm},\\
\delta_{\bsm &3\\2\esm} &= \delta_{\bsm 2\esm}\cdot\delta_{\bsm 3\esm} -
\delta_{\bsm 2\\&3\esm},\\
\delta_{\bsm 1&&3\\&2\esm} &= \delta_{\bsm &2\\1&&3\esm}
+\delta_{\bsm 1\esm}\cdot\delta_{\bsm 2\esm}\cdot\delta_{\bsm 3\esm} -
\delta_{\bsm 1\esm}\cdot\delta_{\bsm 2\\&3\esm} -
\delta_{\bsm 3\esm}\cdot\delta_{\bsm &2\\1\esm}.
\end{align*}
So we wrote all cluster variables as linear combinations of
dual PBW-basis vectors.

\subsection{Non-adaptable Weyl group elements: an example}
Let $\GG$ be the graph
{\small
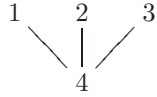
$$
\xymatrix@-0.8pc{
1 \ar@{-}[dr] & 2\ar@{-}[d] & 3 \ar@{-}[dl]\\
&4
}
$$
}\noindent
of type $\D_4$.
Let $w$ be the Weyl group element
$s_3s_4s_2s_1s_4$.
The set of reduced expressions for $w$ is
$R(w) = \{ (3,4,2,1,4), (3,4,1,2,4) \}$.
It follows that $w$ is not adaptable.
Furthermore, an easy calculation shows that
$$
\Delta_w^+ = \left\{ 
\bsm 0&0&0\\&1\esm, \bsm 1&0&0\\&1\esm,\bsm 0&1&0\\&1\esm,
\bsm 1&1&0\\&1\esm,\bsm 1&1&1\\&2\esm
 \right\}.
$$
Let $\ii=(3,4,2,1,4)$.
Using a construction dual to the one in \cite[Section II.2]{BIRS}, we
obtain a subcategory $\CC_w$ of $\md(\LL)$, which by definition has
a generator
$$
T_\ii := T_1 \oplus \cdots \oplus T_5
:=
\bsm 4 \esm \oplus 
\bsm 4\\1 \esm \oplus
\bsm 4\\2 \esm \oplus
\bsm &4\\1&&2\\&4 \esm \oplus 
\bsm &4\\1&&2\\&4\\&3\esm.
$$
It follows that
$$
\CC_w = \add\left(T_\ii \oplus \bsm &4\\1&&2 \esm\right).
$$
We can think of $\CC_w$ as a categorification of a cluster
algebra of type $\A_1$ with four coefficients.

Now let $Q$ be any quiver with $|Q| = \GG$.
By $M_1,\ldots,M_5$ we denote the five indecomposable
$KQ$-modules with
$\Delta_w^+ = \{ \dimv(M_i) \mid 1 \le i \le 5 \}$.
Then the projection $\pi_Q(T_5)$ is not in $\add(M)$.
Note also that the dimension vector of $T_5$ is a root, but
$\pi_Q(T_5)$ is a decomposable $KQ$-module.
Another calculation shows that
$$
\pi_Q^{-1}\left(\pi_Q(\CC_w)\right) \not= \CC_w.
$$

This indicates that our proof of Theorem~\ref{proofmain1}
does not work for non-adaptable Weyl group elements.

\subsection{Some generalities on bases of algebras}\label{generalbase}

We start with the following:

\begin{Lem}\label{multiproj}
Let
$M' = M_1' \oplus M_2'$ be in $\add(M)$ such that
$$
M_2' \cong \bigoplus_{a=0}^{t_i} \tau^a(I_i)
$$
for some $i$.
Then we have
$\delta_{(M',f_{M'})} = \delta_{(M_1',f_{M_1'})} \cdot \delta_{(M_2',f_{M_2'})}$.
\end{Lem}

\begin{proof}
Note that the $\LL$-module $(M_2',f_{M_2'})$ is a projective-injective object
in $\CC_M$.
The claim then follows easily from
\cite[Theorem 1.1]{GLSSemi1} in combination with
explanations in \cite[Section 2.6]{GLSSemi1}.
\end{proof}

This lemma gives rise to the following definition:
A $\C Q$-module $M' \in \add(M)$ is called {\it interval-free} 
\index{interval-free $\C Q$-module}
if
$M'$ does not have a direct summand
isomorphic to $I_i \oplus \tau(I_i) \oplus \cdots \oplus \tau^{t_i}(I_i)$.

Let ${\rm B} := \{ b_i \mid i \ge 1 \}$ be a $K$-basis of a
commutative $K$-algebra $A$.
For some fixed $n \ge 1$ let ${\rm C} := \{ b_1,\ldots,b_n \}$.
A basis vector $b \in {\rm B}$ is called 
${\rm C}$-{\it free}
if $b \notin b_i{\rm B}$ for some $b_i \in {\rm C}$.
Assume that the following hold:
\begin{itemize}

\item[(i)]
For all $b_i \in {\rm C}$ we have
$b_i{\rm B} \subseteq {\rm B}$;

\item[(ii)]
If $b_1^{z_1}\cdots b_n^{z_n}b = b_1^{z_1'}\cdots b_n^{z_n'}b'$
for some $z_i,z_i' \ge 0$ and some ${\rm C}$-free elements 
$b,b' \in {\rm B}$, then
$b=b'$ and $z_i = z_i'$ for all $i$.

\end{itemize}
It follows that
$$
{\rm B} = \left\{ b_1^{z_1} \cdots b_n^{z_n}b \mid b \in {\rm B} 
\text{ is {\rm C}-free}, z_i \ge 0 \right\}.
$$

Define
$$
\underline{A} := A/(b_1-1,\ldots,b_n-1).
$$
For $a \in A$, let $\underline{a}$ be the residue class
of $a$ in $\underline{A}$.
Furthermore, 
let $A_{b_1,\ldots,b_n}$ be the localization of
$A$ at $b_1,\ldots,b_n$.
The following lemma is easy to show:

\begin{Lem}\label{basislemma}
With the notation above, the following hold:
\begin{itemize}

\item[(1)]
The set
$
\underline{\rm B} := \{ \underline{b} \mid b 
\text{ is {\rm C}-free} \}
$
is a $K$-basis of $\underline{A}$;

\item[(2)]
The set
$
{\rm B}_{b_1,\ldots,b_n} := \left\{ b_1^{z_1} \cdots b_n^{z_n}b \mid b 
\text{ is {\rm C}-free}, z_i \in \Z \right\}
$
is a $K$-basis of $A_{b_1,\ldots,b_n}$.

\end{itemize}
\end{Lem}

\subsection{Inverting and specializing coefficients}\label{sectinvertspecial}
One can rewrite the basis $\cS_M^*$ appearing in 
Theorem~\ref{basesthm} as
$$
\cS_M^* = 
\left\{ 
(\delta_{T_{1,[0,t_1]}})^{z_1} \cdots (\delta_{T_{n,[0,t_n]}})^{z_n}
\delta_{(M',f_{M'})} 
\mid 
M' \in \add(M), M' \text{ interval-free}, z_i \ge 0 \right\}.
$$
The next two theorems deal with the situation of invertible
coefficients and specialized coefficients.

\begin{Thm}[Invertible coefficients]\label{basesthm2}
Let $M$ be a terminal $\C Q$-module.
Then 
$$
\widetilde\cS_M^* := 
\left\{ 
(\delta_{T_{1,[0,t_1]}})^{z_1} \cdots (\delta_{T_{n,[0,t_n]}})^{z_n}
\delta_{(M',f_{M'})} 
\mid 
M' \in \add(M), M' \text{ interval-free}, z_i \in \Z \right\}
$$ 
\index{$\widetilde\cS_M^*$}\noindent
is a $\C$-basis of $\widetilde{\RR}(\CC_M)$,
and
$\widetilde\cS_M^*$
contains all cluster monomials of the cluster algebra
$\widetilde{\RR}(\CC_M)$.
\end{Thm}

Next, we specialize all $n$ coefficients $\delta_{T_{i,[0,t_i]}}$
of the cluster algebra $\RR(\CC_M)$ to $1$.
We obtain a new cluster algebra 
$\underline{\RR}(\CC_M)$
\index{$\underline{\RR}(\CC_M)$}
which does not have any coefficients.
The residue class of $\delta_X \in \RR(\CC_M)$ 
is denoted by $\underline{\delta}_X$.

\begin{Thm}[No coefficients]\label{basesthm3}
Let $M$ be a terminal $\C Q$-module.
Then 
$$
\underline\cS_M^* := \left\{ \underline{\delta}_{(M',f_{M'})} \mid 
M' \in \add(M), M' \text{ interval-free} \right\}
$$ 
\index{$\underline\cS_M^*$}\noindent
is a $\C$-basis of $\underline{\RR}(\CC_M)$,
and
$\underline\cS_M^*$
contains all cluster monomials of the cluster algebra
$\underline{\RR}(\CC_M)$.
\end{Thm}

\begin{proof}[Proof of Theorem~\ref{basesthm2} and Theorem~\ref{basesthm3}]
Let ${\rm B} := \{ b_i \mid i \ge 1 \} := \cS_M^*$
be the dual semicanonical basis of $\RR(\CC_M)$.
We can label the $b_i$ such that
$$
\{ b_1,\ldots,b_n \} = \left\{ \delta_{T_{1,[0,t_1]}},\ldots,
\delta_{T_{n,[0,t_n]}} \right\}.
$$
Using Lemma~\ref{multiproj} it is easy to check that
the elements $b_i$ satisfy the properties (i) and (ii) mentioned
in Section~\ref{generalbase}.
Then apply Lemma~\ref{basislemma}.
\end{proof}


{\Large\section{Acyclic cluster algebras}\label{acycliccase}}


In this section we will study the case of acyclic cluster algebras, 
which is of special interest.
Assume that
$M$ is a terminal $\C Q$-module with $t_i(M) = 1$ for all $i$.
Thus $M$ is of the form
$$
M = \bigoplus_{i=1}^n (I_i \oplus \tau(I_i)), 
$$
then $\RR(\CC_M)$ is an acyclic cluster algebra
associated to $Q$ having $n$ coefficients,
whereas $\underline{\RR}(\CC_M)$ is the acyclic
cluster algebra associated to $Q$ having no coefficients.

\begin{Thm}
Let $M = M_1 \oplus \cdots \oplus M_{2n}$ be a terminal $\C Q$-module with
$t_i(M) = 1$ for all~$i$.
Then the following hold:
\begin{itemize}

\item[(i)]
$\RR(\CC_M) = \C[\delta_{(M_1,0)},\ldots,\delta_{(M_{2n},0)}]
= \Span_\C\ebrace{\delta_X \mid X \in \CC_M}$;

\item[(ii)]
$\left\{ \delta_{(M',0)} \mid M' \in \add(M) \right\}$ is a $\C$-basis 
of $\RR(\CC_M)$;

\item[(iii)]
$\left\{ \underline{\delta}_{(M',0)} \mid M' \in \add(M), 
M' \text{ interval-free} \right\}$ 
is a $\C$-basis of $\underline{\RR}(\CC_M)$;

\item[(iv)]
There is an isomorphism of cluster algebras
$\underline{\RR}(\CC_M) \cong \cA_Q$,
where $\cA_Q$ is the coefficient free acyclic 
cluster algebra associated to $Q$.

\end{itemize}
\end{Thm}

\begin{proof}
Part (i) and (ii) were already proved before for $M$ arbitrary.
Part (iv) is clear from our description of the initial seed
(labelled by $T_M$) for the cluster algebra $\RR(\CC_M)$.
It remains to prove (iii):
We have
$$
\RR(\CC_M) = \bigoplus_{d \in \N^n} \RR_d
$$
where $\RR_d$ is the $\C$-vector space with basis
$\left\{ \delta_{(M',f_{M'})} \mid M' \in \add(M) \cap \rep(Q,d) \right\}$.
We know that
$\left\{ \delta_{(M',0)} \mid M' \in \add(M) \cap \rep(Q,d) \right\}$
is a basis of $\RR_d$ as well.
After specializing the coefficients $\delta_{T_{i,[0,1]}}$ to 1,
we get
$$
\underline{\RR}(\CC_M) = \bigoplus_{d \in \N^n} 
\underline{\RR}_d
$$
where $\underline{\RR}_d$ is the $\C$-vector space with
basis
$$
\left\{ \underline{\delta}_{(M',f_{M'})} \mid M' \in \add(M) \cap \rep(Q,d),
M' \text{ interval-free} \right\}.
$$
Now one can use the formula
$$
\delta_{T_{i,[0,1]}} = \delta_{T_{i,[1,1]}} \cdot \delta_{T_{i,[0,0]}} -
\prod_{i \to j} \delta_{T_{j,[1,1]}} \cdot
\prod_{k \to i} \delta_{T_{k,[0,0]}}
$$
(where the products are taken over all arrows of $Q^\op$ which
start and end in $i$, respectively) and induction on the vertices
of $Q$ to show that for every interval-free $M' \in \add(M)$,
the vector $\underline{\delta}_{(M',f_{M'})}$ is a linear combination
of elements of the form $\underline{\delta}_{(M'',0)}$ where $M''$ is 
interval-free in $\add(M)$ and $|\dimv(M'')| \leq |\dimv(M')|$.
For dimension reasons we get that the vectors $\underline{\delta}_{(M'',0)}$
with $M''$ interval-free form a linearly independent set.
This implies (iii).
\end{proof}

It is interesting to compare (iii) to Berenstein, Fomin and
Zelevinsky's construction of a basis for 
the acyclic cluster algebra $\cA_Q$.
Let ${\mathbf y} := \{ y_1,\ldots,y_n \}$ be the initial
cluster whose exchange matrix $B_Q$ is encoded by $Q$, as in
Section~\ref{clustintro}.
Let $\{ y_1^*,\ldots,y_n^* \}$ be the $n$ cluster variables
obtained from $\mathbf{y}$ by mutation in the $n$ possible
directions.  

The $n$ sets $\{y_1,\ldots,y_n\} \setminus \{y_k\} \cup\{y_k^*\}$
are the neighbouring clusters of our initial cluster
$\mathbf{y}$.
Using a simple Gr\"obner basis argument, the following is shown
in \cite{BFZ}:

\begin{Thm}[Berenstein, Fomin, Zelevinsky]\label{BFZbasis}
The monomials
$$
\left\{ y_1^{p_1}(y_1^*)^{q_1} \cdots y_n^{p_n}(y_n^*)^{q_n} \mid
p_i,q_i \ge 0,\; p_iq_i = 0 \right\}
$$
form a $\C$-basis of the acyclic cluster algebra $\cA_Q$.
\end{Thm}

Now assume that the vertices $1,\ldots,n$ of $Q$ are numbered 
in such a way that $1$ is a sink of $Q$, and $k+1$ is a sink
of $\sigma_k \cdots \sigma_2\sigma_1(Q)$ for $1 \le k \le n-1$.
We perform $n$ mutations $\mu_k \cdots \mu_2\mu_1$,
$1 \le k \le n$ of the initial seed $({\mathbf y},B_Q)$.
In each step we obtain a new cluster variable which we denote by
$y_k^\dag$.
Note that $y_1^\dag = y_1^*$, but already $y_2^\dag$ and $y_2^*$
may be different.
Observe that $\mu_n \cdots \mu_2\mu_1(B_Q) = B_Q$.
We get that
$$
(\{y_1^\dag,\ldots,y_n^\dag\},B_Q)
$$
is a seed of the cluster algebra $\cA_Q$
where 
$$
\{ y_1,\ldots,y_n \} \cap \{ y_1^\dag,\ldots,y_n^\dag \} = \varnothing.
$$
Our version of Theorem~\ref{BFZbasis} looks then as follows:
\begin{Thm}
The monomials
$$
\left\{ 
y_1^{p_1}(y_1^\dag)^{q_1} \cdots y_n^{p_n}(y_n^\dag)^{q_n} \mid
p_i,q_i \ge 0, \; p_iq_i = 0 
\right\}
$$
form a $\C$-basis of the acyclic cluster algebra $\cA_Q$.
\end{Thm}
Note that in our setting, 
the initial cluster $\{y_1,\ldots,y_n\}$ comes from $T_M$
and the cluster 
$\{y_1^\dag,\ldots,y_n^\dag\}$ 
comes from
$T_M^\vee$.


{\Large\section{The coordinate rings $\C[N(w)]$ and $\C[N^w]$}\label{coordinaterings}}


\subsection{Nilpotent Lie algebras and unipotent groups} \label{ssec:nilp}
As before, let 
$$
\widehat{\n} := \prod_{\alpha \in \Delta^+} \n_\alpha
$$
\index{$\widehat{\n}$}\noindent
be the completion of the Lie algebra $\n$.
This is a pro-nilpotent pro-Lie algebra.
Let $N$ be the pro-unipotent pro-group
with Lie algebra $\widehat{\n}$.
We refer to Kumar's book \cite[Section 4.4]{Ku} for 
all missing definitions.

We can assume that $N = \widehat{\n}$ as a set and that the
multiplication of $N$ is defined via the Baker-Campbell-Hausdorff formula.
Hence the exponential map ${\rm Exp}\df \widehat{\n} \to N$
is just the identity map.

Put ${\mathcal H} := U(\n)^*_{\rm gr}$.
This is a commutative Hopf algebra.
We can regard ${\mathcal H}$ as the coordinate ring $\C[N]$ of $N$, that is,
we can identify $N$ with the set 
$$
{\rm maxSpec}({\mathcal H}) \equiv \Hom_{\rm alg}({\mathcal H},\C)
$$ 
of $\C$-algebra
homomorphisms ${\mathcal H} \to \C$.
An element $f \in \Hom_{\rm alg}({\mathcal H},\C)$ is determined
by the images $c_i := f(p_i^*)$ for all $i \ge 1$.

It is well known (see e.g. \cite[\S 3.4]{Ab}) that 
$\Hom_{\rm alg}({\mathcal H},\C)$ can also be identified with
the group $G({\mathcal H}^\circ)$ of all group-like elements of
the dual Hopf algebra ${\mathcal H}^\circ$ of ${\mathcal H}$, by mapping 
$f\in \Hom_{\rm alg}({\mathcal H},\C)$ to
\[
y_f = \sum_{\bf m}\left( \prod_i c_i^{m_i}\right) p_{\bf m}\in
G({\mathcal H}^\circ).
\]
Note that the map $f\mapsto y_f$ does not depend on the
choice of the PBW-basis $\{p_{\bf m}\}$.
Note also that $G({\mathcal H}^\circ)$ is contained in the vector space dual 
${\mathcal H}^*$
of ${\mathcal H}$, which
is the completion $\widehat{U(\n)}$ of $U(\n)$ with respect
to its natural grading.
When we use this second identification, an element 
$x\in N = \widehat{\n}$ is simply represented by
the group-like element 
$$
\exp(x) := \sum_{k\ge 0} x^k/k!
$$
in $\widehat{U(\n)}$.
To summarize, we have ${\mathcal H}=U(\n)^*_{\rm gr}\equiv \C[N]$ and
\[
N\equiv {\rm maxSpec}({\mathcal H}) \equiv \Hom_{\rm alg}({\mathcal H},\C)
\equiv 
G({\mathcal H}^\circ)\subset {\mathcal H}^\circ\subset {\mathcal H}^*
=\widehat{U(\n)}.
\]

Let $\Theta$ be a bracket closed subset of $\Delta^+$.
Recall from Section~\ref{nw} the definition of
the Lie subalgebra $\widehat{\n}(\Theta) \subseteq \widehat{\n}$.
Let $N(\Theta) := {\rm Exp}(\widehat{\n}(\Theta))$ be the corresponding
pro-unipotent pro-group.
For example, if $\alpha \in \Delta_{\rm re}^+$, then
$\Theta_\alpha := \{ \alpha \}$ is bracket closed.
In this case, $N(\alpha)$ is called the 
{\it one-parameter subgroup}
\index{one-parameter subgroup}
of $N$ associated to $\alpha$.
We have an isomorphism of groups $N(\alpha) \cong (\C,+)$.

If $\Theta$ is bracket closed and bracket coclosed, then
set $N'(\Theta) := N(\Delta^+ \setminus \Theta)$.
\index{$N'(\Theta)$}
In this case, the multiplication yields a bijection~\cite[Lemma 6.1.2]{Ku}
$$
m\df N(\Theta) \times N'(\Theta) \to N.
$$

For $w \in W$ let
$N(w) := N(\Delta_w^+)$
\index{$N(w)$} 
be the unipotent algebraic group
of dimension $l(w)$ associated to the Lie algebra $\n(w)$.
Similarly, define $N'(w) := N'(\Delta_w^+)$.
\index{$N'(w)$}

\subsection{The coordinate ring $\C[N(w)]$ as a ring of invariants}\label{invariant_ring}

We start by noting that the PBW-basis of $U(\n)$ and the dual PBW-basis 
of $U(\n)^*$
are homogeneous with respect to the (root lattice) $\N^n$-grading  of $U(\n)$. 
We write $\abs{\mm} = d \in \N^n$ in case $p_\mm$ is a homogeneous element
of degree $d\in\N^n$.
Let us denote by
$(\ee_i)_{i\in J}$ the usual coordinate vectors of $\Z^{(J)}_{\phantom{0}}$.
For example, $\abs{\ee_i}=\alpha(i)$ for $1\le i\le r$. 

The multiplication $\mu\df U(\n)\otimes U(\n) \to U(\n)$ is given
by its effect on the PBW-basis, say
$$
p_\mm \cdot p_{\nn} = 
\sum_{\abs{\kk} = \abs{\mm+\nn}} \!c_{\mm,\nn}^{\kk}\, p_{\kk}.
$$

Next, the comultiplication $\mu^*\df \C[N] \to \C[N] \otimes \C[N]$ is a
ring homomorphism, so it is determined by the value on the generators
$p^*_i = p^*_{\ee_i}$. 
By construction, we have
\[
\mu^*(p^*_i) = \sum_{\abs{\mm+\nn}=|\ee_i|} 
c_{\mm,\nn}^{\ee_i}p^*_\mm\otimes p^*_{\nn}
\]

\begin{Lem} \label{lem:coord}
Let $1 \le i \le r$ and $0 \neq \nn \in \N^{(J)}_{\phantom{0}}$ such
that $n_j=0$ for $1\leq j\leq r$. Then $c_{\mm,\nn}^{\ee_i} = 0$.
\end{Lem}

\begin{proof}
Let $\mm=\mm^<+\mm^>$ such that $m^<_j=0$ for $j>r$ and $m^>_j=0$ 
for $1\le j\leq r$, so $p_{\mm}=p_{\mm^<}\cdot p_{\mm^>}$.  
Since $\Delta_w^+$ is bracket closed and coclosed
we have
\[
p_{\mm^>}\cdot p_\nn=\sum_{\abs{\kk'}=\abs{\mm^>+\nn}} c_{\mm^>,\nn}^{\kk'} p_{\kk'}
\]
with $k'_j=0$ for $1\leq j\leq r$.
Thus
\[
p_\mm\cdot p_{\nn}=\sum_{\abs{\kk'}=\abs{\mm^>+\nn}}
c_{\mm^>,\nn}^{\kk'} p_{\kk'+\mm^<}.
\]
Putting $\kk=\kk'+\mm^<$ we get 
$c_{\mm,\nn}^\kk=c_{\mm^>,\nn}^{\kk'}$.
Thus, if in our situation $c_{\mm,\nn}^{\kk}\neq 0$ then
$k_j\neq 0$ for some $k>r$.
\end{proof}

Now, let us turn to the subgroups $N(w)$ and $N'(w)$. 
Consider the
ideals 
\[
I(w):= (p^*_{r+1},p^*_{r+2},\ldots ),\quad I'(w)=(p^*_1,\ldots,p^*_r)
\]
in $\C[N]$. 
Then we have
\begin{align*}
N(w) & =\{\nu\in\Hom_{\text{alg}}(\C[N],\C)\mid \nu(I(w))=0\},\text{ and}\\
N'(w)& =\{\nu'\in\Hom_{\text{alg}}(\C[N],\C)\mid \nu'(I'(w))=0\}.
\end{align*}
In other words we have canonically $\C[N(w)]=\C[N]/I(w)$ and 
$\C[N'(w)]=\C[N]/I'(w)$.

We consider the action of $N'(w)$ on $N$ via {\em right} multiplication. 
By definition, this comes from the 
left action of $N'(w)$ on $\C[N]$ given by
\[
\nu'\cdot f=(\text{id}\otimes\nu')\mu^*(f)
\]
for $f\in\C[N]$ and $\nu'\in N'(w)$. 
(Here we identify $\C[N]\otimes \C \equiv \C[N]$ in the canonical way.)

We denote by $\C[N]^{N'(w)}$ the 
invariant subring for this group action.

\begin{Prop}
Consider the injective ring homomorphism 
$$
\tilde{\pi}^*_w\df \C[N(w)]\to \C[N]
$$
defined by $p^*_i+I(w)\mapsto p^*_i$ for $1\leq i\leq r$. 
The corresponding
morphism (of schemes) $\tilde{\pi}_w\df N\to N(w)$ 
is  $N'(w)$-invariant and is a
retraction for the inclusion of $N(w)$ into $N$. As a consequence,
$\tilde{\pi}^*_w$ identifies $\C[N(w)]$ with 
$\C[N]^{N'(w)}=\C[p_1^*,\ldots,p_r^*]$.
\end{Prop}

\begin{proof}
We have
\[
\mu^*(p^*_i)=1\otimes p^*_i + p^*_i\otimes 1+
\sum_{\abs{\mm+\nn}=\abs{\ee_i}} c_{\mm,\nn}^{\ee_i} p^*_\mm\otimes p^*_\nn
\]
where in the last sum $\abs{\mm}\neq 0\neq\abs{\nn}$.
Thus for $1\leq i\leq r$ and $\nu'\in N'(w)$ we get
\[\nu'\cdot p^*_i = 1\cdot 0 + p^*_i \cdot 1 +
\sum_{\abs{\mm+\nn}=\abs{\ee_i}} c_{\mm,\nn}^{\ee_i} p^*_\mm\cdot\nu'(p^*_\nn)
\]
with the last sum vanishing by Lemma~\ref{lem:coord} and
the definition of $N'(w)$. 
In other words, $p^*_i \in \C[N]^{N'(w)}$ for $1\leq i\leq r$.
Thus, $\tilde{\pi}_w\df N \to N(w)$ is $N'(w)$-equivariant,
that is, $\tilde{\pi}_w(nn')=\tilde{\pi}_w(n)$ for any $n'\in N'(w)$.

Now, since the multiplication map $N(w)\times N'(w)\to N$ is bijective,
each $N'(w)$-orbit on $N$ is of the form $n\cdot N'(w)$ for 
a unique $n\in N(w)$.
We conclude that the inclusion $N(w)\hookrightarrow N$ is a section for
$\tilde{\pi}_w$. Our claim follows.
\end{proof}

By Theorem~\ref{proofmain1} we know that 
$\RR(\CC_M)=\C[p_1^*,\ldots,p_r^*]$.
Therefore we have proved:

\begin{Cor}\label{caN(w)}
Under the identification $U(\n)^*_{\rm gr}\equiv \C[N]$
the cluster algebra $\RR(\CC_M)$ gets identified to the
ring of invariants $\C[N]^{N'(w)}$, which is isomorphic to
$\C[N(w)]$.
\end{Cor}

\subsection{The coordinate ring $\C[N^w]$ as a localization 
of $\C[N]^{N'(w)}$}\label{coordinateunipotent}
We start with some generalities on Kac-Moody groups. 
Let $\Gmin$ be the Kac-Moody group with $\Lie(\Gmin)=\g$ constructed by
Kac-Peterson (see \cite[7.4]{Ku}).
It has a refined Tits system 
\[
(\Gmin, \Norm_\Gmin(H),N\cap\Gmin,N_-,H).
\] 
Write $\Nmin :=\Gmin\cap N$. 
Moreover, $\Gmin$ is
an affine ind-variety in a unique way
\cite[7.4.8]{Ku}.

For any real root $\alpha$ of $\g$, the one-parameter root 
subgroup $N(\alpha)$ is contained in $\Gmin$, 
and the $N(\alpha)$ together with $H$ generate $\Gmin$ as a group.
We have an anti-automorphism $g\mapsto g^T$ of $\Gmin$ which maps
$N(alpha)$ to $N(-\alpha)$ for each real root $\alpha$, and fixes $H$.
We have another anti-automorphism $g\mapsto g^\iota$ which fixes 
$N(\alpha)$ for every real root $\alpha$ and such that $h^\iota=h^{-1}$
for every $h\in H$.

For $i= 1,\ldots,n$ we have a unique homomorphism 
$\vph_i \df SL_2(\C) \to \Gmin$ satisfying
\[
\vph_i
\begin{pmatrix}
1& t \cr 0& 1
\end{pmatrix} 
=
\exp(te_i),
\qquad
\vph_i
\begin{pmatrix}
1& 0 \cr t& 1
\end{pmatrix} 
=
\exp(tf_i),
\qquad
(t\in\C).
\]
We define 
\[
\overline{s}_i
=
\vph_i
\begin{pmatrix}
0 & -1 \cr 1& 0
\end{pmatrix}. 
\]
For $w \in W$, we define 
$\overline{w} = \overline{s}_{i_r}\cdots\overline{s}_{i_1}$,
where $(i_r,\ldots,i_1)$ is a reduced expression for $w$.
Thus, we choose for every $w\in W$ a particular 
representative $\overline{w}$ of $w$ in the normalizer
$\Norm_G(H)$.
  
We have the following analogue of the Gaussian decomposition.

\begin{Prop}\label{Gauss}
Let $G_0$ be the subset $N_-\cdot H\cdot \Nmin$ of $\Gmin$.
\begin{itemize}

\item[(i)]
The subset $G_0$ is dense open in $\Gmin$ and each element $g\in G_0$ 
admits a unique
factorization $g=[g]_- [g]_0 [g]_+$ with $[g]_-\in N_-$, $[g]_0\in H$ and
$[g]_+\in \Nmin$.

\item[(ii)]
The map $g\mapsto [g]_+$ (resp.\@ $g\mapsto [g]_0$)
is a morphism of ind-varieties from $G_0$ to 
$\Nmin$ (resp.\@ to $H$). 

\end{itemize}
\end{Prop}

Part (i) follows from the fundamental properties of a refined Tits 
system~\cite[Theorem 5.2.3]{Ku}.
For part (ii), see \cite[Proposition 7.4.11]{Ku}.

Following Fomin and Zelevinsky~\cite{FZ5} we can now define for each 
$\vpi_j$ a generalized minor 
$\Delta_{\vpi_j,\vpi_j}$ as the regular function on $\Gmin$ such that 
\[
\Delta_{\vpi_j,\vpi_j}(g) = [g]_0^{\vpi_j},\qquad (g \in G_0). 
\]
For $w\in W$, we also define $\Delta_{\vpi_j,w(\vpi_j)}$ by
\[
\Delta_{\vpi_j,w(\vpi_j)}(g) = \Delta_{\vpi_j,\vpi_j}(g\bar{w}).
\]
The generalized minors have the following alternative description.
Let $L(\vpi_j)$ denote the irreducible highest weight $\g$-module
with highest weight $\vpi_j$. Let $v_j$ be a highest weight vector
of $L(\vpi_j)$. This is an integrable module, so it is also a representation
of $\Gmin$.

\begin{Prop}
Let $g\in G$.
The coefficient of $v_j$ in the projection of $gv_j$ on the weight
space $L(\vpi_j)_{\vpi_j}$ is equal to $\Delta_{\vpi_j,\vpi_j}(g)$.
\end{Prop}

\begin{proof}
Let $g=[g]_-[g]_0[g]_+\in G_0$.
We have $[g]_+v_j=v_j$, and $[g]_0v_j=[g]_0^{\vpi_j} v_j$.
The result then follows from the fact that $[g]_-v_j$ is
equal to $v_j$ plus elements in lower weights.
\end{proof}

\begin{Prop}
We have 
$$
G_0 = 
\left\{ g \in \Gmin \mid \Delta_{\vpi_i,\vpi_i}(g) \neq 0 
\text{ for } 1 \le i \le n \right\}.
$$
\end{Prop}

\begin{proof} 
We use the Birkhoff decomposition \cite[Theorem 5.2.3]{Ku}
\[
\Gmin = \bigsqcup_{w\in W} N_-\overline{w}H\Nmin,
\]
where $G_0$ is the subset of the right-hand side corresponding
to $w=e$.
If $g=[g]_-[g]_0[g]_+\in G_0$, then 
$\Delta_{\vpi_j,\vpi_j}(g)=[g]_0^{\vpi_j}\neq 0$.
Conversely,
if $g \not\in G_0$ we have $g= n_-whn$ for some $n_-\in N_-$,
$n\in \Nmin$, $h\in H$ and $w\not = e$.
Then for some $j$ we have $w(\vpi_j)\not = \vpi_j$ and
$\overline{w}hnv_j$ is a multiple of the extremal weight vector 
$\overline{w}v_j$.
Since the projection of $n_-\overline{w}v_j$ on the highest weight space  
of $L(\vpi_j)$ is zero, it follows that $\Delta_{\vpi_j,\vpi_j}(g)=0$.  
Finally, note that for any $j>n$ the minor $\Delta_{\vpi_j,\vpi_j}$ does 
not vanish on $\Gmin$. Indeed, the corresponding highest weight irreducible
module $L(\vpi_j)$ is one-dimensional since $\vpi_j(\alpha_i)=0$
for any $i$.
Hence in the above description of $G_0$, we may omit the 
minors $\Delta_{\vpi_j,\vpi_j}$ with $j>n$.
\end{proof}

Let us now consider the groups $N(w)$ and $N'(w)$ introduced in 
Section~\ref{ssec:nilp}.

\begin{Lem}
We have 
\begin{align*}
N(w) &= N\cap (w^{-1} N_- w),\\
N'(w) &= N \cap (w^{-1}N w),\\
N'(w) \cap \Nmin &= \Nmin \cap (w^{-1}\Nmin w).
\end{align*}
\end{Lem}

\begin{proof}
This follows from \cite[5.2.3]{Ku} and \cite[6.2.8]{Ku}.
\end{proof}

It follows that $\Delta_{\vpi_j,w^{-1}(\vpi_j)}$ is invariant under the 
action of
$N'(w) \cap \Nmin$ on $\Gmin$ via right multiplication. 
Indeed, for $g\in \Gmin$ and $n'\in N'(w) \cap \Nmin$, we have
$n'w^{-1}=w^{-1}n''$ for some $n''\in N'(w) \cap \Nmin$,
hence
\begin{align*}
\Delta_{\vpi_j,w^{-1}(\vpi_j)}(gn') &= 
\Delta_{\vpi_j,\vpi_j}(gn'\overline{w}^{-1})=
\Delta_{\vpi_j,\vpi_j}(g\overline{w}^{-1}n'')\\
&=\Delta_{\vpi_j,\vpi_j}(g\overline{w}^{-1})
=\Delta_{\vpi_j,w^{-1}(\vpi_j)}(g).
\end{align*}

Define
\[
N_w := \left\{ n\in\Nmin\mid \Delta_{\vpi_i,w^{-1}(\vpi_i)}(n)\neq 0 
\text{ for all } i \right\}.
\]
This is the subset of $\Nmin$ consisting of elements $n$ such
that $\overline{w}n^T\in G_0$.
Indeed,
\[
\Delta_{\vpi_i,w^{-1}(\vpi_i)}(n)=
\Delta_{\vpi_i,\vpi_i}(n\overline{w}^{-1})=
\Delta_{\vpi_i,\vpi_i}((n\overline{w}^{-1})^T)=
\Delta_{\vpi_i,\vpi_i}(\overline{w}n),
\]
since $\overline{w}^{-1}=\overline{w}^T$.

Following \cite[Section 5]{BZ1}, we can now define the map  
$\tilde{\eta}_w\df N_w\to \Nmin$
given by 
\[
\tilde{\eta}_w(z)=[\overline{w}z^T]_+.
\]

\begin{Prop}\label{BZprop}
The following properties hold:
\begin{itemize}

\item[(i)] The map $\tilde{\eta}_w$ is a morphism of ind-varieties.

\item[(ii)] The image of $\tilde{\eta}_w$ is $N^w$.

\item[(iii)] $\tilde{\eta}_w(x) = \tilde{\eta}_w(y)$ if and only if 
$x=yn'$ for some $n'\in N'(w)\cap \Nmin$.

\item[(iv)] $\tilde{\eta}_w$ restricts to a bijective 
morphism $N(w)\cap N_w \to N^w$.

\item[(v)] We have $N^w\subset N_w$, and 
$\tilde{\eta}_w$ restricts to a bijection 
$\eta_w\df N^w \to N^w$.

\item[(vi)] The inverse of $\eta_w$ is given by
$\eta_w^{-1}(x) = \eta_{w^{-1}}(x^\iota)^\iota$ for $x \in N^w$.
It follows that $\eta_w$ is an automorphism of $N^w$.

\end{itemize}
\end{Prop}

\begin{proof}
Property (i) follows from Proposition~\ref{Gauss}~(ii).
Next, we have
\[
[\overline{w}z^T]_+=([\overline{w}z^T]_0^{-1}
[\overline{w}z^T]_-^{-1})\overline{w}z^T \in B_-\overline{w}B_{-}.
\]
This shows that the image of $\tilde{\eta}_w$ is contained in $N^w$.
The rest of Property (ii)
and Property~(iii) are proved as in \cite[Proposition 5.1]{BZ1}.
Property~(iv) follows from (ii), (iii), and the decomposition
$\Nmin = N(w) \times (N'(w)\cap \Nmin)$. 
Finally, (v) and (vi) are proved exactly in the same way as
in \cite[Propositions 5.1, 5.2]{BZ1}.  
\end{proof}

\begin{Prop}\label{coordNw}
The map $\tilde{\pi}_w$ restricts to a morphism   
$\pi_w\df N^w\to N_w\cap N(w)$.
This is an isomorphism with inverse
$$
\eta^{-1}_w\tilde{\eta}_w\df N_w\cap N(w)\to N^w.
$$
In particular, $N^w$ is an affine
variety with coordinate ring identified to the localized
ring $\C[N]^{N'(w)}_{\Delta_w}$, 
where 
$$
\Delta_w=\prod_{i=1}^n \Delta_{\vpi_i,w^{-1}\vpi_i}.
$$
\end{Prop}

\begin{proof}
By Proposition~\ref{BZprop}~(iv) and (v), we know that 
$\eta^{-1}_w\tilde{\eta}_w$ is a bijection.
On the other hand $\tilde{\pi}_w(N^w)\subset N(w)\cap N_w$
because $N^w\subset N_w$.
Now, by Proposition~\ref{BZprop}~(iii), we have that
$$
\tilde{\eta}_w(\pi_w(x)) = \tilde{\eta}_w(x) = \eta_w(x)
$$
for every $x\in N^w$. Hence
$\eta^{-1}_w\tilde{\eta}_w\pi_w(x)=x$ for every $x$ in
$N^w$. 
So we have $\eta^{-1}_w\tilde{\eta}_w\pi_w = \id_{{N^w}}$,
and this proves that $\pi_w$ is the inverse of 
$\eta^{-1}_w\tilde{\eta}_w$.

These maps are morphisms of varieties so they induce 
isomorphisms 
\[
\C[N^w] \xrightarrow{\sim} \C[N(w)\cap N_w]=\C[N(w)]_{\Delta_w}
\xrightarrow{\sim} 
\C[N]_{\Delta_w}^{N'(w)}.
\]
\end{proof}

\subsection{Generalized minors arising from $T_M^\vee$}\label{adaptedordering}
For $w\in W$ and $1\le j\le n$, we denote by 
$D_{\vpi_j,w(\vpi_j)}$ the restriction of the
generalized minor $\Delta_{\vpi_j,w(\vpi_j)}$ to $N$.
For example, $D_{\vpi_j,\vpi_j}$ is equal to the constant
function $1$.

Recall the identifications $\MM^* \equiv U(\n)^*_{\rm gr} = \C[N]$.
To every $X \in \nil(\LL)$, we have associated a linear form
$\delta_X \in U(\n)^*_{\rm gr}$. 
We shall also denote $\delta_X$ by $\vph_X$ when we 
regard it as a function on $N$.
For $i=1,\ldots,n$, define 
$$
x_i(t) := \exp(te_i).
$$ 
The following formula shows how to evaluate $\vph_X$
on a product of $x_i(t)$'s.
To state it, we introduce some notation.
Given a sequence $\ii = (i_1,\ldots,i_k)$ and $X \in \nil(\LL)$,
we denote by $\F_{\ii,X}$ the variety of ascending flags of submodules
\[
0 = X_0 \subset X_1 \subset \cdots \subset X_k=X
\]
with $X_j/X_{j-1}\cong S_{i_j}\ (j=1,\ldots , k).$
Equivalently, a point of $\F_{\ii,X}$ can be seen as a 
composition series of $X$ of type $\ii$.

\begin{Prop}\label{phi_form}
Let $X \in \nil(\LL)$. 
We have  
\[
\vph_X(x_{i_1}(t_1)\cdots x_{i_k}(t_k)) =
\sum_{\aaa=(a_1,\ldots,a_k)\in\N^k}\chi_c(\F_{\ii^\aaa,X})
\frac{t_1^{a_1}\cdots t_k^{a_k}}{a_1!\cdots a_k!}.
\] 
Here $\ii^\aaa$ is short for the sequence 
$(i_1,\ldots,i_1,\ldots,i_k,\ldots,i_k)$ consisting of $a_1$ letters
$i_1$ followed by $a_2$ letters $i_2$, etc.
\end{Prop}

\begin{proof}
By Section~\ref{ssec:nilp} we can regard 
$x_{i_1}(t_1)\cdots x_{i_k}(t_k)$ as an element of $\widehat{U(\nn)}$,
namely,
\[
x_{i_1}(t_1)\cdots x_{i_k}(t_k)=\sum_{\aaa=(a_1,\ldots,a_k)\in\N^k}
\frac{t_1^{a_1}\cdots t_k^{a_k}}{a_1!\cdots a_k!}
e_{i_1}^{a_1}\cdots e_{i_k}^{a_k}.
\]
It follows from the identification $\vph_X=\delta_X$ that
\[
\vph_X(x_{i_1}(t_1)\cdots x_{i_k}(t_k))=
\sum_{\aaa=(a_1,\ldots,a_k)\in\N^k}
\frac{t_1^{a_1}\cdots t_k^{a_k}}{a_1!\cdots a_k!}
\delta_X(e_{i_1}^{a_1}\cdots e_{i_k}^{a_k}).
\] 
Now, in the geometric realization $\MM$ of the enveloping 
algebra $U(\n)$ in terms of constructible functions, 
$e_{i_1}^{a_1}\cdots e_{i_k}^{a_k}$ becomes the convolution 
product ${\bf 1}_{i_1}^{a_1} \star \cdots \star {\bf 1}_{i_k}^{a_k}$
and it is easy to see that 
$$
\delta_X({\bf 1}_{i_1}^{a_1} \star \cdots \star {\bf 1}_{i_k}^{a_k})
= \chi_c(\F_{\ii^\aaa,X}).
$$
\end{proof}

\begin{Rem}
{\rm
The formula for $\vph_X$ given in \cite[\S 9]{GLSRigid} 
involves descending flags instead of ascending flags of
submodules of $X$. This is because in the present paper
we have taken a convolution product $\star$ opposite to
that of our previous papers (see Remark~\ref{newconv}).
}
\end{Rem}

We are going to show that all generalized minors 
$D_{\vpi_j,w(\vpi_j)}$ can be expressed as some functions
$\vph_X$ for certain $\LL$-modules $X$.
In order to do this, we need to recall some results on 
Kac-Moody groups.

Let $\C[\Gmin]_{\rm s.r.}$ denote the algebra of strongly regular
functions on $\Gmin$ \cite[\S 2C]{KP1}.
Define the invariant ring
\[
\C[N_-\backslash \Gmin]_{\rm s.r.}
=
\left\{f\in \C[\Gmin]_{\rm s.r.}\mid f(ng)=f(g)
\text{ for all } n\in N_-,\ g\in\Gmin \right\}.
\]
This ring is endowed with the usual left action of $\Gmin$
given by
\[
(g\cdot f)(g') = f(g'g),\qquad (f\in \C[N_-\backslash \Gmin]_{\rm s.r.},
\,g,g'\in \Gmin).
\] 
It was proved by Kac and Peterson \cite[Cor. 2.2]{KP1} that
as a left $\Gmin$-module, it decomposes as follows
\[
\C[N_-\backslash \Gmin]_{\rm s.r.}=\bigoplus_{\lambda\in P^+} L(\lambda).
\]
This is a multiplicity-free decomposition, in which the highest 
weight irreducible module $L(\lambda)$ is carried by the subspace
\[
S(\lambda) = \left\{f\in \C[N_-\backslash\Gmin]_{\rm s.r.}
\mid f(hg) = \Delta_\lambda(h)f(g) 
\text{ for all } h\in H,\ g\in\Gmin \right\},
\]
where we denote 
$$
\Delta_\lambda : = \prod_j \Delta_{\vpi_j,\vpi_j}^{\lambda(\alpha_j^\vee)}.
$$
Clearly, $\Delta_\lambda\in S(\lambda)$ and is a highest weight vector.
Moreover, for any $w\in W$, the one-dimensional extremal
weight space of $S(\lambda)$ with weight $w(\lambda)$ is spanned
by 
$$
\Delta_{w(\lambda)} := 
\prod_j \Delta_{\vpi_j,w(\vpi_j)}^{\lambda(\alpha_j^\vee)}.
$$ 

Now consider the restriction map 
\[
\rho\df \C[N_-\backslash \Gmin]_{\rm s.r.} \to \C[\Nmin]_{\rm s.r.}
\]
given by restriction of functions from $\Gmin$ to $\Nmin$.

\begin{Lem}
For every $\lambda\in P^+$, the restriction 
$$
\rho_\lambda\df S(\lambda) \to \C[N^{\min}]_{\rm s.r.}
$$ 
of $\rho$ to $S(\lambda)$ is injective.
\end{Lem}

\begin{proof}
We have 
$$
\Nmin \subset G_0 = N_-H\Nmin = B_-\Nmin.
$$
It follows that the natural projection from $\Gmin$ onto
$B_-\backslash \Gmin$ restricts to an embedding of $\Nmin$,
with image the open subset of the flag variety 
$\X = B_-\backslash\Gmin$ defined by the non-vanishing of the minors
$\Delta_{\vpi_j,\vpi_j}$.
Now $\C[N_-\backslash \Gmin]_{\rm s.r.}$ can be regarded as the 
multi-homogeneous coordinate ring of $\X$ with homogeneous
components $S(\lambda)\ (\lambda\in P^+)$.
It follows that $\C[\Nmin]$ can be identified with the
subring of degree 0 homogeneous elements of the localized
ring obtained from $\C[N_-\backslash \Gmin]_{\rm s.r.}$
by formally inverting the element 
$$
\Delta=\prod_j\Delta_{\vpi_j,\vpi_j}.
$$
Therefore, the restriction $\rho_\lambda$ of $\rho$ to every 
homogeneous piece $S(\lambda)$ is an embedding.
\end{proof}

It follows that we can transport the $\Gmin$-module structure from
$S(\lambda)$ to $\rho(S(\lambda))$ by setting
\[
g\cdot \vph = \rho(g\cdot \rho_\lambda^{-1}(\vph)),
\qquad
(g\in\Gmin,\ \vph\in \rho(S(\lambda))).
\]
In this way, we can identify the highest weight module $L(\lambda)$
with the subspace $\rho(S(\lambda))$ of $\C[\Nmin]_{\rm s.r.}$. 
The highest weight vector is now $\rho(\Delta_\lambda)=1$, and the 
extremal weight vectors are the minors 
$$
D_{w(\lambda)} := \prod_j D_{\vpi_j,w(\vpi_j)}^{\lambda(\alpha_j^\vee)},
$$
for $w\in W$. 

At this point, we note that a strongly regular function on $\Nmin$
is just the same as an element of $U(\n)^*_{\rm gr}$.
Indeed, the elements of $\C[N^{\min}]_{\rm s.r.}$ are the restrictions
to $\Nmin$ of the linear combinations of matrix coefficients 
of the irreducible integrable representations 
$L(\lambda)\ (\lambda\in P^+)$ of $\Gmin$
(see \cite[Lemma 4.2]{KP1}).
Now, by Theorem~\ref{thm:irred-l}, we can realize every  
$L(\lambda)$ as a subspace of $U(\n)^*_{\rm gr}$, and every 
$f\in U(\n)^*_{\rm gr}$ belongs to such a subspace for $\lambda$
sufficiently dominant. It follows that each element of $U(\n)^*_{\rm gr}$
can be seen as a matrix coefficient for some   
$L(\lambda)$, and vice-versa.
We can therefore identify
\[
\C[\Nmin]_{\rm s.r.} \equiv U(\n)^*_{\rm gr} \equiv \C[N].
\]
Moreover, these two ways of embedding $L(\lambda)$ in $\C[N]$ coincide.

\begin{Lem}\label{embedding2}
Under the identification $U(\n)^*_{\rm gr} \equiv \C[\Nmin]_{\rm s.r.}$,
the embedding of $L(\lambda)$ in $U(\n)^*_{\rm gr}$ given by
Theorem~\ref{thm:irred-l}
coincides
with $\rho(S(\lambda))$.
\end{Lem}

\begin{proof}
The natural right action of $U(\n)$ on $U(\n)^*_{\rm gr}$ defined
before Corollary~\ref{embedding1} coincides with the
right action of $U(\n)$ on $\C[\Nmin]_{\rm s.r.}$ obtained by differentiating
the right regular representation of $\Nmin$:
\[
(f\cdot n)(x) = f(nx),\qquad (x,n\in\Nmin,\ f\in\C[\Nmin]_{\rm s.r.}).
\]
Consider first the case of a fundamental weight  $\lambda=\vpi_j$.
It is easy to check that 
\[
\Delta_{\vpi_j,\vpi_j}(x_i(t)g)=
\left\{
\begin{array}{ll}
\Delta_{\vpi_j,\vpi_j}(g)&
\mbox{if $i\not = j$}, \\
\Delta_{\vpi_j,\vpi_j}(g)
+t\Delta_{s_j(\vpi_j),\vpi_j}(g)&
\mbox{if $i = j$}.
\end{array}
\right.
\]
Now, the subspace $\rho(S(\lambda))$ is spanned by the functions
$n \mapsto \Delta_{\vpi_j,\vpi_j}(ng),\ (n \in N_-, \ g\in\Gmin)$.
By differentiating the previous equation with respect to $t$
and setting $t=0$, we obtain that 
\[
\rho(S(\lambda)) \subset \left\{ f\in \C[\Nmin]_{\rm s.r.}\mid f\cdot e_i = 0 
\text{ for } i\not =j, \ 
f\cdot e_j^2 = 0 \right\}.
\] 
Hence, using Corollary~\ref{embedding1}, we see that 
$\rho(S(\lambda))$ is contained in the embedding of $L(\vpi_j)$
given by the dual Verma module. Since these spaces have
the same graded dimensions, they must coincide.
The case of a general $\lambda$ follows using the fact that
$$
\Delta_\lambda = \prod_j \Delta_{\vpi_j,\vpi_j}^{\lambda(\alpha_j^\vee)}
$$
and that the $e_i$'s act as derivations on $\C[\Nmin]_{\rm s.r.}$.
\end{proof}

Let $M = M_1 \oplus \cdots \oplus M_r$ be a terminal $\C Q$-module.
Let
$w = w(M) = s_{i_r}\cdots s_{i_1}$ be the element of $W$ attached to $M$
with its reduced expression $(i_r,\ldots,i_1)$
obtained from a $\GG_M$-adapted ordering, as in
Lemma~\ref{bracketclosed2}.
Set 
\[
w^{-1}_{\le k} = s_{i_1}\cdots s_{i_k},\qquad (k=1,\ldots ,r).
\]
Finally, let $T_M^\vee = T_1 \oplus \cdots \oplus T_r$ be the 
$\CC_M$-maximal rigid $\LL$-module 
defined in Section~\ref{completerigid}.
Here we number the indecomposable direct summands of $T_M^\vee$
using the same $\GG_M$-adapted ordering.

\begin{Prop}\label{minors}
We have
\[
\vph_{T_k}=
D_{\vpi_{i_k},w^{-1}_{\le k}(\vpi_{i_k})},\qquad (k=1,\ldots ,r).
\]
In particular, we have $D_{\vpi_i,w^{-1}(\vpi_i)}=\vph_{T_{i,[0,t_i]}}$
for every $1\le i \le n$.
\end{Prop}

\begin{proof} 
Using Lemma~\ref{embedding2}, we can realize the fundamental module
$L(\vpi_{i_k})$ as the subspace $\rho(S(\vpi_k))$ of $\C[N]$.
Then using Theorem~\ref{thm:irred-l} and 
arguing as in the proof of \cite[Theorem 2]{GLSAus},
we can check that the function $\vph_{T_k}$ is an extremal vector
of weight $w^{-1}_{\le k}(\vpi_{i_k})$ in $L(\vpi_{i_k})$, hence
it coincides with $D_{\vpi_{i_k},w^{-1}_{\le k}(\vpi_{i_k})}$ up to a scalar.
Moreover, its image under $e_{i_k}^{\rm max}\cdots e_{i_1}^{\rm max}$
is equal to $1$, so the normalizations agree and we have
$\vph_{T_k} = D_{\vpi_{i_k},w^{-1}_{\le k}(\vpi_{i_k})}$. 
\end{proof}

\subsection{Example}\label{sect22_5}
Let $\GG_M$ be the following quiver, which appeared already in
Section~\ref{preinjective}.
{\small
$$
\xymatrix@-0.5pc{ 
(1,2) \ar@<0.5ex>[dr]\ar@<-0.5ex>[dr] &&
(1,1) \ar@<0.5ex>[dr]\ar@<-0.5ex>[dr] && (1,0)
\\
&(2,1) \ar@<0.5ex>[ur]\ar@<-0.5ex>[ur]\ar[dr] && 
(2,0) \ar@<0.5ex>[ur]\ar@<-0.5ex>[ur]
\\
(3,1) \ar[ur] && 
(3,0) \ar[ur]
}
$$
}\noindent
The following is
a $\GG_M$-adapted ordering:
{\small
$$
\xymatrix@-0.5pc{ 
x(6) \ar@<0.5ex>[dr]\ar@<-0.5ex>[dr] &&
x(3) \ar@<0.5ex>[dr]\ar@<-0.5ex>[dr] && x(1)
\\
&x(5) \ar@<0.5ex>[ur]\ar@<-0.5ex>[ur]\ar[dr] && 
x(2) \ar@<0.5ex>[ur]\ar@<-0.5ex>[ur]
\\
x(7) \ar[ur] && 
x(4) \ar[ur]
}
$$
}\noindent
This gives $\ii=(3,1,2,3,1,2,1)$.
The indecomposable direct summands of $T_M^\vee$ are
\begin{align*}
T_1 &= {\bsm1\esm}, &
T_2 &= {\bsm1&&1\\&2\esm}, &
T_3 &= {\bsm1&&1&&1\\&2&&2\\&&1\esm}, \\
T_4 &= {\bsm1&&1\\&2\\&3\esm},&
T_5 &= 
{\bsm&\\&\\1&&1&&1&&1&&1&&1\\&2&&2&&&&2&&2\\&&1&&&3&&&1\\&&&&&2\esm}, \\
T_6 &= 
{\bsm&\\&\\1&&1&&1&&1&&1&&1&&1&&1&&1\\&2&&2&&&&2&&2&&&&2&&2\\
&&1&&&3&&&1&&&3&&&1
\\&&&&&2&&&&&&2\\&&&&&&&&1\esm}, &
T_7 &=
 {\bsm&\\&\\1&&1&&1&&1\\&2&&2&&2\\&&1&&1\\&&&2\\&&&3\esm}.
\end{align*}
Here, the $\LL$-modules are represented by their radical filtration.
The indecomposable $\CC_M$-injective modules are $T_5, T_6$ and $T_7$.
The corresponding functions $\vph_{T_i}$ are given by
\begin{align*}
\vph_{T_1} &= D_{\vpi_1,s_1(\vpi_1)}, &
\vph_{T_2} &= D_{\vpi_2,s_1s_2(\vpi_2)}, &
\vph_{T_3} &= D_{\vpi_1,s_1s_2s_1(\vpi_1)}, \\
\vph_{T_4} &= D_{\vpi_3,s_1s_2s_1s_3(\vpi_3)}, &
\vph_{T_5} &= D_{\vpi_2,s_1s_2s_1s_3s_2(\vpi_2)}, \\
\vph_{T_6} &= D_{\vpi_1,s_1s_2s_1s_3s_2s_1(\vpi_1)}, &
\vph_{T_7} &= D_{\vpi_3,s_1s_2s_1s_3s_2s_1s_3(\vpi_3)}.
\end{align*}

It is also interesting to calculate the expansions of
these minors in terms of the dual PBW-basis 
using our generalized determinantal identity in Theorem~\ref{detthm}.
For this, it is convenient to change notation and to
write 
\begin{align*}
T_1 &= T_{1,[0,0]}, &
T_2 &= T_{2,[0,0]}, &
T_3 &= T_{1,[0,1]}, &
T_4 &= T_{3,[0,0]}, \\
T_5 &= T_{2,[0,1]}, &
T_6 &= T_{1,[0,2]}, &
T_7 &= T_{3,[0,1]}.
\end{align*}
With this notation, our initial seed becomes
{\small
$$
\xymatrix@-0.5pc{ 
T_{1,[0,2]} \ar@<0.5ex>[dr]\ar@<-0.5ex>[dr] &&
T_{1,[0,1]}  \ar@<0.5ex>[dr]\ar@<-0.5ex>[dr] \ar[ll]&&
T_{1,[0,0]} \ar[ll]
\\
&T_{2,[0,1]} \ar@<0.5ex>[ur]\ar@<-0.5ex>[ur]\ar[dr] && 
T_{2,[0,0]} \ar@<0.5ex>[ur]\ar@<-0.5ex>[ur] \ar[ll]
\\
T_{3,[0,1]} \ar[ur] && 
T_{3,[0,0]} \ar[ur]\ar[ll]
}
$$
}\noindent
Now, writing $T_{i,[a,b]}$ instead of $\vph_{T_{i,[a,b]}}$,
we obtain immediately
\begin{align*}
T_{1,[0,1]} &= T_{1,[1,1]}T_{1,[0,0]} - T_{2,[0,0]}^2,\\
T_{2,[0,1]} &= T_{2,[1,1]}T_{2,[0,0]} - T_{1,[1,1]}^2T_{3,[0,0]},\\
T_{3,[0,1]} &= T_{3,[1,1]}T_{3,[0,0]} - T_{2,[1,1]}.
\end{align*}
Again by Theorem~\ref{detthm} and a short calculation we get
\begin{multline*}
T_{1,[0,2]} = T_{1,[2,2]}T_{1,[1,1]}T_{1,[0,0]} -
T_{1,[2,2]}T_{2,[0,0]}^2 - 
T_{2,[1,1]}^2T_{1,[0,0]} +\\
+ 2\, T_{2,[1,1]}T_{2,[0,0]}T_{1,[1,1]}T_{3,[0,0]} -
T_{1,[1,1]}^3T_{3,[0,0]}^2.
\end{multline*}
%

\subsection{Proof of Theorem~\ref{wident}}
Everything is now ready for the proof of Theorem~\ref{wident}.
By Proposition~\ref{coordNw}, we know that $\C[N^w]$ is
the localization of the ring $\C[N(w)]$ with respect to the
minors $D_{\vpi_i,w^{-1}(\vpi_i)}$. 
By Corollary~\ref{caN(w)}, $\C[N(w)]$ is equal to the cluster
algebra $\RR(\CC_M)$. 
By Proposition~\ref{minors}, the minors $D_{\vpi_i,w^{-1}(\vpi_i)}$
coincide with the functions $\vph_X$ where $X$ runs through the 
set of indecomposable $\CC_M$-projective-injective.
In other words, the $D_{\vpi_i,w^{-1}(\vpi_i)}$ coincide with
the generators of the coefficient ring of $\RR(\CC_M)$.
Hence $\C[N^w]$ is equal to the cluster algebra
$\widetilde{\RR}(\CC_M)$.

\subsection{The subcategory $\CC_w$}\label{subcatCw}
In Lemma~\ref{bracketclosed2} we have associated to a terminal 
$\C Q$-module $M$ a $Q^\op$-adaptable element $w$ of the Weyl group.
It is not difficult to see that the map $M \mapsto w$ is a bijection
from the set of isomorphism classes of terminal $\C Q$-modules to the set of
sincere $Q^\op$-adaptable elements $w$ of $W$.
(Here $w$ is called {\it sincere} if for any reduced expression
$(i_t,\ldots,i_1) \in R(w)$ we have $i \in \{ i_1,\ldots,i_t \}$
for all $1 \le i \le n$. Looking at sincere Weyl group elements
is not really a restriction, since we could just pass from $Q$
to a smaller quiver by deleting the vertices which do not
occur in $\{ i_1,\ldots,i_t \}$.)

On the other hand, an element $w$ of $W$ may be adaptable to several
orientations $Q$ of the Dynkin graph of $\g$.
(For example if $\g$ is finite-dimensional, the longest element 
$w_0$ of $W$ is adaptable to every orientation of the Dynkin graph.)
In this case, $w$ is associated with a terminal $\C Q$-module $M=M_Q$
for several orientations $Q$ of the Dynkin graph.

\begin{Lem}
Let $Q$ and $Q'$ be two orientations of the Dynkin graph
of $\g$.
Let $w \in W$ be $Q^\op$-adaptable and ${Q'}^\op$-adaptable,
and let $M \in \md(\C Q)$ and $M' \in \md(\C Q')$ be the
corresponding terminal modules. 
Then, the subcategories $\CC_{M}$ and $\CC_{M'}$ of $\nil(\LL)$
are equal. 
\end{Lem}

\begin{proof}
Recall that, by Proposition~\ref{minors}, the elements of $\C[N]$
attached to the $\CC_M$-injective indecomposable direct summands 
of $T_M^\vee$ are the $D_{\vpi_i,w^{-1}(\vpi_i)}$ for $1\le i\le n$.
In particular, they depend only on $w$ and not on the choice
of a reduced expression.
Now, if $X$ and $Y$ are two rigid modules such that 
$\vph_X=\vph_Y$, we have that $X$ is isomorphic to $Y$.
Otherwise, the closures of their orbits would be different
irreducible components of the nilpotent variety on which they lie, 
and therefore $\vph_X$ and $\vph_Y$ would be different elements
of the dual semicanonical basis.
It follows that 
$\CC_{M}$ and $\CC_{M'}$ have the same generator-cogenerator,
so they are equal.
\end{proof}

Therefore we may define $\CC_w := \CC_M$, and 
we have associated to every \emph{adaptable} element $w\in W$
a subcategory $\CC_w$ of $\nil(\LL)$. 
This subcategory is, up to duality, the same as the one introduced
in a different manner in \cite{BIRS}.
To check this, one only needs to compare our maximal rigid module
$T_M^\vee$
(as described in Proposition~\ref{minors}) with the cluster tilting
object 
$$
\bigoplus _{k=1}^r \LL/I_{s_{i_1}\cdots s_{i_k}}
$$ 
defined in \cite[Theorem II.2.6]{BIRS}. 
More precisely, if in our construction we would replace our 
terminal $\C Q$-module $M$ in the preinjective component by
an ``initial'' $\C Q$-module in the preprojective component,
then we would get exactly the cluster tilting objects and
the subcategories $\CC_w$ of
\cite[Theorem II.2.8]{BIRS}, but only for the adaptable
elements $w$ of $W$.
 
Thus we have proved:

\begin{Thm}
Conjecture~III.3.1 of \cite{BIRS} holds for every adaptable element 
$w\in W$.
\end{Thm}

\subsection{Calculation of Euler characteristics}

We retain the notation of Section~\ref{adaptedordering}. 
In particular $M=M_1\oplus \cdots \oplus M_r$ is a terminal
$\C Q$-module, and $w=w(M)=s_{i_r}\cdots s_{i_1}$ is the corresponding
Weyl group element.
Let $T$ be a $\CC_M$-maximal rigid module in the mutation class of
$T_M$, or equivalently in the mutation class of 
$T_M^\vee=T_1\oplus\cdots\oplus T_r$.
Let $X$ be an indecomposable direct summand of $T$
and let $\jj=(j_1,\ldots,j_d)$. 
By Proposition~\ref{phi_form},
the Euler characteristic
$\chi_c(\F_{\jj,X})$ is equal to the coefficient of $t_1\cdots t_d$ in
$\vph_X(x_{j_1}(t_1)\cdots x_{j_d}(t_d))$.
Using mutations, we can express algorithmically $\vph_X$
as a Laurent polynomial in the functions $\vph_{T_i}\ (i=1,\ldots,r)$.
Again, by Proposition~\ref{phi_form}, to evaluate $\vph_{T_i}$
on $x_{j_1}(t_1)\cdots x_{j_d}(t_d)$, we only need to know
the Euler characteristic $\chi_c(\F_{\kk,T_i})$ for arbitrary
types $\kk$ of composition series.
These Euler characteristics can in turn be calculated via
a simple algorithm that we shall now describe.

To this end, it will be convenient to embed $U(\n)^*_{\rm gr}\cong\C[N]$
in the shuffle algebra $F^*$, as explained in \cite[\S 2.8]{LeMZ}.
As a $\C$-vector space, $F^*$ has a basis consisting of all words
\[
w[\kk]=w[k_1,k_2,\ldots,k_s]:=w_{k_1}w_{k_2}\cdots w_{k_s},
\qquad (1\le k_1,\ldots,k_s\le n,\ s\in\N),
\]
in the letters $w_1,\ldots, w_n$.
The multiplication in $F^*$ is the classical 
commutative shuffle product $\shuffle$ of words 
with unit the empty word $w[]$ (see e.g. 
\cite[\S 2.5]{LeMZ}).
By \cite[Prop. 9, 10]{LeMZ}, for any $X\in\nil(\LL)$ the image of $\vph_X$ in 
this embedding is just the generating function 
\[
g_X := \sum_{\kk} \chi_c(\F_{\kk,X}) w[\kk]
\]
of the Euler characteristics $\chi_c(\F_{\kk,X})$ for all types
$\kk$ of composition series. 

Let $\la\in P^+$ and $1\le i \le n$. 
Define endomorphisms $\rho_\la(e_i), \rho_\la(f_i)$ of the vector space $F^*$ by 
\begin{eqnarray*}
\rho_\la(e_i) (w[j_1,\ldots ,j_k]) &=& \delta_{i,j_k} w[j_1, \ldots
,j_{k-1}],\label{eq2}\\[2mm]
\rho_\la(f_i) (w[j_1,\ldots ,j_k]) &=& \sum_{r=0}^k 
(\la-\alpha_{j_1}-\cdots
-\alpha_{j_r})(\alpha_i^\vee)\,w[j_1,\ldots,j_r,i,j_{r+1},\ldots , j_k].
\end{eqnarray*}
\begin{Prop}\label{prop1}
The formulas above extend to a linear representation  
$\rho_\la$ of $\g$ on $F^*$. This turns $F^*$ into a $U(\g)$-module.
The image of $\C[N]$ in its embedding in $F^*$
is a $U(\g)$-submodule 
isomorphic to the dual Verma module $M(\la)_{\rm low}^*$
(see Section~\ref{dualVerma}).
In particular the set
\[
\{\rho_\la(f_{i_1}\cdots f_{i_s})(w[]) \mid s\in\N, 1\le i_1,\ldots ,i_s
\le n\}
\]
spans a copy of the irreducible module $L(\la)$. 
\end{Prop}
The above formulas for $\rho_\la(e_i)$ and $\rho_\la(f_i)$
can be obtained by specializing $q$ to 1 in the formulas of the
proof of \cite[Prop. 50]{LeMZ}. We omit the details.

By Proposition~\ref{minors}, we have
\[
\vph_{T_k}=
D_{\vpi_{i_k},w^{-1}_{\le k}(\vpi_{i_k})},\qquad (k=1,\ldots ,r),
\]
that is, $\vph_{T_k}$ is the (suitably normalized) extremal weight vector
of $L(\vpi_{i_k})$ with weight $s_{i_1}\cdots s_{i_k}(\vpi_{i_k})$. 
This implies that $\vph_{T_k}$ is obtained by acting on the highest weight vector
of $L(\vpi_{i_k})$ with the product 
$f_{i_1}^{(b_1)}\cdots f_{i_k}^{(b_k)}$ of divided powers of the Chevalley 
generators, where $b_k=\vpi_{i_k}(\alpha_{i_k}^\vee)=1$ and   
$b_j=(s_{i_{j+1}}\cdots s_{i_k}(\vpi_{i_k}))(\alpha_{i_j}^\vee)
\ (j=1,\ldots, k-1)$.
Therefore we have
\begin{equation}\label{calculEuler}
g_{T_k}=\rho_{\vpi_{i_k}}\left(f_{i_1}^{(b_1)}\cdots f_{i_k}^{(b_k)}\right)(w[]).
\end{equation}
Hence to calculate the generating function $g_{T_k}$ one only needs
to apply $b_1+\cdots + b_k=\dim T_k$ times the above combinatorial formula 
for $\rho_{\varpi_{i_k}}(f_i)$.

Thus we have obtained an algorithm for calculating the 
Euler characteristics $\chi_c(\F_{\kk,T})$ 
for any rigid module $T$ in the mutation class of $T_M^\vee$.
This applies in particular to every summand $M_i$ of the
terminal $\C Q$-module $M$. Hence, by varying $M$, we see that the
Euler characteristics $\chi_c(\F_{\kk,L})$ are (in principle)
computable for any preinjective $\C Q$-module $L$.

\subsection{Example}

We continue the example of Section~\ref{sect22_5}.
Clearly, we have 
\[
g_{T_1}=\rho_{\vpi_1}(f_1)(w[])= \vpi_1(\alpha_1^\vee)w[1]=w[1].
\]
Similarly
\[
g_{T_2}=\rho_{\vpi_2}(f_1^{(2)}f_2)(w[]).
\]
Now we calculate successively
\[
\begin{array}{rcl}
\rho_{\vpi_2}(f_2)(w[])&=&\vpi_2(\alpha_2^\vee)\,w[2]\ =\ w[2],\\
\rho_{\vpi_2}(f_1)(w[2])&=&\vpi_2(\alpha_1^\vee)\,w[1,2]
+(\vpi_2-\alpha_2)(\alpha_1^\vee)\,w[2,1]
\ =\ 2\,w[2,1],\\
\rho_{\vpi_2}(f_1)(2\,w[2,1])&=&
2(\vpi_2(\alpha_1^\vee)\,w[1,2,1]
+(\vpi_2-\alpha_2)(\alpha_1^\vee)\,w[2,1,1]\\
&&+(\vpi_2-\alpha_2-\alpha_1)(\alpha_1^\vee)\,w[2,1,1])
\\
&=&4\,w[2,1,1].
\end{array}
\]
Hence, taking into account that $f_1^{(2)}=f_1^{2}{/}2$, we get
\[
g_{T_2}=2\,w[2,1,1].
\]
Similar applications of formula~(\ref{calculEuler}) yield the following results
\[
\begin{array}{rcl}
g_{T_3}&=& \rho_{\vpi_1}\left(f_1^{(3)}f_2^{(2)}f_1\right)(w[]) 
= 4\,w[1,2,1,2,1,1]+12\,w[1,2,2,1,1,1],\\
g_{T_4}&=& \rho_{\vpi_3}\left(f_1^{(2)}f_2f_3\right)(w[])
= 2\,w[3,2,1,1],\\
g_{T_7}&=&\rho_{\vpi_3}\left(f_1^{(4)}f_2^{(3)}f_1^{(2)}f_2f_3\right)(w[])\\
&=& 288\,w[3,2,1,1,2,2,2,1,1,1,1]+144\,w[3,2,1,1,2,2,1,2,1,1,1]\\
       && +96\,w[3,2,1,2,1,2,2,1,1,1,1]+48\,w[3,2,1,1,2,2,1,1,2,1,1]\\
       && +48\,w[3,2,1,2,1,1,2,2,1,1,1]+48\,w[3,2,1,2,1,2,1,2,1,1,1]\\
       && +48\,w[3,2,1,1,2,1,2,2,1,1,1]+16\,w[3,2,1,2,1,2,1,1,2,1,1]\\
       && +16\,w[3,2,1,2,1,1,2,1,2,1,1]+16\,w[3,2,1,1,2,1,2,1,2,1,1].
\end{array}
\]
The generating functions $g_{T_5}$ and $g_{T_6}$ are too large to be included
here. For example $g_{T_5}$ is a linear combination of 402 words.

\subsection{Coordinate rings of unipotent radicals}

In this section, we assume that $Q$ is of finite Dynkin type $\A, \D, \E$.
We first recall some standard notation (we refer the reader to 
\cite{GLSFlag} for more details).
The group $G$ is now a simple complex algebraic group of the same type as $Q$.
Let $J$ be a subset of the set of vertices of $Q$, and let $K$ be
the complementary subset. To $K$ one can attach a standard parabolic
subgroup $B_K$ containing the Borel subgroup $B=HN$. We denote by
$N_K$ the unipotent radical of $B_K$. This is a subgroup of~$N$.
Let $W_K$ be the subgroup of the Weyl group $W$ generated by the 
reflexions $s_k\ (k\in K)$. This is a finite Coxeter group and
we denote by $w_0^K$ its longest element. The longest element
of $W$ is denoted by $w_0$.

In finite type, the preprojective algebra $\LL$ is finite-dimensional
and self-injective.
In agreement with \cite{GLSFlag}, we shall denote by $P_i$ the
indecomposable projective $\LL$-module with top $S_i$ and by
$Q_i$ the indecomposable injective module with socle $S_i$.
We write 
$$
Q_J = \bigoplus_{j\in J} Q_j 
\text{\;\;\; and \;\;\;}
P_J = \bigoplus_{j\in J} P_j.
$$

In \cite{GLSFlag}, we have shown that $\C[N_K]$ is naturally
isomorphic to the subalgebra
\[
R_K := \Span_\C\ebrace{\vph_X \mid X\in\Cogen(Q_J)}
\]
of $\C[N]$.
As before, $\Cogen(Q_J)$ is the full subcategory of $\md(\LL)$ whose objects
are submodules of direct sums of finitely many copies of $Q_J$.
This allowed us to introduce a cluster algebra $\cA_J \subseteq R_K$,
whose cluster monomials are of the form $\vph_X$ with $X$ a rigid
object of $\Cogen(Q_J)$. We conjectured that in fact $\cA_J = R_K$
\cite[Conj. 9.6]{GLSFlag}.

We are going to prove that this conjecture follows from the results
of this paper if $w_0w_0^K$ is adaptable.

\begin{Lem}
We have $N_K = N'(w_0^K) = N(w_0w_0^K)$.
\end{Lem}

\begin{proof}
We know that $N'(w_0^K)$ is the subgroup of $N$ 
generated by the one-parameter subgroups $N(\alpha)$ with $\alpha>0$
and $w_0^K(\alpha)>0$. These are exactly the one-parameter subgroups
of $N$ which do not belong to the Levi subgroup of $B_K$, hence the
first equality follows.
Now, since $N=w_0N_-w_0$, we have 
\[
 N'(w_0^K)=N\cap w_0^KNw_0^K = N\cap w_0^Kw_0N_{-}w_0w_0^K=N(w_0w_0^K).
\]
\end{proof}

As before, let 
$\Gen(P_J)$ be
the subcategory of $\md(\LL)$ whose objects
are factor modules of direct sums of finitely many copies of $P_J$.

\begin{Lem}
We have $\CC_{w_0w_0^K} = \Gen(P_J)$.
\end{Lem}

\begin{proof}
By Proposition~\ref{minors}, we know that the indecomposable 
projective-injective object $I_i$ of $\CC_{w_0w_0^K}$ with socle $S_i$
satisfies 
\[
\vph_{I_i} = D_{\vpi_i,w_0^Kw_0(\vpi_i)},\qquad (i\in I).
\]
By \cite[\S 6.2]{GLSFlag}, it follows that $I_i={\mathcal E}_{w_0^K}Q_i$,
where ${\mathcal E}_{w_0^K}$ is the functor defined in \cite[\S 5.2]{GLSFlag}.
It readily follows that $I_i$ is the projective-injective indecomposable 
object of $\Gen(P_J)$ with simple socle $S_i$. 
Hence $\CC_{w_0w_0^K}$ and $\Gen(P_J)$ have the same projective-injective
generator.
\end{proof}

Let $S$ denote the self-duality of $\md(\LL)$ induced 
by the involution $a \mapsto a^*$
mapping an arrow $a$ of $\overline{Q}$ to its opposite arrow $a^*$
(see \cite[\S 1.7]{GLSAus}).
This restricts to an
anti-equivalence of categories 
$\Gen(P_J) \to \Cogen(Q_J)$, that we shall again denote by $S$.

\begin{Lem}\label{lem:inverse}
For every $X\in\md(\LL)$ and every $n\in N$ we have
\[
\vph_{X}(n^{-1})=(-1)^{\dm X}\vph_{S(X)}(n).
\]
\end{Lem}

\begin{proof}
We know that $N$ is generated by the one-parameter subgroups $x_i(t)$
attached to the simple positive roots.
By Proposition~\ref{phi_form} we have
\[
\vph_X(x_{i_1}(t_1)\cdots x_{i_k}(t_k))=
\sum_{\aaa=(a_1,\ldots,a_k)\in\N^k}\chi_c(\F_{\ii^\aaa,X})
\frac{t_1^{a_1}\cdots t_k^{a_k}}{a_1!\cdots a_k!}.
\] 
Now, if $n=x_{i_1}(t_1)\cdots x_{i_k}(t_k)$, we have
$n^{-1}=x_{i_k}(-t_k)\cdots x_{i_1}(-t_1)$ and the result 
follows from the fact that 
$\F(\ii^\aaa,X) \cong \F(\ii_\op^{\aaa_\op},S(X))$,
where $\ii_\op$ and $\aaa_\op$ denote the sequences 
obtained by reading $\ii$ and $\aaa$ from right to left.
\end{proof}

We can now prove:

\begin{Thm}
Conjecture  9.6 of \cite{GLSFlag} holds if
$w_0w_0^K$ is adaptable.
\end{Thm}

\begin{proof}
Suppose that $w_0w_0^K$ is $Q$-adapted. Let $\CC_M=\CC_{w_0w_0^K}$
be the corresponding subcategory of $\md(\LL)$.
The cluster algebra $\RR(\CC_M)= \RR(\Gen(P_J))$ is 
isomorphic to $\cA_J$ via the map $\vph_X\mapsto\vph_{S(X)}$.
This comes from the fact that $S\df \Gen(P_J) \to \Cogen(Q_J)$
is an anti-equivalence which maps the maximal rigid module
$T_M^\vee$ used to define the initial seed of $\RR(\CC_M)$
to the maximal rigid module
$U_\ii$ of \cite[\S 9.2]{GLSFlag} used to define the initial seed of $\cA_J$. 
(Here we assume that $\ii$ is the reduced word of $w_0^Kw_0$
obtained by reading the $Q$-adapted reduced word of $w_0w_0^K$
from right to left.)
In particular the cluster variables $\vph_{M_i}$ which,
by Theorem~\ref{proofmain1}, generate
$\RR(\Gen(P_J))=\C[N(w_0w_0^K)]$ are mapped to cluster variables
$\vph_{S(M_i)}$ of $\cA_J$. They also form a system of generators of the 
polynomial algebra $\C[N(w_0w_0^K)]=\C[N_K]$ by Lemma~\ref{lem:inverse},
because the map $n\mapsto n^{-1}$ is a biregular automorphism of $N_K$.
Hence $\cA_J=\C[N_K]$.  
\end{proof}

\begin{Rem}\label{finrem}
{\rm
The previous discussion shows that we obtain two different cluster
algebra structures on $\C[N_K]$, coming from the two different subcategories
$\Gen(P_J)$ and $\Cogen(Q_J)$.

When using $\Gen(P_J)=\CC_{w_0w_0^K}$, we regard $\C[N_K]$ as the subring
of $N'(w_0w_0^K)$-invariant functions of $\C[N]$ for the action of
$N'(w_0w_0^K)$ on $N$ by \emph{right} translations 
(see Section~\ref{invariant_ring}).  
This allows us to relate the first cluster structure to the cluster structure
of the unipotent cell $\C[N^{w_0w_0^K}]$ (see Proposition~\ref{coordNw}).

When using $\Cogen(Q_J)$, we regard $\C[N_K]$ as the subring
of $N'(w_0w_0^K)$-invariant functions of $\C[N]$ for the action of
$N'(w_0w_0^K)=N(w_0^K)$ on $N$ by \emph{left} translations.
These functions can then be ``lifted'' to $B_K^{-}$-invariant functions
on $G$ for the action of $B_K^{-}$ on $G$ by left translations.
This allows us to ``lift'' the second cluster structure
to a cluster structure on $\C[B_K^-\backslash G]$ (see \cite[\S 10]{GLSFlag}).
}
\end{Rem}


{\Large\section{Open problems}
\label{openproblems}}


\subsection{}
It is known that the dual canonical basis ${\mathcal B}^*$ 
and the dual semicanonical
basis $\cS^*$ of $\MM^* \equiv U(\n)^*_{\rm gr}$ 
do not coincide, see \cite{GLSSemi1}.
But one might at least hope that both bases have some interesting
elements in common:

\begin{Conj}[Open Orbit Conjecture]\label{finalconj2}
Let $Z$ be an irreducible component of $\LL_d$,
and let $b_Z$ and $\rho_Z$ be the associated dual canonical
and dual semicanonical basis vectors of $\MM^*$.
If $Z$ contains an open $\GL_d$-orbit, then
$b_Z = \rho_Z$.
\end{Conj}

We know that each cluster monomial of the cluster algebra 
$\cA(\CC_M)$ is of the form $\rho_Z$, where
$Z$ contains an open $\GL_d$-orbit.
So if the conjecture is true, then all cluster monomials belong
to the dual canonical basis.

\subsection{}
Finally, we would like to ask the following question.
Is it possible to realize every element of
the dual canonical basis of $\MM^*$ as 
a $\delta$-function?
We know several examples of elements $b$ of ${\mathcal B}^*$ which
do not belong to $\cS^*$. In all these examples,
$b$ is however equal to $\delta_X$ for a non-generic $\Lambda$-module $X$.

Let us look at an example.
Let $Q$ be the quiver 
{\small
$$
\xymatrix{
1 & 2 \ar@<0.5ex>[l]\ar@<-0.5ex>[l]
}
$$
}\noindent
and let $\LL$ be the associated preprojective algebra.
For $\lambda \in \C^*$ we define representations $M(\lambda,1)$ and
$M(\lambda,2)$ of $Q$ as follows:
$$
M(\lambda,1) := \;
\xymatrix@+1.0pc{
\C & \C \ar@<0.5ex>[l]^{\left(\bsm \lambda \esm\right)}
\ar@<-0.5ex>[l]_{\left(\bsm 1 \esm\right)}
}
\text{\;\; and \;\;}
M(\lambda,2) := \;
\xymatrix@+1.0pc{
\C^2 & \C^2 \ar@<0.5ex>[l]^{\left(\bsm \lambda & 1\\0&\lambda\esm\right)}
\ar@<-0.5ex>[l]_{\left(\bsm 1 &0\\0&1\esm\right)}
}
$$
Let $\iota\df \rep(Q,(2,2)) \to \LL_{(2,2)}$ be the canonical
embedding, defined by $M' \mapsto (M',0)$.
Clearly, the image of $\iota$ is an irreducible component
of $\LL_{(2,2)}$, which we denote by $Z_Q$.
It is not difficult to check that the set
$$
\left\{ (M(\lambda,1) \oplus M(\mu,1),0) \mid \lambda,\mu \in \C^* \right\}
$$
is a dense subset of $Z_Q$.
It follows that 
$$
\delta_{(M(\lambda,1) \oplus M(\mu,1),0)} = \rho_{Z_Q}
$$
is an element of the dual semicanonical basis $\cS^*$.
It is easy to check that 
$$
\delta_{(M(\lambda,2),0)}\not = \delta_{(M(\lambda,1) \oplus M(\mu,1),0)}.
$$
Indeed the variety of ascending flags of type $(1,2,1,2)$ 
has Euler characteristic $2$ for $(M(\lambda,1) \oplus M(\mu,1),0)$
and Euler characteristic $1$ for $(M(\lambda,2),0)$.
Furthermore, one can show that
$$
\delta_{(M(\lambda,2),0)} = b_{Z_Q}
$$ belongs to the dual canonical basis 
${\mathcal B}^*$ of $\MM^*$.

Note that the functions $\delta_{(M(\lambda,1) \oplus M(\mu,1),0)}$
and $\delta_{(M(\lambda,2),0)}$ do not belong to any of the
subalgebras $\RR(\CC_M)$, since 
$M(\lambda,1)$ and $M(\lambda,2)$ are regular representations
of the quiver $Q$ for all $\lambda$.


\bigskip
{\parindent0cm \bf Acknowledgements.}\,
We are grateful to {\O}yvind Solberg for answering our
questions on relative homology theory.
We thank Shrawan Kumar for his kind help concerning
Kac-Moody groups.
It is a pleasure to thank the Mathematisches Forschungsinstitut
Oberwolfach (MFO) for two weeks of hospitality  
in August/September 2006, where part of this work was done.
Furthermore, the first and second authors like to thank the 
Max-Planck Institute for Mathematics in Bonn for a research stay in
September - December 2006 and October 2006, respectively.
Finally, the three authors are grateful to the Mathematical
Sciences Research Institute in Berkeley (MSRI) for 
an invitation in Spring 2008 during which this work was finalized. 
The first author was partially supported by PAPIIT grant IN103507-2.




\begin{thebibliography}{999}


\bibitem[Ab]{Ab}
{\it E. Abe},
Hopf algebras.
Cambridge Tracts in Mathematics 74.
Cambridge University Press, 1980.

\bibitem[ASS]{ASS}
{\it I. Assem, D. Simson, A. Skowro{\'n}ski},
Elements of the representation theory of associative algebras. 
Vol. 1. Techniques of representation theory. 
London Mathematical Society Student Texts, 65. 
Cambridge University Press, Cambridge, 2006. x+458 pp.

\bibitem[A]{A}
{\it M. Auslander},
Representation theory of artin algebras. I. 
Comm. Alg. 1 (1974), 177--268. 

\bibitem[APR]{APR}
{\it M. Auslander, M. Platzeck, I. Reiten},
Coxeter functors without diagrams.
Trans. Amer. Math. Soc. 250 (1979), 1--46.

\bibitem[ARS]{ARS}
{\it M. Auslander, I. Reiten, S. Smal\o}, 
Representation theory of Artin algebras. 
Corrected reprint of the 1995 original. 
Cambridge Studies in Advanced Mathematics, 36. Cambridge University 
Press, Cambridge, 1997. xiv+425 pp. 

\bibitem[AS1]{AS1}
{\it M. Auslander, \O. Solberg},
Relative homology and representation theory. I.
Relative homology and homologically finite subcategories.
Comm. Algebra 21 (1993), no. 9, 2995--3031.

\bibitem[AS2]{AS2}
{\it M. Auslander, \O. Solberg},
Relative homology and representation theory. II.
Relative cotilting theory.
Comm. Algebra 21 (1993), no. 9, 3033--3079.


\bibitem[BZ]{BZ1}
{\it A. Berenstein, A. Zelevinsky},
Total positivity in Schubert varieties.
Comment. Math. Helv. 72 (1997), 128--166.

\bibitem[Be]{Be}
{\it R. B\'edard},
On commutation classes of reduced words in Weyl groups.
Europ. J. Combinatorics (1999) 20, 483--505.

\bibitem[BFZ]{BFZ}
{\it A. Berenstein, S. Fomin, A. Zelevinsky},
Cluster algebras III: Upper bounds and double Bruhat
cells.
Duke Math. J. 126 (2005), no. 1, 1--52.

\bibitem[Bo1]{Bo1}
{\it K. Bongartz}, 
Algebras and quadratic forms. 
J. London Math. Soc. (2) 28 (1983), no. 3, 461--469. 

\bibitem[Bo2]{Bo2}
{\it K. Bongartz}, 
Tilted algebras.
Representations of algebras (Puebla, 1980), pp. 26--38,
Lecture Notes in Math., 903, Springer, Berlin-New York, 1981.

\bibitem[BMRRT]{BMRRT}
{\it A. Buan, R. Marsh, M. Reineke, I. Reiten, G. Todorov},
Tilting theory and cluster combinatorics.  
Adv. Math. 204 (2006), no. 2, 572--618.

\bibitem[BM]{BM}
{\it A. Buan, R. Marsh},
Cluster-tilting theory.
In: Trends in Representation Theory of Algebras and Related Topics, 
Contemp. Math. 406 (2006), 1--30.

\bibitem[BIRS]{BIRS}
{\it A. Buan, O. Iyama, I. Reiten, J. Scott},
Cluster structures for 2-Calabi-Yau categories and
unipotent groups.
arXiv:math.RT/0701557.

\bibitem[CK1]{CK}
{\it P. Caldero, B. Keller},
From triangulated categories to cluster algebras.
Invent. Math. 172 (2008), 169--211.

\bibitem[CK2]{CK2}
{\it P. Caldero, B. Keller},
From triangulated categories to cluster algebras II.
Ann. Sci. {\'E}cole Norm. Sup. (4) 39 (2006), no. 6, 983--1009. 

\bibitem[CR]{CR}
{\it P. Caldero, M. Reineke},
On the quiver Grassmannian in the acyclic case.
J. Pure Appl. Algebra (to appear).

\bibitem[CPS]{CPS}
{\it E. Cline, B. Parshall, L. Scott},
Finite-dimensional algebras and highest weight categories.  
J. Reine Angew. Math. 391 (1988), 85--99. 

\bibitem[CB1]{CB}
{\it W. Crawley-Boevey}, 
On the exceptional fibres of Kleinian singularities.
Amer. J. Math. 122 (2000), 1027--1037.

\bibitem[CB2]{CB2}
{\it W. Crawley-Boevey}, 
Geometry of the moment map for representations of quivers.
Compositio Math. 126 (2001), no. 3, 257--293.

\bibitem[DFK]{DFK}
{\it P. Di Francesco, R. Kedem},
Q-systems as cluster algebras II: Cartan matrix of finite type and 
the polynomial property.
arXiv:0803.0362. 

\bibitem[D]{D}
{\it J. Dixmier},
Enveloping algebras.
North-Holland Mathematical Library, Vol. 14. Translated from the French.
North-Holland Publishing Co., Amsterdam-Now York-Oxford, 1977. xvi+375 pp.

\bibitem[FZ1]{FZ5}
{\it S. Fomin, A. Zelevinsky}, 
Double Bruhat cells and total positivity.
J. Amer. Math. Soc. 12 (1999), no. 2, 335--380.

\bibitem[FZ2]{FZ1}
{\it S. Fomin, A. Zelevinsky}, 
Cluster algebras. I. Foundations. 
J. Amer. Math. Soc. 15 (2002), no. 2, 497--529.

\bibitem[FZ3]{FZ2}
{\it S. Fomin, A. Zelevinsky}, 
Cluster algebras. II. Finite type classification.
Invent. Math. 154 (2003), no. 1, 63--121. 

\bibitem[FZ4]{FZSurv}
{\it S. Fomin, A. Zelevinsky}, 
Cluster algebras: notes for the CDM-03 conference.  
Current developments in mathematics, 2003,  1--34, Int. Press,
Somerville, MA, 2003.

\bibitem[FZ5]{FZ3}
{\it S. Fomin, A. Zelevinsky}, 
Cluster algebras. IV. Coefficients.
Compositio Math. 143 (2007), 112--164.

\bibitem[GR]{GR}
{\it P. Gabriel, A.V. Roiter},
Representations of finite-dimensional algebras.
Translated from the Russian.
With a chapter by B. Keller.
Reprint of the 1992 English translation.
Springer-Verlag, Berlin, 1997. iv+177 pp.

\bibitem[GLS1]{GLSSemi1}
{\it C. Gei{\ss}, B. Leclerc, J. Schr\"oer},
Semicanonical bases and preprojective algebras.
Ann. Sci. {\'E}cole Norm. Sup. (4)  38  (2005),  no. 2, 193--253.

\bibitem[GLS2]{GLSAus}
{\it C. Gei{\ss}, B. Leclerc, J. Schr\"oer},
Auslander algebras and initial seeds for cluster algebras.
J. London Math. Soc. (2) 75 (2007), 718--740.

\bibitem[GLS3]{GLSVerma}
{\it C. Gei{\ss}, B. Leclerc, J. Schr\"oer},
Verma modules and preprojective algebras.
Nagoya Math. J. 182 (2006), 241--258.

\bibitem[GLS4]{GLSSemi2}
{\it C. Gei{\ss}, B. Leclerc, J. Schr\"oer},
Semicanonical bases and preprojective algebras II: A multiplication
formula.
Compositio Math. 143 (2007), 1313--1334.

\bibitem[GLS5]{GLSRigid}
{\it C. Gei{\ss}, B. Leclerc, J. Schr\"oer},
Rigid modules over preprojective algebras.
Invent. Math. 165 (2006), no. 3, 589--632.

\bibitem[GLS6]{GLSFlag}
{\it C. Gei{\ss}, B. Leclerc, J. Schr\"oer},
Partial flag varieties and preprojective algebras.
arXiv:math.RT/0609138, 
Ann. Inst. Fourier (Grenoble) (to appear).

\bibitem[GS]{GSExt}
{\it C. Gei{\ss}, J. Schr\"oer},
Extension-orthogonal components of preprojective varieties.
Trans. Amer. Math. Soc. 357 (2005), 1953--1962. 

\bibitem[H1]{H1}
{\it D. Happel},
On the derived category of a finite-dimensional algebra.  
Comment. Math. Helv.  62 (1987),  no. 3, 339--389. 

\bibitem[H2]{H2}
{\it D. Happel},
Triangulated categories in the representation theory of 
finite-dimensional algebras.
London Mathematical Society Lecture Note Series, 119. 
Cambridge University Press, Cambridge, 1988. x+208 pp.

\bibitem[H3]{H3}
{\it D. Happel},
Partial tilting modules and recollement.
Proceedings of the International Conference on Algebra, Part 2 
(Novosibirsk, 1989), 345--361, Contemp. Math., 131, Part 2, Amer. Math. Soc.,
Providence, RI, 1992.

\bibitem[HR]{HR}
{\it D. Happel, C.M. Ringel},
Tilted algebras.  
Trans. Amer. Math. Soc.  274 (1982), no. 2, 399--443. 

\bibitem[HU1]{HU4}
{\it D. Happel, L. Unger},
On the quiver of tilting modules.
J. Algebra 284 (2005), 857--868.

\bibitem[HU2]{HU3}
{\it D. Happel, L. Unger},
On a partial order of tilting modules.
Algebras and Representation Theory 8 (2005), 147--156.

\bibitem[Ig]{Ig}
{\it K. Igusa}, Notes on the no loops conjecture.
J. Pure Appl. Algebra 69 (1990), 161--176.


\bibitem[Iy1]{I1}
{\it O. Iyama},
Higher dimensional Auslander-Reiten theory on maximal orthogonal subcategories.
Adv. Math. 210 (2007), no. 1, 22--50. 

\bibitem[Iy2]{I2}
{\it O. Iyama},
Auslander correspondence. 
Adv. Math. 210 (2007),  no. 1, 51--82.

\bibitem[K1]{K2}
{\it V. Kac},
Infinite root systems, representations of graphs and invariant theory.
Invent. Math. 56 (1980), no. 1, 57--92.

\bibitem[K2]{K}
{\it V. Kac},
Infinite-dimensional Lie algebras. 
Third edition. 
Cambridge University Press, Cambridge, 1990. xxii+400 pp. 

\bibitem[KP]{KP1}
{\it V. Kac, D. Peterson},
Regular functions on certain infinite-dimensional groups.
Arithmetic and geometry, Vol II, 141--166, Progr. Math., 36,
Birkh\"auser Boston, Boston, MA, 1983.

\bibitem[KS]{KSa}
{\it M. Kashiwara, Y. Saito},
Geometric construction of crystal bases.
Duke Math. J. 89 (1997), 9--36.

\bibitem[K]{Ke}
{\it B. Keller},
On triangulated orbit categories.  
Doc. Math. 10 (2005), 551--581 (electronic).

\bibitem[KR]{KR}
{\it B. Keller, I. Reiten},
Acyclic Calabi-Yau categories. 
arXiv:math.RT/0610594.

\bibitem[Ku]{Ku}
{\it S. Kumar},
Kac-Moody groups, their flag varieties and representation theory. 
Progress in Mathematics, 204. Birkh\"auser Boston, Inc., Boston, MA, 2002. 

\bibitem[Le]{LeMZ}
{\it B. Leclerc},
Dual canonical bases, quantum shuffles and $q$-characters.  
Math. Z.  246  (2004), 691--732. 

\bibitem[L1]{Lu1}
{\it G. Lusztig},
Quivers, perverse sheaves, and quantized enveloping algebras. 
J. Amer. Math. Soc. 4 (1991), no. 2, 365--421.

\bibitem[L2]{Lu3}
{\it G. Lusztig},
Constructible functions on the Steinberg variety. 
Adv. Math. 130 (1997), no. 2, 287--310.

\bibitem[L3]{Lu2}
{\it G. Lusztig},
Semicanonical bases arising from enveloping algebras. 
Adv. Math. 151 (2000), no. 2, 129--139.


\bibitem[R]{R}
{\it N. Richmond},
A stratification for varieties of modules.
Bull. London Math. Soc. 33 (2001), no. 5, 565--577.

\bibitem[RS]{RS2}
{\it C. Riedtmann, A. Schofield},
On a simplicial complex associated with tilting modules. 
Comment. Math. Helv. 66 (1991), no. 1, 70--78. 


\bibitem[Ri1]{Ri1}
{\it C.M. Ringel},
Tame algebras and integral quadratic forms. 
Lecture Notes in Mathematics, 1099. 
Springer-Verlag, Berlin, 1984. xiii+376 pp. 

\bibitem[Ri2]{Ri3}
{\it C.M. Ringel},
Representation theory of finite-dimensional algebras.
In: Representations of algebras (Durham, 1985), 7--79, 
London Math. Soc. Lecture Note Ser., 116, 
Cambridge Univ. Press, Cambridge, 1986.

\bibitem[Ri3]{Ri5}
{\it C.M. Ringel},
The category of modules with good filtrations over a quasi-hereditary 
algebra has almost split sequences.  
Math. Z.  208  (1991),  no. 2, 209--223.

\bibitem[Ri4]{Ri4}
{\it C.M. Ringel},
The category of good modules over a quasi-hereditary algebra.  
Proceedings of the Sixth International Conference on Representations of 
Algebras (Ottawa, ON, 1992),  17 pp., 
Carleton-Ottawa Math. Lecture Note Ser., 14, Carleton Univ., Ottawa, ON, 1992. 

\bibitem[Ri5]{Ri2}
{\it C.M. Ringel},
The preprojective algebra of a quiver.  
Algebras and modules, II (Geiranger, 1996),  467--480, CMS Conf. Proc., 24,
Amer. Math. Soc., Providence, RI, 1998.

\bibitem[Ri6]{Ri6}
{\it C.M. Ringel},
PBW-bases of quantum groups.
J. Reine Angew. Math. 470 (1996), 51--88.

\bibitem[S]{S}
{\it A. Schofield},
Quivers and Kac-Moody Lie algebras.
Unpublished Preprint.


\bibitem[SZ]{SZ}
{\it P. Sherman, A. Zelevinsky},
Positivity and canonical bases in rank 2 cluster algebras of
finite and affine type.
Moscow Math. J. 4 (2004), no. 4, 947--974.

\end{thebibliography}
\end{document}